\documentclass[12pt,a4paper,twoside,openright]{report}

\usepackage[dvipsnames]{xcolor}
\usepackage[unicode,psdextra,final]{hyperref}
\hypersetup{
  colorlinks=true,
  linkcolor=Purple,
  citecolor=cyan,  
}
\usepackage{mystyle}
\usepackage{times}
\usepackage{sectsty}
\allsectionsfont{\scshape}
\usepackage{geometry}
\usepackage{fancyhdr}

\begin{document}
\newgeometry{margin=1in}
\begin{titlepage}
   \includegraphics[width=0.4\textwidth]{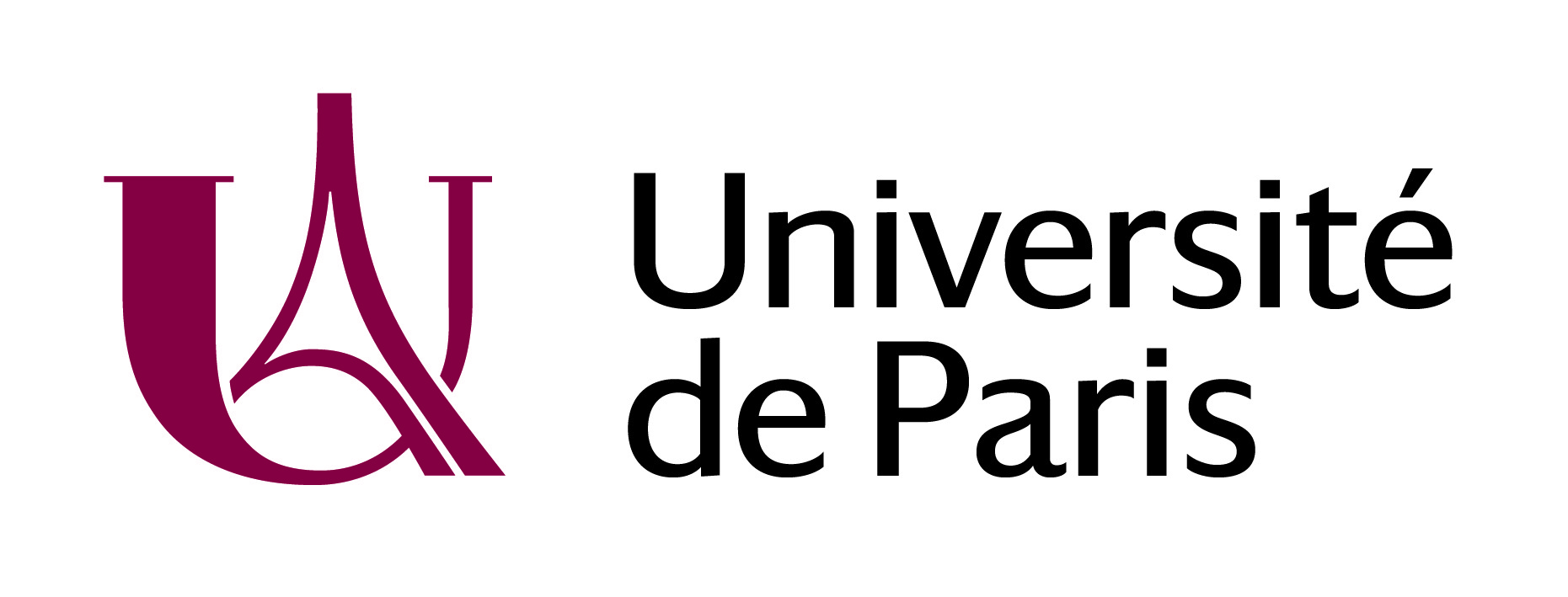}
    
       \vspace*{1cm}
   \raggedleft
   {
   \scshape{\textbf{Université de Paris}} \par}
 \vspace{0.25cm}
 {\large École de Sciences Mathématiques Paris Centre (ED 386)}
 
 {Institut de Recherche Fondamentale en Informatique (IRIF)}

 \vspace*{0.25cm}

 {\Large
   \scshape
   \textbf{Thèse de doctorat en Mathématiques}}

      {\large Dirigée par François Métayer \\
        et par Clemens Berger \par}

      \vspace*{2cm}
      \hrulefill
      \vspace*{0.25cm}
      
      \centering
        {\Huge \scshape \textbf{Homology of strict $\omega$-categories}\par}

        \vspace*{0.25cm}
        \hrulefill
        
        \vspace{1cm}

        {\Large
          
        \textbf{par Léonard Guetta}}

   \vspace*{2cm}
   {\slshape Présentée et soutenue publiquement en ligne le 28 janvier 2021 \\
     devant un jury composé de : \par}
    \vspace{0.8cm}

    \begin{tabular}{ll}
      \textbf{M. Benoît Fresse} (Prof., Université de Lille) & Rapporteur \\
      \textbf{Mme Eugenia Cheng} (Prof., School of the Art Institute) & Examinatrice \\
      \textbf{M. Carlos Simpson} (DR CNRS, Université de Nice) & Président du jury \\
      \textbf{M. Samuel Mimram} (Prof., École Polytechnique) & Examinateur \\
      \textbf{M. Georges Maltsiniotis} (DR CNRS, Université de Paris) & Examinateur \\
      \textbf{M. François Métayer} (MCF, Université de Paris Nanterre) & Directeur de Thèse \\
      \textbf{M. Clemens Berger} (MCF, Université de Nice) & Co--directeur de Thèse \\
      \textbf{M. Dimitri Ara} (MCF, Université Aix-Marseille) & Membre invité 
    \end{tabular}
    
      \vspace{0.8cm}
      {\slshape Rapporteur non présent à la soutenance : \par}
           \vspace{0.8cm}
    \begin{tabular}{ll}
      \textbf{M. Richard Garner} (Senior Lecturer, Macquarie University) & \\
    \end{tabular}

  
    \end{titlepage}

\selectlanguage{french}
\begin{abstract}
Dans cette thèse, on compare l'homologie \og classique \fg{} d'une
$\oo$\nbd{}catégorie (définie comme l'homologie de son nerf de Street) avec
son homologie polygraphique. Plus précisément, on prouve que les deux
homologies ne coïncident pas en général et qualifions d'\emph{homologiquement
  cohérentes} les $\oo$\nbd{}catégories particulières pour lesquelles l'homologie
polygraphique coïncide effectivement avec l'homologie du nerf. Le but poursuivi
est de trouver des critères abstraits et concrets permettant de détecter les
$\oo$\nbd{}catégories homologiquement cohérentes. Par exemple, on démontre que
toutes les (petites) catégories, que l'on considère comme des
$\oo$\nbd{}catégories strictes dont toutes les cellules au-delà de la dimension
$1$ sont des unités, sont homologiquement cohérentes. On introduit également la
notion de $2$\nbd{}catégorie \emph{sans bulles} et on conjecture qu'une
$2$\nbd{}catégorie cofibrante est homologiquement cohérente si et seulement si
elle est sans bulles. On démontre également des résultats importants concernant
les $\oo$\nbd{}catégories strictes qui sont libres sur un polygraphe, comme
le fait que si $F : C \to D$ est un $\oo$\nbd{}foncteur de Conduché discret et
si $D$ est libre sur un polygraphe alors $C$ l'est aussi. Dans son ensemble,
cette thèse établit un cadre général dans lequel étudier l'homologie des
$\oo$\nbd{}catégories en faisant appel à des outils d'algèbre homotopique
abstraite, tels que la théorie des catégories de modèles de Quillen ou la théorie
des dérivateurs de Grothendieck.

\bigskip

\noindent\textbf{Mots-clés : } Catégories supérieures, $\oo$\nbd{}catégories, homologie,
théorie de l'homotopie, polygraphes.
\end{abstract}
\selectlanguage{english}
\frenchspacing
\begin{abstract}
  In this dissertation, we compare the ``classical''
  homology of an $\oo$\nbd{}category (defined as the homology of its Street
  nerve) with its polygraphic homology. More precisely, we prove that both
  homologies generally do not coincide and call \emph{homologically coherent} the
  particular strict $\oo$\nbd{}categories for which polygraphic homology and
  homology of the nerve do coincide. The goal pursued is to find abstract
  and concrete criteria to detect homologically coherent $\oo$\nbd{}categories. For
  example, we prove that all (small) categories, considered as strict
  $\oo$\nbd{}categories with unit cells above dimension $1$, are homologically
  coherent. We also introduce the notion of \emph{bubble-free} $2$\nbd{}category
  and conjecture that a cofibrant $2$\nbd{}category is homologically
  coherent if and only if it is bubble-free.
  We also prove important results concerning free strict
  $\oo$\nbd{}categories on polygraphs (also known as computads), such as the
  fact that if $F : C \to D$ is a discrete Conduché $\oo$\nbd{}functor and $D$
  is a free strict $\oo$\nbd{}category on a polygraph, then so is $C$.
  Overall, this thesis achieves to build a general framework in which to study the
  homology of strict $\oo$\nbd{}categories using tools of abstract homotopical
  algebra such as Quillen's theory of model categories or Grothendieck's theory
  of derivators.
  
\bigskip

  
\end{abstract}
\restoregeometry
\pagestyle{empty}
\clearpage\mbox{}\clearpage 

\pagenumbering{roman}
\pagestyle{plain}
\chapter*{Remerciements}
\addcontentsline{toc}{chapter}{Remerciements} 
\selectlanguage{french}
Avec toute thèse, viennent les habituels remerciements. L'écriture de ces derniers
semble presque aussi délicate que l'écriture de la thèse elle-même car, il faut
bien le dire, la plupart des lecteurs ne s'aventurent pas au-delà (et je
les comprends !). Je n'ai d'ailleurs jamais réussi à expliquer à mes proches
les mathématiques que je fais, qui sont pourtant vachement chouettes. 

Merci à François Métayer et Clemens Berger d'avoir dirigé ma thèse avec autant de bienveillance
et de gentillesse. Vous avez su me guider tout en me laissant libre de suivre
mes intuitions,
parfois de façon têtue, je dois le reconnaître. Je voudrais aussi remercier Georges Maltsiniotis
au même titre que mes directeurs de thèse. Véritable mentor pour moi, j'ose
orgueilleusement me penser son élève par adoption. Il m'est impossible de
résumer en quelques mots, ni même quelques phrases, la portée de son influence
mathématique sur moi et je ne compte plus le nombre de fois où ses explications ont rendu
limpides des concepts qui me paraissaient mystérieux. 

Je dois également remercier chaleureusement Alain Prouté, dont la rencontre a
été déterminante pour cette thèse. C'est avec lui que j'ai appris la topologie
algébrique et la théorie des catégories, ce qui résume essentiellement les
mathématiques que je fais aujourd'hui. C'est également grâce à lui que j'ai
rencontré François, dans le cadre de mon mémoire de M2 qu'ils ont co-encadré.

Je voudrais remercier Benoît Fresse et Richard Garner d'avoir accepté de
rapporter ma thèse.  Richard, I am
indebted to you for having read my dissertation so carefully and thoroughly. I
sincerely hope that we will have the occasion to meet and discuss in the future.
Je voudrais également remercier Eugenia Cheng, Carlos
Simpson, Samuel Mimram, Georges Maltsiniotis et Dimitri Ara d'avoir accepté de
faire partie de mon jury. J'aurais beaucoup apprécié pouvoir échanger de vive
voix avec vous lors de la soutenance et je suis très attristé que le climat
actuel ne nous permette pas de nous réunir physiquement.

Merci à Albert Burroni, Mathieu Anel, Eric Finster, Yves Guiraud, Pierre-Louis
Curien, Bérénice Delcroix-Oger, Yves Lafont et tous les autres habitués du groupe de travail ``catégories
supérieures, polygraphes et homotopie'' du vendredi après-midi, avec qui il a
été si agréable d'échanger.

Merci à tous les doctorants que j'ai croisés à l'IRIF avec qui j'ai eu la chance de partager un
bout de chemin et le pain au déjeuner : Cyrille, Maxime, Pierre V., Pierre C.,
Jovana, Cédric, Zeinab, Léo, Rémi, Costia, Nicolas, Victor, Loïc, El-Mehdi, Alen, Loïc,
Farzad, Gianluca, Mehdi, Simon, Kenji,... Il y a parmi eux plusieurs à qui je dois une mention
spéciale. Antoine, avec qui j'ai tant apprécié discuter. Chaitanya qui, en plus d'être un ami cher, a eu la patience de
corriger mes nombreuses fautes d'anglais dans cette thèse, et
pour ça je lui suis extrêmement redevable. Enfin Jules, qui fut mon
veritable compagnon de route ces dernières années. Merci pour cette \og colla-bureaution\fg{} qui a résulté en une sincère amitié. Je ris
d'avance en pensant au jour où les mémoires auront oublié qui furent les
responsables de la décoration du bureau 3055.
Et puis bien sûr, il y a Daniel; de toutes les choses qui le rendent si
formidablement singulier, je crois
que ce que j'admire le plus c'est sa liberté d'esprit. Je compte bien te faire
prendre exceptionnellement des vacances pour venir me rendre visite quand je
serai parti de l'IRIF.

Merci à mes s\oe{}urs Margot, Naomi et Oriane d'êtres présentes dans ma vie. Vous
occupez une place dans mon c\oe{}ur bien plus grande que vous ne le soupçonnez
certainement. Merci Maman de m'avoir toujours poussé et encouragé dans
toutes mes entreprises. Il ne fait aucun doute que je n'en serais pas là sans
toi. Merci Papa d'avoir éveillé en moi le goût pour les
énigmes et la fascination pour la science. Le fait que tu te sois remis à
apprendre les mathématiques après plus de quarante ans est absolument
fantastique et admirable; j'ai hâte de pouvoir partager avec toi ce que je fais.
Merci à mes grands-parents Pierrot, Michèle et Francine d'avoir fait de mon
enfance une époque heureuse. Merci aussi à la famille Moatti de m'avoir acceuilli affectueusement. Une
grande partie de cette thèse a été écrite chez vous en Provence et la
chaleur de votre entourage a été essentielle.

Une pensée également à David, Willie, Ariel, Samuel, Paul et tous ceux qui
remplissent ma vie d'amitié.

Merci aussi à Georges Brassens et Django Reinhardt, sans qui j'aurais terminé
cette thèse dans les temps. 

Et puis bien évidemment, il y a Emma. Le rayonnement que tu apportes dans ma vie
rend toute tentative de remerciement d'une fadeur ingrate. Indirectement, cette
thèse doit autant à toi qu'à moi et c'est pourquoi elle t'est dédiée...

\selectlanguage{english}
\frenchspacing

\tableofcontents

\newpage
\pagenumbering{arabic}
\pagestyle{fancy}
\fancyhf{}
\fancyfoot[C]{\thepage}
\fancyhead[RO,LE]{INTRODUCTION}
\chapter*{Introduction}
\addcontentsline{toc}{chapter}{Introduction} 

The general framework in which this dissertation takes place is the
\emph{homotopy theory of strict $\oo$\nbd{}categories}, and, as the title
suggests, its focus is on homological aspects of this theory. The goal is to
study and compare two different homological invariants for strict
$\oo$\nbd{}categories; that is to say, two different
functors \[\mathbf{Str}\oo\Cat \to \ho(\Ch)\] from the category of strict
$\oo$\nbd{}categories to the homotopy category of non-negatively graded chain
complexes (i.e.\ the localization of the category of non-negatively graded chain
complexes with respect to the quasi\nbd{}isomorphisms).

Before we enter into the heart of the subject, let us emphasize that, with the
sole exception of the end of this introduction, all the $\oo$\nbd{}categories
that we consider are strict $\oo$\nbd{}categories. Hence, we drop the adjective
``strict'' and simply say \emph{$\oo$\nbd{}category} instead of \emph{strict
  $\oo$\nbd{}category} and we write $\oo\Cat$ instead of $\mathbf{Str}\oo\Cat$
for the category of (strict) $\oo$\nbd{}categories.

\begin{named}[Background: $\oo$-categories as spaces] The homotopy theory of
  $\oo$\nbd{}categories begins with the nerve functor introduced by Street in
  \cite{street1987algebra}
  \[
    N_{\omega} : \oo\Cat \to \Psh{\Delta}
  \]
  that associates to every $\oo$\nbd{}category $C$ a simplicial set $N_{\oo}(C)$
  called the \emph{nerve of $C$}, generalizing the usual nerve of (small)
  categories. Using this functor, we can transfer the homotopy theory of
  simplicial sets to $\oo$\nbd{}categories, as it is done for example in the
  articles
  \cite{ara2014vers,ara2018theoreme,gagna2018strict,ara2019quillen,ara2020theoreme,ara2020comparaison}.
  Following the terminology of these articles, a morphism $f : C \to D$ of
  $\oo\Cat$ is a \emph{Thomason equivalence} if $N_{\omega}(f)$ is a Kan--Quillen
  weak equivalence of simplicial sets. By definition, the nerve functor induces
  a functor at the level of homotopy categories
  \[
    \overline{N_{\omega}} : \ho(\oo\Cat^{\Th}) \to \ho(\Psh{\Delta}),
  \]
  where $\ho(\oo\Cat^{\Th})$ is the localization of $\oo\Cat$ with respect to
  the Thomason equivalences and $\ho(\Psh{\Delta})$ is the localization of
  $\Psh{\Delta}$ with respect to the Kan--Quillen weak equivalences of simplicial
  sets. As it so happens, the functor $\overline{N_{\omega}}$ is an equivalence
  of categories, as proved by Gagna in \cite{gagna2018strict}. In other words,
  the homotopy theory of $\oo$\nbd{}categories induced by Thomason equivalences
  is the same as the homotopy theory of spaces. Gagna's result is in fact a
  generalization of the analogous result for the usual nerve of small
  categories, which is attributed to Quillen in \cite{illusie1972complexe}. In
  the case of small categories, Thomason even showed the existence of a model
  structure whose weak equivalences are the ones induced by the nerve functor
  \cite{thomason1980cat}. The analogous result for $\oo\Cat$ is conjectured but
  not yet established \cite{ara2014vers}.
\end{named}
\begin{named}[Two homologies for $\oo$-categories]
  Keeping in mind the nerve functor of Street, a natural thing to do is to
  define the \emph{$k$-th homology group of an $\oo$\nbd{}category $C$} as the
  $k$\nbd{}th homology group of the nerve of $C$. In light of Gagna's result, these
  homology groups are just another way of looking at the homology groups of
  spaces. In order to explicitly avoid future confusion, we shall now use the
  name \emph{singular homology groups} of $C$ for these homology groups and the
  notation $H^{\sing}_k(C)$.

  On the other hand, Métayer gives a definition in \cite{metayer2003resolutions}
  of other homology groups for $\oo$\nbd{}categories. This definition is based
  on the notion of \emph{$\oo$\nbd{}categories free on a polygraph} (also known
  as \emph{$\oo$\nbd{}categories free on a computad}), which are
  $\oo$\nbd{}categories that are obtained from the empty category by recursively
  freely adjoining cells. From now on, we simply say \emph{free
    $\oo$\nbd{}category}. Métayer observed that every $\oo$\nbd{}category $C$ admits
  what we call a \emph{polygraphic resolution}, which means that there exists a
  free $\oo$\nbd{}category $P$ and a morphism of $\oo\Cat$
  \[
    f : P \to C
  \]
  that satisfies properties formally resembling those of trivial fibrations of
  topological spaces (or of simplicial sets). Furthermore, every free
  $\oo$\nbd{}category $P$ can be ``abelianized'' to a chain complex $\lambda(P)$ and
  Métayer proved that for two different polygraphic resolutions $P \to C$ and
  $P' \to C$ of the same
  $\oo$\nbd{}category $C$, the chain complexes
  $\lambda(P)$ and $\lambda(P')$ are quasi-isomorphic. Hence, we can define the
  \emph{$k$-th polygraphic homology group} of $C$, denoted by $H_k^{\pol}(C)$,
  as the $k$-th homology group of $\lambda(P)$ for any polygraphic resolution $P
  \to C$.

  One is then led to the following question:
  \begin{center}
    Do we have $H_{\bullet}^{\pol}(C) \simeq H_{\bullet}^{\sing}(C)$ for every
    $\oo$\nbd{}category $C$?
  \end{center}
  A first partial answer to this question is given by Lafont and Métayer in
  \cite{lafont2009polygraphic}: for a monoid $M$ (seen as category with one
  object and hence as
  an $\oo$\nbd{}category), we have $H_{\bullet}^{\pol}(M) \simeq
  H_{\bullet}^{\sing}(M)$. In fact, the original motivation for polygraphic
  homology was the homology of monoids and is part of a program that generalizes
  to higher dimension the results of Squier on the rewriting theory of monoids
  \cite{guiraud2006termination,lafont2007algebra,guiraud2009higher,guiraud2018polygraphs}. However, interestingly
  enough, the general answer to the above question is \emph{no}. A
  counterexample was found by Maltsiniotis and Ara. Let $B$ be the commutative
  monoid $(\mathbb{N},+)$, seen as a $2$\nbd{}category with only one $0$\nbd{}cell and no
  non-trivial $1$-cells. This $2$\nbd{}category is free (as an $\oo$\nbd{}category)
  and a quick computation shows that:
  \[
    H_k^{\pol}(B)=\begin{cases} \mathbb{Z} &\text{ if } k=0,2 \\ 0 &\text{
        otherwise. }\end{cases}
  \]
  On the other hand, it is shown in \cite[Theorem 4.9 and Example
  4.10]{ara2019quillen} that the nerve of $B$ is a $K(\mathbb{Z},2)$; hence, it
  has non-trivial homology groups in all even dimension.

  A question that still remains is:
  \begin{center}
    \textbf{(Q)} Which are the $\oo$\nbd{}categories $C$ such that
    $H_{\bullet}^{\pol}(C) \simeq H_{\bullet}^{\sing}(C)$ ?
  \end{center}
  This is precisely the question around which this dissertation revolves.
  Nevertheless, the reader will also find several new notions and results within
  this document that, although primarily motivated by the above question, are of
  interest in the theory of $\oo$\nbd{}categories and whose \emph{raisons
    d'être} go beyond the above considerations.
\end{named}
\begin{named}[Another formulation of the problem] One of the achievements of the
  present work is a more abstract reformulation of the question of comparison of
  singular and polygraphic homology of
  $\oo$\nbd{}categories. 

  In order to do so, recall first that by a variation of the Dold--Kan
  equivalence (see for example \cite{bourn1990another}), the category of abelian
  group objects in $\oo\Cat$ is equivalent to the category of non-negatively
  graded chain complexes
  \[
    \Ab(\oo\Cat) \simeq \Ch.
  \]
  Hence, we have a forgetful functor $\Ch \simeq \Ab(\oo\Cat) \to \oo\Cat$,
  which has a left adjoint
  \[
    \lambda : \oo\Cat \to \Ch.
  \]
  Moreover, for a \emph{free} $\oo$\nbd{}category $C$, the chain complex
  $\lambda(C)$ is exactly the one obtained by the ``abelianization'' process
  considered in Métayer's definition of polygraphic homology.

  Now, the category $\oo\Cat$ admits a model structure, known as the \emph{folk
    model structure} \cite{lafont2010folk}, whose weak equivalences are the
  \emph{equivalences of $\oo$\nbd{}categories} (a generalization of the usual
  notion of equivalence of categories) and whose cofibrant objects are exactly
  the free $\oo$\nbd{}categories \cite{metayer2008cofibrant}. Polygraphic
  resolutions are then nothing but cofibrant replacements in this model
  category. As the definition of polygraphic homology groups strongly suggests,
  the functor $\lambda$ is left Quillen with respect to this model structure. In
  particular, it admits a left derived functor
  \[
    \LL \lambda^{\folk} : \ho(\oo\Cat^{\folk}) \to \ho(\Ch)
  \]
  and we tautologically have that $H_k^{\pol}(C) = H_k(\LL \lambda^{\folk}(C))$
  for every $\oo$\nbd{}category $C$ and every $k \geq 0$. From now on, we set
  \[
    \sH^{\pol}(C):=\LL \lambda^{\folk}(C).
  \]
  This way of understanding polygraphic homology as a left derived functor has
  been around in the folklore for some time and I claim absolutely no
  originality for
  it. 

  On the other hand, $\lambda$ is also left derivable when $\oo\Cat$ is equipped
  with Thomason equivalences, yielding a left derived functor
  \[
    \LL \lambda^{\Th} : \ho(\oo\Cat^{\Th}) \to \ho(\Ch).
  \]
  This left derived functor being such that $ H_k^{\sing}(C) = H_k(\LL
  \lambda^{\Th}(C))$ for every $\oo$\nbd{}category $C$ and every $k \geq 0$.
  Contrary to the ``folk'' case, this result is new and first appears within
  this document (at least to my knowledge). Note that since, as mentioned
  earlier, the existence of a Thomason-like model structure on $\oo\Cat$ is
  still conjectural, usual tools from Quillen's theory of model categories were
  unavailable to prove the left derivability of $\lambda$ and the difficulty was
  to find a workaround solution.

  From now on, we set
  \[
    \sH^{\sing}(C):=\LL \lambda^{\Th}(C).
  \]

  Finally, it can be shown that every equivalence of $\oo$\nbd{}categories is a
  Thomason equivalence. Hence, the identity functor of $\oo\Cat$ induces a
  functor $\J$ at the level of homotopy categories
  \[
    \J : \ho(\oo\Cat^{\folk}) \to \ho(\oo\Cat^{\Th}),
  \]
  and altogether we have a triangle
  \[
    \begin{tikzcd}
      \ho(\oo\Cat^{\folk}) \ar[rd,"\LL \lambda^{\folk}"'] \ar[r,"\J"] & \ho(\oo\Cat^{\Th}) \ar[d,"\LL \lambda^{\Th}"] \\
      & \ho(\Ch).
    \end{tikzcd}
  \]
  This triangle is \emph{not} commutative (even up to isomorphism), since this
  would imply that the singular and polygraphic homology groups coincide for every
  $\oo$\nbd{}category. However, since both functors $\LL \lambda^{\folk}$ and
  $\LL \lambda^{\Th}$ are left derived functors of the same functor $\lambda$,
  the existence of a natural transformation $\pi : \LL \lambda^{\Th} \circ \J
  \Rightarrow \LL \lambda^{\folk}$ follows by universal property. Since $\J$ is
  the identity on objects, for every $\oo$\nbd{}category $C$, this natural
  transformation yields a map
  \[
    \pi_C : \sH^{\sing}(C) \to \sH^{\pol}(C),
  \]
  which we refer to as the \emph{canonical comparison map}. Let us say that $C$
  is \emph{homologically coherent} if $\pi_C$ is an isomorphism (which means
  exactly that for every $k\geq 0$, the induced map $H^{\sing}_k(C) \to
  H_k^{\pol}(C)$ is an isomorphism). The question of study then becomes:
  \begin{center}
    \textbf{(Q')} Which $\oo$\nbd{}categories are homologically coherent ?
  \end{center}
  Note that, in theory, question \textbf{(Q')} is more precise than question
  \textbf{(Q)} since we impose which morphism has to be an isomorphism in the
  comparison of homology groups. However, for all the concrete examples that we shall
  meet in practice, it is always question \textbf{(Q')} that will be answered.

  As will be explained in this thesis, a formal consequence of the above is that
  polygraphic homology is \emph{not} invariant under Thomason equivalence. This
  means that there exists at least one Thomason equivalence $f : C \to D$ such
  that the induced map
  \[
    \sH^{\pol}(C) \to \sH^{\pol}(D)
  \]
  is \emph{not} an isomorphism. 
  In other words, if we think of $\oo\Cat$ as a model of homotopy types (via the
  localization by the Thomason equivalences), then the polygraphic homology is
  \emph{not} a well-defined invariant. Another point of view would be to
  consider that the polygraphic homology is an intrinsic invariant of
  $\oo$\nbd{}categories (and not up to Thomason equivalence) and in that way is
  finer than singular homology. This is not the point of view adopted here, and
  the reason will be motivated at the end of this introduction. The slogan to
  retain is:
  \begin{center}
    Polygraphic homology is a way of computing singular homology groups of a
    homologically coherent $\oo$\nbd{}category.
  \end{center}
  The point is that given a \emph{free} $\oo$\nbd{}category $P$ (which is thus
  its own polygraphic resolution), the chain complex $\lambda(P)$ is much
  ``smaller'' than the chain complex associated to the nerve of $P$ and hence
  the polygraphic homology groups of $P$ are much easier to compute than its
  singular homology groups. The situation is comparable to using cellular
  homology for computing singular homology of a CW-complex. The difference is
  that in this last case, such a thing is always possible while in the case of
  $\oo$\nbd{}categories, one must ensure that the (free) $\oo$\nbd{}category is
  homologically
  coherent. 
\end{named}
\begin{named}[Finding homologically coherent $\oo$-categories]
  One of the main results presented in this dissertation is:
  \begin{center}
    Every (small) category $C$ is homologically coherent.
  \end{center}
  In order for this result to make sense, one has to consider categories as
  $\oo$\nbd{}categories with only unit cells above dimension $1$. Beware that
  this does not make the result trivial because given a polygraphic resolution
  $P \to C$ of a small category $C$, the $\oo$\nbd{}category $P$ need \emph{not}
  have only unit cells above dimension $1$.
  
  As such, this result is only a small generalization of Lafont and Métayer's
  result concerning monoids (although this new result, even restricted to
  monoids, is more precise because it means that the \emph{canonical comparison
    map} is an isomorphism). But the true novelty lies in the proof which is
  more conceptual that the one of Lafont and Métayer. It requires the
  development of several new concepts and results which in the end combine
  together smoothly to yield the desired result. This dissertation has been
  written so that all the elements needed to prove this result are spread over
  several chapters; a more condensed version of it is the object of the article
  \cite{guetta2020homology}. Among the new notions developed along the way, that
  of discrete Conduché $\oo$\nbd{}functor is probably the most significant. An
  $\oo$\nbd{}functor $f : C \to D$ is a \emph{discrete Conduché
    $\oo$\nbd{}functor} when for every cell $x$ of $C$, if $f(x)$ can be written
  as
  \[
    f(x)=y'\comp_k y'',
  \]
  then there exists a unique pair $(x',x'')$ of cells of $C$ that are
  $k$\nbd{}composable and such that
  \[
    f(x')=y',\, f(x'')=y'' \text{ and } x=x'\comp_k x''.
  \]
  The main result that we prove concerning discrete Conduché $\oo$\nbd{}functors
  is that for a discrete Conduché $\oo$\nbd{}functor $f : C \to D$, if the
  $\oo$\nbd{}category $D$ is free, then $C$ is also free. The proof of this
  result is long and tedious, though conceptually not extremely hard, and first
  appears in the paper \cite{guetta2020polygraphs}, which is dedicated to it.

  After having settled the case of ($1$\nbd{})categories, it is natural to move
  on to $2$\nbd{}categories. Contrary to the case of ($1$\nbd{})categories, not
  all $2$\nbd{}categories are \good{} and the situation seems to be much harder
  to understand. As a simplification, one can focus on $2$\nbd{}categories which
  are free (as $\oo$\nbd{}categories). This is what is done in this
  dissertation. With this extra hypothesis, the problem of characterization of
  \good{} free $2$\nbd{}categories may be reduced to the following question:
  given a cocartesian square of the form
  \[
    \begin{tikzcd}
      \sS_1 \ar[r] \ar[d] & P \ar[d]\\
      \sD_2 \ar[r] & P', \ar[from=1-1,to=2-2,"\ulcorner",phantom,very near end]
    \end{tikzcd}
  \]
  where $P$ is a free $2$\nbd{}category, when is it \emph{homotopy cocartesian}
  with respect to the Thomason equivalences? As a consequence, a substantial part of
  the work presented here consists in developing tools to detect homotopy
  cocartesian squares of $2$\nbd{}categories with respect to the Thomason
  equivalences. While it appears that these tools do not allow to completely
  answer the above question, they still make it possible to detect such homotopy
  cocartesian squares in many concrete situations. In fact, a whole section of
  the thesis is dedicated to giving examples of (free) $2$\nbd{}categories and
  computing the homotopy type of their nerve using these tools. Among all these
  examples, a particular class of well-behaved $2$\nbd{}categories, which I have
  coined ``bubble-free $2$\nbd{}categories'', seems to stand out. This class is
  easily characterized as follows. Given a $2$\nbd{}category, let us call
  \emph{bubble} a non-trivial $2$\nbd{}cell whose source and target are units on
  a $0$\nbd{}cell (necessarily the same). A \emph{bubble-free $2$\nbd{}category}
  is then nothing but a $2$\nbd{}category that has no bubbles. The archetypal
  example of a $2$\nbd{}category that is \emph{not} bubble-free is the
  $2$\nbd{}category $B$ introduced earlier (which is the commutative monoid
  $(\mathbb{N},+)$ seen as a $2$\nbd{}category). As already said, this
  $2$\nbd{}category is not \good{} and this does not seem to be a coincidence.
  It is indeed remarkable that of all the many examples of $2$\nbd{}categories
  studied in this work, the only ones that are not \good{} are exactly the ones
  that are \emph{not} bubble-free. This leads to the conjecture below, which
  stands as a conclusion of the thesis.
  \begin{center}
    \textbf{(Conjecture)} A free $2$\nbd{}category is \good{} if and only if it
    is bubble-free.
  \end{center}
\end{named}
\begin{named}[The big picture]
  Let us end this introduction with another point of view on the comparison of
  singular and polygraphic homologies. This point of view is highly conjectural
  and is not addressed at all in the rest of the dissertation. It should be
  thought of as a guideline for future work.

  In the same way that (strict) $2$\nbd{}categories are particular cases of
  bicategories, strict $\oo$\nbd{}categories are in fact particular cases of
  what are usually called \emph{weak $\oo$\nbd{}categories}. Such mathematical
  objects have been defined, for example, by Batanin using globular operads
  \cite{batanin1998monoidal} or by Maltsiniotis following ideas of Grothendieck
  \cite{maltsiniotis2010grothendieck}. Similarly to the fact that the theory of
  quasi-categories (which is a homotopical model for the theory of weak
  $\oo$\nbd{}categories whose cells are invertible above dimension $1$) may be
  expressed using the same language as the theory of usual categories, it is
  generally believed that all ``intrinsic'' notions (in a precise sense to be
  defined) of the theory of strict $\oo$\nbd{}categories have weak counterparts.
  For example, it is believed that there should be a folk model structure on the
  category of weak $\oo$\nbd{}categories and that there should be a good notion
  of free weak $\oo$\nbd{}category. In fact, this last notion should be defined
  as weak $\oo$\nbd{}categories that are recursively obtained from the empty
  $\oo$\nbd{}category by freely
  adjoining cells, which is the formal analogue of the strict version but in the
  weak context. The important point here is that a free strict
  $\oo$\nbd{}category is \emph{never} free as a weak $\oo$\nbd{}category (except
  for the empty $\oo$\nbd{}category).
  Moreover, there are good candidates for the polygraphic homology of weak
  $\oo$\nbd{}categories obtained by mimicking the definition in the strict case.
  But in general the polygraphic homology of a strict $\oo$\nbd{}category need
  not be the same as its ``weak polygraphic homology''. Indeed, since free
  strict $\oo$\nbd{}categories are not free as weak $\oo$\nbd{}categories,
  taking a ``weak polygraphic resolution'' of a strict $\oo$\nbd{}category is
  not the same as taking a polygraphic resolution. In fact, when trying to
  compute the weak polygraphic homology of $B$, it would seem that it gives the
  homology groups of a $K(\mathbb{Z},2)$, which is what we would have expected
  of its polygraphic homology in the first place. From this observation, it is
  tempting to make the following conjecture:
  \begin{center}
    The weak polygraphic homology of a strict $\oo$\nbd{}category coincides
    with its singular homology.
  \end{center}
  In other words, we conjecture that the fact that polygraphic and singular
  homologies of strict $\oo$\nbd{}categories do not coincide is a defect due to
  working in too narrow a setting. The ``good'' definition of polygraphic
  homology ought to be the weak one.

  We can go even further and conjecture the same thing for weak
  $\oo$\nbd{}categories. In order to do so, we need a definition of
  singular homology for weak $\oo$\nbd{}categories. This is
  conjecturally done as follows. To every weak $\oo$\nbd{}category
  $C$, one can associate a weak $\oo$\nbd{}groupoid $L(C)$ by formally
  inverting all the cells of $C$. Then, if we believe in
  Grothendieck's conjecture (see \cite{grothendieck1983pursuing} and
  \cite[Section 2]{maltsiniotis2010grothendieck}), the category of
  weak $\oo$\nbd{}groupoids equipped with the weak equivalences of weak
  $\oo$\nbd{}groupoids (see
  \cite[Paragraph 2.2]{maltsiniotis2010grothendieck}) is a model for the homotopy
  theory of spaces. In particular, every weak $\oo$\nbd{}groupoid has
  homology groups and we can define the singular homology groups of a
  weak $\oo$\nbd{}category $C$ as the homology groups of $L(C)$.

\end{named}
\begin{named}[Organization of the thesis]
  In the first chapter, we review some aspects of the theory of
  $\oo$\nbd{}categories. In particular, we study with great care free
  $\oo$\nbd{}categories, which are at the heart of the present work. It is the
  only chapter of the thesis that does not contain any reference to homotopy
  theory whatsoever. It is also there that we introduce the notion of discrete
  Conduché $\oo$\nbd{}functor and study their relation with free
  $\oo$\nbd{}categories. The culminating point of the chapter is Theorem
  \ref{thm:conduche}, which states that given a discrete Conduché
  $\oo$\nbd{}functor $F : C \to D$, if $D$ is free, then so is $C$. The proof of
  this theorem is long and technical and is broke down into several distinct
  parts.

  The second chapter is devoted to recalling some tools of homotopical algebra.
  More precisely, basic aspects of the theory of homotopy colimits using the
  formalism of Grothendieck's derivators are quickly presented. Note that this
  chapter does \emph{not} contain any original result and can be skipped at
  first reading. It is only intended to give the reader a summary of useful
  results on homotopy colimits that are used in the rest of the dissertation.

  In the third chapter, we delve into the homotopy theory of
  $\oo$\nbd{}categories. It is there that we define the different notions of
  weak equivalences for $\oo$\nbd{}categories and compare them. The two most
  significant new results to be found in this chapter are probably Proposition
  \ref{prop:folkisthom}, which states that every equivalence of
  $\oo$\nbd{}categories is a Thomason equivalence, and Theorem
  \ref{thm:folkthmA}, which states that equivalences of $\oo$\nbd{}categories
  satisfy a property reminiscent of Quillen's Theorem $A$ \cite[Theorem
  A]{quillen1973higher} and its $\oo$\nbd{}categorical generalization by Ara and
  Maltsiniotis \cite{ara2018theoreme,ara2020theoreme}.

  In the fourth chapter, we define the polygraphic and singular homologies of
  $\oo$\nbd{}categories and properly formulate the problem of their comparison.
  Up to Section \ref{section:polygraphichmlgy} included, all the results were
  known prior to this thesis (at least in the folklore), but starting from
  Section \ref{section:singhmlgyderived} all the results are original. Three
  fundamental results of this chapter are: Theorem \ref{thm:hmlgyderived},
  which states that singular homology is obtained as a derived functor of an
  abelianization function, Proposition \ref{prop:criteriongoodcat}, which gives
  an abstract criterion to detect \good{} $\oo$\nbd{}categories, and Proposition
  \ref{prop:comphmlgylowdimension}, which states that low-dimensional singular
  and polygraphic homology groups always coincide.

  The fifth chapter is mainly geared towards the fundamental Theorem
  \ref{thm:categoriesaregood}, which states that every category is \good{}. To
  prove this theorem, we first focus on a particular class of
  $\oo$\nbd{}categories, which we call \emph{contractible
    $\oo$\nbd{}categories}, and show that every contractible
  $\oo$\nbd{}category is \good{} (Proposition \ref{prop:contractibleisgood}).

  Finally, the sixth and last chapter of the thesis revolves around the homology
  of free $2$\nbd{}categories. The goal pursued is to try to understand which
  free $2$\nbd{}categories are \good{}. In order to do so, we give a criterion
  to detect homotopy cocartesian square with respect to Thomason equivalences
  (Proposition \ref{prop:critverthorThomhmtpysquare}) based on the homotopy
  theory of bisimplicial sets. Then, we apply this criterion and some other \emph{ad
  hoc} techniques to compute many examples of homotopy type of free
  $2$\nbd{}categories. The conclusion of the chapter is Conjecture
  \ref{conjecture:bubblefree}, which states that a free $2$\nbd{}category is
  \good{} if and only if it is bubble-free.
\end{named}

\fancyhf{}
\fancyfoot[C]{\thepage}
\fancyhead[RO,LE]{INTRODUCTION (FRANÇAIS)}
\chapter*{Introduction (français)}
\addcontentsline{toc}{chapter}{Introduction (français)} 
\selectlanguage{french}
Cette thèse a pour cadre général la \emph{théorie de l'homotopie des
$\oo$\nbd{}catégories strictes}, et, comme son titre le suggère, ce sont les
aspects homologiques de cette théorie qui sont traîtés. Le but est d'étudier et de
comparer deux invariants homologiques différents associés aux
$\oo$\nbd{}catégories strictes; c'est-à-dire, deux foncteurs différents
\[\mathbf{Str}\oo\Cat \to \ho(\Ch)\]
de la catégorie des $\oo$\nbd{}catégories strictes vers la catégorie homotopique
des complexes de chaînes en degré positif (i.e.\ la localisation de la catégorie
des complexes de chaînes en degré positif relativement aux
quasi\nbd{}isomorphismes).

Avant d'entrer dans le vif du sujet, précisons sans plus tarder qu'à
l'unique exception de la toute fin de cette introduction, toutes les
$\oo$\nbd{}catégories que nous considérerons sont des $\oo$\nbd{}catégories
\emph{strictes}. C'est pourquoi nous omettrons l'adjectif \og
strict\fg{} et parlerons simplement de \emph{$\oo$\nbd{}catégorie} plutôt
que de \emph{$\oo$\nbd{}catégorie stricte}. De même, nous noterons $\oo\Cat$ plutôt que
$\mathbf{Str}\oo\Cat$ la catégorie des $\oo$\nbd{}catégories (strictes).

\begin{named}[Le contexte : les $\oo$-categories en tant qu'espaces] L'étude de la théorie
  de l'homotopie des $\oo$\nbd{}catégories commence avec le foncteur nerf
  introduit par Street \cite{street1987algebra}
  \[
    N_{\omega} : \oo\Cat \to \Psh{\Delta}
  \]
  qui associe à toute $\oo$\nbd{}catégorie $C$ un ensemble simplicial
  $N_{\oo}(C)$, appelé le \emph{nerf de $C$}, généralisant le nerf usuel des
  (petites) catégories. En utilisant ce foncteur, il est possible de transférer
  la théorie de l'homotopie des ensembles simpliciaux aux $\oo$\nbd{}catégories,
  comme cela est fait par exemple dans les articles
  \cite{ara2014vers,ara2018theoreme,gagna2018strict,ara2019quillen,ara2020theoreme,ara2020comparaison}.
  Suivant la terminologie de ces derniers, un morphisme $f : C \to D$ de la catégorie
  $\oo\Cat$ est une \emph{équivalence de Thomason} si $N_{\oo}(f)$ est une
  équivalence faible de Kan--Quillen d'ensembles simpliciaux. Par définition, le
  foncteur nerf induit un foncteur au niveau des catégories homotopiques
    \[
    \overline{N_{\omega}} : \ho(\oo\Cat^{\Th}) \to \ho(\Psh{\Delta}),
  \]
  où $\ho(\oo\Cat^{\Th})$ est la localisation de $\oo\Cat$ relativement aux
  équivalences de Thomason et $\ho(\Psh{\Delta})$ est la localisation de
  $\Psh{\Delta}$ relativement aux équivalences faibles de Kan--Quillen d'ensembles
  simpliciaux. Comme l'a démontré Gagna \cite{gagna2018strict}, ce dernier foncteur est en
  fait une équivalence de catégories. Autrement dit, la théorie de l'homotopie
  des $\oo$\nbd{}catégories induite par les équivalences de Thomason est la
  même que la théorie de l'homotopie des espaces. Le résultat de Gagna est une
  généralisation du résultat analogue pour le nerf usuel des petites catégories,
  attribué à Quillen dans \cite{illusie1972complexe}. Dans le cas des petites
  catégories, Thomason a même démontré l'existence d'une structure de catégorie de modèles dont les équivalences
  faibles sont celles induites par le nerf \cite{thomason1980cat}. Le résultat analogue pour $\oo\Cat$ est toujours une
  conjecture \cite{ara2014vers}. 
\end{named}  

\begin{named}[Deux homologies pour les $\oo$-catégories]
  Armé du foncteur nerf de Street, il est naturel de définir le
  \emph{$k$\nbd{}ème groupe d'homologie d'une $\oo$\nbd{}catégorie $C$} comme
  étant le $k$\nbd{}ème groupe d'homologie du nerf de $C$. À la lumière du
  résultat de Gagna, ces groupes d'homologies sont simplement les groupes
  d'homologie des espaces vus sous un autre angle. Afin d'éviter de potentielles
  confusions à venir, nous appellerons désormais ces groupes d'homologie les
  \emph{groupes d'homologie singulière} de $C$ et nous utiliserons la notation
  $H^{\sing}_k(C)$.

  D'autres groupes d'homologie pour les $\oo$\nbd{}catégories ont aussi été
  définis par Métayer dans \cite{metayer2003resolutions}. La définition de
  ceux-ci repose sur la notion de \emph{$\oo$\nbd{}catégorie libre sur un
    polygraphe} (aussi connue sous le nom de \emph{computade}), c'est-à-dire
  de $\oo$\nbd{}catégorie obtenue de manière récursive à partir de la catégorie
  vide en attachant librement des cellules. Désormais, nous dirons simplement
  \emph{$\oo$\nbd{}catégorie libre}. Il a été observé par Métayer que toute
  $\oo$\nbd{}catégorie $C$ admet une \emph{résolution polygraphique},
  c'est-à-dire qu'il existe une $\oo$\nbd{}catégorie libre $P$ et un morphisme
  de $\oo\Cat$
  \[
    f : P \to C,
  \]
  satisfaisant des propriétés analogues à celles des fibrations triviales
  d'espaces topologiques (ou d'ensemples simpliciaux). De plus, toute
  $\oo$\nbd{}catégorie libre $P$ peut être \og abélianisée \fg{} en un complexe
  de chaînes $\lambda(P)$ et il a été démontré par Métayer que pour deux
  résolutions polygraphiques $P \to C$ et $P' \to C$ d'une même
  $\oo$\nbd{}catégorie, les complexes de chaînes $\lambda(P)$ et $\lambda(P')$
  sont quasi-isomorphes. Ainsi, on peut définir le \emph{$k$\nbd{}ème groupe
    d'homologie polygraphique} de $C$, noté $H_k^{\pol}(C)$, comme étant le
  $k$\nbd{}ème groupe d'homologie de $\lambda(P)$ pour n'importe quelle
  résolution polygraphique $P \to C$.

  Ces considérations invitent à se poser la question suivante :
  \begin{center}
    A-t-on $H_{\bullet}^{\pol}(C) \simeq H_{\bullet}^{\sing}(C)$ pour toute
    $\oo$\nbd{}catégorie $C$ ?
  \end{center}
  Une première réponse partielle a été donnée par Lafont et
  Métayer \cite{lafont2009polygraphic} : pour tout monoïde $M$ (vu comme une $\oo$\nbd{}catégorie
  à un seul objet et dont toutes les cellules de dimension supérieure à $1$ sont
  des unités), on a $H_{\bullet}^{\pol}(M) \simeq H_{\bullet}^{\sing}(M)$.
  Mentionnons au passage que l'homologie polygraphique a été conçue originellement pour étudier l'homologie des monoïdes et fait partie d'un programme
  dont le but est de généraliser en dimension supérieure les travaux de Squier
  sur la théorie de la réécriture des monoïdes \cite{guiraud2006termination,
    lafont2007algebra, guiraud2009higher, guiraud2018polygraphs}. Malgré le
  résultat de Lafont et Métayer, la réponse générale à la
  question précédente est \emph{non}. Un contre-exemple a été découvert par
  Maltsiniotis et Ara. Soit $B$ le monoïde commutatif $(\mathbb{N},+)$, vu
  comme une $2$\nbd{}catégorie avec une seule $0$\nbd{}cellule et pas de
  $1$\nbd{}cellule non-triviale. Cette $2$\nbd{}catégorie est libre (en tant
  que $\oo$\nbd{}catégorie) et un calcul rapide montre que:
  \[
    H_k^{\pol}(B)=\begin{cases} \mathbb{Z} &\text{ pour } k=0,2 \\ 0 &\text{
        sinon. }\end{cases}
  \]
  D'autre part, Ara a démontré \cite[Theorem 4.9 et Example
  4.10]{ara2019quillen} que le nerf de $B$ est un $K(\mathbb{Z},2)$, qui a donc
  des groupes d'homologie non triviaux en toute dimension paire.

  La question devient ainsi :
  \begin{center}
    \textbf{(Q)} Quelles sont les $\oo$\nbd{}catégories $C$ pour lesquelles
    $H_{\bullet}^{\pol}(C) \simeq H_{\bullet}^{\sing}(C)$?
  \end{center}
  C'est précisément à cette question que tente de répondre cette thèse.
  Néanmoins, le lecteur trouvera dans ce document plusieurs notions nouvelles et
  résultats qui, bien qu'originellement motivés par la question ci-dessus, sont
  intrinsèquement intéressants pour la théorie des $\oo$\nbd{}catégories et
  dont la portée dépasse les considérations homologiques précédentes.
\end{named}
\begin{named}[Une autre formulation du problème] Un des accomplissements
  du travail présenté ici est l'établissement d'un cadre conceptuel qui permet
  une reformulation plus abstraite et plus satisfaisante de la question
  de comparaison de l'homologie polygraphique et de l'homologie singulière des
  $\oo$\nbd{}catégories.

  Pour cela, rappelons tout d'abord qu'une variation de l'équivalence de
  Dold--Kan (voir par exemple \cite{bourn1990another}) permet d'affirmer que la
  catégorie des objets en groupes abéliens dans la catégorie $\oo\Cat$ est
  équivalente à la catégorie des complexes de chaînes en degré positif
  \[
    \Ab(\oo\Cat) \simeq \Ch.
  \]
  Ainsi, on a un foncteur d'oubli $\Ch \simeq \Ab(\oo\Cat) \to \oo\Cat$, qui
  admet un adjoint à gauche
  \[
    \lambda : \oo\Cat \to \Ch.
  \]
  En outre, pour une $\oo$\nbd{}catégorie \emph{libre} $C$, le complexe de
  chaînes $\lambda(C)$ est exactement celui obtenu par le processus
  d'\og abélianisation \fg{}
  que Métayer utilise dans sa définition d'homologie polygraphique.
  
  Par ailleurs, la catégorie $\oo\Cat$ peut être munie d'une structure de
  catégorie de modèles, communément appelée \emph{la structure de catégorie de modèles folk}
  \cite{lafont2010folk}, dont les équivalences faibles sont les
  \emph{équivalences de $\oo$\nbd{}catégories} (notion généralisant celle d'équivalence de
  catégories) et dont les objets cofibrants sont les $\oo$\nbd{}catégories
  libres \cite{metayer2008cofibrant}. Les résolutions polygraphiques ne sont alors
  rien d'autre que des remplacements cofibrants pour cette structure de catégorie de modèles.
  Comme la définition des groupes d'homologie polygraphique le laissait deviner,
  le foncteur $\lambda$ est Quillen à gauche relativement à cette structure de
  catégorie de modèles. En particulier, ce foncteur admet un foncteur dérivé à gauche
  \[
    \LL \lambda^{\folk} : \ho(\oo\Cat^{\folk}) \to \ho(\Ch)
  \]
  et on a tautologiquement $H_k^{\pol}(C) = H_k(\LL \lambda^{\folk}(C))$ pour
  toute $\oo$\nbd{}catégorie $C$ et pour tout $k\geq 0$. Désormais, on posera même
  \[
    \sH^{\pol}(C):=\LL \lambda^{\folk}(C).
  \]
  Je précise que cette façon de comprendre l'homologie polygraphique comme foncteur dérivé à
  gauche fait partie du folklore depuis un certain temps et je ne prétends à
  aucune originalité concernant ce point précis.

  D'autre part, le foncteur $\lambda$ est aussi dérivable à gauche quand
  $\oo\Cat$ est munie des équivalences de Thomason, ce qui permet d'obtenir un
  foncteur dérivé à gauche
    \[
    \LL \lambda^{\Th} : \ho(\oo\Cat^{\Th}) \to \ho(\Ch).
  \]
  En outre, ce foncteur est tel que  $ H_k^{\sing}(C) = H_k(\LL
  \lambda^{\Th}(C))$ pour toute $\oo$\nbd{}catégorie $C$ et pour tout $k \geq
  0$. Contrairement au cas \og folk \fg{}, ce résultat est
  complètement nouveau et apparaît pour la première fois dans ce manuscrit (à ma
  connaissance, du moins). Notons également que, puisque l'existence d'une
  structure de catégorie de modèles \og à la Thomason \fg{} sur $\oo\Cat$
  est toujours conjecturale, les outils habituels de la théorie de Quillen des
  catégories de modèles sont inutilisables pour démontrer que $\lambda$ est
  dérivable à gauche. La difficulté fut de trouver un moyen de
  contourner ce problème.

  Désormais, on posera
  \[
    \sH^{\sing}(C):=\LL \lambda^{\Th}(C).
  \]

  Finalement, on peut montrer que toute équivalence de $\oo$\nbd{}catégories
  est une équivalence de Thomason. Ainsi, le foncteur identité de $\oo\Cat$
  induit formellement un foncteur $\J$ au niveau des catégories homotopiques
  \[
    \J : \ho(\oo\Cat^{\folk}) \to \ho(\oo\Cat^{\Th}),
  \]
  et on obtient donc un triangle
  \[
    \begin{tikzcd}
      \ho(\oo\Cat^{\folk}) \ar[rd,"\LL \lambda^{\folk}"'] \ar[r,"\J"] & \ho(\oo\Cat^{\Th}) \ar[d,"\LL \lambda^{\Th}"] \\
      & \ho(\Ch).
    \end{tikzcd}
  \]
  Ce triangle n'est \emph{pas} commutatif (même à un isomorphisme près), car
  cela impliquerait que les groupes d'homologie singulière et les groupes
  d'homologie polygraphique coïncident pour toute $\oo$\nbd{}catégorie.
  Néanmoins, puisque les foncteurs $\LL \lambda^{\folk}$ et $\LL \lambda^{\Th}$
  sont tous les deux des foncteurs dérivés à gauche du même foncteur, l'existence
  d'une transformation naturelle $\pi : \LL \lambda^{\Th} \circ \J
  \Rightarrow \LL \lambda^{\folk}$ découle formellement par propriété
  universelle. De plus, $\J$ étant l'identité sur les objets, cette
  transformation naturelle fournit pour toute $\oo$\nbd{}catégorie $C$ un
  morphisme
    \[
    \pi_C : \sH^{\sing}(C) \to \sH^{\pol}(C),
  \]
  que nous appellerons le \emph{morphisme de comparaison canonique}. Une $\oo$\nbd{}catégorie $C$ est qualifiée d'\emph{homologiquement cohérente} si $\pi_C$ est
  un isomorphisme (ce qui signifie exactement que le
  morphisme induit $H^{\sing}_k(C) \to
  H_k^{\pol}(C)$ est un isomorphisme pour tout $k \geq 0$). La question devient alors : 
  \begin{center}
    \textbf{(Q')} Quelles $\oo$\nbd{}catégories sont homologiquement cohérentes ?
  \end{center}
  Notons au passage que la question \textbf{(Q')} est théoriquement plus précise
  que la question \textbf{(Q)}. Cependant, dans tous les exemples concrets
  que nous rencontrerons, c'est toujours à la question \textbf{(Q')} que nous répondrons.

  Comme il sera expliqué dans la thèse, de cette
  reformulation en terme de foncteurs dérivés, il est également possible de déduire formellement  que l'homologie polygraphique n'est \emph{pas} invariante
  relativement aux équivalences de Thomason. Cela signifie qu'il existe au moins
  une équivalence de Thomason $f : C \to D$ telle que le morphisme induit en
  homologie polygraphique
  \[
    \sH^{\pol}(C) \to \sH^{\pol}(D)
  \]
  n'est \emph{pas} un isomorphisme. En d'autres termes, si nous voyons les
  $\oo$\nbd{}catégories comme modèles des types d'homotopie (via la
  localisation de $\oo\Cat$ relativement aux équivalences de Thomason), alors
  l'homologie polygraphique n'est \emph{pas} un invariant bien défini. Un autre
  point de vue possible serait de considérer que l'homologie polygraphique est
  un invariant intrinsèque des $\oo$\nbd{}catégories (et non pas à équivalence
  de Thomason près) et, de cette façon, est un invariant plus fin que
  l'homologie singulière. Ce n'est pas le point de vue adopté dans cette thèse
  et ce choix sera motivé à la fin de cette introduction. Le
  slogan à retenir est :
  \begin{center}
    L'homologie polygraphique est un moyen de calculer les groupes d'homologie
    singulière des $\oo$\nbd{}catégories homologiquement cohérentes.
  \end{center}
  L'idée étant que pour une $\oo$\nbd{}catégorie libre $P$ (qui est donc sa propre résolution polygraphique), le complexe de chaînes
  $\lambda(P)$ est beaucoup moins \og gros \fg{} que le
  complexe de chaînes associé au nerf de $P$, et ainsi les groupes d'homologie
  polygraphique de $P$ sont beaucoup plus faciles à calculer que les groupes
  d'homologie singulière. La situation est comparable à l'utilisation de
  l'homologie cellulaire afin de calculer les groupes d'homologie singulière
  d'un CW\nbd{}complexe. La différence étant que dans ce dernier cas il est toujours
  possible de procéder ainsi, alors que dans le cas des $\oo$\nbd{}catégories,
  on doit d'abord s'assurer que la $\oo$\nbd{}catégorie (libre) en question est
  bien homologiquement cohérente.
\end{named}
\begin{named}[Détecter les $\oo$-catégories qui sont homologiquement cohérentes]
  Un des principaux résultats de cette thèse est le suivant :
  \begin{center}
    Toute (petite) catégorie est homologiquement cohérente.
  \end{center}
  Afin de donner du sens à cette assertion, il faut considérer les catégories
  comme des $\oo$\nbd{}catégories dont les cellules au delà de la dimension $1$
  sont des unités. Le résultat ci-dessus n'est pas pour autant trivial car pour
  une résolution polygraphique $P \to C$ d'une petite catégorie $C$, la
  $\oo$\nbd{}catégorie $P$, elle, n'a pas forcément que des cellules unités au
  delà de la dimension 1.

  En tant que tel, ce résultat est seulement une petite généralisation du
  résultat de Lafont et Métayer sur les monoïdes (bien qu'il soit plus précis,
  même restreint aux monoïdes, car il dit que c'est le \emph{morphisme de
    comparaison canonique} qui est un isomorphisme). Mais la véritable nouveauté
  du résultat en est sa démonstration qui est plus conceptuelle que celle de
  Lafont et Métayer. Elle repose sur l'introduction de nouvelles notions et le
  développement de nouveaux résultats; le tout s'assemblant élégamment pour
  finalement produire le résultat voulu. Cette thèse a été écrite de telle façon
  que tous les élements nécessaires à la démonstration du résultat précédent
  sont répartis sur plusieurs chapitres; une version plus condensée de celle-ci
  se trouve dans l'article \cite{guetta2020homology}. Parmi les nouvelles notions
  développées, la plus significative est probablement celle
  de $\oo$\nbd{}foncteur de Conduché discret. Un $\oo$\nbd{}foncteur $f : C \to
  D$ est un \emph{$\oo$\nbd{}foncteur de Conduché discret} quand pour toute
  cellule $x$ de $C$, si $f(x)$ peut être décomposé en
  \[
    f(x)=y'\comp_k y'',
  \]
  alors il existe une unique paire $(x',x'')$ de cellules de $C$ qui sont
  $k$\nbd{}composables et telles que
  \[
    f(x')=y',\, f(x'')=y'' \text{ and } x=x'\comp_k x''.
  \]
  Le résultat principal démontré concernant cette notion est que pour tout
  $\oo$\nbd{}foncteur de Conduché discret $f : C \to D$, si $D$ est libre alors
  $C$ est aussi libre. La démonstration est longue et fastidieuse, bien que
  relativement facile d'un point de vue conceptuel, et apparaît pour la première
  fois dans le papier \cite{guetta2020polygraphs}, qui lui est dédié.

  Une fois le cas de l'homologie des ($1$\nbd{})catégories complètement résolu,
  il est naturel de s'intéresser aux $2$\nbd{}catégories. Contrairement au cas
  des (1\nbd{})catégories, les 2\nbd{}catégories ne sont pas toutes
  homologiquement cohérentes et la situation est beaucoup plus compliquée. En
  premier lieu, on peut se restreindre aux $2$\nbd{}catégories qui sont libres
  (en tant que $\oo$\nbd{}catégories) et c'est ce qui est fait dans cette thèse.
  Avec cette hypothèse simplificatrice, le problème de caractérisation des
  $2$\nbd{}catégories libres homologiquement cohérentes peut être réduit à la
  question suivante : soit un carré cocartésien de la forme
  \[
    \begin{tikzcd}
      \sS_1 \ar[r] \ar[d] & P \ar[d]\\
      \sD_2 \ar[r] & P', \ar[from=1-1,to=2-2,"\ulcorner",phantom,very near end]
    \end{tikzcd}
  \]
  où $P$ est une $2$\nbd{}catégorie libre, quand est-il \emph{homotopiquement
    cocartésien} relativement aux équivalences de Thomason? En conséquence, une
  partie substantielle du travail présenté ici consiste à développer des outils
  permettant de reconnaître les carrés homotopiquement cocartésiens de
  $2$\nbd{}catégories relativement aux équivalences de Thomason. Bien que les
  outils qui seront présentés ne permettent pas de répondre entièrement à la
  question ci-dessus, ils permettent tout de même de détecter de tels carrés
  cocartésiens dans beaucoup de situations concrètes. Il y a même une section
  entière de la thèse qui consiste uniquement en une liste d'exemples détaillés
  de calculs du type d'homotopie de $2$\nbd{}catégories libres en utilisant ces
  outils. De ces exemples se dégage très clairement une classe particulière de
  $2$\nbd{}catégories, que j'ai nommées les \og $2$\nbd{}catégories sans bulles
  \fg{} et qui sont caractérisées comme suit. Pour une $2$\nbd{}catégorie,
  appelons \emph{bulle} une $2$\nbd{}cellule non-triviale dont la source et le
  but sont des unités sur une $0$\nbd{}cellule (nécessairement la même). Une
  \emph{$2$\nbd{}catégorie sans bulles} est simplement une $2$\nbd{}catégorie
  qui n'a aucune bulle. L'archétype de la $2$\nbd{}catégorie qui n'est
  \emph{pas} sans bulles est la $2$\nbd{}catégorie $B$ que nous avons déjà
  rencontrée plus haut (c'est-à-dire le monoïde commutatif $(\mathbb{N},+)$ vu
  comme une $2$\nbd{}catégorie). Comme dit précédemment, cette
  $2$\nbd{}catégorie n'est pas homologiquement cohérente et cela ne semble pas
  être une coïncidence. Il est tout à fait remarquable que de toutes les $2$\nbd{}catégories étudiées dans cette thèse, les seules qui ne
  sont pas homologiquement cohérentes sont exactement celles qui ne sont
  \emph{pas} sans bulles. Cela conduit à la conjecture ci-dessous, qui est le
  point d'orgue de la thèse.
  \begin{center}
    \textbf{(Conjecture)} Une $2$\nbd{}catégorie libre est homologie cohérente
    si et seulement si elle est sans bulles.
  \end{center}
\end{named}
  \begin{named}[Une vue d'ensemble]
    Terminons cette introduction par un autre point
    de vue sur la comparaison des homologies polygraphique et singulière.
    Précisons immédiatement que ce point de vue est hautement conjectural et
    n'est pas du tout abordé dans le reste de la thèse. Il s'agit plus d'un
    guide pour des travaux futurs qu'autre chose.

    De la même façon que les $2$\nbd{}catégories (strictes) sont des cas
    particuliers de bicatégories, les $\oo$\nbd{}catégories strictes sont en
    réalité des cas particuliers de ce qui est communément appelé des
    \emph{$\oo$\nbd{}catégories faibles}. Ces objets mathématiques ont été
    définis, par exemple, par Batanin en utilisant le formalisme des opérades
    globulaires \cite{batanin1998monoidal} ou par Maltsiniotis en suivant des
    idées de Grothendieck \cite{maltsiniotis2010grothendieck}. Tout comme la
    théorie des quasi-catégories (qui est un modèle homotopique pour la théorie
    des $\oo$\nbd{}catégories faibles dont les cellules sont toutes inversibles
    au delà de la dimension 1) s'exprime avec le même langage que la théorie des
    catégories usuelle, il est attendu que toutes les notions \og
    intrinsèques \fg{} (dans un sens précis à définir) de la théorie des
    $\oo$\nbd{}catégories strictes ont des analogues \emph{faibles}. Par
    exemple, il est attendu qu'il y ait une structure de catégorie de modèles folk sur la
    catégorie des $\oo$\nbd{}catégories faibles et qu'il y ait une bonne notion
    de $\oo$\nbd{}catégorie faible libre. En fait, ces dernières seraient
    certainement définies comme les $\oo$\nbd{}catégories faibles qui sont
    récursivement obtenues à partir de la $\oo$\nbd{}catégorie vide en attachant librement
    des cellules, ce qui est l'analogue formelle du cas strict. Le point clé ici
    est qu'une $\oo$\nbd{}catégorie stricte libre n'est \emph{jamais} libre en
    tant que $\oo$\nbd{}catégorie faible (excepté la $\oo$\nbd{}catégorie vide).
    Par ailleurs, il existe de bons candidats pour l'homologie polygraphique des
    $\oo$\nbd{}catégories faibles qui sont obtenus par mimétisme de la
    définition du cas strict. Mais il n'y aucune raison en général que
    l'homologie polygraphique d'une $\oo$\nbd{}catégorie stricte soit la même
    que son \og homologie polygraphique faible \fg{}. En effet, puisque les
    $\oo$\nbd{}catégories strictes libres ne sont pas libres en tant
    que $\oo$\nbd{}catégories faibles, prendre une \og résolution polygraphique
    faible \fg{} d'une $\oo$\nbd{}catégorie libre ne revient pas à
    prendre une résolution polygraphique. De fait, lorsqu'on essaye de calculer
    l'homologie polygraphique faible de $B$, il semblerait que cela donne les groupes
    d'homologie d'un $K(\mathbb{Z},2)$, ce qui aurait été attendu de l'homologie
    polygraphique au départ. De cette observation, il est tentant de faire la
    conjecture suivante :
    \begin{center}
      L'homologie polygraphique faible d'une $\oo$\nbd{}categorie stricte
        coïncide avec son homologie singulière. 
    \end{center}
    En d'autres termes, nous conjecturons que le fait que l'homologie
    polygraphique et l'homologie singulière d'une $\oo$\nbd{}catégorie stricte
    ne coïncident pas est un défaut dû à un cadre de travail trop étroit. La \og
    bonne \fg{} définition de l'homologie polygraphique devrait être la faible.

    Nous pourrions même aller plus loin et conjecturer la même chose pour les
    $\oo$\nbd{}catégories faibles. Pour cela, il est nécessaire de disposer
    d'une définition de l'homologie singulière des $\oo$\nbd{}catégories
    faibles. Conjecturellement, on procède de la manière suivante. À toute
    $\oo$\nbd{}catégorie faible $C$, on peut associer un $\oo$\nbd{}groupoïde
    faible $L(C)$ en inversant formellement toutes les cellules de $C$. Puis, si
    on en croit la conjecture de Grothendieck (voir
    \cite{grothendieck1983pursuing} et \cite[Section
    2]{maltsiniotis2010grothendieck}), la catégorie des $\oo$\nbd{}groupoïdes
    faibles munie des équivalences de $\oo$\nbd{}groupoïdes faibles (voir
    \cite[Paragraph 2.2]{maltsiniotis2010grothendieck}) est un modèle de la
    théorie homotopique des espaces. En particulier, chaque $\oo$\nbd{}groupoïde
    a des groupes d'homologie et on peut définir les groupes d'homologie
    singulière d'une $\oo$\nbd{}catégorie faible $C$ comme les groupes
    d'homologie de $L(C)$.
  \end{named}
  \begin{named}[Organisation de la thèse]
    Dans le premier chapitre, nous passerons en revue quelques aspects de la
    théorie des $\oo$\nbd{}catégories. En particulier, nous étudierons avec
    grand soin les $\oo$\nbd{}catégories libres, qui sont au c\oe{}ur de cette
    thèse. C'est le seul chapitre de la thèse qui ne contient aucune référence à
    la théorie de l'homotopie. C'est également dans ce chapitre que nous introduirons la notion
    de $\oo$\nbd{}foncteur de Conduché discret et que nous étudierons leur
    relation avec les $\oo$\nbd{}catégories libres. Le point
    culminant du chapitre étant le théorème \ref{thm:conduche}, qui dit que pour
    un $\oo$\nbd{}foncteur de Conduché discret $F : C \to D$, si $D$ est libre,
    alors $C$ l'est aussi. La démonstration de ce théorème est longue et
    technique et est décomposée en plusieurs parties distinctes.

    Le second chapitre a pour but de rappeler quelques outils d'algèbre
    homotopique. En particulier, les aspects élémentaires de la
    théorie des colimites homotopiques en utilisant le formalisme de Grothendieck des
    dérivateurs y sont rapidement présentés. Notons au passage que ce chapitre
    ne contient \emph{aucun} résultat original et peut être omis en première
    lecture. Son unique raison d'être est de donner au lecteur un catalogue de
    résultats concernant les colimites homotopiques qui seront utilisés par la
    suite.

    Dans le troisième chapitre, nous aborderons enfin la théorie de l'homotopie
    des $\oo$\nbd{}catégories. C'est là que seront définies et comparées entre elles les différentes
    notions d'équivalences faibles pour les $\oo$\nbd{}catégories. Les deux
    résultats les plus significatifs de ce chapitre sont probablement
    la proposition \ref{prop:folkisthom}, qui dit que toute équivalence de
    $\oo$\nbd{}catégorie est une équivalence de Thomason, et le théorème
    \ref{thm:folkthmA}, qui dit que les équivalences de $\oo$\nbd{}catégories
    satisfont une propriété réminiscente du théorème A de Quillen \cite[Theorem
    A]{quillen1973higher} et sa généralisation $\oo$\nbd{}catégorique par Ara
    et Maltsiniotis \cite{ara2018theoreme,ara2020theoreme}.

    Dans le quatrième chapitre, nous définirons les homologies polygraphique et
    singulière des $\oo$\nbd{}catégories et fomulerons précisément le problème
    de leur comparaison. Jusqu'à la section \ref{section:polygraphichmlgy}
    incluse, tous les résultats étaient connus avant cette thèse (au moins dans
    le folklore), mais à partir de la section \ref{section:singhmlgyderived}
    tous les résultats sont orignaux. Trois résultats fondamentaux de ce
    chapitre sont les suivants : le théorème \ref{thm:hmlgyderived}, qui dit que
    l'homologie singulière s'obtient comme le foncteur dérivé d'un foncteur
    d'abélianisation, la proposition \ref{prop:criteriongoodcat}, qui donne un
    critère abstrait pour détecter les $\oo$\nbd{}catégories homologiquement
    cohérentes, et la proposition \ref{prop:comphmlgylowdimension}, qui dit que les
    groupes d'homologies polygraphique et singulière coïncident toujours en basse
    dimension.

    Le cinquième chapitre a pour but de démontrer le théorème fondamental
    \ref{thm:categoriesaregood}, qui dit que toute catégorie est homologiquement
    cohérente. Pour cela, nous nous intéresserons en premier lieu
    à une classe particulière de $\oo$\nbd{}catégories, dites
    \emph{contractiles}, et nous montrerons que toute $\oo$\nbd{}catégorie
    contractile est homologiquement cohérente (Proposition
    \ref{prop:contractibleisgood}).

    Enfin, le sixième et dernier chapitre de la thèse s'intéresse à l'homologie
    des $2$\nbd{}catégories libres. Le but est d'essayer de comprendre quelles
    sont les $2$\nbd{}catégories libres qui sont homologiquement cohérentes.
    Pour cela, un critère pour détecter les carrés homotopiquement cocartésiens
    relativement aux équivalences de Thomason y est donné (Proposition
    \ref{prop:critverthorThomhmtpysquare}). Ce critère est fondé sur la théorie
    de l'homotopie des ensembles bisimpliciaux. Ensuite, nous appliquerons ce
    critère ainsi que d'autres techniques \emph{ad hoc} au calcul du type d'homotopie d'un grand
    nombre de $2$\nbd{}catégories libres. La conclusion du chapitre est la
    conjecture \ref{conjecture:bubblefree}, qui énonce qu'une $2$\nbd{}catégorie
    libre est homologiquement cohérente si et seulement si elle est sans bulles.
    \end{named}
\selectlanguage{english}
\frenchspacing

\fancyhf{}
\fancyfoot[C]{\thepage}
\fancyhead[LE]{\leftmark}
\fancyhead[RO]{\rightmark}
\chapter{Yoga of \texorpdfstring{$\oo$}{\textomega}-categories}
\section{\texorpdfstring{$\oo$}{ω}-graphs, \texorpdfstring{$\oo$}{ω}-magmas and \texorpdfstring{$\oo$}{ω}-categories}
\begin{paragr}\label{paragr:defoograph}
  An \emph{$\oo$-graph} $X$ consists of an infinite sequence of sets $(X_n)_{n \in \mathbb{N}}$  together with maps 
 \[ \begin{tikzcd}
    X_{n} &\ar[l,"\src",shift left] \ar[l,"\trgt"',shift right] X_{n+1}
  \end{tikzcd}
  \]
  for every $n \in \mathbb{N}$, subject to the \emph{globular identities}:
  \begin{equation*}
  \left\{
  \begin{aligned}
    \src \circ \src &= \src \circ \trgt, \\
    \trgt \circ \trgt &= \trgt \circ \src.
  \end{aligned}
  \right.
  \end{equation*}
  Elements of $X_n$ are called \emph{$n$\nbd{}cells} or \emph{$n$\nbd{}arrows} or \emph{cells of dimension $n$}. For $n=0$, elements of $X_0$ are also called \emph{objects}. For $x$ an $n$\nbd{}cell with $n>0$, $\src(x)$ is the \emph{source} of $x$ and $\trgt(x)$ is the \emph{target} of $x$. We use the notation \[x : a \to b\] to say that $a$ is the source of $x$ and $b$ is the target of $x$.

  More generally, for $0\leq k <n$, we define maps $\src_k : X_n \to X_k$ and $\trgt_k : X_n \to X_k$ as
  \[
  \src_k = \underbrace{\src\circ \dots \circ \src}_{n-k \text{ times}}
  \]
  and
  \[
  \trgt_k = \underbrace{\trgt\circ \dots \circ \trgt}_{n-k \text{ times}}.
  \]
  For an $n$\nbd{}cell $x$, the $k$\nbd{}cells $\src_k(x)$ and $\trgt_k(x)$ are respectively the \emph{$k$\nbd{}source} and the \emph{$k$\nbd{}target} of $x$.

  Two $n$\nbd{}cells $x$ and $y$ are \emph{parallel} if
  \[
  n=0
  \]
  or
  \[
  n>0 \text{ and }\src(x)=\src(y) \text{ and } \trgt(x)=\trgt(y).
  \]

  Let $0 \leq k <n$. Two $n$\nbd{}cells $x$ and $y$ are \emph{$k$\nbd{}composable} if
  \[
  \src_k(x)=\trgt_k(y).
  \]
  Note that the expression \og $x$ and $y$ are $k$\nbd{}composable\fg{} is \emph{not} symmetric in $x$ and $y$ and we \emph{should} rather speak of a ``$k$\nbd{}composable pair $(x,y)$'', although we won't always do it. The set of pairs of $k$\nbd{}composable $n$\nbd{}cells is denoted by $X_n\underset{X_k}{\times}X_n$, and is characterized as the following fibred product
      \[
    \begin{tikzcd}
      X_n\underset{X_k}{\times}X_n \ar[r] \ar[dr,phantom,"\lrcorner", very near start] \ar[d] &X_n \ar[d,"\trgt_k"]\\
      X_n \ar[r,"\src_k"] & X_k.
      \end{tikzcd}
    \]
    A \emph{morphism of $\oo$-graphs} $f : X \to Y$ is a sequence $(f_n : X_n \to Y_n)_{n \in \mathbb{N}}$ of maps that is compatible with source and target, i.e.\ for every $n$\nbd{}cell $x$ of $X$ with $n >0$, we have
    \[
    f_{n-1}(\src(x))=\src(f_n(x)) \text{ and } f_{n-1}(\trgt(x))=\trgt(f_n(x)).
    \]
    For an $n$\nbd{}cell $x$ of $X$, we often write $f(x)$ instead of $f_n(x)$.

    The category of $\oo$-graphs and morphisms of $\oo$-graphs is denoted by $\oo\Grph$.
\end{paragr}
\begin{paragr}\label{paragr:defoomagma}
  An \emph{$\oo$-magma} consists of an $\oo$-graph $X$ together with maps
\begin{align*}
     1_{(\shortminus)}: X_n &\to X_{n+1}\\
    x &\mapsto 1_x
  \end{align*}
  for every $n\geq 0$, and maps
  \begin{align*}
   (\shortminus)\comp_k(\shortminus): X_n \underset{X_k}{\times}X_n &\to X_n\\
    (x,y)&\mapsto x\comp_k y
  \end{align*}
  for all $0 \leq k <n$, subject to the following axioms:
\begin{enumerate}[label=(\alph*)]
    \item For every $n\geq 0$ and every $n$\nbd{}cell $x$,
    \[
    \src(1_x)=\src(1_x)=x.
    \]
  \item For all  $0 \leq k< n$ and all $k$\nbd{}composable $n$\nbd{}cells $x$ and $y$,
    \[
    \src(x\comp_k y) =
    \begin{cases}
      \src(y) &\text{ when }k=n-1,\\
      \src(x)\comp_k \src(y) &\text{ otherwise,}
      \end{cases}
    \]
    and
        \[
    \trgt(x\comp_k y) =
    \begin{cases}
      \trgt(x) & \text{ when }k=n-1,\\
      \trgt(x)\comp_k \trgt(y) &\text{ otherwise.}
      \end{cases}
    \]

  \end{enumerate}
We will use the same letter to denote an $\oo$\nbd{}magma and its underlying $\oo$\nbd{}graph.

For an $n$\nbd{}cell $x$, the $(n+1)$\nbd{}cell $1_{x}$ is referred to as the \emph{unit on $x$}. More generally, for all $0 \leq k < n$, we define maps $\1^n_{(\shortminus)} : C_k \to C_n$ as
\[
\1^{n}_{(\shortminus)} := \underbrace{1_{(\shortminus)} \circ \dots \circ 1_{(\shortminus)}}_{n-k \text{ times }} : C_k \to C_n.
\]
For a $k$\nbd{}cell $x$ and $n>k$, the $n$\nbd{}cell $\1^n_x$ is referred to as
the \emph{$n$\nbd{}dimensional unit on $x$}. For consistency, we also set
\[\1^k_x:=x\]
for every $k$\nbd{}cell $x$. An $n$\nbd{}cell that is a unit on a strictly
lower dimensional cell is sometimes referred to as a \emph{trivial $n$\nbd{}cell}.

For two $k$\nbd{}composable $n$\nbd{}cells $x$ and $y$, the $n$\nbd{}cell $x \comp_k y$ is referred to as the \emph{$k$\nbd{}composition} of $x$ and $y$.

More generally, we extend the notion of $k$\nbd{}composition for cells of different dimension in the following way. Let $x$ be an $n$\nbd{}cell, $y$ be an $m$-cell with $m \neq n$ and $k~<~\min\{m,n\}$. The cells $x$ and $y$ are \emph{$k$\nbd{}composable} if $\src_k(x)=\trgt_k(y)$, in which case we define the cell $x\comp_k y$ of dimension $\max\{m,n\}$ as
\[
x\comp_k y :=
\begin{cases}
  1^n_x \comp_k y &\text{ if } m<n\\
  x \comp_k 1^m_y  &\text{ if } m>n\\
  \end{cases}
\]
We also follow the convention that if $n<m$, then $\comp_n$ has priority over $\comp_m$. This means that 
\[
x\comp_n y \comp_m z = (x \comp_n y) \comp_m z \text{ and } x \comp_m y \comp_n z = x \comp_m (y \comp_n z)
\]
whenever these equations make sense.

A \emph{morphism of $\oo$\nbd{}magmas} $f : X \to Y$ is a morphism of underlying $\oo$\nbd{}graphs that is compatible with units and compositions, i.e.\ for every $n$\nbd{}cell $x$, we have 
    \[
    f(1_x)=1_{f(x)},
    \]
    and for every $k$\nbd{}composable $n$\nbd{}cells $x$ and $y$, we have
    \[
    f(x\comp_k y)= f(x) \comp_k f(y).
    \]
    We write $\oo\Mag$ for the category of $\oo$\nbd{}magmas and morphisms of $\oo$\nbd{}magmas.
\end{paragr}
\begin{paragr}\label{paragr:defoocat}
  An \emph{$\oo$-category} is an $\oo$-magma $X$ that satisfies the following axioms:
  \begin{description}
      \item[Units:] for all $k<n$, for every $n$\nbd{}cell $x$, we have
    \[
    \1^n_{\trgt_k(x)}\comp_k x =x= x \comp_k\1^n_{\src_k(x)},
    \]
    \item[Functoriality of units:] for all $k<n$ and for all $k$\nbd{}composable $n$\nbd{}cells $x$ and $y$, we have 
    \[
    1_{x\comp_k y}=1_{x}\comp_k 1_{y},
    \]

  \item[Associativity:] for all  $k<n$, for all $n$\nbd{}cells $x, y$ and $z$ such that $x$ and $y$ are $k$\nbd{}composable, and $y$ and $z$ are $k$\nbd{}composable, we have
    \[
    (x\comp_{k}y)\comp_{k}z=x\comp_k(y\comp_kz),
    \]
    \item[Exchange rule:] for all $k,l,n \in \mathbb{N}$ with $k<l<n$, for all $n$\nbd{}cells $x,x',y$ and $y'$  such that
    \begin{itemize}
    \item[-] $x$ and $y$ are $l$-composable, $x'$ and $y'$ are $l$-composable,
    \item[-] $x$ and $x'$ are $k$\nbd{}composable, $y$ and $y'$ are $k$\nbd{}composable,
    \end{itemize}
    we have
    \[
    (x \comp_k x')\comp_l (y \comp_k y')=(x \comp_l y)\comp_k (x' \comp_l y').
    \]
  \end{description}
  We will use the same letter to denote an $\oo$\nbd{}category and its underlying $\oo$\nbd{}magma. A \emph{morphism of $\oo$\nbd{}categories} (or \emph{$\oo$\nbd{}functor}) $f : X \to Y$ is simply a morphism of the underlying $\oo$\nbd{}magmas. We denote by $\oo\Cat$ the category of $\oo$\nbd{}categories and morphisms of $\oo$\nbd{}categories. This category is locally presentable.
\end{paragr}
\begin{paragr}\label{paragr:defncat}
  For $n \in \mathbb{N}$, the notions of \emph{$n$\nbd{}graph}, \emph{$n$\nbd{}magma} and
  \emph{$n$\nbd{}category} are defined as truncated version of $\oo$\nbd{}graph,
  $\oo$\nbd{}magma and $\oo$\nbd{}category in an obvious way. For example, a
  $0$-category is a set and a $1$-category is nothing but a usual (small)
  category. The category of $n$\nbd{}categories and morphisms of $n$\nbd{}categories (or
  $n$\nbd{}functors) is denoted by $n\Cat$. When $n=0$ and $n=1$, we almost always use the notation $\Set$ and $\Cat$ instead of $0\Cat$ and $1\Cat$.

  For every $n\geq 0$, there is a canonical functor
  \[
  \tau_{\leq n}^s : \oo\Cat \to n\Cat
  \]
  that simply discards all the cells of dimension strictly higher than $n$. This functor has a left adjoint
  \[
  \iota_n : n\Cat \to \oo\Cat,
  \]
  where for an $n$\nbd{}category $C$, the $\oo$\nbd{}category $\iota_n(C)$ has the same $k$\nbd{}cells as $C$ for $k\leq n$ and only unit cells in dimension strictly higher than $n$. This functor itself has a left adjoint
  \[
  \tau_{\leq n }^i : \oo\Cat \to n\Cat,
  \]
  where for an $\oo$\nbd{}category $C$, the $n$\nbd{}category $\tau_{\leq n}^i(C)$ has the same $k$\nbd{}cells as $C$ for $k < n$ and whose set of $n$\nbd{}cells is the quotient of $C_n$ under the equivalence relation $\sim$ generated by
  \[
  x \sim y \text{ if there exists } z \in C_{n+1} \text{ of the form } z : x \to y.
  \]
  The functor $\tau_{\leq n}^s$ also have a right adjoint
  \[
  \kappa_n : n\Cat \to \oo\Cat,
  \]
  where for an $n$\nbd{}category $C$, the $\oo$\nbd{}category $\kappa_n(C)$ has
  the same $k$\nbd{}cells as $C$ for $k \leq n$ and has exactly one
  $k$\nbd{}cell $x \to y$ for every pair of parallel $(k-1)$\nbd{}cells $(x,y)$
  for $k>n$.
  
  The sequence of adjunctions
  \[
   \tau^{i}_{\leq n} \dashv \iota_n \dashv \tau^s_{\leq n} \dashv \kappa_n
  \]
  is maximal in that $\kappa_n$ doesn't have a right adjoint and $\tau^{i}_{\leq
    n}$ doesn't have a left adjoint.

  The functors $\tau^{s}_{\leq n}$ and $\tau^{i}_{\leq n}$ are respectively referred to as the \emph{stupid truncation functor} and the \emph{intelligent truncation functor}.

  The functor $\iota_n$ is fully faithful and preserves both limits and
  colimits; in regards to these properties, we often identify $n\Cat$ with the
  essential image of $\iota_n$, which is the full subcategory of $\oo\Cat$
  spanned by the $\oo$\nbd{}categories whose $k$\nbd{}cells for $k >n$ are all units.
\end{paragr}
\begin{paragr}
  For $n \geq 0$, we define the $n$\nbd{}skeleton functor $\sk_n : \oo\Cat \to \oo\Cat$ as
  \[
  \sk_n := \iota_n \circ \tau^{s}_{\leq n}.
  \]
  This functor preserves both limits and colimits. For an $\oo$\nbd{}category
  $C$, $\sk_n(C)$ is the sub-$\oo$\nbd{}category of $C$ generated by the
  $k$\nbd{}cells of $C$ with $k\leq n$ in an obvious sense. It is also
  convenient to define $\sk_{-1}(C)$ to be the empty $\oo$\nbd{}category
  \[
    \sk_{-1}(C)=\emptyset
  \]
  for every $\oo$\nbd{}category $C$. Note that the functor $\sk_{-1} : \oo\Cat
  \to \oo\Cat$ preserves colimits but does not preserve limits.
  The inclusion induces a canonical filtration
  \[
  \emptyset=\sk_{-1}(C) \hookrightarrow \sk_{0}(C) \hookrightarrow \sk_{1}(C) \hookrightarrow \cdots \hookrightarrow\sk_{n}(C) \hookrightarrow\cdots,
  \]
  and we leave the proof of the following lemma as an easy exercise for the reader.
\end{paragr}
\begin{lemma}\label{lemma:filtration}
  Let $C$ be an $\oo$\nbd{}category. The colimit of the canonical filtration
    \[
  \sk_{-1}(C) \hookrightarrow \sk_{0}(C) \hookrightarrow \sk_{1}(C) \hookrightarrow \cdots \hookrightarrow\sk_{n}(C) \hookrightarrow\cdots
  \]
  is $C$ and for $n \geq 0$ the universal arrow $\sk_{n}(C) \to C$ is given by the co-unit of the adjunction $\tau^{s}_{\leq n} \dashv \iota_n$.
\end{lemma}
\begin{paragr}\label{paragr:defglobe}
  For $n \in \mathbb{N}$, the \emph{$n$\nbd{}globe} $\sD_n$ is the $n$\nbd{}category that has:
  \begin{itemize}[label=-]
  \item exactly one non-trivial $n$\nbd{}cell, which we refer to as the \emph{principal $n$\nbd{}cell} of $\sD_n$, and which we denote by $e_n$,
  \item exactly two non-trivial $k$\nbd{}cells for every $k<n$; these $k$\nbd{}cells being parallel and given by the $k$\nbd{}source and the $k$\nbd{}target of $e_n$.
  \end{itemize}
  This completely describes the $n$\nbd{}category $\sD_n$ as no non-trivial composition can occur. Here are pictures in low dimension:
      \[
  \sD_0= \begin{tikzcd}\bullet,\end{tikzcd}
  \]
  \[
  \sD_1 = \begin{tikzcd} \bullet \ar[r] &\bullet, \end{tikzcd}
  \]
  \[
  \sD_2 = \begin{tikzcd}
    \bullet \ar[r,bend left=50,""{name = U,below}] \ar[r,bend right=50,""{name=D}]&\bullet, \ar[Rightarrow, from=U,to=D]
  \end{tikzcd}
  \]
    \[
    \sD_3 = \begin{tikzcd}
         \bullet \ar[r,bend left=50,""{name = U,below,near start},""{name = V,below,near end}] \ar[r,bend right=50,""{name=D,near start},""{name = E,near end}]&\bullet. \ar[Rightarrow, from=U,to=D, bend right,""{name= L,above}]\ar[Rightarrow, from=V,to=E, bend left,""{name= R,above}]
    \arrow[phantom,"\Rrightarrow",from=L,to=R]
  \end{tikzcd}
    \]
    For every $\oo$\nbd{}category $C$, the map
    \begin{align*}
      \Hom_{\oo\Cat}(\sD_n,C) &\to C_n \\
      F &\mapsto F(e_n)
    \end{align*}
    is a bijection natural in $C$. In other words, the $n$\nbd{}globe represents the functor
    \begin{align*}
      \oo\Cat &\to \Set\\
      C &\mapsto C_n.
    \end{align*}
    For an $n$\nbd{}cell $x$ of $C$, we denote by
    \[
    \langle x \rangle : \sD_n \to C
    \]
    the canonically associated $\oo$\nbd{}functor.
\end{paragr}
\begin{paragr}\label{paragr:inclusionsphereglobe}
  For $n \in \mathbb{N}$, the $n$\nbd{}sphere $\sS_n$ is the $n$\nbd{}category that has
  exactly two parallel non-trivial $k$\nbd{}cells for every $k\leq n$. In other words, we have
  \[
  \sS_{n}=\sk_n(\sD_{n+1}),
  \]
  and in particular, we have a canonical inclusion functor
  \[
  i_{n+1} : \sS_{n} \to \sD_{n+1}.
  \]
  It is also customary to define $\sS_{-1}$ to be the empty $\oo$\nbd{}category and $i_{-1}$ to be the unique $\oo$\nbd{}functor
  \[
  \emptyset \to \sD_0.
  \]
  Notice that for every $n\geq 0$, the following commutative square
  \[
  \begin{tikzcd}
    \sS_{n-1} \ar[r,"i_n"] \ar[d,"i_n"] & \sD_n \ar[d,"j_n^+"]\\
    \sD_n \ar[r,"j_n^-"'] & \sS_{n},
  \end{tikzcd}
  \]
  where we wrote $j_n^+$ (resp.\ $j_n^-$) for the morphism $\langle \trgt(e_{n+1}) \rangle : \sD_n \to \sS_n$ (resp.\ $\langle \src(e_{n+1}) \rangle : \sD_n \to \sS_n$), is cocartesian.

  Here are some pictures of the $n$\nbd{}spheres in low dimension:
         \[
  \sS_0= \begin{tikzcd}\bullet & \bullet, \end{tikzcd}
  \]
  \[
  \sS_1 = \begin{tikzcd}     \bullet \ar[r,bend left=50,""{name = U,below}] \ar[r,bend right=50,""{name=D}]&\bullet,  \end{tikzcd}
  \]
  \[
  \sS_2 = \begin{tikzcd}
         \bullet \ar[r,bend left=50,""{name = U,below,near start},""{name = V,below,near end}] \ar[r,bend right=50,""{name=D,near start},""{name = E,near end}]&\bullet. \ar[Rightarrow, from=U,to=D, bend right,""{name= L,above}]\ar[Rightarrow, from=V,to=E, bend left,""{name= R,above}]
    \end{tikzcd}
  \]
  
  For an $\oo$\nbd{}category $C$ and $n\geq 0$, an $\oo$\nbd{}functor
  \[
  \sS_n \to C
  \]
  amounts to the data of two parallel $n$\nbd{}cells of $C$. In other words, $\sS_{n}$ represents the functor $\oo\Cat \to \Set$ that sends an $\oo$\nbd{}category to the set of its parallel $n$\nbd{}cells. For $(x,y)$ a pair of parallel $n$\nbd{}cells of $C$, we denote by
  \[
  \langle x,y \rangle : \sS_{n} \to C
  \]
  the canonically associated $\oo$\nbd{}functor. For example, the $\oo$\nbd{}functor $i_n$ is nothing but
  \[
  \langle \src(e_{n+1}),\trgt(e_{n+1}) \rangle : \sS_{n} \to \sD_{n+1}.
  \]
  \end{paragr}
\section{Free \texorpdfstring{$\oo$}{ω}-categories}
\begin{definition}\label{def:nbasis}
  Let $C$ be an $\oo$\nbd{}category and $n \geq 0$. A subset $E \subseteq C_n$ of the $n$\nbd{}cells of $C$ is an \emph{$n$\nbd{}basis of $C$} if the commutative square
  \[
  \begin{tikzcd}[column sep=huge, row sep=huge]
   \displaystyle \coprod_{ x \in E}\ar[d,"\displaystyle\coprod_{x \in E} i_n"'] \sS_{n-1} \ar[r,"{\langle \src(x),\trgt(x)\rangle_{x \in E}}"] & \sk_{n-1}(C) \ar[d,hook] \\
   \displaystyle \coprod_{x \in E} \sD_n \ar[r,"\langle x \rangle_{x \in E}"] & \sk_{n}(C) 
  \end{tikzcd}
  \]
  is cocartesian.
\end{definition}
\begin{remark}\label{remark:nbasisncat}
  Note that since for all $n<m$, we have $\sk_n \circ \sk_m = \sk_n$, an $\oo$\nbd{}category $C$ has an $n$\nbd{}basis if and only if $\sk_n(C)$ has an $n$\nbd{}basis. 
\end{remark}
\begin{paragr}\label{paragr:defnbasisdetailed}
  Let us unfold Definition \ref{def:nbasis}. For $n=0$, $E$ is an $0$\nbd{}basis
  of $C$ if $E=C_0$. For $n>0$,
  $E$ is $n$\nbd{}basis of $C$ if for every $n$\nbd{}category $D$, for every $(n-1)$-functor
  \[
  F : \tau_{\leq n-1}^{s}(C) \to \tau_{\leq n-1}^{s}(D),
  \]
  and for every map
  \[
  f : E \to D_n
  \]
  such that for every $x \in E$,
  \[
  \src(f(x))=F(\src(x)) \text{ and } \trgt(f(x))=F(\trgt(x)),
  \]
 there exists a \emph{unique} $n$\nbd{}functor
  \[
  \tilde{F} : \tau^{s}_{\leq n}(C) \to D
  \]
  such that $\tilde{F}_k = F_k$ for every $k<n$ and $\tilde{F}_n(x) = f(x)$ for every $x \in E$.

  Intuitively speaking, this means that $\sk_{n}(C)$ has been obtained from $\sk_{n-1}(C)$ by freely adjoining the cells in $E$.
\end{paragr}

\begin{example}\label{dummyexample}
  An $n$\nbd{}category (seen as an $\oo$\nbd{}category) always has a $k$\nbd{}basis for every $k>n$, namely the empty set.
\end{example}
 Less trivial examples will come along soon.
\begin{definition}\label{def:freeoocat}
  An $\oo$\nbd{}category is \emph{free}\footnote{Other common terminology for ``free $\oo$\nbd{}category'' is ``$\oo$\nbd{}category free on a polygraph'' \cite{burroni1993higher} or ``$\oo$\nbd{}category free on a computad'' \cite{street1976limits,makkai2005word}.} if it has $n$\nbd{}basis for every $n \geq 0$.
\end{definition}
\begin{paragr}\label{paragr:freencat}
  By considering $n\Cat$ as a subcategory of $\oo\Cat$, the previous definition also works for $n$\nbd{}categories. It follows from Example \ref{dummyexample} that an $n$\nbd{}category is free if and only if it has a $k$\nbd{}basis for every $0 \leq k \leq n$.
\end{paragr}
We now wish to recall an important result due to Makkai concerning the
uniqueness of the $n$\nbd{}basis for a free $\oo$\nbd{}category. First we need the following definition.
\begin{definition}\label{def:indecomposable}
  Let $C$ be an $\oo$\nbd{}category. For $n>0$, an $n$\nbd{}cell $x$ of $C$ is \emph{indecomposable} if both following conditions are satisfied:
  \begin{enumerate}[label=(\alph*)]
  \item $x$ is not a unit on a lower dimensional cell,
  \item if $x$ is of the form
    \[
    x=a\comp_k b
    \]
    with $k<n$, then either
    \[
    a=\1^n_{\trgt_k(x)},
    \]
    or
    \[
    b=\1^n_{\src_k(x)}.
    \]
  \end{enumerate}
  For $n=0$, all $0$-cells are, by convention, indecomposable.
\end{definition}
We can now state the promised result, whose proof can be found in \cite[Section 4, Proposition 8.3]{makkai2005word}.
\begin{proposition}[Makkai]\label{prop:uniquebasis}
  Let $C$ be a free $\oo$\nbd{}category. For every $n \in \mathbb{N}$, $C$ has a \emph{unique} $n$\nbd{}basis. The cells of this $n$\nbd{}basis are exactly the indecomposable $n$\nbd{}cells of $C$.
\end{proposition}
\begin{remark}
  Note that there is a subtlety in the previous proposition. It is not true in
  general that if an $\oo$\nbd{}category $C$ has an $n$\nbd{}basis then it is
  unique. The point is that we need the existence of the $k$\nbd{}bases for
  $k<n$ in order to prove the uniqueness of the $n$\nbd{}basis. (See the paper of Makkai cited previously for details.)
\end{remark}
\begin{paragr}
  Proposition \ref{prop:uniquebasis} allows us to speak of \emph{the} $k$\nbd{}basis of a free $\oo$\nbd{}category $C$ and more generally of the \emph{basis} of $C$ for the sequence
  \[
  (\Sigma_k)_{k \in \mathbb{N}}
  \]
  where each $\Sigma_k$ is the $k$\nbd{}basis of $C$. In the case that $C$ is a free $n$\nbd{}category with $n$ finite and in light of Example \ref{dummyexample}, we will also speak of \emph{the basis of $C$} as the finite sequence
  \[
  (\Sigma_k)_{0 \leq k \leq n}.
  \]
  
  We often refer to the elements of the $n$\nbd{}basis of a free $\oo$\nbd{}category as the \emph{generating $n$\nbd{}cells}. This sometimes leads to use the alternative terminology \emph{set of generating $n$\nbd{}cells} instead of \emph{$n$\nbd{}basis}.      
\end{paragr}
\begin{definition}\label{def:rigidmorphism}
  Let $C$ and $D$ be two free $\oo$\nbd{}categories. An $\oo$\nbd{}functor $f : C \to D$ is \emph{rigid} if for every $n\geq 0$ and every generating $n$\nbd{}cell $x$ of $C$, $f(x)$ is a generating $n$\nbd{}cell of $D$.
\end{definition}
So far, we have not yet seen examples of free $\oo$\nbd{}categories. In order to
do so, we will explain in a further section a recursive way of constructing free
$\oo$\nbd{}categories; but let us first take a little detour.

\section{Suspension of monoids and counting generators}\label{sec:suspmonoids}
\begin{paragr}\label{paragr:suspmonoid}
  Let $M$ be a monoid. For every $n >0$, let $B^{n}M$ be the $n$\nbd{}magma such that:
  \begin{itemize}[label=-]
  \item it has only one object $\star$,
    \item it has only one $k$\nbd{}cell for $0 < k <n$, which is $\1^k_{\star}$,
  \item the set of $n$\nbd{}cells is (the underlying set of) $M$,
  \item for every $k<n$, the $k$\nbd{}composition of $n$\nbd{}cells is given by the composition law of the monoid (which makes sense since all $n$\nbd{}cells are $k$\nbd{}composable) and the only unital $n$\nbd{}cell is given by the neutral element of the monoid.
  \end{itemize}
  It is sometimes useful to extend the above construction to the case $n=0$ by saying that $B^0M$ is the underlying set of the monoid $M$.
  
  For $n=1$, $B^1M$ is nothing but the monoid $M$ seen as a $1$\nbd{}category with one object.

  For $n>1$, while it is clear that all first three axioms for $n$\nbd{}categories (units, functoriality of units and associativity) hold, it is not always true that the exchange rule is satisfied. If $\ast$ denotes the composition law of the monoid, this axiom states that for all $a,b,c,d \in M$, we must have
  \[
  (a \ast b) \ast (c \ast d) = (a \ast c ) \ast (b \ast d).
  \]
  It is straightforward to see that this equation holds if and only if $M$ is commutative. Hence, we have proved the following lemma.
\end{paragr}
\begin{lemma}
  Let $M$ be a monoid and $n \in \mathbb{N}$. Then:
  \begin{itemize}[label=-]
  \item if $n=1$, $B^1M$ is a $1$-category,
  \item if $n>1$, the $n$\nbd{}magma $B^nM$ is an $n$\nbd{}category if and only if $M$ is commutative.
    \end{itemize}
\end{lemma}
This construction will turn out to be of great use many times in this dissertation and we now explore a few of its properties.
\begin{lemma}\label{lemma:nfunctortomonoid}
  Let $C$ be an $n$\nbd{}category with $n\geq 1$ and let $M=(M,\ast,1)$ be a monoid (commutative if $n>1$). The map
  \begin{align*}
    \Hom_{n\Cat}(C,B^nM) &\to \Hom_{\Set}(C_n,M)\\
    F &\mapsto F_n
  \end{align*}
  is injective and its image consists exactly of those functions $f : C_n \to M$ such that:
  \begin{itemize}[label=-]
  \item for every $0 \leq k <n$ and every pair $(x,y)$ of $k$\nbd{}composable $n$\nbd{}cells of $C$, we have
    \[
    f(x\comp_ky)=f(x)\ast f(y),
    \]
  \item for every $x \in C_{n-1}$, we have
    \[
    f(1_x)=1.
    \]
    \end{itemize}
\end{lemma}
\begin{proof}
  The injectivity part follows from the fact that $(B^nM)_k$ is a singleton set for every $k<n$ and hence, an $n$\nbd{}functor $F : C \to B^nM$ is entirely determined by its restriction to the $n$\nbd{}cells $F_n : C_n \to M$.

  The characterization of the image is immediate once noted that the requirements are only the reformulation of the axioms of $n$\nbd{}functors in this particular case.
\end{proof}
\begin{lemma}\label{lemma:freencattomonoid}
  Let $C$ be an $n$\nbd{}category with $n \geq 1$ and $M$ a monoid (commutative
  if $n>1$). If $C$ has an $n$\nbd{}basis $E$, then the map
  \begin{align*}
    \Hom_{n\Cat}(C,B^nM) &\to \Hom_{\Set}(E,M)\\
    F &\mapsto F_n\vert_{E}
  \end{align*}
  is bijective.
\end{lemma}
\begin{proof}
This is an immediate consequence of the universal property of $n$\nbd{}bases (as explained in Paragraph \ref{paragr:defnbasisdetailed})
\end{proof}
We can now prove the important proposition below.
\begin{proposition}\label{prop:countingfunction}
  Let $C$ be an $\oo$\nbd{}category and suppose that $C$ has an $n$\nbd{}basis $E$ with $n\geq 0$. For every $\alpha \in E$, there exists a unique function
  \[
  w_{\alpha} : C_n \to \mathbb{N}
  \]
  such that:
  \begin{enumerate}[label=(\alph*)]
  \item $w_{\alpha}(\alpha)=1$,
  \item $w_{\alpha}(\beta)=0$ for every $\beta \in E$ such that $\beta \neq \alpha$,
  \item for every $0 \leq k<n$ and every pair $(x,y)$ of $k$\nbd{}composable $n$\nbd{}cells of $C$, we have
    \[
    w_{\alpha}(x\comp_k y)=w_{\alpha}(x)+w_{\alpha}(y).
    \]
  \end{enumerate}
\end{proposition}
\begin{proof}
  Notice first that $C$ has an $n$\nbd{}basis if and only if $\sk_n(C)$ has an
  $n$\nbd{}basis (Remark \ref{remark:nbasisncat}). Hence we can suppose that $C$ is an $n$\nbd{}category.

  For $n=0$, conditions (c) is vacuous and the assertion is trivial.

  Now let $n>0$ and consider the monoid  $\mathbb{N}=(\mathbb{N},+,0)$. The existence of a function $C_n \to \mathbb{N}$ satisfying conditions (a) and (b) follows from Lemma \ref{lemma:freencattomonoid} and the fact that it satisfies (c) follows Lemma \ref{lemma:nfunctortomonoid}.

  For the uniqueness, notice that for every $x \in C_{n-1}$ we have $1_x =1_x\comp_{n-1} 1_x$ and thus condition (c) implies that  
  \[
  w_{\alpha}(1_x)=0.
  \]
  Hence, we can apply Lemma \ref{lemma:nfunctortomonoid} which shows that $w_{\alpha}$ necessarily comes from an $n$\nbd{}functor $C \to B^n\mathbb{N}$. Then, the uniqueness follows from conditions (a) and (b) and Lemma \ref{lemma:freencattomonoid}.
\end{proof}
\begin{paragr}\label{paragr:weight}
  Let $C$ be an $n$\nbd{}category with an $n$\nbd{}basis $E$. For an $n$\nbd{}cell $x$ of $C$,
  we refer to the integer $w_{\alpha}(x)$ as the \emph{weight of $\alpha$ in
    $x$}. The reason for such a name will become clearer after Remark
  \ref{remark:weightexplicitly} where we give an
  explicit construction of $w_{\alpha}$ as a function that ``counts the number
  of occurrences of $\alpha$ in an $n$\nbd{}cell''.

  For later reference, let us also highlight the fact that in the proof of the previous proposition, we have shown the important property that if $n>0$, then for $y \in C_{n-1}$, we have
  \[
  w_{\alpha}(1_y)=0.
  \]
  This implies that for $n>1$, there might be $n$\nbd{}cells $x$ such that
  \[
  x \neq \alpha \text{ and }w_{\alpha}(x)=1.
  \]
 Indeed, suppose that there exists a $k$\nbd{}cell $z$ with $0<k<n-1$ which is not a unit on a lower dimensional cell and such that $\trgt_{k-1}(z)=\src_{k-1}(\alpha)$, then we have
  \[
  w_{\alpha}(\alpha\comp_{k-1}\1^n_z)=w_{\alpha}(\alpha)+w_{\alpha}(\1^n_z)=1.
  \]
\end{paragr}

\section{Recursive construction of free \texorpdfstring{$\oo$}{ω}-categories}\label{section:freeoocataspolygraph}
\begin{definition}\label{def:cellularextension}
  Let $n \in \mathbb{N}$. An \emph{$n$\nbd{}cellular extension} is a quadruplet $\E=(C,\Sigma,\sigma,\tau)$ where:
  \begin{itemize}[label=-]
  \item $C$ is an $n$\nbd{}category,
  \item $\Sigma$ is a set, whose elements are referred to as the \emph{indeterminates} of $\E$,
    \item $\sigma$ and $\tau$ are maps $\Sigma \to C_n$ such that for every element $x \in \Sigma$, the $n$\nbd{}cells $\sigma(x)$ and $\tau(x)$ are parallel.
  \end{itemize}
\end{definition}
\begin{paragr}
  If we are given an $n$\nbd{}category $C$, then we also say that an $n$\nbd{}cellular extension $\E$ is a \emph{cellular extension of $C$} if it is of the form $\E=(C,\Sigma,\sigma,\tau)$.
  
  Intuitively speaking, the indeterminates are formal extra $(n+1)$\nbd{}cells attached to $C$ via $\sigma$ and $\tau$. For every $x \in \Sigma$, the $n$\nbd{}cells $\sigma(x)$ and $\tau(x)$ are understood respectively the source and target of $x$ (which makes sense since these two $n$\nbd{}cells are parallel). Consequently, we often adopt the notation
  \[
  x : a \to b
  \]
  for an indeterminate such that $\sigma(x)=a$ and $\tau(x)=b$.
  \end{paragr}
\begin{definition}
  Let $\E=(C,\Sigma,\sigma,\tau)$ and $\E'=(C',\Sigma',\sigma',\tau')$ be two $n$\nbd{}cellular extensions. A morphism of $n$\nbd{}cellular extensions $\E \to \E'$ consists of a pair $(F,\varphi)$ where:
  \begin{itemize}[label=-]
  \item F is an $n$\nbd{}functor $C \to C'$,
    \item $\varphi$ is a map $\Sigma \to \Sigma'$,
  \end{itemize}
  such that for every $x \in \Sigma$, we have
  \[
  \sigma'(\varphi(x))=F(\sigma(x)) \text{ and } \tau'(\varphi(x))=F(\tau(x)).
  \]
  \end{definition}
\begin{paragr}\label{paragr:freecext}
   For $n\geq 0$, we denote by $n\Cat^{+}$ the category of $n$\nbd{}cellular extensions and morphisms of $n$\nbd{}cellular extensions. Every $(n+1)$\nbd{}category $C$ canonically defines an $n$\nbd{}cellular extension $(\tau^s_{\leq n }(C),C_{n+1},\src,\trgt)$ where $\src,\trgt : C_{n+1} \to C_n$ are the source and target maps of $C$. This defines a functor
  \begin{align*}
    U_n : (n+1)\Cat &\to n\Cat^+\\
    C &\mapsto (\tau^s_{\leq n }(C),C_{n+1},\src,\trgt).
  \end{align*}
  On the other hand, every $n$\nbd{}cellular extension $\E=(D,\Sigma,\sigma,\tau)$ yields an $(n+1)$\nbd{}category $\E^*$ defined as the following amalgamated sum:
   \begin{equation}\label{squarefreecext}
  \begin{tikzcd}[column sep=huge, row sep=huge]
    \displaystyle\coprod_{x \in \Sigma}\sS_n \ar[d,"\displaystyle\coprod_{x \in \Sigma}i_{n+1}"']\ar[r,"{\langle \sigma(x),\tau(x)\rangle_{x \in \Sigma}}"] & D \ar[d] \\
    \displaystyle\coprod_{x \in \Sigma}\sD_{n+1}\ar[r]&\E^*.
    \ar[from=1-1,to=2-2,very near end,phantom,"\ulcorner"]
  \end{tikzcd}
  \end{equation}
  This defines a functor
  \begin{align*}
    n\Cat^+ &\to (n+1)\Cat \\
    \E &\mapsto \E^*,
  \end{align*}
  which is easily checked to be left adjoint to $U_n$.

   Now let $\phi : \coprod_{x \in \Sigma} \sD_{n} \to \E^*$ be the bottom map of square \eqref{squarefreecext}. It induces a canonical map
  \begin{align*}
    j: \Sigma &\to (\E^*)_{n+1}\\
    x &\mapsto \phi_x(e_{n+1}),
  \end{align*}
  where $e_{n+1}$ is the principal $(n+1)$\nbd{}cell of $\sD_{n+1}$ (\ref{paragr:defglobe}). Notice that this map is natural in that,  for every morphism of $n$\nbd{}cellular extensions
  \[G=(F,\varphi) : \E \to \E',\] the square
  \[
  \begin{tikzcd}
    \Sigma \ar[d,"j"] \ar[r,"\varphi"] & \Sigma' \ar[d,"j'"] \\
    (\E^*)_{n+1} \ar[r,"(G^*)_{n+1}"]&(\E'^*)_{n+1}
  \end{tikzcd}
  \]
  is commutative. Notice also that $j$ is compatible with source and target in the sense that for every $x \in \Sigma$, we have
  \[
  \src(j(x))=\sigma(x) \text{ and } \trgt(j(x))=\tau(x).
  \]
\end{paragr}
\begin{lemma}\label{lemma:basisfreecext}
  Let $\E=(C,\Sigma,\sigma,\tau)$ be an $n$\nbd{}cellular extension. The canonical map
  \[
   j: \Sigma \to (\E^*)_{n+1}
   \]
   is injective.
\end{lemma}
\begin{proof}
   A thorough reading of the techniques used in the proofs of Lemma \ref{lemma:nfunctortomonoid}, Lemma \ref{lemma:freencattomonoid} and Proposition \ref{prop:countingfunction} shows that the universal property defining $\E^*$ as the amalgamated sum \eqref{squarefreecext} is sufficient enough to prove the existence, for each $x \in \Sigma$, of a function
  \[
  w_{x} : (\E^*)_{n+1} \to \mathbb{N}
  \]
  such that $w_{x}(j(x))=1$ and $w_{x}(j(y))=0$ for every $y \in \Sigma$ with $y\neq x$. In particular, this implies that $j$ is injective.
\end{proof}
\begin{paragr}
  In consequence of the previous lemma, we will often identify $\Sigma$ with a subset of $(\E^*)_{n+1}$. When we do so, it will \emph{always} be via the map $j$. This identification is compatible with source and target in the sense that the source (resp. target) of $x \in \Sigma$, seen as an $(n+1)$\nbd{}cell of $\E^*$, is exactly $\sigma(x)$ (resp. $\tau(x)$).
\end{paragr}
We can now prove the following proposition, which is the key result of this section. It is slightly less trivial than it appears.
\begin{proposition}\label{prop:fromcexttocat}
  For  every $n$\nbd{}cellular extension $\E=(C,\Sigma,\sigma,\tau)$, the subset $\Sigma \subseteq (\E^{*})_{n+1}$ is an $(n+1)$\nbd{}basis of $\E^*$.
  \end{proposition}

\begin{proof}
  Notice first that since the map $i_{n+1} : \sS_n \to \sD_{n+1}$ is nothing but
  the canonical inclusion $\sk_{n}(\sD_{n+1}) \to
  \sk_{n+1}(\sD_{n+1})=\sD_{n+1}$, it follows easily from square
  \eqref{squarefreecext} and the fact that the skeleton functors preserve colimits, that $C$ is canonically isomorphic to $\sk_n(\E^*)$ and that the map $C \to \E^*$ can be identified with the canonical inclusion $\sk_n(\E^*) \to \sk_{n+1}(\E^*)=\E^*$. Hence, cocartesian square \eqref{squarefreecext} can be identified with
    \[
  \begin{tikzcd}[column sep=huge, row sep=huge]
    \displaystyle\coprod_{x \in \Sigma}\sS_n \ar[d,"\displaystyle\coprod_{x \in \Sigma}i_{n+1}"']\ar[r,"{\langle \src(x),\trgt(x)\rangle_{x \in \Sigma}}"] & \sk_{n}(\E^*) \ar[d,hook] \\
    \displaystyle\coprod_{x \in \Sigma}\sD_{n+1}\ar[r]&\sk_{n+1}(\E^*).
    \ar[from=1-1,to=2-2,very near end,phantom,"\ulcorner"]
  \end{tikzcd}
  \]
  Since we have identified $\Sigma$ to a subset of the $n$\nbd{}cells of $\E^*$ via $j$, the above cocartesian square means exactly that $\Sigma$ is an $(n+1)$\nbd{}basis of $\E^*$. 
\end{proof}
\begin{paragr}\label{paragr:cextfromsubset}
  Let $C$ be an $(n+1)$\nbd{}category and $E$ be a subset $E \subseteq C_{n+1}$. This defines an $n$\nbd{}cellular extension
  \[
  \E_E =(\tau_{\leq n}^s(C),E,\src,\trgt),
  \]
  where $\src$ and $\trgt$ are simply the restriction to $E$ of the source and target maps $C_{n+1} \to C_n$. The canonical inclusion $E \hookrightarrow C_{n+1}$ induces a morphism of $n$\nbd{}cellular extensions
  \[
  \E_E \to U_n(C),
  \]
  and then, by adjunction, an $(n+1)$\nbd{}functor
  \[
  \E_E^* \to C.
  \]
\end{paragr}
\begin{proposition}\label{prop:criterionnbasis}
  Let $C$ be an $(n+1)$\nbd{}category. A subset $E \subseteq C_{n+1}$ is an ${(n+1)}$\nbd{}basis of $C$ if and only if the canonical $(n+1)$\nbd{}functor
  \[
  \E_E^* \to C
  \]
  is an isomorphism. 
\end{proposition}
\begin{proof}
  It is clear that the canonical $(n+1)$\nbd{}functor $\E^*_E \to C$ sends $E$, seen as a subset of $(\E^*_E)_{n+1}$, to $E$, seen as a subset of $C_{n+1}$. Hence, it follows from Proposition \ref{prop:fromcexttocat} that if this $(n+1)$\nbd{}functor is an isomorphism, then $E$ is an $(n+1)$\nbd{}base of $C$.

  Conversely, if $E$ is an $(n+1)$\nbd{}base of $C$, then we can define an $(n+1)$\nbd{}functor $C \to \E_E^*$ that sends $E$, seen as a subset of $C_{n+1}$, to $E$, seen as a subset of $(\E^*_E)_{n+1}$ (and which is obviously the identity on cells of dimension strictly lower than $n+1$). The fact that $C$ and $\E^*$ have $E$ as an $(n+1)$\nbd{}base implies that this $(n+1)$\nbd{}functor $C \to \E^*$ is the inverse of the canonical one $\E^* \to C$.
\end{proof}
\begin{paragr}\label{paragr:cextlowdimension}
  We extend the definitions and the results from \ref{def:cellularextension} to
  \ref{prop:criterionnbasis} to the case $n=-1$ by saying that a $(-1)$-cellular
  extension is simply a set $\Sigma$ (which is the set of indeterminates) and $(-1)\Cat^+$ is the category of sets. Since a $0\Cat$ is also the category of sets, it makes sense to define the functors
  \[
  U_{-1} : 0\Cat \to (-1)\Cat^+
  \]
  and
  \[
  (-)^* : (-1)\Cat^+ \to 0\Cat
  \]
  to be both the identity functor on $\Set$. 
\end{paragr}
\begin{proposition}\label{prop:freeonpolygraph}
  Let $(\E^{(n)})_{n \geq -1}$ be a sequence where:
  \begin{itemize}[label=-]
  \item $\E^{(-1)}$ is a $(-1)$-cellular extension,
    \item for every $n\geq 0$, $\E^{(n)}$ is a cellular extension of the $n$\nbd{}category $(\E^{(n-1)})^*$.
  \end{itemize}
  Then, the $\oo$\nbd{}category defined as the colimit of the canonical diagram
  \[
  (\E^{(-1)})^* \to (\E^{(0)})^* \to \cdots \to (\E^{(n)})^* \to \cdots
  \]
  is free and for every $n \in \mathbb{N}$, its $n$\nbd{}basis is (canonically isomorphic to) the set of indeterminates of $\E^{(n+1)}$.

  Moreover, suppose we are given another sequence $(\E'^{(n)})_{n\geq -1}$ as above and a sequence 
  \[
  (G^{(n)}=(F^{(n)},\varphi^{(n)}) : \E^{(n)} \to \E'^{(n)})_{n \geq -1}
  \]
  where each $G^{(n)}$ is a morphism of $n$\nbd{}cellular extensions such that for every $n \geq 0$
  \[
  F^{(n)}=(G^{(n-1)})^*.
  \]
  Then, the $\oo$\nbd{}functor
  \[
  \colim_{n\geq -1 }(\E^{(n)})^* \to \colim_{n\geq -1 }(\E'^{(n)})^*
  \]
  induced by colimit is rigid.
  
  Conversely, every free $\oo$\nbd{}category and every rigid $\oo$\nbd{}functor arise this way.
\end{proposition}
\begin{proof}
  From Proposition \ref{prop:fromcexttocat}, we know that each $(\E^{(n)})^*$ has an $(n+1)$\nbd{}basis, which is canonically isomorphic to
  the set of indeterminates of $\E^{(n)}$. Besides, since for every $n\geq 0$,
  $\E^{(n)}$ is a cellular extension of $(\E^{(n-1)})^*$, we have
  \[\sk_{n-1}((\E^{(n)})^*)=(\E^{(n-1)})^*\] by definition. Hence, by a straightforward induction, each $\E^{(n)}$ is a free $(n+1)$\nbd{}category and its $k$\nbd{}basis for $0 \leq k \leq n+1$ is (canonically isomorphic to) the set of indeterminates of $\E^{(k-1)}$.
  Now let $C :=\colim_{n \geq -1}(\E^{(n)})^*$. Since for every $k\geq 0$, $\sk_k$ preserves colimits and since $\sk_k((\E^{(n)})^*)=(\E^{(k-1)})^*$ for all $0\leq k <n$, we have that
  \[
  \sk_k(C)=(\E^{(k-1)})^*
  \]
  for every $k \geq 0$. Altogether, this proves that $C$ is free and its $k$\nbd{}basis is the set of indeterminates of $\E^{(k-1)}$ for every $k \geq 0$.

  The fact that a sequence of morphisms of cellular extensions that satisfy the hypothesis given in the statement of the proposition induces a rigid $\oo$\nbd{}functor is proven in a similar fashion.
  
  For the converse part, notice that a free $\oo$\nbd{}category $C$, whose basis is denoted by $(\Sigma_k)_{k \in \mathbb{N}}$, induces a sequence of cellular extensions:
  \[
  \E_C^{(-1)}:=\Sigma_0
  \]
  and
  \[
  \E_C^{(n)}:=(\sk_n(C),\Sigma_{n+1},\src,\trgt) \text{ for } n\geq 0.
  \]
  It follows from Proposition \ref{prop:criterionnbasis} that $\sk_n(C)\simeq (\E_C^{(n-1)})^*$ and, then, from Lemma \ref{lemma:filtration} that $C \simeq \colim_{n \geq -1}\E_C^{(n)}$.

  Finally, notice that the construction $C \mapsto (\E_C^{(n)})_{n \geq -1}$ described above is obviously functorial with respect to rigid $\oo$\nbd{}functors and the isomorphism $\sk_n(C)\simeq (\E_C^{(n)})^*$ is natural with respect to rigid $\oo$\nbd{}functors. Since the statement of Lemma \ref{lemma:filtration} is also natural in $C$, this easily implies that every rigid $\oo$\nbd{}functor arises as the colimit of sequence of morphisms of cellular extensions as described in the statement of the proposition.
\end{proof}
\begin{remark}
  The previous proposition admits an obvious truncated version for free $n$\nbd{}categories with $n$ finite. In that case, we only need a finite sequence $(\E^{(k)}))_{ -1 \leq k \leq n-1}$ of cellular extensions.
\end{remark}
\begin{remark}
  The data of a sequence $(\E^{(n)})_{n \geq -1}$ as in Proposition \ref{prop:freeonpolygraph} is commonly referred to in the literature of the field as a \emph{computad} \cite{street1976limits} or \emph{polygraph} \cite{burroni1993higher}; consequently a $\oo$\nbd{}category which is free in the sense of definition \ref{def:freeoocat} is sometimes referred to as \emph{free on a computad}-or-\emph{polygraph}. Note that the underlying polygraph of a free $\oo$\nbd{}category is uniquely determined by the free $\oo$\nbd{}category itself (a straightforward consequence of Proposition \ref{prop:uniquebasis}), and this is why we chose the shorter terminology \emph{free $\oo$\nbd{}category}. 
\end{remark}
\begin{paragr}
  Concretely, Proposition \ref{prop:freeonpolygraph} gives us a recipe to construct free $\oo$\nbd{}categories. It suffices to give a formal list of generating cells of the form:
  \begin{itemize}[label=-]
  \item generating $0$-cells : $x^0, y^0, \dots$
  \item generating $1$-cells : $x^1 : \sigma(x^1) \to \tau(x^1),  y^1 : \sigma(y^1) \to \tau(y^1), \dots$
  \item generating $2$-cells : $x^2 : \sigma(x^2) \to \tau(x^2),  y^2 : \sigma(y^2) \to \tau(y^2), \dots$
    \item $\dots$,
  \end{itemize}
  where for a generating $k$\nbd{}cell $x$ with $k>0$, $\sigma(x)$ and $\tau(x)$ are parallel $(k-1)$\nbd{}cells of the free $(k-1)$\nbd{}category recursively generated by the generating cells of dimension strictly lower than $k$.
\end{paragr}
\begin{example}
  The data of a $1$\nbd{}cellular extension $\E$ is nothing but the data of a graph $G$ (or $1$\nbd{}graph  in the terminology of \ref{paragr:defncat}), and in that case it is not hard to see that $\E^*$ is nothing but the free category on $G$. That is to say, the category whose objects are those of $G$ and whose arrows are strings of composable arrows of $G$; the composition being given by concatenation of strings. Hence, from Proposition \ref{prop:freeonpolygraph}, a ($1$-)category is free in the sense of Definition \ref{def:freeoocat} if and only if it is (isomorphic to) a free category on a graph.
\end{example}
\begin{example}
  The notion of free category on a graph is easily generalized to the notion of free $n$\nbd{}category on an $n$\nbd{}graph (with $n \in \mathbb{N}\cup\{\infty\}$). As in the previous example, every free $n$\nbd{}category on an $n$\nbd{}graph is free in the sense of Definition \ref{def:freeoocat}. \emph{However}, the converse is not true for $n>1$. The point is that in a free $n$\nbd{}category on an $n$\nbd{}graph, the source and target of a $k$\nbd{}generating cell must be $(k-1)$\nbd{}generating cells; whereas for a free $n$\nbd{}category, they can be any $(k-1)$\nbd{}cells (not necessarily generating). For example, the free $2$-category described as
  \begin{itemize}[label=-]
  \item generating $0$-cells : $A,B,C,D$
  \item generating $1$-cells : $f:A \to B$, $g : B \to C$, $h : A \to D$, $i : D \to C$
    \item generating $2$-cells : $\alpha : i\comp_0 h \Rightarrow g \comp_0 f$,
  \end{itemize}
  which can be pictured as
  \[
  \begin{tikzcd}
    & B\ar[rd,"g"] & \\
    A\ar[ru,"f"] \ar[rd,"h"'] & & C \\
    & D\ar[ru,"i"'] &
    \ar[from=3-2,to=1-2,Rightarrow,shorten <= 1.5em, shorten >= 1.5em, "\alpha"]
  \end{tikzcd}
  \]
  is not free on a $2$-graph. The reason is that the source (resp. the target) of $\alpha$ is $g \comp_0 f$ (resp. $i \comp_0 h$) which are not generating $1$-cells. 
\end{example}
\begin{example}\label{example:freemonoid}
  Let $n \geq 1$ and $M$ be a monoid (commutative if $n>1$). The $n$\nbd{}category $B^nM$ is free if and only if $M$ is a free monoid (free commutative monoid if $n>1$). If so, it has exactly one generating cell of dimension $0$, no generating cells of dimension $0 <k< n$, and the set of generators of the monoid (which is unique) as generating $n$\nbd{}cells. 
\end{example}
\section{Cells of free \texorpdfstring{$\oo$}{ω}-categories as words}\label{section:cellsaswords}

In this section, we undertake to give a more explicit construction of the $(n+1)$\nbd{}category $\E^*$ generated by an $n$\nbd{}cellular extension $\E=(C,\Sigma,\sigma,\tau)$. By definition of $\E^*$, this amounts to give an explicit description of a particular type of colimit in $\oo\Cat$. Note also that since $\tau_{\leq n}(\E^*)=C$, all we need to do is to describe the $(n+1)$\nbd{}cells of $\E^*$. This will take place in two steps: first we construct what ought to be called the \emph{free $(n+1)$\nbd{}magma generated by $\E$}, for which the $(n+1)$\nbd{}cells are really easy to describe, and then we quotient these cells as to obtain an $(n+1)$\nbd{}category, which will be $\E^*$.

Recall that an $n$\nbd{}category is a particular case of $n$\nbd{}magma.
\begin{paragr}\label{paragr:defwords}
  Let $\E=(C,\Sigma,\sigma,\tau)$ be an $n$\nbd{}cellular extension. We denote by $\W[\E]$ the set of finite words on the alphabet that has:
  \begin{itemize}[label=-]
  \item a symbol $\cc_{\alpha}$ for each $\alpha \in \Sigma$,
  \item a symbol $\ii_{x}$ for each $x \in C_n$,
  \item a symbol $\fcomp_k$ for each $0 \leq k \leq n$,
  \item a symbol of opening parenthesis $($,
  \item a symbol of closing parenthesis $)$.
  \end{itemize}
  If $w$ and $w'$ are two elements of $\W[\E]$, we write $ww'$ for their concatenation. We define the subset $\T[\E] \subseteq \W[\E]$ of \emph{well formed words} (or \emph{terms}) on the previous alphabet together with two maps $\src, \trgt : \T[\E] \to C_n$ in the following recursive way:
  \begin{itemize}[label=-]
  \item for every $\alpha \in \Sigma$, the word $(\cc_{\alpha})$ is well formed and we have
    \[
    \src((\cc_{\alpha}))=\sigma(\alpha) \text{ and } \trgt((\cc_{\alpha}))=\tau(\alpha),
    \]
  \item for every $x \in C_n$, the word $(\ii_{x})$ is well formed and we have
    \[
    \src((\ii_{x}))=\trgt((\ii_{x}))=x,
    \]
  \item if $w$ and $w'$ are well formed words such that $\src(w)=\trgt(w')$, then the word $(w\fcomp_n w')$ is well formed and we have
    \[
    \src((w\fcomp_n w'))=s(w') \text{ and } \trgt((w\fcomp_n w'))=t(w),
    \]
  \item if $w$ and $w'$ are well formed words such that $\src_k(\src(w))=\trgt_k(\trgt(w))$ with $0 \leq k < n$, then the word $(w\fcomp_k w')$ is well formed and we have
    \[
    \src((w \fcomp_k w'))=\src(w)\comp_k s(w') \text{ and } \trgt((w \fcomp_k w'))=\trgt(w) \comp_k \trgt(w').
    \]
  \end{itemize}
  As usual, for $0\leq k<n$, we define $\src_k,\trgt_k : \T[\E] \to C_k$ to be respectively the iterated source and target (and we set $\src_n=\src$ and $\trgt_n=\trgt$ for consistency). 
\end{paragr}
\begin{definition}\label{def:sizeword}
The \emph{size} of a well-formed word $w$, denoted by $\vert w \vert$, is the number of symbols $\fcomp_k$ for any $0 \leq k \leq n$ that appear in $w$.
  \end{definition}
\begin{paragr}
  Let $\E$ be an $n$\nbd{}cellular extension and let us write $\E^{+}$ for the $(n+1)$\nbd{}magma defined in the following fashion:
  \begin{itemize}[label=-]
  \item for every $0 \leq k \leq n$, we have $(\E^+)_k:=C_k$; the source, target, compositions of $k$\nbd{}cells for $0 < k \leq n$ and units on $k$\nbd{}cells for $0 \leq k <n$ are those of $C$, 
  \item $(\E^{+})_{n+1}=\T[\E]$,
  \item the source and target maps $\src,\trgt : (\E^{+})_{n+1} \to (\E^{+})_{n}$ are the ones defined in \ref{paragr:defwords},
  \item for every $n$\nbd{}cell $x$, the unit on $x$ is given by the word
    \[
    (\ii_x),
    \]
  \item for $0 \leq k \leq n$, the $k$\nbd{}composition of two $(n+1)$\nbd{}cells $w$ and $w'$ such that $\src_k(w)=\trgt_k(w')$ is given by the word
    \[
    (w\fcomp_k w').
    \]
  \end{itemize}
  By definition, $\E^+$ satisfy all the axioms for $\oo$\nbd{}categories up to
  dimension $n$. But on the other hand, the $(n+1)$\nbd{}cells of $\E^+$ make it as far as possible from being an ${(n+1)}$\nbd{}category as \emph{none} of the axioms of $\oo$\nbd{}categories are satisfied for cells of dimension $n+1$. 
\end{paragr}
\begin{remark}\label{remark:formalunformal}
  Notice that, by definition, we have
  \[
  w \comp_k w' := (w\fcomp_k w')
  \]
  for every pair $(w,w')$ of $k$\nbd{}composable $(n+1)$\nbd{}cells of $\E^+$, and we ought to be careful not to confuse the ``real'' composition symbol ``$\; \comp_k \;$''
  with the ``formal'' composition symbol ``$\; \fcomp_k\; $''. As a rule of thumb, it is better not to use both symbols in the same expression. Note also that, since we use the usual symbols ``(`` and ``)'' as formal symbols of opening and closing parenthesis, things can get really messy if we don't apply the previous rule because it would be hard to distinguish a formal parenthesis from an ``non-formal'' one.
\end{remark}
In the following definition, we consider that a binary relation $\R$ on a set $E$ is nothing but a subset of $E \times E$, and we write $x \; \R \; x'$ to say $(x,x') \in \R$.
\begin{definition}\label{def:congruence}
  Let $n\geq 1$. A \emph{congruence} on an $n$\nbd{}magma $X$ is a binary relation $\R$ on the set of $n$\nbd{}cells $X_n$ such that:
  \begin{enumerate}[label=(\alph*)]
  \item $\R$ is an equivalence relation,
  \item  if $x\; \R\; x'$ then $x$ and $x'$ are parallel,
  \item if $x \;\R\; y$ and $x' \;\R\; y'$, and if $x$ and $x'$ are $k$\nbd{}composable for some $0 \leq k <n$ then
    \[
    x \comp_k x' \; \R \;  y\comp_k y'
    \]
    (which makes sense since $y$ and $y'$ are $k$\nbd{}composable by the previous axiom).
    \end{enumerate}
\end{definition}
\begin{remark}
  Beware that in the previous definition, the relation $\R$ is \emph{only} on the set of cells of dimension $n$. 
\end{remark}
\begin{example}
  Let $F : X \to Y$ be a morphism of $n$\nbd{}magmas with $n>1$. The binary relation $\R$ on $X_n$ defined as
  \[
  x \; \R \;y \text{ if } F(x)=F(y)
  \]
  is a congruence. 
  \end{example}
\begin{paragr}
  Let $X$ be an $n$\nbd{}magma with $n \geq 1$ and $\R$ a congruence on $X$. By the first axiom of Definition \ref{def:congruence}, $\R$ is an equivalence relation and we can consider the quotient set $X_n/\R$. We write $[x]$ for the equivalence class of an $n$\nbd{}cell $x$ of $X$. From the second axiom of Definition \ref{def:congruence}, we can define unambiguously
  \[
  \src([x]):=\src(x) \text{ and } \trgt([x]):=\trgt(x),
  \]
  for $x \in X_n$ and from the third axiom, we can define unambiguously
  \[
  [x]\comp_k[y] := [x \comp_k y]
  \]
  for every pair $(x,y)$ of $k$\nbd{}composable $n$\nbd{}cells of $X$. Altogether, this defines an $n$\nbd{}magma, which we denote by $X/{\R}$, whose set $k$\nbd{}cells is $X_k$ for $0 \leq k < n$, and  $X_n/{\R}$ for $k=n$. The composition, source, target and units of cells of dimension strictly lower than $n$ being those of $X$ and the composition, source and target of $n$\nbd{}cells being given by the above formulas. 
\end{paragr}
\begin{definition}\label{def:categoricalcongruence}
  Let $\R$ be a congruence on an $n$\nbd{}magma $X$ with $n \geq 1$. We say that $\R$ is \emph{categorical} if it satisfies all four following axioms:
  \begin{enumerate}
  \item for every $k<n$ and every $n$\nbd{}cell $x$ of $X$, we have
    \[
    \1^{n}_{\trgt_k(x)}\comp_kx\; \R \; x \text{ and } x \: \R \; x \comp_k\1^n_{\src_k(x)},
    \]
    \item for every $k<n$ and for all $k$\nbd{}composable $n$\nbd{}cells $x$ and $y$ of $X$, we have 
    \[
    1_{x\comp_k y} \; \R \; 1_{x}\comp_k 1_{y},
    \]
  \item  for every $k<n$, for all $n$\nbd{}cells $x, y$ and $z$ of $X$ such that $x$ and $y$ are $k$\nbd{}composable, and $y$ and $z$ are $k$\nbd{}composable, we have
    \[
    (x\comp_{k}y)\comp_{k}z \; \R \; x\comp_k(y\comp_kz),
    \]
    \item  for all $k,l \in \mathbb{N}$ with $k<l<n$, for all $n$\nbd{}cells $x,x',y$ and $y'$ of $X$ such that
    \begin{itemize}
    \item[-] $x$ and $y$ are $l$-composable, $x'$ and $y'$ are $l$-composable,
    \item[-] $x$ and $x'$ are $k$\nbd{}composable, $y$ and $y'$ are $k$\nbd{}composable,
    \end{itemize}
    we have
    \[
    (x \comp_k x')\comp_l (y \comp_k y')\; \R \; (x \comp_l y)\comp_k (x' \comp_l y').
    \]
    \end{enumerate}
\end{definition}
\begin{example} Let $C$ be an $n$\nbd{}category with $n>1$, which we consider as an $n$\nbd{}magma. The equality on the set of $n$\nbd{}cells of $C$ is, by definition, categorical. 
\end{example}
\begin{example}\label{example:categoricalcongruence}
  Let $F : X \to Y$ be a morphism of $n$\nbd{}magmas with $n>1$ and suppose that $Y$ is an $n$\nbd{}category. Then the binary relation $\R$ on $X_n$ defined as
  \[
  x\; \R \;y \text{ if } F(x)=F(y)
  \]
  is obviously a \emph{categorical} congruence. 
\end{example}
Similarly to the ``stupid'' truncation of $\oo$\nbd{}categories (\ref{paragr:defncat}), given an $(n+1)$\nbd{}magma $X$, we write $\tau_{\leq n}^s(X)$ for the $n$\nbd{}magma obtained by simply forgetting the cells of dimension $(n+1)$.

The following lemma is trivial but nonetheless important. Its immediate proof is omitted.
\begin{lemma}\label{lemma:quotientmagma}
  Let $X$ be an $n$\nbd{}magma with $n>1$ and $\R$ a congruence on $X$. If $\tau_{\leq n}^s(X)$ is an $(n-1)$-category and $\R$ is categorical, then $X/{\R}$ is an $n$\nbd{}category.
\end{lemma}
We now wish to see how to prove the existence of a congruence defined with a condition such as ``the smallest congruence that contains a given binary relation on the $(n+1)$\nbd{}cells''. 
\begin{lemma}\label{lemma:intersectioncongruence}
  Let $X$ be an $n$\nbd{}magma with $n \geq 1$ and $(\R_i)_{i \in I}$ a \emph{non-empty} family of congruences on $X$ (i.e.\ $I$ is not empty). Then, the binary relation 
  \[
  \R:=\bigcap_{i \in I}\R_i
  \]
  is a congruence. 
\end{lemma}
\begin{proof}
  The fact that $\R$ satisfies that first and third axiom of Definition \ref{def:congruence} is immediate and do not even require that $I$ be non empty. The only left thing to prove is that if $x \; \R \; x'$, then $x$ and $x'$ are parallel. To see that, notice that since $I$ is not empty, we can choose $i \in I$. Then, by definition, we have $x \; \R_i \; x'$ and thus, $x$ and $x'$ are parallel. 
\end{proof}
\begin{remark}
  The hypothesis that $I$ be non empty in the previous lemma cannot be omitted because, in this case, we would have
  \[
  \R=\bigcap_{\emptyset}\R_i = X_n\times X_n,
  \]
  which means that \emph{all} $n$\nbd{}cells would be equivalent under $\R$. This certainly does not guarantee that the second axiom of Definition \ref{def:congruence}, i.e.\ that equivalent cells are parallel, is satisfied.
\end{remark}
\begin{lemma}\label{lemma:congruencegenerated}
  Let $X$ be an $n$\nbd{}magma with $n\geq 1$ and $E$ a set of pairs of parallel $n$\nbd{}cells of $X$. There exists a smallest congruence $\R$ on $X$ such that for every $(x,y)\in E$, we have $x \; \R \;y$.
\end{lemma}
\begin{proof}
  Let $I$ be the set of congruences $\mathcal{S}$ on $X$ such that for every $(x,y) \in E$, we have $x\; \mathcal{S} \; y$. All we have to prove is that $I$ is not empty, since in that case, we can apply Lemma \ref{lemma:intersectioncongruence} to the binary relation
  \[
  \R:=\bigcap_{\mathcal{S} \in I}\mathcal{S},
  \]
  which will obviously be the smallest congruence satisfying the desired condition. To see that $I$ is not empty, it suffices to notice that the binary relation ``being parallel $n$\nbd{}cells'' is a congruence, which obviously is in $I$. 
\end{proof}
\begin{proposition}\label{prop:smallestcategoricalcongruence}
  Let $X$ be an $n$\nbd{}magma with $n\geq 1$. There exists a smallest categorical congruence on $X$. 
\end{proposition}
\begin{proof}
  Each of the four axioms of Definition \ref{def:categoricalcongruence} says that some pairs of parallel $n$\nbd{}cells must be equivalent under a congruence $\R$ for it to be categorical. The result follows then from Lemma \ref{lemma:congruencegenerated}.
\end{proof}
\begin{paragr}
  In particular, for every $n$\nbd{}cellular extension $\E$ there exists a smallest categorical congruence on $\E^+$, which we denote here by $\equiv$. Applying Lemma \ref{lemma:quotientmagma} gives us that $\E^+/{\equiv}$ is an $(n+1)$\nbd{}category. The construction \[\E \mapsto \E^+/{\equiv}\] clearly defines a functor $\Cat_n^+ \to (n+1)\Cat$.
\end{paragr}
We can now prove the expected result.
  \begin{proposition}\label{prop:Estarasquotient}
    Let $\E$ be an $n$\nbd{}cellular extension and let $\equiv$ be the smallest categorical congruence on $\E^+$. The $(n+1)$\nbd{}category $\E^+/{\equiv}$ is naturally isomorphic to $\E^*$.
  \end{proposition}
  \begin{proof}
    The strategy of the proof is to show that the functor $\E \mapsto \E^+/{\equiv}$ is left adjoint to $U_n : (n+1)\Cat \to \Cat_n^+$. The result will follow then from the uniqueness (up to a natural isomorphism) of left adjoints.

    Let $\E = (C,\Sigma,\sigma,\tau)$ be an $n$\nbd{}cellular extension, let $D$ be an $(n+1)$\nbd{}category and let
    \[
    G=(F,\varphi) : \E \to U_n(D)
    \]
    be a morphism of $n$\nbd{}cellular extensions. We recursively define a map
  \[
  \overline{\varphi} : \T[\E] \to D_{n+1}
  \]
  as:
  \begin{itemize}[label=-]
  \item $\overline{\varphi}((\cc_{\alpha})) = \varphi(\alpha)$ for $\alpha \in \Sigma$,
  \item $\overline{\varphi}((\ii_x))=1_{F(x)}$ for $x \in C_n$,
    \item $\overline{\varphi}((w \fcomp_k w'))=\overline{\varphi}(w)\comp_k\overline{\varphi}(w')$.
  \end{itemize}
  A straightforward induction shows that $\overline{\varphi}$ is compatible with source, target and units. Hence, we can define a morphism of $(n+1)$\nbd{}magmas
  \[
  \overline{G} : \E^+ \to D
  \]
  as
  \[
  \overline{G}_{n+1}=\overline{\varphi}
  \]
  and
  \[
  \overline{G}_k=F_k \text{ for }0 \leq k \leq n.
  \]
  Since $D$ is an $(n+1)$\nbd{}category, the binary relation $\R$ on $\T[\E]$ defined as
  \[
  x\; \R \; y \text{ if } \overline{\varphi}(x)=\overline{\varphi}(y)
  \]
  is a \emph{categorical} congruence on $\E^+$ (see Example
  \ref{example:categoricalcongruence}). Since $\equiv$ is by definition the smallest categorical congruence on $\E^+$, $\equiv$ is included in $\R$. In particular, the map $\overline{\varphi}$ induces a map
  \[
  \widehat{\varphi} : \T[\E]/{\equiv} \to D_{n+1}.
  \]
  Since $\equiv$ is a congruence, this map is compatible with source and target, and it is straightforward to check that it is also compatible with units. Hence, we have an $(n+1)$\nbd{}functor
  \[
  \widehat{G} : \E^+/{\equiv} \to D
  \]
  defined as
  \[
  \widehat{G}_{n+1} = \widehat{\varphi}
  \]
  and
  \[
  \tau_{\leq n }^s(\widehat{G})=F.
  \]
  Altogether, we have constructed a map
  \begin{align*}
    \Hom_{\Cat_n^+}(\E,U_n(D)) &\to \Hom_{(n+1)\Cat}(\E^+/{\equiv},D)\\
    G &\mapsto \widehat{G}
  \end{align*}
  which is clearly natural in $\E$ and $D$. What we want to prove is that this map is a bijection.

  Let us start with the surjectivity. Let $H : \E^+/{\equiv} \to D$ be a $(n+1)$\nbd{}functor. We define a map $\varphi : \Sigma \to D_{n+1}$ as
  \begin{align*}
    \varphi : \Sigma &\to D_{n+1}\\
    \alpha &\mapsto H_{n+1}([(\cc_{\alpha})]),
  \end{align*}
  where $[w]$ is the equivalence class under $\equiv$ of an element $w \in \T[\E]$. All we need to show is that
  \[
  \widehat{\varphi}=H_{n+1}.
  \]
  Let $z$ be an element of $\T[\E]/{\equiv}$ and let us choose $w \in \T[\E]$ such that $z=[w]$. We proceed to show that $\widehat{\varphi}(z)=H_{n+1}(z)$ by induction on the size of $w$ (Definition \ref{def:sizeword}).
  If $|w|=0$, then either $w=(\cc_{\alpha})$ for some $\alpha \in \Sigma$ or $w = (\ii_x)$ for some $x \in C_n$. In the first case, we have
  \[
  \widehat{\varphi}([(\cc_{\alpha})])=\overline{\varphi}((\cc_{\alpha}))=\varphi(\alpha)=H_{n+1}([(\cc_{\alpha})]),
  \]
  and in the second case we have,
  \[
  \widehat{\varphi}([(\ii_x)])=\overline{\varphi}((\ii_x))=1_{H_n(x)}=H_{n+1}([(\ii_x)])
  \]
  where for the last equality, we used the fact that $[\ii_x]$ is the unit on $x$ in $\E^+/{\equiv}$.

  Now if $|w| = n+1$ with $n \geq 0$, then $w=(w'\fcomp_k w'')$ for some $w', w'' \in \T[\E]$ that are $k$\nbd{}composable with $k\leq n$. Hence, using the induction hypothesis, we have
  \begin{align*}
  \widehat{\varphi}([w])=\widehat{\varphi}([(w' \fcomp_k w'')])=\widehat{\varphi}([w']\comp_k [w''])&=\widehat{\varphi}([w'])\comp_k\widehat{\varphi}([w''])\\&= H_{n+1}([w'])\comp_kH_{n+1}([w''])\\&=H_{n+1}([w']\comp_k [w''])\\&=H_{n+1}([(w' \fcomp_k w'')])\\&=H_{n+1}([w]).
  \end{align*}
  We now turn to the injectivity. Let $G=(F,\varphi)$ and $G'=(F',\varphi')$ be
  two morphisms of $n$\nbd{}cellular extensions $\E \to U_n(D)$ such that
  $\widehat{G}=\widehat{G'}$. Since we have \[F=\tau_{\leq n}^s(\widehat{G})=\tau_{\leq n}^s(\widehat{G'})=F',\] all we have to show is that
  \[
  \varphi=\varphi'.
  \]
  But, by definition, for every $\alpha \in \Sigma$ we have
  \[
  \varphi(\alpha)=\widehat{\varphi}([(\cc_{\alpha})])=\widehat{\varphi'}([(\cc_{\alpha})])=\varphi'(\alpha).\qedhere
  \]
  \end{proof}
  \begin{paragr}
    In other words, the previous proposition tells us that for an $n$\nbd{}cellular extension $\E$, the $(n+1)$\nbd{}cells of $\E^*$ can be represented as equivalence classes of (well formed) words made up of the indeterminates of $\E$ and units on $n$\nbd{}cells. Note, however, that the equivalences classes are enormous. For example, for an indeterminate $\alpha$, all three following words are in the equivalence class of $(\cc_{\alpha})$:
    \[
   ((\cc_{\alpha})\fcomp_n(\ii_{\sigma(\alpha)})) \quad (((\cc_{\alpha})\fcomp_n(\ii_{\sigma(\alpha)}))\fcomp_n (\ii_{\sigma(\alpha)})) \quad ((\cc_{\alpha})\fcomp_n((\ii_{\sigma(\alpha)})\fcomp_n (\ii_{\sigma(\alpha)}))).
    \]
    Note also that when $n=0$, we have already said that a $0$-cellular
    extension $\E$ is nothing but a graph and that $\E^*$ is the free category
    on the graph. In particular, the $1$\nbd{}cells of $\E^*$ can simply be
    encoded as strings of composable arrows of the graph, which is a much
    simpler description that the one obtained from the previous proposition
    since no equivalence relation is involved. A natural question to ask is
    whether it would be possible to obtain a similar simple description for
    $n>0$. While it is certainly possible to have a simpler description of the
    $(n+1)$\nbd{}cells of $\E^*$ than the one obtained in the present section
    (for example by getting rid of some ``pathological'' well formed words and
    reducing the size of the equivalence classes), it seems not possible to
    completely avoid the use of an equivalence relation, at least not in a
    canonical way. Indeed, we have already seen in Example
    \ref{example:freemonoid} that a free \emph{commutative} monoid $M$ gives
    rise to a free $n$\nbd{}category $B^nM$ whose generating $n$\nbd{}cells are
    in bijection with the generators of the monoid. And, as soon as there are at
    least two generators of the monoid, say $a$ and $b$, there is no canonical way of representing elements of the monoid by a unique word on the generators, for we have
    \[
    ab = ba.
    \]
  \end{paragr}
 
  \section{Discrete Conduché \texorpdfstring{$\oo$}{ω}-functors}\label{section:conduche}

  \begin{paragr}
    Recall that given a category $\C$ and a class $M$ of arrows of $\C$, we say that an arrow $f : X \to Y$ of $\C$ is \emph{left orthogonal} to $M$ if for every $m : A \to B$ in $M$ and every solid arrow square
    \[
    \begin{tikzcd}
      X \ar[d,"f"'] \ar[r] & A \ar[d,"m"] \\
      Y \ar[r] & B,
      \ar[from=2-1,to=1-2,dashed,"u"]
    \end{tikzcd}
    \]
    there exists a \emph{unique} arrow $u : Y \to A$ making the whole diagram commute. 
  \end{paragr}
 \emph{In the two following paragraphs, we freely use the notations from \ref{paragr:defglobe} and advise the reader to refer to it if needed.}
  \begin{paragr}
    For all $0 \leq k < n$, we denote by $\sD_n \coprod_{\sD_k} \sD_n$ the $n$\nbd{}category defined as the following amalgamated sum:
    \[
    \begin{tikzcd}[column sep = large, row sep=large]
      \sD_k \ar[r,"\langle \src_k(e_n) \rangle"] \ar[d,"\langle \trgt_k(e_n) \rangle"'] & \sD_n \ar[d] \\
      \sD_n \ar[r] & \displaystyle\sD_n \coprod_{\sD_k} \sD_n. \ar[from=1-1,to=2-2,phantom,"\ulcorner",very near end]
    \end{tikzcd}
    \]
    For an $\oo$\nbd{}category $C$, an arrow $\sD_n \coprod_{\sD_k} \sD_n \to C$ amounts to the same data as a pair of $k$\nbd{}composable $n$\nbd{}cells of $C$. Hence, $\sD_n \coprod_{\sD_k} \sD_n$ represents the functor
    \begin{align*}
      \oo\Cat &\to \Set \\
      C &\mapsto C_n \times_{C_k} C_n.
    \end{align*}
    Given a pair of $k$\nbd{}composable $n$\nbd{}cells $(x,y)$ of an $\oo$\nbd{}category $C$, we write
    \[
    \langle x,y \rangle : \sD_n \coprod_{\sD_k} \sD_n \to C
    \]
    for the associated $\oo$\nbd{}functor.
    
    By the Yoneda Lemma, the $k$\nbd{}composition operation
    \[(\shortminus) \comp_k (\shortminus) : C_n \times_{C_k} C_n \to C_n,\] which is obviously natural in $C$, yields an $n$\nbd{}functor
    \[
    \nabla^n_k : \sD_n \to \sD_n \coprod_{\sD_k} \sD_n.
    \]
    The $n$\nbd{}cell $\nabla^n_k$ is nothing but the $k$\nbd{}composition of the two generating $n$\nbd{}cells of $\sD_n \coprod_{\sD_k} \sD_n$. In more practical terms, a commutative triangle
    \[
    \begin{tikzcd}
      \sD_n \ar[r,"\langle x \rangle"] \ar[d,"\nabla^n_k"'] & C \\
       \sD_n \coprod_{\sD_k} \sD_n \ar[ru,"{\langle x',x'' \rangle}"'] & 
    \end{tikzcd}
    \]
    means exactly that $x = x' \comp_k x''$.
  \end{paragr}
  \begin{paragr}
    For all $0 \leq k < n$, we denote by $\kappa^n_k : \sD_n \to \sD_k$, the $\oo$\nbd{}functor
    \[
   \kappa^n_k := \langle \1^{n}_{e_k}\rangle : \sD_n \to \sD_k.
   \]
   More conceptually, for every $\oo$\nbd{}category $C$ the unit map
   \[\1^{n}_{(\shortminus)} : C_k \to C_n,\]
   is obviously natural in $C$ and the map $\kappa^n_k$ is the one induced by the Yoneda Lemma. In practical terms, a commutative triangle
   \[
   \begin{tikzcd}
      \sD_n \ar[r,"\langle x \rangle"] \ar[d,"\kappa^n_k"'] & C \\
      \sD_k \ar[ru,"{\langle y \rangle}"'] &
      \end{tikzcd}
   \]
   means exactly that $x=1^{n}_y$.
    \end{paragr}
  \begin{definition}\label{def:conduche}
    An $\oo$\nbd{}functor $F : C \to D$ is \emph{discrete Conduché} if it is right orthogonal to the $\oo$\nbd{}functors
        \[
   \kappa^n_k : \sD_n \to \sD_k
   \]
   and
    \[
    \nabla^n_k : \sD_n \to \sD_n \coprod_{\sD_k} \sD_n
    \]
   for all $k,n \in \mathbb{N}$ such that $0 \leq k < n$.
  \end{definition}
  Before unfolding the previous definition, let us put here for later reference the immediate following lemma.
  \begin{lemma}\label{lemma:pullbackconduche}
    The class of discrete Conduché $\oo$\nbd{}functors is stable by pullback. This means that for every cartesian square in $\oo\Cat$
    \[
    \begin{tikzcd}
      C' \ar[d,"F'"'] \ar[r] & C \ar[d,"F"] \\
      D' \ar[r] & D,\ar[from=1-1,to=2-2,"\lrcorner",very near start,phantom]
      \end{tikzcd}
    \]
    if $F$ is a discrete Conduché $\oo$\nbd{}functor, then so is $F'$.
    \end{lemma}
  \begin{proof}
    This is a standard fact about right orthogonal classes in a category that admits pullbacks. See for example the dual of \cite[Proposition 2.1.1(e)]{freyd1972categories}.
    \end{proof}
  \begin{paragr}
    Let us now give a practical version of Definition \ref{def:conduche}. The right orthogonality to $\kappa^n_k$ means that for every $n$\nbd{}cell $x$ of $C$ and every $k$\nbd{}cell $y$ of $D$ such that
    \[
    F(x)=\1^{n}_y,
    \]
    there exists a unique\footnote{The map $z \mapsto 1^n_k(z)$ being injective, the uniqueness actually comes for free.} $k$\nbd{}cell $z$ of $C$ such that
    \[
    x=\1^n_z \text{ and } F(z)=y.
    \]
    Similarly, the right orthogonality to $\nabla^n_k$ means that for every $k$\nbd{}cell $x$ of $C$, if
    \[
    f(x) = y' \comp_k y''
    \]
    with $(y',y'')$ a pair of $k$\nbd{}composable $n$\nbd{}cells of $D$, then there exists a \emph{unique} pair $(x',x'')$ of $k$\nbd{}composable $n$\nbd{}cells of $C$ such that
    \begin{enumerate}
    \item $x=x'\comp_k x''$,
    \item $F(x')=y'$ and $F(x'')=y''$.
    \end{enumerate}
  \end{paragr}
    As it turns out, the definition we gave of discrete Conduché $\oo$\nbd{}functor is highly redundant.
  \begin{lemma}\label{lemma:redundencies}
    Let $F : C \to D$ be an $\oo$\nbd{}functor and let $k<n \in \mathbb{N}$. If $F$ is right orthogonal to $\nabla^n_k$ then it is right orthogonal to $\kappa^n_k$.
  \end{lemma}
  \begin{proof}
    Let $x$ be an $n$\nbd{}cell of $C$ and suppose that $F(x)=\1^n_y$ with $y$ a $k$\nbd{}cell of $D$. Notice that
    \[
    F(x)=\1^n_y \comp_k \1^n_y
    \]
    and
    \[
    x=x\comp_k \1^n_{\src_k(x)}=\1^n_{\trgt_k(x)}\comp_k x
    \]
    and
    \[
    F(\1^n_{\src_k(x)})=\1^n_{\src_k(F(x))}=\1^n_y=\1^n_{\trgt_k(F(x))}=F(\1^n_{\trgt_k(x)}).
    \]
    Using the uniqueness part of the right orthogonality to $\nabla^n_k$, we deduce that $x=\1^n_{\src_k(x)}=\1^n_{\trgt_k(x)}$. Thus, if we set $z=\src_k(x)=\trgt_k(x)$, we have $x=\1^n_z$ and $F(z)=y$, which is what we needed to prove.
  \end{proof}
  \begin{remark}
    In light of the previous lemma, the reader might wonder why the right orthogonality to $\kappa^n_k$  was included in Definition \ref{def:conduche}. The motivation for such a choice is that it \emph{should} be possible to apply this definition \emph{mutatis mutandis} to weak $\oo$\nbd{}categories where Lemma \ref{lemma:redundencies} might not hold anymore.
  \end{remark}
  Somewhat related to Lemma \ref{lemma:redundencies} is the following lemma.
  \begin{lemma}\label{lemma:technicalorthogonality}
    Let $F : C \to D$ be an $\oo$\nbd{}functor and $k<m<n \in \mathbb{N}$. If $F$ is right orthogonal to $\nabla^n_k$ and $\kappa^n_m$, then it is right orthogonal to $\nabla^m_k$.
  \end{lemma}
  \begin{proof}
    Let $x$ be an $m$-cell of $C$ and suppose that
    \[
    F(x)=y'\comp_k y''
    \]
    with $(y',y'')$ a pair of $k$\nbd{}composable $m$-cells of $D$. Then, we have
    \[
    F(\1^n_x)=\1^n_{y'} \comp_k \1^n_{y''}.
    \]
    From the right orthogonality to $\nabla^n_k$, we know that there exist $z$ and $z''$ $k$\nbd{}composable $n$\nbd{}cells of $C$ such that $F(z')=\1^n_{y'}$, $F(z'')=\1^n_{y''}$ and
    \[
    \1^n_x=z' \comp_k z''.
    \]
    From the right orthogonality to $\kappa^n_m$, we know that there exist $x'$ and $x''$ $m$\nbd{}cells of $C$ such that $z'=\1^n_{x'}$, $z''=\1^n_{x''}$, $F(x')=y'$ and $F(x'')=y''$. It follows that $\src_k(x')=\src_k(x'')$ and
    \[
    \1^n_x=\1^n_{x'} \comp_k \1^n_{x''}=\1^n_{x'\comp_k x''},
    \]
    hence $x=x'\comp_k x''$. This proves the existence part of the right orthogonality to $\nabla^m_k$.

    Now suppose that there are two pairs $(x_1',x_1'')$ and $(x_2',x_2'')$ that lift the pair $(y',y'')$ in the required way. It follows that the pairs $(\1^n_{x_1'},\1^n_{x_1''})$ and $(\1^n_{x_2'},\1^n_{x_2''})$ lift the pair $(\1^n_{y'},\1^n_{y''})$ in the required way.

    From the uniqueness part of the right orthogonality to $\nabla^n_k$, we deduce that $\1^n_{x_1'}=\1^n_{x_2'}$ and $\1^n_{x_1''}=\1^n_{x_2''}$, hence $x'_1=x'_2$ and $x''_1=x''_2$.
  \end{proof}
  Definition \ref{def:conduche} admits an obvious truncated version as follows.
  \begin{definition}\label{def:mconduche}
 Let $m \in \mathbb{N}$. An $m$\nbd{}functor $F : C \to D$ is a \emph{discrete Conduché} if it is right orthogonal to $\kappa^n_k$ and $\nabla^n_k$ for all $k,n \in \mathbb{N}$ such that $0 \leq k < n \leq m$. 
  \end{definition}
   \begin{paragr}
    Let us dwell on a subtlety here. Since we have considered $n\Cat$ as a full
    subcategory of $\oo\Cat$ for every $n \in \mathbb{N}$ (see
    \ref{paragr:defncat}), an $n$\nbd{}functor is a particular case of
    $\oo$\nbd{}functor. Hence, it also makes sense to call an $n$\nbd{}functor
    \emph{discrete Conduché} when, seen as an
    $\oo$\nbd{}functor, it is discrete Conduché in the
    sense of Definition \ref{def:conduche}. This might have been conflicting
    with Definition \ref{def:mconduche} but the following lemma tells us that,
    in fact, the two notions are equivalent. Consequently, there is no
    distinction to make. 
  \end{paragr}
  \begin{lemma}
    Let $n \in \mathbb{N}$. An $n$\nbd{}functor $F : C \to D$ is discrete
    Conduché in the sense of Definition \ref{def:mconduche} if and only if, seen as an $\oo$\nbd{}functor, it is discrete Conduché in the sense of Definition \ref{def:conduche}.
  \end{lemma}
  \begin{proof}
    We only have to prove that if $F$ is an discrete Conduché $n$\nbd{}functor then, when seen as an $\oo$\nbd{}functor, it is a discrete Conduché $\oo$\nbd{}functor; the other implication being trivial. From Lemma \ref{lemma:redundencies}, this amounts to showing that if $F$ is right orthogonal to $\nabla^m_k$ for all $k,m \in \mathbb{N}$ such that $0 \leq k < m \leq n$, then it is also right orthogonal to $\nabla^m_k$ for all $k,m \in \mathbb{N}$ such that $0 \leq k < m$ and $m > n$. This follows easily from the fact that for every $m$\nbd{}cell $x$  in an $n$\nbd{}category (seen as an $\oo$\nbd{}category) with $m > n$, there exists a unique $n$\nbd{}cell $x'$ such that $x = \1^m_{x'}$. Details are left to the reader. 
  \end{proof}
In practice, we will use the following criterion to detect discrete Conduché $n$\nbd{}functors. 
  \begin{proposition}\label{prop:conduchepractical}
    Let $n \in \mathbb{N}$. An $n$\nbd{}functor $F : C \to D$ is a discrete Conduché $n$\nbd{}functor if and only if it is right orthogonal to $\nabla^n_k$ for every $k \in \mathbb{N}$ such that $k<n$.
    \end{proposition}
  \begin{proof}
     Immediate consequence of Lemmas \ref{lemma:redundencies} and \ref{lemma:technicalorthogonality}.
  \end{proof}
  \begin{remark}
    In the case $n=1$, we recover the usual definition of what's commonly referred to as discrete Conduché fibration (see for example \cite{johnstone1999note}).
    \end{remark}
  As we shall now see, discrete Conduché $\oo$\nbd{}functors have a deep connection with free $\oo$\nbd{}categories. We begin with an easy property.
  \begin{lemma}\label{lemma:conducheindecomposable}
    Let $F : C \to D$ be a discrete Conduché $\oo$\nbd{}functor and $n \in \mathbb{N}$. An $n$\nbd{}cell $x$ of $C$ is indecomposable (Definition \ref{def:indecomposable}) if and only if $F(x)$ is indecomposable.
  \end{lemma}
  \begin{proof}
    When $n=0$, there is nothing to prove since every $0$-cell is indecomposable. We assume now that $n>0$.

    Suppose that $x$ is indecomposable. Since $x$ is not a unit on a lower dimensional cell, the right orthogonality to $\kappa^n_k$ for every $0 \leq k <n$ implies that $F(x)$ is also not a unit on a lower dimensional cell. Now, if we have
    \[
    F(x)=y'\comp_k y'',
    \]
    with $(y',y'')$ a pair of $k$\nbd{}composable $n$\nbd{}cells of $D$, then the right orthogonality to $\nabla^n_k$ implies that
    \[
    x=x'\comp_k x''
    \]
    with $F(x')=y'$ and $F(x'')=y''$. Since $x$ is indecomposable, $x'$ or $x''$ has to be of the form $\1^n_z$ with $z$ and $k$\nbd{}cell of $C$, and thus $y'$ or $y''$ is of the form $\1^n_{F(z)}$. This proves that $F(x)$ is indecomposable.

    Conversely, suppose that $F(x)$ is indecomposable. The $n$\nbd{}cell $x$ cannot be a unit on a lower dimensional cell as if it were, then $F(x)$ would also be a unit on a lower dimensional cell, which is impossible since it is indecomposable. Now, if
    \[
    x=x'\comp_k x''
    \]
    with $(x',x'')$ a pair of $k$\nbd{}composable $n$\nbd{}cells, then
    \[F(x)=F(x')\comp_k F(x'').\]
    Hence, either $F(x')$ or $F(x'')$ is a unit on a lower dimensional cell and from the right orthogonality to $\kappa^n_k$ we deduce that either $x'$ or $x''$ is a unit on a lower dimensional cell. This proves that $x$ is indecomposable. 
  \end{proof}
  From the previous lemma and Proposition \ref{prop:uniquebasis}, we deduce the following proposition.
  \begin{proposition}
    Let $F : C \to D$ be an $\oo$\nbd{}functor with $C$ and $D$ free $\oo$\nbd{}categories. If $F$ is a discrete Conduché $\oo$\nbd{}functor then it is rigid. 
  \end{proposition}

  \begin{remark}
    The converse of the above proposition does not hold. For details, see \cite[Appendix A]{guetta2020polygraphs}.
    \end{remark}
  We now turn to the main result concerning the relation between discrete Conduché $\oo$\nbd{}functors and free $\oo$\nbd{}categories. As a matter of fact, a significant amount of results and definitions that we have seen so far were geared towards the statement and proof of the following theorem. 
  \begin{theorem}\label{thm:conduche}\label{THM:CONDUCHE}
    Let $F : C \to D$ be a discrete Conduché $\oo$\nbd{}functor.
    \begin{enumerate}
    \item If $D$ is free then so is $C$.
    \item If $C$ is free and if $F_n : C_n \to D_n$ is surjective for every $n \in \mathbb{N}$, then $D$ is also free.
      \end{enumerate}
  \end{theorem}
  The proof of the previous theorem is long and technical and the next three sections are devoted to it. Before that, let us make Theorem \ref{thm:conduche} a little more precise.
  \begin{paragr}
    Let $F : C \to D$ be a discrete Conduché $\oo$\nbd{}functor.
    \begin{enumerate}
    \item If $D$ is free, we know from the previous theorem that $C$ is also free. Now, if we write $\Sigma^D_n$ for the $n$\nbd{}basis of $D$, then it follows from Proposition \ref{prop:uniquebasis} and Lemma \ref{lemma:conducheindecomposable} that the $n$\nbd{}basis of $C$ is given by
      \[
      \Sigma^C_n :=\left\{ x \in C_n \vert F(x) \in \Sigma^D_n\right\}.
      \]

    \item Similarly, if $C$ is free and if $F_k : C_k \to D_k$ is surjective for every $k \in \mathbb{N}$, we know that $D$ is also free. If we write $\Sigma^C_n$ for the $n$\nbd{}basis of $C$, then it follows once again from Proposition \ref{prop:uniquebasis} and Lemma \ref{lemma:conducheindecomposable} that the $n$\nbd{}basis of $D$ is given by
      \[
       \Sigma^D_n :=\left\{ F(x) \in D_n \vert  x \in \Sigma_n^C \right\}.
      \]
    \end{enumerate}
  \end{paragr}
  \section{Proof of Theorem \ref{THM:CONDUCHE}: part I}
  This first part of the proof of Theorem \ref{thm:conduche} consists of several technical results on words. They lay a preliminary foundation on which the key arguments of the proof will later rely.
  
  \emph{For the whole section, we fix a cellular $n$\nbd{}extension $\E=(C,\Sigma,\sigma,\tau)$. A ``word'' always means an element of $\W[\E]$ and a ``well formed word'' always means an element of $\T[\E]$.}
\begin{paragr}
  For a word $w$, the number of symbols that appear in $w$ is referred to as the \emph{length of $w$} and denoted by $\mathcal{L}(w)$. Recall also from Definition \ref{def:sizeword} that when $w$ is well formed, the number of occurrences of symbols $\fcomp_k$ for any $0 \leq k \leq n$ that appear in $w$ is referred to as the \emph{size of $w$} and denoted by $\vert w \vert$.
\end{paragr}
\begin{definition}\label{def:subword}A word $v$ is a \emph{subword} of a word $w$ if there exist words $a$ and $b$ such that $w$ can be written as
  \[w=avb.\]
\end{definition}
\begin{remark}
  Beware that in the previous definition, none of the words were supposed to be well formed. In particular, a subword of a well formed word is not necessarily well formed. 
\end{remark}
\begin{paragr}\label{paragr:lengthandsize} Since a word $w$ is a finite sequence of symbols, it makes sense to write $w(i)$ for the symbol at position $i$ of $w$, with $0 \leq i \leq \mathcal{L}(w)-1$.

  For every $0 \leq i \leq \mathcal{L}(w)-1$, define $P_{w}(i)$ to be the number of opening parentheses in $w$ with position $ \leq i$ minus the number of closing parentheses in $w$ with position $\leq i$. This defines a function
  \[
    P_{w} : \{0,\dots,\mathcal{L}(w)-1\} \to \mathbb{Z}.
  \]
\end{paragr}
\begin{remark}
  Such a counting function is standard in the literature about formal languages. For example see \cite[chapter 1, exercice 1.4]{hopcroft1979introduction}.
  \end{remark}
\begin{definition}
  A word $w$ is \emph{well parenthesized} if:
  \begin{enumerate}
    \item it is not empty,
  \item $P_w(i) \geq 0$ for every $0 \leq i \leq \mathcal{L}(w)-1$,
    \item $P_w(i) = 0$ if and only if $i = \mathcal{L}(w)-1$.
    \end{enumerate}
\end{definition}
\begin{paragr}
  It follows from the previous definition that the first letter of a well parenthesized word is necessarily an opening parenthesis and that the last letter is necessarily a closing parenthesis. Thus, the length of a well parenthesized word is not less than 2.

  Moreover, it is immediate that if $w_1$ and $w_2$ are well parenthesized words, then for every $0 \leq k \leq n$, the word
  \[
  (w_1 \fcomp_k w_2)
  \]
  is well parenthesized.
  \end{paragr}
\begin{lemma}\label{lemmawellformedispseudo}
  A well formed word is  well parenthesized.
\end{lemma}
\begin{proof}
  Let $w$ be a well formed word. We proceed by induction on $|w|$. If $|w|=0$, then $w$ is either of the form
  \[
  (\cc_{\alpha})
  \]
  or of the form
  \[
  (\ii_{x}).
  \]
  In either case, the assertion is trivial.
  Now suppose that $|w|>0$, we know by definition that
  \[
  w=(w_1 \fcomp_k w_2)
  \]
  with $w_1,w_2$ well formed words such that $|w_1|, |w_2| < |w|$. The desired properties follow easily from the induction hypothesis. Details are left to the reader.
\end{proof}
The converse of the previous lemma is obviously not true. However, Lemma \ref{lemmapartialconverse} below is a partial converse.
\begin{lemma}\label{lemmasubwordpseudo}
  Let $w$ be a well parenthesized word of the form
  \[
  w=(w_1 \fcomp_k w_2)
  \]
  with $w_1$ and $w_2$ well parenthesized words and $0 \leq k \leq n$, and let $v$ be a subword of $w$. If $v$ is well parenthesized then one the following holds:
  \begin{enumerate}
  \item $v=w$,
  \item $v$ is a subword of $w_1$,
    \item $v$ is a subword of $w_2$.
    \end{enumerate}
\end{lemma}
\begin{proof}
  Let $a$ and $b$ be words such that
  \[
  avb=w=(w_1\fcomp_k w_2).
  \]
  Let $l_1, l_2, l, l_a, l_b, l_v$ respectively be the lengths of $w_1,w_2,w,a,b,v$. Notice that
  \[
  l_a+l_v+l_b=l=l_1+l_2+3.
  \]
  Notice also that since $v$ is well parenthesized, the following cases are forbidden:
  \begin{enumerate}
  \item $l_1 \leq l_a \leq l_1 +1$,
  \item $l_2 \leq l_b \leq l_2 +1$,
  \item $l_a\geq l-1$,
    \item $l_b\geq l-1$.
  \end{enumerate}
  Indeed, the first case would imply that the first letter of $v$ is a closing parenthesis or the symbol $\fcomp_k$. Similarly, the second case would imply that the last letter of $v$ is an opening parenthesis or the symbol $\fcomp_k$. The third and fourth cases would imply that $l_v<2$ which is also impossible.

  This leaves us with the following cases:
  \begin{enumerate}
  \item $l_a=0$,
  \item $l_b=0$,
  \item $0<l_a<l_1$ and $0<l_b < l_2$,
  \item $0<l_a<l_1$ and $l_b > l_2+1$,
  \item $l_1 +1 < l_a$ and $0 < l_b < l_2$.
  \end{enumerate}
  If we are in the first case, then
  \[
  P_w(j)=P_v(j)
  \]
  for $0 \leq j \leq l_v-1$. This implies that $P_w(l_v-1)=0$ which means that $l=l_v$, hence $w=v$.
  
  By a similar argument left to the reader, we can show that the second case implies that $w=v$.

  If we are in the fourth (resp.\ fifth) case, then it is clear that $v$ is a subword of $w_1$ (resp.\ $w_2$).

  Suppose now that we are in the third case. Intuitively, it means that the first letter of $v$ is inside $w_1$ and the last letter of $v$ is inside $w_2$. Notice first that
  \begin{equation}\label{usefulinequality}\tag{$\star$}
  l_a<l_1 <l_a+l_v-3,
  \end{equation}
   where the inequality on the right comes from the fact that $l_v \geq 2$ (because $v$ is well formed).
  
  Besides, by definition of $P_w$,
  \[
  P_w(j)=P_v(j-l_a)+P_w(l_a)
  \]
  for $l_a \leq j < l_v+l_a$. In particular, we have
  \[
  1=P_{w_1}(l_1-1)+1=P_w(l_1-1)=P_v(l_1-1)+P_w(l_a).
  \]
  From \eqref{usefulinequality} and since $v$ is well parenthesized, we deduce that
  \[
  P_v(l_1-1)>0.
  \]
  Hence $P_w(l_a)\leq 0$, which is impossible since $w$ is well formed and since ${l_a < l-1}$.
  \end{proof}
  \begin{lemma}\label{lemmapartialconverse}
    Let $w$ be a well formed word. Any well parenthesized subword of $w$ is also
    well formed.
\end{lemma}
\begin{proof}
  We proceed by induction on $|w|$. If $|w|=0$, then $w$ is either of the form
  \[
  (\cc_{\alpha})
  \]
  or of the form
  \[
  (\ii_x).
  \]
  In both cases, the only well parenthesized subword of $w$ is $w$ itself, which
  is well formed by hypothesis.

  Suppose now that $|w|>0$ and let $u$ be a well parenthesized subword of $w$. By definition,
  \[
  w=(w_1 \fcomp_k w_2)
  \]
  with $|w_1|,|w_2| < |w|$. By Lemmas \ref{lemmawellformedispseudo} and \ref{lemmasubwordpseudo}, we have that either:
  \begin{itemize}
  \item[-] $u=w$ in which case $u$ is well formed by hypothesis,
  \item[-] $u$ is a subword of $w_1$ and from the induction hypothesis we deduce that $u$ is well formed,
    \item[-] $u$ is a subword of $w_2$ which is similar to the previous case.\qedhere
    \end{itemize}
\end{proof}
\begin{proposition}\label{prop:subwords} Let $w$ be a well formed word of the form
  \[
  w = (w_1 \fcomp_k w_2)
  \]
  with $w_1$ and $w_2$ well formed words and $0 \leq k \leq n$, and let $v$ be a subword of $w$. If $v$ is well formed, then we are in one of the following cases:
  \begin{enumerate}
  \item $v=w$,
  \item $v$ is a subword of $w_1$,
    \item $v$ is a subword of $w_2$.
    \end{enumerate}
\end{proposition}
\begin{proof}
  This follows immediately from Lemma \ref{lemmawellformedispseudo} and Lemma \ref{lemmasubwordpseudo}.
  \end{proof}
\begin{proposition}\label{prop:subwordsubstitution}
  Let $u$ be a well formed word of the form
  \[
  vew 
  \]
  with $v$, $w$ and $e$ words and such that $e$ is well formed. If $e'$ is a well formed word that is parallel to $e$, then the word
  \[
  ve'w 
  \]
is also well formed.
\end{proposition}
\begin{proof}
  We proceed by induction on $|u|$.
  \begin{description}
    \item[Base case] If $|u|=0$, then necessarily $v$ and $w$ are both the empty word and the assertion is trivial.

  \item[Inductive step] If $|u| \geq 1$, then
  \[
  u=(u_1 \fcomp_k u_2)
  \]
  with $u_1$ and $u_2$ well formed words such that $|u_1|,|u_2| < |u|$. By hypothesis, $e$ is a subword of $u$ and from Proposition \ref{prop:subwords}, we are in one of the following cases: 
  \begin{itemize}
  \item[-] $u=e$ in which case the assertion is trivial.
  \item[-] $e$ is a subword of $u_1$, which means that there exist words $\tilde{v}, \tilde{w}$ such that
    \[
    u_1=\tilde{v}e\tilde{w}.
    \]
    Moreover, we have \[v=(\tilde{v}\] and \[w=\tilde{w}\fcomp_k u_2).\] By induction hypothesis, the word
    \[
    \tilde{v}e'\tilde{w}
    \]
    is well formed and thus
    \[
    (\tilde{v}e'\tilde{w}\fcomp_k u_2)=vew
    \]
    is well formed.
    \item[-] $e$ is a subword of $u_2$, which is symmetric to the previous case.\qedhere
  \end{itemize}
  \end{description}
\end{proof}
\begin{lemma}\label{lemmaunicitydecompositionpseudo}
  Let $w_1,w_2,w_1',w_2'$ be well parenthesized words, and let $0 \leq k \leq n$ and $0 \leq k' \leq n$ be such that
  \[
  (w_1\fcomp_k w_2) = (w_1' \fcomp_{k'} w_2').
  \]
  Then, $w_1=w_1'$, $w_2=w_2'$ and $k=k'$.
\end{lemma}
\begin{proof}
  Let us define $l:=\mathrm{min}(\mathcal{L}(w_1),\mathcal{L}(w_1'))$. Notice that
  \[
  P_w(j)=P_{w_1}(j-1)+1=P_{w_1'}(j-1)+1
  \]
  for $0 < j \leq l$, hence
  \[
  P_{w_1}(l-1)=P_{w_1'}(l-1).
  \]
Since $w_1$ and $w_1'$ are well parenthesized, one of the members of the last equality (and thus both) is equal to $0$. This implies that $\mathcal{L}(w_1)=\mathcal{L}(w_1')$ and the desired properties follow immediately from that.
\end{proof}

\begin{proposition}\label{prop:unicitydecomposition}
  Let $w_1,w_1',w_2,w_2'$ be well formed words, and let $0 \leq k \leq n$ and $0 \leq k' \leq n$ be such that $(w_1\fcomp_k w_2)$ and $(w_1'\fcomp_{k'} w_2')$ are well formed.
  If
  \[
  (w_1 \fcomp_k w_2) = (w_1' \fcomp_{k'} w_2'),
  \]
  then
  \[
  w_1=w_1', w_2=w_2' \text{ and }k=k'.
  \]
  \end{proposition}
\begin{proof}
This follows from Lemma \ref{lemmawellformedispseudo} and Lemma \ref{lemmaunicitydecompositionpseudo}.
\end{proof}
\begin{corollary}\label{cor:compatiblesubwords}
  Let $w$ be a well formed word and suppose that it can be written as
  \[w=(w_1 \fcomp_k w_2)\]
with $w_1$ and $w_2$ well formed words and $0\leq k \leq n$. Then $\src_k(w_1)=\trgt_k(w_2)$.
\end{corollary}
\begin{proof}
  By hypothesis, $|w| \geq 1$. From the definition of well formed words, we know that $w$ is of the form
  \[
  (w_1'\fcomp_{k'}w_2')
  \]
  with $w_1'$ and $w_2'$ well formed words and $0 \leq k' \leq n$ such that
  \[
  \src_{k'}(w_1')=\trgt_{k'}(w_2').
  \]
  From Proposition \ref{prop:unicitydecomposition}, we have that $w_1'=w_1$, $w_2'=w_2$ and $k=k'$.
\end{proof}

\section{Proof of Theorem \ref{thm:conduche}: part II}
Recall that we have seen in Proposition \ref{prop:smallestcategoricalcongruence}
that there exists a smallest categorical congruence on every $n$\nbd{}magma $X$
with $n \geq 1$. However, the description of this congruence that we used when
proving its existence was rather abstract. The main goal of this section is to
give a more concrete description of the smallest categorical congruence in the
case that $X= \E^+$ for an $n$\nbd{}cellular extension $\E$. This description
will turn out to be crucial for the last part of the proof of Theorem
\ref{thm:conduche}.

\begin{definition}\label{def:elementarymove}
  Let $n \geq 0$, $\E$ be an $n$\nbd{}cellular extension
  and $u,u' \in \T[\E]$. An \emph{elementary move} from $u$ to $u'$ is a
  quadruple $\mu=(v,w,e,e')$ with $v,w \in \W[\E]$ and $e,e' \in \T[\E]$ such
  that
  \[
    u=vew,
  \]
  \[
    u'=ve'w,
  \]
  and such that one of the following holds:
  \begin{enumerate}[label=(\arabic*)]
  
  \item $e$ is of the form
    \[((x\fcomp_ky)\fcomp_kz)\] and $e'$ is of the form
    \[
      (x\fcomp_k(y\fcomp_kz))\] with $x,y,z \in \T[\E]$ and $0\leq k \leq n$,
  \item $e$ is of the form \[((\ii_c)\fcomp_kx)\] and $e'$ is of the form
    \[x\] with $x \in \T[\E]$, $0 \leq k \leq n$ and $c=\1^n_{\trgt_k(x)}$,
  \item $e$ is of the form
    \[(x\fcomp_k(\ii_{c}))\] and $e'$ is of the form
    \[x\] with $x \in \T[\E]$, $0 \leq k \leq n$ and $c=\1^n_{\src_k(x)}$,
  \item $e$ is of the form \[((\ii_{c})\fcomp_k(\ii_{d}))\] and $e'$ is of the
    form \[(\ii_{c\comp_kd})\]with $(c,d)$ a pair of $k$\nbd{}composable
    $n$\nbd{}cells of $C$ with $0 \leq k < n$,
  \item $e$ is of the form \[((x\fcomp_ky)\fcomp_l(z\fcomp_kt))\] and $e'$ is of
    the form
    \[((x\fcomp_lz)\fcomp_k(y\fcomp_lt))\] with $x,y,z,t \in \T[\E]$ and $0\leq
    l < k \leq n$.
  \end{enumerate}
\end{definition}
\begin{paragr}
  We will use the notation
  \[
    \mu : u \rightarrow u'
  \]
  to say that $\mu$ is an elementary move from $u$ to $u'$.

  We write $\G[\E]$ for the graph (or $1$-graph in the terminology of
  \ref{paragr:defncat}) defined as:
  \begin{itemize}
  \item[-] the set of objects of $\G[\E]$ is $\T[\E]$,
  \item[-] for all $u,u'$ in $\T[\E]$, an arrow of $\G[\E]$ from $u$ to $u'$ is
    an elementary move from $u$ to $u'$.
  \end{itemize}
  We will use the categorical notation
  \[
    \G[\E](u,u')
  \]
  for the set of arrows from $u$ to $u'$.

  Finally, we will also sometimes write
  \[
    u \leftrightarrow u'
  \]
  to say that there exists an elementary move from $u$ to $u'$ or from $u'$ to
  $u$.
\end{paragr}
\begin{definition}\label{definitionequivalence}
  Let $\E$ be an $n$\nbd{}cellular extension and $u,u'
  \in \T[\E]$. We say that the well formed words $u$ and $u'$ are
  \emph{equivalent} and write
  \[
    u \sim u'
  \]
  if they are in the same connected component of $\G[\E]$. More precisely, this
  means that there exists a finite sequence $(u_j)_{0\leq j \leq N}$ of well
  formed words with $u_0=u$, $u_N=u'$ and $u_j \leftrightarrow u_{j+1}$ for $0
  \leq j < N$. The equivalence class of a well formed word $u$ will be denoted
  by $[u]$.
\end{definition}

\begin{lemma}\label{lemma:equivparallel}Let $u,u' \in \T[\E]$. If $u \sim u'$,
  then $u$ and $u'$ are parallel.
\end{lemma}
\begin{proof}
  Let
  \[
    \mu =(v,w,e,e') : u \to u'
  \]
  be an elementary move from $u$ to $u'$. Let us show that $\src(u)=\src(u')$
  and $\trgt(u)=\trgt(u')$ with an induction on $\mathcal{L}(v)+\mathcal{L}(w)$.
  Notice first that, by definition of elementary moves, $|u|\geq 1$ and thus
  \[
    u=(u_1 \fcomp_k u_2)
  \]
  with $u_1,u_2 \in \T[\E]$.
  \begin{description}
    
  \item[Base case] If $\mathcal{L}(v)+\mathcal{L}(w)=0$, it means that $v$ and
    $w$ are both the empty word. It is then straightforward to check the desired
    property using Definition \ref{def:elementarymove}.
  \item[Inductive step] Suppose now that $\mathcal{L}(v)+\mathcal{L}(w)\geq 0$.
    Since $e$ is a subword of $u$ and is well formed, we know from Proposition
    \ref{prop:subwords} that we are in one of the following cases:
    \begin{itemize}
    \item[-] $e=u$, which is exactly the base case.
    \item[-] $e$ is a subword of $u_1$, which means that there exist
      $\tilde{v},\tilde{w} \in \T[\E]$ such that
      \[
        u_1 = \tilde{v}e\tilde{w}.
      \]
      Moreover, we have
      \[
        v=(\tilde{v}
      \]
      and
      \[
        w=\tilde{w}\fcomp_ku_2).
      \]
      From Proposition \ref{prop:subwordsubstitution}, the word
      \[
        u_1':=\tilde{v}e'\tilde{w}
      \]
      is well formed. Since $\mathcal{L}(\tilde{v})+\mathcal{L}(\tilde{w}) <
      \mathcal{L}(v)+\mathcal{L}(w)$, we can use the induction hypothesis on
      \[
        \tilde{\mu} :=(\tilde{v},\tilde{w},e,e') : u_1 \to u_1'.
      \]
      This shows that $\src(u_1)=\src(u_1')$ and $\trgt(u_1)=\trgt(u_1')$ and
      since
      \[
        u=(u_1\fcomp_k u_2) \text{ and } u'=(u_1' \fcomp_k u_2)
      \]
      it follows easily that $\src(u)=\src(u')$ and $\trgt(u)=\trgt(u')$.
    \item[-] $e$ is a subword of $u_2$, which is symmetric to the previous case.
    \end{itemize}

  \end{description}
  By definition of $\sim$, this suffices to show the desired property.
\end{proof}
\begin{lemma}\label{lemma:equivcong}
  Let $v_1, v_2, v_1', v_2' \in \T[\E]$ and $0 \leq k \leq n$ such that $v_1$
  and $v_2$ are $k$\nbd{}composable, and $v_1'$ and $v_2'$ are
  $k$\nbd{}composable. If $v_1 \sim v_2$ and $v_1' \sim v_2'$, then
  \[
    (v_1 \fcomp_k v_2) \sim (v_1' \fcomp_k v_2').
  \]
\end{lemma}
\begin{proof}
  Let
  \[
    \mu = (v,w,e,e') : v_1 \to v_1'
  \]
  be an elementary move. Set
  \[
    \tilde{v} := (v
  \]
  and
  \[
    \tilde{w} := w\fcomp_k v_2).
  \]
  Hence, $(\tilde{v},\tilde{w},e,e')$ is an elementary move from $(v_1 \fcomp_k
  v_2)$ to $(v_1' \fcomp_k v_2)$. Similarly, if we have an elementary move from
  $v_2$ to $v_2'$, we obtain an elementary move from $(v_1 \fcomp_k v_2)$ to
  $(v_1 \fcomp_k v_2')$. By definition of $\sim$, this suffices to show the desired property.
\end{proof}
\begin{lemma}\label{lemma:equivsmallestcong}
  Let $\R$ be a categorical congruence on $\E^+$ and $u,u' \in \T[\E]$. If $u
  \sim u'$, then $u \; \R \; u'$.
\end{lemma}
\begin{proof}
  The proof is very similar to the proof of Lemma \ref{lemma:equivparallel}. Let
  \[
    \mu =(v,w,e,e') : u \to u'
  \]
  be an elementary move from $u$ to $u'$. Let us prove that $u \; \R \; u'$ with
  an induction on $\mathcal{L}(v)+\mathcal{L}(w)$. Once again, by definition of
  elementary moves, we have $|u|\geq 1$ and thus
  \[
    u=(u_1 \fcomp_k u_2)
  \]
  with $u_1,u_2 \in \T[\E]$.
  \begin{description}
  \item[Base case] If $\mathcal{L}(v)+\mathcal{L}(w)=0$, then we have $u=e$ and
    $u'=e'$. In this case, all the different cases of elementary moves from
    Definition \ref{def:elementarymove} correspond to the different axioms of
    categorical congruences (Definition \ref{def:categoricalcongruence}). Hence,
    we have $u \; \R u'$.
  \item[Inductive step] Suppose now that $\mathcal{L}(v) + \mathcal{L}(w) \geq
    0$. Since $e$ is a subword of $u$, we know from Proposition
    \ref{prop:subwords} that we are in one of the following cases:
    \begin{itemize}[label=-]
    \item $e=u$, which is exactly the base case.
    \item $e$ is a subword of $u_1$, which means that there exist
      $\tilde{v},\tilde{w} \in \T[\E]$ such that
      \[
        u_1=\tilde{v}e\tilde{w}.
      \]
      Moreover, we have
      \[
        v=(\tilde{v}
      \]
      and
      \[
        w=\tilde{w}\fcomp_k u_2).
      \]
      From Proposition \ref{prop:subwordsubstitution}, the word
      \[
        u_1':=\tilde{v}e'\tilde{w}
      \]
      is well formed and since
      $\mathcal{L}(\tilde{v})+\mathcal{L}(\tilde{w})<\mathcal{L}(v)+\mathcal{L}(w)$,
      we can apply the induction hypothesis on
      \[
        \tilde{\mu} :=(\tilde{v},\tilde{w},e,e') : u_1 \to u_1'.
      \]
      This proves that $u_1 \; \R \; u_1'$. Now, since we have
      \[
        u =(u_1 \fcomp_k u_2) =u_1 \comp_k u_2
      \]
      and
      \[
        u' =(u_1' \fcomp_k u_2) =u_1' \comp_k u_2
      \]
      (see Remark \ref{remark:formalunformal} for the distinction between the
      symbol ``$\fcomp_k$'' and the symbol ``$\comp_k$'') and since $\R$ is a
      congruence, we have
      \[
        u \; \R \; u'.
      \]
    \item $e$ is a subword of $u_2$, which is symmetric to the previous
      case.\qedhere
    \end{itemize}
  \end{description}
\end{proof}
Altogether, Lemmas \ref{lemma:equivparallel}, \ref{lemma:equivcong} and
\ref{lemma:equivsmallestcong} prove the following proposition, which is what we
were aiming for.
\begin{proposition}\label{prop:equivsmallestcategoricalcongruence}
  Let $\E$ be an $n$\nbd{}cellular extension. The
  equivalence relation $\sim$ on $\T[\E]$ is the smallest categorical congruence
  on $\E^+$.
\end{proposition}
We end this section with yet another characterisation of $n$\nbd{}bases of
$\oo$\nbd{}categories.
\begin{paragr}\label{paragr:defrho}
  Let $C$ be an $\oo$\nbd{}category, $n \geq 0$ and $\Sigma$ a subset of
  $C_{n+1}$, and consider the $n$\nbd{}cellular extension
  \[
    \E_{\Sigma}=(\tau^{s}_{\leq n }(C),\Sigma,\src,\trgt)
  \]
  (see Paragraph \ref{paragr:cextfromsubset}). When there is no ambiguity on the
  rest of the data, we allow ourselves to write $\W[\Sigma]$, $\T[\Sigma]$ and
  $\G[\Sigma]$ instead of $\W[\E_{\Sigma}]$, $\T[\E_{\Sigma}]$ and
  $\G[\E_{\Sigma}]$.

  We recursively define a map $\rho_{\Sigma} : \T[\Sigma] \to C_{n+1}$ as
  \begin{itemize}[label=-]
  \item $\rho_{\Sigma}((\cc_{\alpha}))=\alpha$ for every $\alpha \in \Sigma$,
  \item $\rho_{\Sigma}((\ii_x))=1_x$ for every $x \in C_n$,
  \item $\rho_{\Sigma} ((w\fcomp_kw'))=\rho(w)\comp_k\rho(w')$ for every $0 \leq
    k \leq n$ and every pair $(w,w')$ of $k$\nbd{}composable elements of
    $\T[\Sigma]$.
  \end{itemize}
  Intuitively speaking, $\rho_{\Sigma}$ is to be understood as an ``evaluation
  map'': given $w$ a well formed expression on units and formal cells of
  $\Sigma$, $\rho_{\Sigma}(w)$ is the evaluation of $w$ as an $(n+1)$\nbd{}cell
  of $C$.
\end{paragr}
\begin{lemma}\label{lemma:rhosourcetarget}
  The map $\rho_{\Sigma}$ is compatible with source and target, i.e.\ for $w \in
  \T[\Sigma]$, we have
  \[
    \src(\rho_{\Sigma}(w))=\src(w) \text{ and }
    \trgt(\rho_{\Sigma}(w))=\trgt(w).
  \]
\end{lemma}
\begin{proof}
  This is proved by an immediate induction left to the reader.
\end{proof}
\begin{lemma}\label{lemma:rhoequiv}
  Let $w$ and $w'$ be elements of $\T[\Sigma]$. If $w \sim w'$, then
  \[
    \rho_{\Sigma}(w)=\rho_{\Sigma}(w').
  \]
\end{lemma}
\begin{proof}
  From Lemma \ref{lemma:rhosourcetarget} we know that $\rho_{\Sigma}$ is
  compatible with source and target and, by definition, $\rho_{\Sigma}$ is
  compatible with composition and units in an obvious sense. Hence, we have a
  morphism of $(n+1)$\nbd{}magmas
  \[
    F : \E_{\Sigma}^+ \to \tau_{\leq n+1}^s(C)
  \]
  defined as:
  \begin{itemize}[label=-]
  \item $F_{n+1}=\rho_{\Sigma}$,
  \item $F_k$ is the identity on $C_k$ for every $0 \leq k \leq n$.
  \end{itemize}
  Since $\tau_{\leq n+1}^s(C)$ is an $(n+1)$\nbd{}category, the binary relation
  $\R$ on $\T[\Sigma]$ defined as $w\; \R \; w'$ if
  $\rho_{\Sigma}(w)=\rho_{\Sigma}(w')$ is a categorical congruence on
  $\E_{\Sigma}$ (see Example \ref{example:categoricalcongruence}). The result
  follows then from \ref{prop:equivsmallestcategoricalcongruence}
\end{proof}
\begin{paragr}\label{paragr:rhofiber}
  Let $C$ be an $\oo$\nbd{}category, $\Sigma \subseteq C_{n+1}$ and let $a$ be
  an $(n+1)$\nbd{}cell of $C$. We define $\T[\Sigma]_a$ to be the set
  \[
    \T[\Sigma]_a := \{w\in\T[\Sigma]  \, \vert \,  \rho_{\Sigma}(w)=a\}.
  \]
  Lemma \ref{lemma:rhoequiv} implies that if $v \in \T[\Sigma]_a$ and $v \sim w$
  then $w \in \T[\Sigma]_a$.

  We define $\G[\Sigma]_a$ to be the full subgraph of $\G[\Sigma]$ whose set of
  objects is $\T[\Sigma]_a$.
\end{paragr}
\begin{proposition}\label{prop:basis0connected}
  Let $C$ be an $\oo$\nbd{}category, $n\geq 0$ and $\Sigma \subseteq C_{n+1}$.
  The set $\Sigma$ is an $(n+1)$\nbd{}basis of $C$ is and only if for every $a
  \in C_{n+1}$, the graph $\G[\Sigma]_a$ is $0$-connected (i.e.\ non-empty and
  connected).
  
  More precisely, this means that for every $a \in C_{n+1}$:
  \begin{itemize}[label=-]
  \item there exists $w \in \T[\Sigma]$ such that $\rho_{\Sigma}(w)=a$,
  \item for all $v,w \in \T[\Sigma]$, if $\rho_{\Sigma}(v)=a=\rho_{\Sigma}(w)$,
    then $v \sim w$.
  \end{itemize}
\end{proposition}
\begin{proof}
  From Lemma \ref{lemma:rhoequiv}, the map $\rho_{\Sigma}$ induces a map
  \[
    \overline{\rho_{\Sigma}} : \T[\Sigma]/{\sim} \to C_{n+1}.
  \]
  In fact, since $\sim$ is a categorical congruence, this map is the
  $(n+1)$\nbd{}dimensional part of the $(n+1)$\nbd{}functor
  \[
    \overline{F} : \E_{\Sigma}^+/{\sim} \to \tau_{\leq n+1}^s(C),
  \]
  where $F$ is the morphism of $(n+1)$\nbd{}magmas from the proof of Lemma
  \ref{lemma:rhoequiv}. Since $\sim$ is the smallest categorical congruence on
  $\E_{\Sigma}^+$, it follows from Proposition \ref{prop:Estarasquotient} that
  $\E_{\Sigma}^+/{\sim}$ is (canonically isomorphic to) $\E_{\Sigma}^*$ and it
  is easily seen that the functor $\overline{F}$ is nothing but the
  $(n+1)$\nbd{}functor $\E_{\Sigma}^* \to \tau_{\leq n+1}^s(C)$ as obtained in
  \ref{paragr:cextfromsubset}. Hence, from Proposition
  \ref{prop:criterionnbasis}, and the fact that $\overline{F}$ is the identity
  on cells of dimension not greater than $n$, we deduce that $\Sigma$ is an
  $(n+1)$\nbd{}basis of $C$ if and only if $\rho_{\Sigma}$ is an isomorphism.
  The result follows immediately from the fact that equivalence classes of
  $\sim$ are in bijection with maximal connected components of $\G[\Sigma]$.
\end{proof}
\begin{remark}\label{remark:weightexplicitly}
  Let $n>0$, $C$ be an $n$\nbd{}category with an $n$\nbd{}basis $\Sigma$ and
  $\alpha \in \Sigma$. It is immediate to check from the definition of
  elementary moves that for two equivalent well formed words $u \sim u'$ of
  $\T[\Sigma]$, the number of occurrences of $\cc_{\alpha}$ in $u$ and $u'$ are
  the same. In particular, for every $a \in C_n$ we can define the integer
  $w_{\alpha}(a)$ to be the number of occurrences of $\cc_{\alpha}$ in any well
  formed word $u$ such that $\rho_{\sigma}(u)=a$. An immediate induction using
  the properties of $\rho_{\Sigma}$ shows that this function $w_{\alpha} : C_n
  \to \mathbb{N}$ is the same as the one whose existence was established in
  Proposition \ref{prop:countingfunction}.
\end{remark}
\section{Proof of Theorem \ref{thm:conduche}: part III}
In this section, we finally go back to Conduché $\oo$\nbd{}functors. The first two parts might be considered as preliminaries and the key points of the proof of Theorem \ref{thm:conduche} lie within this third and last part.

\begin{paragr}
Let $F : C \to D$ be an $\omega$-functor, $n>0$, $\Sigma^C \subseteq C_n$ and $\Sigma^D \subseteq D_n$ such that $F_n(\Sigma^C)\subseteq \Sigma^D$. We recursively define a map 
\[
\widetilde{F} : \W[\Sigma^C] \to \W[\Sigma^D]
\]
with
\begin{itemize}
\item[-] $\wt{F}(\cc_{\alpha})=\cc_{F(\alpha)}$ for $\alpha \in \Sigma^C$,
\item[-] $\wt{F}(\ii_x)=\ii_{F(x)}$ for $x \in C_n$,
\item[-] $\wt{F}(\fcomp_k)=\fcomp_k$ for $0 \leq k < n$,
\item[-] $\wt{F}(\,(\,)=($,
  \item[-] $\wt{F}(\,)\,)=)$.
\end{itemize}
Notice that for every word $w \in \W[\Sigma^C]$, we have
\[\vert \wt{F}(w) \vert=\vert w \vert\text{ and }\mathcal{L}(\wt{F}(w))=\mathcal{L}(w).\]
\end{paragr}
\begin{lemma}\label{lemmamapinducedonwords}
 Let $F : C \to D$ be an $\omega$-functor and let $\Sigma^C \subseteq C_{n}$ and
 $\Sigma^D \subseteq D_{n}$ be such that $F_{n}(\Sigma^C)\subseteq \Sigma^D$. For every $u \in \W[\Sigma^C]$:
  \begin{enumerate}
  \item if $u$ is well formed then $\wt{F}(u)$ is well formed,
  \item if $\wt{F}(u)$ is well formed and if $u$ is a subword (\ref{def:subword}) of a well formed word then it is also well formed.
    \end{enumerate}
\end{lemma}
\begin{proof}
  The first part is proved with a short induction left to the reader. For the second part, first notice that the map
  \[
  \wt{F} : \W[\Sigma^C] \to \W[\Sigma^D]
  \]
  satisfies the following property: 
  \begin{center}
    For every $w \in \W[\Sigma^C]$, $w$ is well parenthesized if and only if $\wt{F}(w)$ is well parenthesized.
    \end{center}
It suffices then to apply Lemma \ref{lemmawellformedispseudo} and then Lemma \ref{lemmapartialconverse}.
\end{proof}
\begin{paragr}
The first part of Lemma \ref{lemmamapinducedonwords} shows that $\wt{F}$ induces a map 
\[
\wt{F} : \T[\Sigma^C] \to \T[\Sigma^D].
\]
Moreover, we have a commutative square
\[
\begin{tikzcd}
  \T[\Sigma^C] \ar[r,"\rho_{C}"] \ar[d,"\wt{F}"] &C_{n} \ar[d,"F_{n}"] \\
  \T[\Sigma^D] \ar[r,"\rho_{D}"] & D_{n}
  \end{tikzcd}
\]
where $\rho_C$ and $\rho_D$ respectively stand for $\rho_{\Sigma^C}$ and $\rho_{\Sigma^D}$.

Thus, for every $a \in C_{n}$ we can define a map:
\[
\begin{aligned}
  \wt{F}_a : \T[\Sigma^C]_a &\to \T[\Sigma^D]_{F(a)} \\
  w &\mapsto \wt{F}(w).
\end{aligned}
\]
\end{paragr}
Recall from Proposition \ref{prop:conduchepractical} that for an $\oo$\nbd{}functor $F : C \to D$ and $n\geq 0$, $\tau^s_{\leq n}(F)$ is a discrete Conduché $n$\nbd{}functor if and only if $F$ is right orthogonal to $\nabla^n_k$ for every $k \in \mathbb{N}$ such that $k<n$.
\begin{proposition}\label{propositionequivalentconditionsconduche}
  
  Let $F : C \to D$ be an $\omega$-functor and $n>0$. The following conditions are equivalent:
  \begin{enumerate}
    \item $\tau^s_{\leq n}(F): \tau^s_{\leq n}(C) \to \tau^s_{\leq n}(D)$ is a discrete Conduché $n$\nbd{}functor,
    \item for every $\Sigma^D \subseteq D_{n}$ and $\Sigma^C:=F^{-1}(\Sigma^D)$ and for every $a \in C_{n}$ the map
  \[
    \wt{F}_a : \T[\Sigma^C]_a \to \T[\Sigma^D]_{F(a)}
    \]
    defined above is bijective.

    \end{enumerate}
\end{proposition}
\begin{proof}
  We begin with $1 \Rightarrow 2$.
  \begin{description}
  \item[Surjectivity] Let us prove the following assertion:
    \[
    \begin{aligned}
      \forall l \in \mathbb{N}, \forall a \in C_{n}, &\forall w \in \T[\Sigma^D]_{F(a)} \text{ such that }|w|\leq l \\
      &\exists v \in \T[\Sigma^C]_a \text{ such that }\wt{F}_a(v)=w.
      \end{aligned}
    \]
    We proceed by induction on $l$.

    Suppose first that $l=0$. We are necessarily in one of two cases:
    \begin{enumerate}
    \item $w=(\cc_{\beta})$ with $\beta \in \Sigma^D$. By hypothesis,
      $\rho_D(w)=F(a)$
      and by definition of $\rho_{D}$,
      $\rho_D(w)=\beta$
      thus $F(a)=\beta$.
      By definition of $\Sigma^C$, $a\in \Sigma^C$ and we can choose
      $v =(\cc_a) $.
    \item $w=(\ii_{y})$ with $y \in D_{n-1}$.
      By hypothesis, $\rho_D(w)=F(a)$ and by definition of $\rho_D$, $\rho_D(w)=1_y$ thus $F(a)=1_y$. Since $\tau^s_{\leq n}(F)$ is a discrete Conduché $n$\nbd{}functor, $F$ is right orthogonal to $\kappa^{n}_{n-1}$. Hence, there exists $x \in C_{n-1}$ such that $a=1_x$ and $F(x)=y$. We can then choose $v=(\ii_x) \in \T[\Sigma^C]_a$.
      \end{enumerate}
    Now suppose that the assertion is true for a fixed $l \in \mathbb{N}$ and let $w \in \T[\Sigma^D]_{F(a)}$ be such that $|w|=l+1$.

    By definition of well formed words, we have
    \[
    w=(w_1\fcomp_k w_2)
    \]
    with $0 \leq k < n$ and  $w_1, w_2 \in \T[\Sigma^D]$ such that $|w_1|\leq l$ and $|w_2| \leq l$.

    By hypothesis, $\rho_D(w)=F(a)$ and by definition of $\rho_D$, \[\rho_D(w)=\rho_D(w_1)\comp_k \rho_D(w_2)\]
    and thus,
    \[\rho_D(w_1)\comp_k\rho_D(w_2)=F(a).\]
    Since by hypothesis $F$ is right orthogonal to $\nabla^n_k$, we know that
    there exist $a_1 \in C_{n}$ and $a_2 \in C_{n}$ that are $k$\nbd{}composable
    and such that
    \[a=a_1\comp_ka_2,\, F(a_1)=\rho_D(w_1) \text{ and }F(a_2)=\rho_D(w_2).\]

    Since $|w_1|\leq l$ and $|w_2| \leq l$, we can apply the induction
    hypothesis. Hence, there exist $v_1 \in \T[\Sigma^C]_{a_1}$ and $v_2 \in
    \T[\Sigma^C]_{a_2}$ such that
    \[
      \wt{F}_{a_1}(v_1)=\wt{F}(v_1)=w_1\text{ and }\wt{F}_{a_2}(v_2)=\wt{F}(v_2)=w_2.
    \]
    Since $\rho_C$ commutes with source and target by Lemma \ref{lemma:rhosourcetarget}, $v_1$ and $v_2$ are $k$\nbd{}composable and the word $(v_1\fcomp_kv_2)$ is well formed. By definition of $\rho_C$, we have
    \[\rho_C((v_1 \fcomp_k v_2))=\rho_C(v_1)\comp_k\rho_C(v_2)=a_1\comp_k a_2=a.\]
    Thus, $(v_1 \fcomp_k v_2) \in \T[\Sigma^C]_a$ and
    \[\wt{F}_a((v_1\fcomp_kv_2)=\wt{F}((v_1 \fcomp_k v_2))=(\wt{F}(v_1)\fcomp_k \wt{F}(v_2))=(w_1\fcomp_k w_2)=w.\]
  \item[Injectivity] Let us prove the following assertion:
\[
\forall l \in \mathbb{N},\forall v \in \T[\Sigma^C]_a, \forall w \in \T[\Sigma^C]_a \text { such that } |v|=|w|\leq l\]
\[
    \wt{F}_a(v)=\wt{F}_a(w) \Rightarrow v=w
\]
    We proceed by induction on $l$.

    Suppose first that $l=0$. We are necessarily in one of four cases:
    \begin{enumerate}
    \item $v=(\cc_{\alpha})$ and $w=(\cc_{\beta})$ with $\alpha$ and $\beta$ in
      $\Sigma^C$. By definition of $\rho_C$, we have $\alpha=\rho_C(v)=a=\rho_C(w)=\beta$. Hence, $v=w$.
    \item $v=(\ii_x)$ and $w=(\ii_y)$ with $x$ and $y$ in $C_{n-1}$. By
      hypothesis, we have  $\rho_C(v)=a=\rho_C(w)$ and by definition of $\rho_C$, $1_x=\rho_C(v)=a=\rho_C(w)=1_y$, thus $x=y$ and $v=w$.
    \item $v=(\cc_{\alpha})$ and $w=(\ii_x)$ with $\alpha \in \Sigma^C$ and $x \in C_{n-1}$. By hypothesis, $(\cc_{F(\alpha)})=\wt{F}(v)=\wt{F}(w)=(\ii_{F(x)})$ which is impossible.
      \item $v=(\ii_x)$ and $w=(\cc_{\alpha})$ with $\alpha \in \Sigma^C$ and $x \in C_{n-1}$, which is symmetric to the previous case.
    \end{enumerate}
    Now suppose that the assertion is true for a fixed $l\in \mathbb{N}$ and let $v,w \in \T[\Sigma^C]$ such that $|v|=|w|=l+1$ and $\wt{F}(v)=\wt{F}(w)$.
    By definition of well formed words, we have
    \[
    v=(v_1 \fcomp_k v_2)
    \]
    and
    \[
    w=(w_1 \fcomp_{k'} w_2)
    \]
    with $|v_1|,|v_2|,|w_1|,|w_2|\leq l$.
    
    By hypothesis, we have \[
    (\wt{F}(v_1)\fcomp_k\wt{F}(v_2))=\wt{F}(v)=\wt{F}(w)=(\wt{F}(w_1)\fcomp_{k'}\wt{F}(w_2)).
    \]
    From Proposition \ref{prop:unicitydecomposition}, we deduce that $\fcomp_k=\fcomp_{k'}$ and \[\wt{F}(v_j)=\wt{F}(w_j)\] for $j \in \{1,2\}$.

    In order to apply the induction hypothesis, we need to show that $\rho_C(v_j)=\rho_C(w_j)$ for $j \in \{1,2\}$.

    By hypothesis,
    \[\rho_C(v_1)\comp_k\rho_C(v_2)=\rho_C(v)=a=\rho_C(w)=\rho_C(w_1)\comp_k\rho_C(w_2).\]
    Hence,
    \[F(\rho_C(v_1))\comp_kF(\rho_C(v_2))=F(a)=F(\rho_C(w_1))\comp_kF(\rho_C(w_2)).\]
    Besides, $F(\rho_C(v_j))=\rho_D(\wt{F}(v_j))=\rho_D(\wt{F}(w_j))=F(\rho_C(w_j))$. We deduce from the fact that $F$ is right orthogonal to $\nabla^n_k$ that
    \[\rho_C(v_j)=\rho_C(w_j)\]
    for $j \in \{1,2\}$.

    From the induction hypothesis we have $v_j=w_j$ for $j \in \{1,2\}$, hence $v=w$.
  \end{description}
  Now let us prove $2 \Rightarrow 1$. 

  Let $a \in C_{n}$ and suppose that $F(a)=b_1\comp_k b_2$. We set
  $\Sigma^D=\{b_1,b_2\}$. By definition,  $((\cc_{b_1})\fcomp_k (\cc_{b_2})) \in \T[\Sigma^D]_{F(a)}$ and by hypothesis there exists a unique $v \in \T[\Sigma^C]_{a}$ such that $\wt{F}_a(v)=((\cc_{b_1})\fcomp_k(\cc_{b_2}))$. Since $|\wt{F}_a(v)|=|v|=1$, we have
  \[
  v=(v_1 \fcomp_{k'} v_2)
  \]
  with $|v_1|=|v_2|=0$, $\src_{k'}(v_1)=\trgt_{k'}(v_2)$ and $0 \leq k' < n$. Thus,
  \[
  (\wt{F}(v_1)\fcomp_{k'}\wt{F}(v_2))=\wt{F}(v)=((\cc_{b_1})\fcomp_k(\cc_{b_2})).
  \]
  Using Proposition \ref{prop:unicitydecomposition}, we deduce that $k=k'$ and $\wt{F}(v_j)=(\cc_{b_j})$ for $j\in \{1,2\}$.
  
  We set $a_1=\rho_C(v_1)$, $a_2=\rho_C(v_2)$ and we have $\src_k(a_1)=\trgt_k(a_2)$, \[a=\rho_C(v)=\rho_C(v_1)\comp_k \rho_C(v_2)=a_1\comp_k a_2\] and \[F(a_j)=F(\rho_C(v_j))=\rho_D(\wt{F}(v_j))=\rho_D(\cc_{b_j})=b_j\] for $j \in \{1,2\}$, which proves the existence part of the right orthogonality to $\nabla^n_k$.

  Now suppose that we have $a_1,a_1',a_2,a_2' \in C_{n}$
  with \[\src_k(a_1)=\trgt_k(a_2)\text{ and }\src_k(a_1')=\trgt_k(a_2'),\]
    \[a_1\comp_k a_2 =
    a_1'\comp_k a_2'=a,\] \[F(a_1)=F(a_1')=b_1 \text{ and }F(a_2)=F(a_2')=b_2.\]

  By definition of $\Sigma^C=F^{-1}(\Sigma^D)$, we have $a_1,a_1',a_2,a_2' \in
  \Sigma^C$. Let us define $w$ and $w'$ as
  \[w=((\cc_{a_1})\fcomp_k(\cc_{a_2}))\text{ and }w'=((\cc_{a_1'})\fcomp_k(\cc_{a_2'})).\] We have $\rho_C(w)=\rho_C(w')=a$ and $\wt{F}(w)=((\cc_{b_1})\fcomp_k(\cc_{b_2}))=\wt{F}(w')$. The injectivity of $\wt{F}_a$ implies that $w=w'$, hence $a_1=a_1'$ and $a_2=a_2'$, which proves the uniqueness part of the right orthogonality to $\nabla^n_k$. 
\end{proof}
\begin{paragr}\label{paragr:imageelemmove}
  Let $F : C \to D$ be an $\omega$-functor, $n > 0$, $\Sigma^C \subseteq C_{n}$ and $\Sigma^D \subseteq D_{n}$ such that $F(\Sigma^C)\subseteq \Sigma^D$. It follows from the definition of $\wt{F} : \T[\Sigma^C] \to \T[\Sigma^D]$ and the definition of elementary move (\ref{def:elementarymove}) that for an elementary move
  \[
  \mu = (v,w,e,e') : u \to u'
  \]
  with $u,u' \in \T[\Sigma^C]$, the quadruple
  \[
  (\wt{F}(v),\wt{F}(w),\wt{F}(e),\wt{F}(e'))
  \]
  is an elementary move from $\wt{F}(u)$ to $\wt{F}(u')$. Thus, we have defined a map
  \[
  \G[\Sigma^C](u,u') \to \G[\Sigma^D](\wt{F}(u),\wt{F}(u')).
  \]
  Together with the map $\wt{F} : \T[\Sigma^C] \to \T[\Sigma^D]$, this defines a morphism of graphs
  \[
  \wt{F} : \G[\Sigma^C] \to \G[\Sigma^D]
  \]
  and, by restriction, a morphism of graphs
    \[
  \wt{F}_a : \G[\Sigma^C]_a \to \G[\Sigma^D]_{F(a)}
  \]
  for every $a \in C_{n}$.
\end{paragr}
\begin{lemma}\label{lemmafaithful} With the notations of the above paragraph, the map
  \[
  \G[\Sigma^C](u,u') \to \G[\Sigma^D](\wt{F}(u),\wt{F}(u'))
  \]
  is injective.
\end{lemma}
  \begin{proof}
    Let $(v_1,w_1,e_1,e_1')$ and $(v_2,w_2,e_2,e_2')$ be two elementary moves from $u$ to $u'$ such that
  \[
   (\wt{F}(v_1),\wt{F}(w_1),\wt{F}(e_1),\wt{F}(e_1'))= (\wt{F}(v_2),\wt{F}(w_2),\wt{F}(e_2),\wt{F}(e_2')).
  \]
 In particular, we have
  \[
  \mathcal{L}(v_1)=\mathcal{L}(v_2) \text{ , }  \mathcal{L}(w_1)=\mathcal{L}(w_2) \text{ , }  \mathcal{L}(e_1)=\mathcal{L}(e_2) \text{ , }  \mathcal{L}(e_1')=\mathcal{L}(e_2'). 
  \]
  Since
  \[
  v_1e_1w_1=u=v_2e_2w_2 \text{ and } v_1e_1'w_1=u'=v_2e_2'w_2,
  \]
  we have
  \[
    v_1=v_2 \text{ , }  w_1=w_2 \text{ , }  e_1=e_2 \text{ , }  e_1'=e_2'. \qedhere
  \]
  \end{proof}
  \begin{lemma}\label{technicallemma}
    With the notations of paragraph \ref{paragr:imageelemmove}, suppose that $\tau^s_{\leq n-1}(F)$ is a discrete Conduché $(n-1)$-functor. Let
  \[
  \mu : v \to v'
  \]
  be an elementary move in $\T[\Sigma^D]$. If there exists $u \in \T[\Sigma^C]$ such that
  \[
  \wt{F}(u) = v,
  \]
  then there exists $u' \in \T[\Sigma^C]$ and an elementary move
  \[
  \lambda : u \to u'
  \]
  such that
  \[
  \wt{F}(u')=v' \text{ and } \wt{F}(\lambda)=\mu.
  \]
\end{lemma}
\begin{proof}
  The proof is long and tedious as we have to check all the different cases of elementary moves. For the sake of clarity, we first outline a sketch of the proof that is common to all the cases of elementary moves and then we proceed to fill in the blanks successively for each case.

 Let
  \[
  \mu =(v_1,v_2,e,e') : v \to v'
  \]
  be an elementary move. Since, by definition,
  \[\wt{F}(u)= v = v_1ev_2,\]
  $u$ is necessarily of the form
  \[
  u=u_1\overline{e}u_2
  \]
  with $\overline{e},u_1,u_2 \in \W[\Sigma^C]$ such that
  \[
  \wt{F}(\overline{e})=e
  \]
    and
    \[
     \wt{F}(u_j)=v_j
     \]
     for $j \in \{1,2\}$.
     From the second part of Lemma \ref{lemmamapinducedonwords}, we deduce that $\overline{e}$ is well formed. In each different case, we will prove the existence of a well formed word $\overline{e'}$ parallel to $\overline{e}$ such that
     \[
     \wt{F}(\overline{e'})=e'.
     \]
     From Proposition \ref{prop:subwordsubstitution}, we deduce that the word
     \[
     \overline{u'}:=u_1\overline{e'}u_2
     \]
     is well formed. By definition, we have
     \[
     \wt{F}(u')=v'.
     \]
     Moreover, in each case, it will be immediate that the pair $(\overline{e},\overline{e'})$ is such that the quadruple
     \[
     \lambda := (u_1,u_2,\overline{e},\overline{e'})
     \]
    is an elementary move and that
     \[
     \wt{F}(\lambda)=\mu.
     \]

  All that is left to prove now is the existence of $\overline{e'}$ with the desired properties.

    \begin{description}
    \item[First case:] The word $e$ is of the form
         \[((x\fcomp_k y )\fcomp_k z)\]  
      and the word $e'$ is of the form
       \[(x\fcomp_k (y \fcomp_k z))\]
      with $x,y,z \in \T[\Sigma^D]$.
      The word $\overline{e}$ is then necessarily of the form
      \[
      ((\overline{x}\fcomp_k \overline{y}) \fcomp_k \overline{z}).
      \]
      Since $\wt{F}(\overline{e})=e$, we deduce from Lemma \ref{lemmamapinducedonwords} that $\overline{x}$, $\overline{y}$, $\overline{z}$ and $(\overline{x} \fcomp_k \overline{y})$ are well formed. From Corollary \ref{cor:compatiblesubwords}, we deduce that
           \[
     \src_k(\overline{x})=\trgt_k(\overline{y})
     \]
     and
          \[
     \src_k(\overline{y})=\trgt_k(\overline{z}).
     \]
     Thus, the word
     \[
     \overline{e'}:=(\overline{x}\fcomp_k(\overline{y}\fcomp_k\overline{z}))
     \]
    is well formed and it satisfies the desired properties.
  \item[Second case:]
    The word $e$ is of the form
    \[
    (x\fcomp_k(\ii_{\1^{n-1}_z}))
    \]
    and the word $e'$ is of the form
    \[
    x
    \]
    with $x \in \T[\Sigma^D]$, $0\leq k < n$ and $z=\src_k(x)$.
    
    Necessarily $\overline{e}$ is of the form
    \[
     (\overline{x}\fcomp_k(\ii_{y}))
    \]
    with $\overline{x} \in \T[\Sigma^C]$ (from  Lemma \ref{lemmamapinducedonwords} again) and $y \in C_{n-1}$ such that
    \[
    \wt{F}(\overline{x})=x
    \]
    and
    \[
    \wt{F}(y)=\1^{n-1}_z.
    \]
    Then, we set
    \[
    \overline{e'}:=\overline{x}.
    \]
    The only thing left to show is that
    \[y=\1^{n-1}_{\src_k(\overline{x})}.\] If $k=n-1$, this follows from Corollary \ref{cor:compatiblesubwords} and the fact that $\overline{e}$ is well formed.
If $k <n-1$, we need first to use the fact that $f$ is right orthogonal to $\kappa^{n-1}_k$ to deduce that 
    \[
    y=\1^{n-1}_{\overline{z}}
    \]
    for some $\overline{z} \in C_k$ such that $F(\overline{z})=z$ and then use Corollary \ref{cor:compatiblesubwords} and the fact that $\overline{e}$ is well formed.
    \item[Third case:] Similar to the second one with the unit on the left.
  \item[Fourth case:] The word $e$ is of the form
    \[
    ((\ii_x) \fcomp_k (\ii_y))
    \]
    and the word $e'$ is of the form
    \[
    (\ii_{x\comp_k y})
    \]
    with $x,y \in D_{n-1}$ such that $\src_k(x)=\trgt_k(y)$.
    Necessarily, $\overline{e}$ is of the form
    \[
    ((\ii_{\overline{x}})\fcomp_k(\ii_{\overline{y}}))
      \]
      with $\overline{x},\overline{y} \in C_{n-1}$ such that
      \[
      F(\overline{x})=x \text{ and } F(\overline{y})=y.
      \]
      Using Corollary \ref{cor:compatiblesubwords} and the fact that $e$ is well formed, we deduce that $\src_k(\overline{x})=\trgt_k(\overline{y})$. Thus, the word
      \[
      \overline{e'}:=(\ii_{\overline{x} \comp_k \overline{y}})
      \]
      is well formed. It satisfies all the desired properties.
    \item[Fifth case:] The word $e$ is of the form
      \[
      ((x\fcomp_ky)\fcomp_l(z \fcomp_k t))
      \]
      and the word $e'$ is of the form
            \[
      ((x\fcomp_lz)\fcomp_k(y \fcomp_l t))
      \]
      with $x,y,z,t \in \T[\Sigma^D]$ and $0 \leq l  < k < n$ such that all the compatibilities of sources and targets needed are satisfied.  

      Necessarily, $\overline{e}$ is of the form
      \[
      ((\overline{x}\fcomp_k \overline{y})\fcomp_l (\overline{z} \fcomp_k \overline{t}))
      \]
      with $\overline{x},\overline{y},\overline{z},\overline{t} \in \W[\Sigma^C]$ such that
      \[
          \wt{F}(\overline{x})=x,
    \]
    \[
    \wt{F}(\overline{y})=y,
    \]
    \[
    \wt{F}(\overline{z})=z,
    \]
    \[
    \wt{F}(\overline{t})=t.
    \]
    From Lemma \ref{lemmamapinducedonwords} and the fact that $\overline{e}$ is
    well formed, we deduce that the words $\overline{x}$, $\overline{y}$,
    $\overline{z}$, $\overline{t}$,$(\overline{x}\fcomp_k \overline{y})$ and $(\overline{z} \fcomp_k \overline{t})$ are well formed and from Corollary \ref{cor:compatiblesubwords}, we deduce that
         \[
     \src_k(\overline{x})=\trgt_k(\overline{y}),
     \]
     \[
     \src_k(\overline{z})=\trgt_k(\overline{t})
     \]
     and
     \[
     \src_l((\overline{x}\fcomp_k\overline{y}))=\trgt_l((\overline{z}\fcomp_k\overline{t})).
     \]
     Since $l<k$, we deduce from this last equality that
     \[
     \src_l(\overline{x})=\src_l(\overline{y})=\trgt_l(\overline{z})=\trgt_l(\overline{t}).
     \]
     Thus, the word
     \[
     \overline{e'}:=((\overline{x}\fcomp_l\overline{z})\fcomp_k(\overline{y} \fcomp_l \overline{t}))
     \]
     is well formed. It satisfies all the desired properties.\qedhere
    \end{description}
  \end{proof}
  \begin{remark}
    In the proof of the previous theorem, we have only used the hypothesis that $F$ is right orthogonal to $\kappa_k^n$ for every $k$ such that $0 \leq k < n-1$.
    \end{remark}
\begin{lemma}\label{lemma:isomorphismgraphs}
  Let $F : C \to D$ be an $\omega$-functor, $n>0$, $\Sigma^D \subseteq D_{n}$ and $\Sigma^C := F^{-1}(\Sigma^D)$. If $\tau^s_{\leq n}(F)$ is a discrete Conduché $n$\nbd{}functor, then for every $a \in C_{n}$
  \[
  \wt{F}_a : \G[\Sigma^C]_a\to \G[\Sigma^D]_{F(a)}
  \]
  is an isomorphism of graphs.
\end{lemma}
\begin{proof}
 Proposition \ref{propositionequivalentconditionsconduche} exactly says that the map
    \[
  \wt{F}_a : \G[\Sigma^C]_a\to \G[\Sigma^D]_{F(a)}
  \]
  is an isomorphism on objects and we know from Lemma \ref{lemmafaithful} that it is a faithful morphism of graphs (same definition as for functors). All that is left to show is that it is also full.

  In other words, we have to show that for all $u,u' \in \T[\Sigma^C]_a$ the map
  \[
  \G[\Sigma^C](u,u') \to \G[\Sigma^D](\wt{F}(u),\wt{F}(u'))
  \]
  is surjective.

  Let $\mu : \wt{F}(u) \to \wt{F}(u')$ be an element of the codomain. From Lemma \ref{technicallemma} we know that there exists
  \[\lambda : u \to v\]
  in $\G[\Sigma^C]$ such that
  \[\wt{F}(\lambda)=\mu.\]
  In particular, we have
  \[
  \wt{F}(v)=\wt{F}(u').
  \]
  Since we have an elementary move from $u$ to $v$ and by hypothesis ${u \in \T[\Sigma^C]_a}$, we also have $v \in \T[\Sigma^C]_a$ (see \ref{paragr:rhofiber}). Using the injectivity of the map
  \[
  \wt{F}_a : \T[\Sigma^C]_a \to \T[\Sigma^D]_{F(a)}
  \]
  we conclude that $v=u'$.
\end{proof}
\begin{proposition}\label{prop:conduchenbasis}
  Let $F : C \to D$ be an $\omega$-functor, $n \in \mathbb{N}$, $\Sigma^D
  \subseteq D_{n}$ and define $\Sigma^C := F^{-1}(\Sigma^D)$. If $\tau_{\leq n}^s(F)$ is a discrete Conduché $n$\nbd{}functor, then:
  \begin{enumerate}
  \item if $\Sigma^D$ is an $n$\nbd{}basis then so is $\Sigma^C$,
    \item if $F_{n} : C_{n}\to D_{n}$ is surjective and $\Sigma^C$ is an $n$\nbd{}basis then so is $\Sigma^D$.
    \end{enumerate}
\end{proposition}
\begin{proof}The case $n=0$ is trivial. We now suppose that $n>0$.
  From Lemma \ref{lemma:isomorphismgraphs} we have that for every $a \in C_n$, the map
    \[
    \wt{F}_a : \G[\Sigma^C]_a \to \G[\Sigma^D]_{F(a)}
    \]
    is an isomorphism of graphs. In particular, $\G[\Sigma^C]_a$ is $0$-connected if and only if $\G[\Sigma^D]_{F(a)}$ is $0$-connected. We conclude with Proposition \ref{prop:basis0connected}.
    \end{proof}
Putting all the pieces together, we finally have the awaited proof.
\begin{paragr}[Proof of Theorem \ref{thm:conduche}]
    Let $F : C \to D$ be a discrete Conduché $\oo$\nbd{}functor.
    \begin{enumerate}
    \item In the case that $D$ is free, it follows immediately from the first part of the Proposition \ref{prop:conduchenbasis} that $C$ is free.
    \item In the case that $C$ is free and $F_n : C_n \to D_n$ is surjective for every $n \in \mathbb{N}$, let us write $\Sigma_n^C$ for the $n$\nbd{}basis of $C$. It follows from Proposition \ref{prop:uniquebasis} and Lemma \ref{lemma:conducheindecomposable} that
      \[
      F^{-1}(F(\Sigma^C_n))=\Sigma_n^C.
      \]
      Hence, we can apply the second part of Proposition \ref{prop:conduchenbasis} and $C$ is free.
      \end{enumerate}
  \end{paragr}



\chapter{Homotopical algebra}
The present chapter stands out from the others as it contains no original
results. Its goal is simply to introduce the language and tools of homotopical
algebra that we shall need in the rest of the dissertation. Consequently, most
of the results are simply asserted and the reader will find references to the
literature for the proofs. The main notion of homotopical algebra we aim for is
the one of \emph{homotopy colimits} and our language of choice is that of
Grothendieck's theory of \emph{derivators} \cite{grothendieckderivators}. We do
not assume that the reader is familiar with this theory and will quickly recall
the basics. If needed, gentle introductions can be found in
\cite{maltsiniotis2001introduction} and in a letter from Grothendieck to
Thomason \cite{grothendieck1991letter}; more detailed introductions can be found
in \cite{groth2013derivators} and in the first section of
\cite{cisinski2003images}; finally, a rather complete (yet unfinished and
unpublished) textbook on the subject is \cite{groth2013book}.

\section{Localization, derivation}

\begin{paragr}\label{paragr:loc}
  A \emph{localizer} is a pair $(\C,\W)$ where $\C$ is a category and $\W$ is a
  class of arrows of $\C$, which we usually refer to as the \emph{weak
    equivalences}. We denote by $\ho^{\W}(\C)$, or simply $\ho(\C)$ when there
  is no ambiguity, the localization of $\C$ with respect to $\W$ and by
  \[
    \gamma : \C \to \ho(\C)
  \]
  the localization functor \cite[1.1]{gabriel1967calculus}. Recall the universal
  property of the localization: for every category $\D$, the functor induced by
  pre-composition
  \[
    \gamma^* : \underline{\Hom}(\ho(\C),\D) \to \underline{\Hom}(\C,\D)
  \]
  is fully faithful and its essential image consists of those functors $F~:~\C~\to~\D$
  that send the morphisms of $\W$ to isomorphisms of $\D$.
  
  We shall always consider that $\C$ and $\ho(\C)$ have the same class of
  objects and implicitly use the equality
  \[
    \gamma(X)=X
  \]
  for every object $X$ of $\C$.

  The class of arrows $\W$ is said to be \emph{saturated} when we have the
  property:
  \[
    f \in \W \text{ if and only if } \gamma(f) \text{ is an isomorphism. }
  \]
\end{paragr}
For later reference, we put here the following definition.
\begin{definition}\label{def:couniversalwe}
  Let $(\C,\W)$ be a localizer such that $\C$ has amalgamated sums. A morphism
  $f : X \to Y$ in $\W$ is a \emph{co-universal weak equivalence} if for every
  cocartesian square of the form
  \[
    \begin{tikzcd}
      X \ar[r] \ar[d,"f"] & X' \ar[d,"f'"] \\
      Y \ar[r] & Y', \ar[from=1-1,to=2-2,phantom,"\ulcorner",very near end]
    \end{tikzcd}
  \]
  the morphism $f'$ is also a weak equivalence.
\end{definition}
\begin{paragr}
  A \emph{morphism of localizers} $F : (\C,\W) \to (\C',\W')$ is a functor
  $F:\C\to\C'$ that preserves weak equivalences, i.e.\ such that $F(\W)
  \subseteq \W'$. The universal property of the localization implies that $F$
  induces a canonical functor
  \[
    \overline{F} : \ho(\C) \to \ho(\C')
  \]
  such that the square
  \[
    \begin{tikzcd}
      \C \ar[r,"F"] \ar[d,"\gamma"] & \C' \ar[d,"\gamma'"]\\
      \ho(\C) \ar[r,"\overline{F}"] & \ho(\C').
    \end{tikzcd}
  \]
  is commutative. Let $G : (\C,\W) \to (\C',\W')$ be another morphism of
  localizers. A \emph{$2$\nbd{}morphism of localizers} from $F$ to $G$ is simply
  a natural transformation ${\alpha : F \Rightarrow G}$. The universal property of
  the localization implies that there exists a unique natural transformation
  \[
    \begin{tikzcd} \ho(\C) \ar[r,bend left,"\overline{F}",""{name=A,below}]
      \ar[r,bend right,"\overline{G}"',""{name=B,above}] & \ho(\C')
      \ar[from=A,to=B,Rightarrow,"\overline{\alpha}"]\end{tikzcd}
  \]
  such that the $2$\nbd{}diagram
  \[
    \begin{tikzcd}[row sep=huge]
      \C\ar[d,"\gamma"] \ar[r,bend left,"F",""{name=A,below}] \ar[r,bend right,"G"',""{name=B,above}] & \C'\ar[d,"\gamma'"] \ar[from=A,to=B,Rightarrow,"\alpha"] \\
      \ho(\C) \ar[r,bend left,"\overline{F}",""{name=A,below}] \ar[r,bend
      right,"\overline{G}"',""{name=B,above}] & \ho(\C')
      \ar[from=A,to=B,Rightarrow,"\overline{\alpha}"]
    \end{tikzcd}
  \]
  is commutative in an obvious sense.
\end{paragr}
\begin{remark}\label{remark:localizedfunctorobjects}
  Since we always consider that for every localizer $(\C,\W)$ the categories
  $\C$ and $\ho(\C)$ have the same class of objects and the localization functor
  is the identity on objects, it follows that for a morphism of localizers $F \colon
    (\C,\W) \to (\C',\W')$, we tautologically have
  \[
    \overline{F}(X)=F(X)
  \]
  for every object $X$ of $\C$.
\end{remark}
\begin{paragr}\label{paragr:defleftderived}
  Let $(\C,\W)$ and $(\C',\W')$ be two localizers. A functor $F : \C \to \C'$ is
  \emph{totally left derivable} when there exists a functor
  \[
    \LL F : \ho(\C) \to \ho(\C')
  \]
  and a natural transformation
  \[
    \alpha : \LL F \circ \gamma \Rightarrow \gamma'\circ F
  \]
  that makes $\LL F$ the \emph{right} Kan extension of $\gamma' \circ F$ along
  $\gamma$:
  \[
    \begin{tikzcd}
      \C \ar[r,"F"] \ar[d,"\gamma"] & \C' \ar[d,"\gamma'"]\\
      \ho(\C) \ar[r,"\LL F"'] & \ho(\C'). \arrow[from=2-1,
      to=1-2,"\alpha",Rightarrow]
    \end{tikzcd}
  \]
  When this right Kan extension is \emph{absolute}, we say that $F$ is
  \emph{absolutely totally left derivable}. When a functor $F$ is totally left
  derivable, the pair $(\LL F,\alpha)$ is unique up to a unique natural
  isomorphism and is referred to as \emph{the total left derived functor of
    $F$}. Often we will abusively discard $\alpha$ and simply refer to $\LL F$
  as the total left derived functor of $F$.

  The notion of \emph{total right derivable functor} is defined dually and
  denoted by $\RR F$ when it exists.
\end{paragr}
\begin{example}\label{rem:homotopicalisder}
  Let $F : (\C,\W) \to (\C',\W')$ a morphism of localizers. The universal
  property of the localization implies that $F$ is
  absolutely totally left and right derivable and $\LL F \simeq \RR F \simeq
  \overline{F}$.
\end{example}
To end this section, we recall a derivability criterion due to Gonzalez, which
we shall use in the sequel.
\begin{paragr}\label{paragr:prelimgonzalez}
  Let $(\C,\W)$ and $(\C',\W')$ be two localizers and let $\begin{tikzcd} F : \C
    \ar[r,shift left] & \C' \ar[l,shift left] : G \end{tikzcd}$ be an adjunction
  whose unit is denoted by $\eta$. Suppose that the functor $G$ is totally right derivable
  with $(\RR G,\beta)$ its total right derived functor and suppose that $\RR G$
  has a left adjoint $F' : \ho(\C) \to \ho(\C')$; the co-unit of this last
  adjunction being denoted by $\epsilon'$. All this data induces a natural
  transformation $\alpha : F' \circ \gamma \Rightarrow \gamma' \circ F$ defined
  as the following composition
  \[
    \begin{tikzcd}
      & \C' \ar[rr,"\gamma'"]\ar[rd,"G"] & &\ho(\C') \ar[rr,"\mathrm{id}",""{name=B,below}]\ar[rd,"\RR G"'] & &\ho(\C') \\
      \C\ar[ru,"F"] \ar[rr,"\mathrm{id}"',""{name=A,above}] && \C
      \ar[rr,"\gamma"'] &&\ho(\C)\ar[ru,"F'"'] &.
      \ar[from=A,to=1-2,"\eta",Rightarrow, shorten <= 0.5em, shorten >= 0.5em]
      \ar[from=2-3,to=1-4,Rightarrow,"\beta",shorten <= 1em, shorten >= 1em]
      \ar[from=2-5,to=B,Rightarrow,"\epsilon'"',shorten <= 0.5em, shorten >=
      0.5em]
    \end{tikzcd}
  \]
\end{paragr}
\begin{proposition}[{\cite[Theorem
    3.1]{gonzalez2012derivability}}]\label{prop:gonz}
  Let $(\C,\W)$ and $(\C',\W')$ be two localizers and
  \[\begin{tikzcd} F : \C \ar[r,shift left] & \C' \ar[l,shift left] :
      G \end{tikzcd}\] be an adjunction. If $G$ is absolutely totally right
  derivable with $(\RR G,\beta)$ its left derived functor and if $\RR G$ has a
  left adjoint $F'$
  \[\begin{tikzcd} F' : \ho(\C) \ar[r,shift left] & \ho(\C') \ar[l,shift left] :
      \RR G, \end{tikzcd}\] then $F$ is absolutely totally left derivable and
  the pair $(F', \alpha)$, with $\alpha$ defined as in the previous paragraph,
  is its left derived functor.
\end{proposition}

\section{(op-)Derivators and homotopy colimits}
\begin{notation}We denote by $\CCat$ the $2$\nbd{}category of small categories
  and $\CCAT$ the $2$\nbd{}category of big categories. For a $2$\nbd{}category
  $\underline{A}$, the $2$\nbd{}category obtained from $\underline{A}$ by
  switching the source and targets of $1$-cells is denoted by
  $\underline{A}^{\op}$.

  The terminal category, i.e.\ the category with only one object and no
  non-trivial arrows, is canonically denoted by $e$. For a (small) category $A$,
  the unique functor from $A$ to $e$ is denoted by
  \[
    p_A : A \to e.
  \]
\end{notation}
\begin{definition}
  An \emph{op\nbd{}prederivator} is a (strict) $2$\nbd{}functor
  \[\sD : \CCat^{\op} \to \CCAT.\]
  More explicitly, an op\nbd{}prederivator consists of the data of:
  \begin{itemize}[label=-]
  \item a big category $\sD(A)$ for every small category $A$,
  \item a functor $u^* : \sD(B) \to \sD(A)$ for every functor $u : A \to B$
    between small categories,
  \item a natural transformation
    \[
      \begin{tikzcd}
        \sD(B)\ar[r,bend left,"u^*",""{name=U,below}] \ar[r,bend
        right,"v^*"',""{name=D,above}] & \sD(A)
        \ar[from=U,to=D,Rightarrow,"\alpha^*"]
      \end{tikzcd}
    \]
    for every natural transformation
    \[
      \begin{tikzcd}
        A \ar[r,bend left,"u",""{name=U,below}] \ar[r,bend
        right,"v"',""{name=D,above}] & B \ar[from=U,to=D,Rightarrow,"\alpha"]
      \end{tikzcd}
    \]
    with $A$ and $B$ small categories,
  \end{itemize}
  compatible with compositions and units (in a strict sense).
\end{definition}
\begin{remark}
  Note that some authors call \emph{prederivator} what we have called
  \emph{op\nbd{}prederivator}. The terminology we chose in the above definition
  is compatible with the original one of Grothendieck, who called
  \emph{prederivator} a $2$\nbd{}functor from $\CCat$ to $\CCAT$ that is
  contravariant at the level of $1$-cells \emph{and} at the level of
  $2$\nbd{}cells.
\end{remark}
\begin{example}\label{ex:repder}
  Let $\C$ be a category. For a small category $A$, we use the notation $\C(A)$
  for the category $\underline{\Hom}(A,\C)$ of functors $A \to \C$ and natural
  transformations between them. The correspondence $A \mapsto \C(A)$ is
  $2$\nbd{}functorial in an obvious sense and thus defines an
  op\nbd{}prederivator
  \begin{align*}
    \C : \CCat^{\op} &\to \CCAT \\
    A &\mapsto \C(A)
  \end{align*}
  which we call the op\nbd{}prederivator \emph{represented by $\C$}. For $u : A
  \to B$ in $\CCat$,
  \[
    u^* : \C(B) \to \C(A)
  \]
  is simply the functor induced from $u$ by pre-composition.
\end{example}

We now turn to the most important way of obtaining op\nbd{}prederivators.
\begin{paragr}\label{paragr:homder}
  Let $(\C,\W)$ be a localizer. For every small category $A$, we write $\W_A$
  the class of \emph{pointwise weak equivalences} of the category $\C(A)$, i.e.\
  the class of arrows $\alpha : d \to d'$ of $\C(A)$ such that $\alpha_a : d(a)
  \to d'(a)$ belongs to $\W$ for every $a \in \Ob(A)$. This defines a localizer
  $(\C(A),\W_A)$. The correspondence $A \mapsto (\C(A),\W_A)$ is
  $2$\nbd{}functorial in that every $u : A \to B$ induces by pre-composition a
  morphism of localizers
  \[
    u^* : (\C(B),\W_B) \to (\C(A),\W_A)
  \]
  and every $\begin{tikzcd}A \ar[r,bend left,"u",""{name=A,below}] \ar[r,bend
    right, "v"',""{name=B,above}] & B
    \ar[from=A,to=B,Rightarrow,"\alpha"]\end{tikzcd}$ induces by pre-composition
  a $2$\nbd{}morphism of localizers
  \[
    \begin{tikzcd}
      (\C(B),\W_B) \ar[r,bend left,"u^*",""{name=A,below}] \ar[r,bend right,
      "v^*"',""{name=B,above}] & (\C(A),\W_A).
      \ar[from=A,to=B,Rightarrow,"\alpha^*"]
    \end{tikzcd}
  \]
  (This last property is trivial since a $2$\nbd{}morphism of localizers is
  simply a natural transformation between the underlying functors.) Then, by the
  universal property of the localization, every morphism $u : A \to B$ of $\Cat$
  induces a functor, again denoted by $u^*$,
  \[
    u^* : \ho(\C(B)) \to \ho(\C(A))
  \]
  and every natural transformation $\begin{tikzcd}A \ar[r,bend
    left,"u",""{name=A,below}] \ar[r,bend right, "v"',""{name=B,above}] & B
    \ar[from=A,to=B,Rightarrow,"\alpha"]\end{tikzcd}$ induces a natural
  transformation, again denoted by $\alpha^*$,
  \[
    \begin{tikzcd}
      \ho(\C(B)) \ar[r,bend left,"u^*",""{name=A,below}] \ar[r,bend right,
      "v^*"',""{name=B,above}] & \ho(\C(A)).
      \ar[from=A,to=B,Rightarrow,"\alpha^*"]
    \end{tikzcd}
  \]
  Altogether, this defines an op\nbd{}prederivator
  \begin{align*}
    \Ho^{\W}(\C) : \CCat^{\op} &\to \CCAT\\
    A &\mapsto \ho(\C(A)),
  \end{align*}
  which we call the \emph{homotopy op\nbd{}prederivator of $(\C,\W)$}. When
  there is no risk of confusion we will simply write $\Ho(\C)$ instead of
  $\Ho^{\W}(\C)$. All the op\nbd{}prederivators we shall work with arise this
  way. Notice that for the terminal category $e$, we have a canonical
  isomorphism
  \[
    \Ho(\C)(e)\simeq \ho(\C),
  \]
  which we shall use without further reference.
\end{paragr}
\begin{definition}
  An op\nbd{}prederivator $\sD$ has \emph{left Kan extensions} if for every $u :
  A \to B$ in $\Cat$, the functor $ u^* : \sD(B) \to \sD(A)$ has a left adjoint
  \[
    u_! : \sD(A) \to \sD(B).
  \]
\end{definition}
\begin{example}
  Let $\C$ be a category. The op\nbd{}prederivator represented by $\C$ has left
  Kan extensions if and only if the category $\C$ has left Kan extensions along
  every morphism $u : A \to B$ of $\Cat$ in the usual sense. By a standard
  categorical argument, this means that the op\nbd{}prederivator represented by
  $\C$ has left Kan extensions if and only if $\C$ is cocomplete. Note that for
  every small category $A$, the functor
  \[
    p_A^* : \C \simeq \C(e) \to \C(A)
  \]
  is nothing but the diagonal functor that sends an object $X$ of $\C$ to the
  constant diagram with value $X$. Hence, the functor $p_{A!}$ is nothing but
  the usual colimit functor of $A$-shaped diagrams
  \[
    p_{A!} = \colim_A : \C(A) \to \C(e) \simeq \C.
  \]
\end{example}
\begin{paragr}
  We say that a localizer $(\C,\W)$ has \emph{homotopy left Kan extensions} when
  the homotopy op\nbd{}prederivator of $(\C,\W)$ has left Kan extensions. In
  this case, for every small category $A$, the \emph{homotopy colimit functor of
    $A$-shaped diagrams} is defined as
  \[
    \hocolim_A := p_{A!} : \ho(\C(A)) \to \ho(\C).
  \]
  For an object $X$ of $\ho(\C(A))$ (which is nothing but a diagram $X : A \to
  \C$ seen ``up to weak equivalence''), the object of $\ho(\C)$
  \[
    \hocolim_A(X)
  \]
  is the \emph{homotopy colimit of $X$}. For consistency, we also use the
  notation
  \[
    \hocolim_{a \in A}X(a).
  \]
  When $\C$ is also cocomplete (which will always be the case in practice), it
  follows from Remark \ref{rem:homotopicalisder} and Proposition \ref{prop:gonz}
  that the functor \[ \colim_A : \C(A) \to \C
  \]
  is left derivable and $\hocolim_A$ is the left derived functor of $\colim_A$:
  \[
    \LL \colim_A \simeq \hocolim_A.
  \]
  In particular, for every $A$-shaped diagram $X : A \to \C$, there is a
  canonical morphism of $\ho(\C)$
  \[
    \hocolim_A(X) \to \colim_A(X).
  \]
  This canonical morphism will be of great importance in the sequel.
\end{paragr}
\begin{paragr}
  Let
  \[
    \begin{tikzcd}
      A \ar[r,"f"] \ar[d,"u"'] & B \ar[d,"v"]\\
      C \ar[r,"g"'] & D \ar[from=1-2,to=2-1,Rightarrow,"\alpha"]
    \end{tikzcd}
  \]
  be a $2$\nbd{}square in $\CCat$. Every op\nbd{}prederivator $\sD$ induces a
  $2$\nbd{}square:
  \[
    \begin{tikzcd}
      \sD(A) & \sD(B)  \ar[l,"f^*"'] \\
      \sD(C) \ar[u,"u^*"] & \sD(D). \ar[u,"v^*"']
      \ar[l,"g^*"]\ar[from=1-2,to=2-1,Rightarrow,"\alpha^*"]
    \end{tikzcd}
  \]
  If $\sD$ has left Kan extensions, we obtain a canonical natural transformation
  \[
    u_!f^* \Rightarrow g^*v_!
  \]
  referred to as the \emph{homological base change morphism induced by $\alpha$}
  and defined as the following composition:
  \[
    \begin{tikzcd}
      \sD(C) & \ar[l,"u_!"'] \sD(A) & \sD(B)  \ar[l,"f^*"'] \\
      & \sD(C) \ar[ul,"\mathrm{id}",""{name=A,above}] \ar[u,"u^*"] & \sD(D)
      \ar[u,"v^*"'] \ar[l,"g^*"]& \sD(B).
      \ar[l,"v_!"]\ar[ul,"\mathrm{id}"',""{name=B,below}]
      \ar[from=1-3,to=2-2,Rightarrow,"\alpha^*"]
      \ar[from=1-2,to=A,Rightarrow,"\epsilon"']
      \ar[from=B,to=2-3,Rightarrow,"\eta"]
    \end{tikzcd}
  \]
 
  In particular, let $u : A \to B$ be a morphism of $\CCat$ and $b$ an object of
  $B$ seen as a morphism $b :e \to B$. We have a square
  \[
    \begin{tikzcd}
      A/b \ar[r,"k"] \ar[d,"p"']& A \ar[d,"u"] \\
      e \ar[r,"b"'] & B \ar[from=1-2,to=2-1,Rightarrow,"\phi"]
    \end{tikzcd}
  \]
  where :
  \begin{itemize}[label=-]
  \item $A/b$ is the category whose objects are pairs $(a, f : u(a) \to b)$ with
    $a$ an object of $A$ and $f$ an arrow of $B$, and morphisms $(a,f) \to
    (a',f')$ are arrows $g : a \to a'$ of $A$ such that $f'\circ u(g) = f$,
  \item $k : A/b \to A$ is the functor $(a,p) \mapsto a$,
  \item $\phi$ is the natural transformation defined by $\phi_{(a,f)}:= f : u(a)
    \to b$.
  \end{itemize}
  Hence, we have a homological base change morphism:
  \[
    p_!\, k^* \Rightarrow b^*u_!.
  \]
  In the case that $\sD$ is the homotopy op\nbd{}prederivator of a localizer
  $(\C,\W)$, for every object $X$ of $\sD(A)$ the above morphism reads
  \[
    \hocolim_{A/b}(X\vert_{A/b}) \rightarrow u_!(X)_b
  \]
  where we use the notation $X\vert_{A/b}$ for $k^*(X)$ and $u_!(F)_b$ for
  $b^*(u_!(X))$. Note that this morphism is reminiscent of the formula that
  computes pointwise left Kan extensions in the ``classical'' sense (see for
  example \cite[chapter X, section
  3]{mac2013categories}). 
\end{paragr}
\begin{definition}[Grothendieck]
  A \emph{right op-derivator} is an op\nbd{}prederivator $\sD$ such that the
  following axioms are satisfied:
  \begin{description}
  \item[Der 1)] For every finite family $(A_i)_{i \in I}$ of small categories,
    the canonical functor
    \[
      \sD(\amalg_{i \in I}A_i) \to \prod_{i \in I}\sD(A_i)
    \]
    is an equivalence of categories. In particular, $\sD(\emptyset)$ is equivalent
    to the terminal category.
  \item[Der 2)]\label{der2} For every small category $A$, the functor
    \[
      \sD(A) \rightarrow \prod_{a \in \Ob(A)}\sD(e)
    \]
    induced by the functors $a^* : \sD(A)\to \sD(e)$ for all $a \in \Ob(A)$
    (seen as morphisms $a : e \to A$), is conservative.
  \item[Der 3d)] $\sD$ admits left Kan extensions.
  \item[Der 4d)] For every $u : A \to B$ in $\CCat$ and $b$ object of $B$, the
    homological base change morphism
    \[
      p_!\, k^* \Rightarrow b^*u_!
    \]
    induced by the square
    \[
      \begin{tikzcd}[column sep=small, row sep=small]
        A/b \ar[r,"k"] \ar[d,"p"']& A \ar[d,"u"] \\
        e \ar[r,"b"'] & B \ar[from=1-2,to=2-1,Rightarrow,"\phi"]
      \end{tikzcd}
    \]
    is an isomorphism.
  \end{description}
\end{definition}
\begin{paragr}
  Let us comment each of the axioms of the previous definition. Axiom
  \textbf{Der 1} ensures that $\sD(A)$ ``looks like'' a category of $A$-shaped
  diagrams. Axiom \textbf{Der 2} says that isomorphisms in $\sD(A)$ can be
  tested ``pointwise''. We have already seen that axiom \textbf{Der 3d} ensures
  the existence of left Kan extensions. Finally, axiom \textbf{Der 4d}
  intuitively says that ``Kan extensions are computed pointwise''.
\end{paragr}
  
  \begin{example}
    Let $\C$ be a category. The op\nbd{}prederivator represented by $\C$ always
    satisfy axioms \textbf{Der 1} and \textbf{Der 2}. We have already seen that
    axioms \textbf{Der 3d} means exactly that $\C$ admits left Kan extensions in
    the classical sense, in which case axiom \textbf{Der 4d} is automatically
    satisfied. Hence, the op\nbd{}prederivator represented by $\C$ is a right
    op\nbd{}prederivator if and only if $\C$ is cocomplete.
  \end{example}
  \begin{remark}
    Beware not to generalize the previous example too hastily. It is not true in
    general that axiom \textbf{Der 3d} implies axiom \textbf{Der 4d}; even in
    the case of the homotopy op\nbd{}prederivator of a localizer.
  \end{remark}
  This motivates the following definition.
  \begin{definition}\label{def:cocompletelocalizer}
    A localizer $(\C,\W)$ is \emph{homotopy cocomplete} if the
    op\nbd{}prederivator $\Ho(\C)$ is a right op-derivator.
  \end{definition}
  \begin{paragr}
    Axioms \textbf{Der 3d} and \textbf{Der 4d} can be dualized to obtain axioms
    \textbf{Der 3g} and \textbf{Der 4g}, which informally say that the
    op\nbd{}prederivator has right Kan extensions and that they are computed
    pointwise. An op\nbd{}prederivator satisfying axioms \textbf{Der 1},
    \textbf{Der 2}, \textbf{Der 3g} and \textbf{Der 4g} is a \emph{left
      op-derivator}. In fact, an op\nbd{}prederivator $\sD$ is a left
    op-derivator if and only if the op\nbd{}prederivator
    \begin{align*}
      \CCat &\to \CCAT \\
      A &\mapsto (\sD(A^{\op}))^{\op}
    \end{align*}
    is a right op\nbd{}prederivator. An op\nbd{}prederivator which is both a
    left and right op-derivator is an \emph{op-derivator}. For details, the
    reader can refer to any of the references on derivators previously cited.
  \end{paragr}
  \section{Morphisms of op-derivators, preservation of homotopy colimits}
  \sectionmark{Morphisms of op-derivators}
  We refer to \cite{leinster1998basic} for the precise definitions of
  pseudo-natural transformation (called strong transformation there) and
  modification.
  \begin{paragr}
    Let $\sD$ and $\sD'$ be two op\nbd{}prederivators. A \emph{morphism of
      op\nbd{}prederivators} $F : \sD \to \sD'$ is a pseudo-natural
    transformation from $\sD$ to $\sD'$. This means that $F$ consists of:
    \begin{itemize}[label=-]
    \item a functor $F_A : \sD(A) \to \sD'(A)$ for every small category $A$,
    \item an isomorphism of functors $F_u: F_A u^* \overset{\sim}{\Rightarrow}
      u^* F_B$,
      \[
        \begin{tikzcd}
          \sD(B) \ar[d,"u^*"] \ar[r,"F_B"] & \sD'(B) \ar[d,"u^*"]\\
          \sD(A) \ar[r,"F_A"'] & \sD'(A),
          \ar[from=2-1,to=1-2,Rightarrow,"F_u","\sim"']
        \end{tikzcd}
      \]
      for every $u : A \to B$ in $\CCat$.
    \end{itemize}
    These data being compatible with compositions and units. The morphism is
    \emph{strict} when $F_u$ is an identity for every $u : A \to B$.

    Let $F : \sD \to \sD'$ and $G : \sD \to \sD'$ be morphisms of
    op\nbd{}prederivators. A \emph{$2$\nbd{}morphism $\phi : F \Rightarrow G$}
    is a modification from $F$ to $G$. This means that $F$ consists of a natural
    transformation $\phi_A : F_A \Rightarrow G_A$ for every small category $A$,
    and is subject to a coherence axiom similar to the one for natural
    transformations.

    We denote by $\PPder$ the $2$\nbd{}category of op\nbd{}prederivators,
    morphisms of op\nbd{}prederivators and $2$\nbd{}morphisms of
    op\nbd{}prederivators.
  \end{paragr}

  \begin{example}
    Let $F : \C \to \C'$ be a functor. It induces a strict morphism at the level
    of op\nbd{}prederivators, again denoted by $F$, where for every small
    category $A$, the functor $F_A : \C(A) \to \C'(A)$ is induced by
    post-composition. Similarly, every natural transformation induces a
    $2$\nbd{}morphism at the level of represented op\nbd{}prederivators.
  \end{example}
  \begin{example}
    Let $F : (\C,\W) \to (\C',\W')$ be a morphism of localizers. For every small
    category $A$, the functor $F_A : \C(A) \to \C'(A)$ preserves weak
    equivalences and the universal property of the localization yields a
    functor \[\overline{F}_A : \ho(\C(A)) \to \ho(\C'(A)).\] This defines a
    strict morphism of op\nbd{}prederivators:
    \[
      \overline{F} : \Ho(\C) \to \Ho(\C').
    \]
    Similarly, every $2$\nbd{}morphism of localizers
    \[
      \begin{tikzcd}
        (\C,\W) \ar[r,bend left,"F",""{name=A,below}] \ar[r,bend right,
        "G"',""{name=B,above}] &(\C',\W') \ar[from=A,to=B,"\alpha",Rightarrow]
      \end{tikzcd}
    \]
    induces a $2$\nbd{}morphism $\overline{\alpha} : \overline{F} \Rightarrow
    \overline{G}$. Altogether, we have defined a $2$\nbd{}functor
    \begin{align*}
      \Loc &\to \PPder\\
      (\C,\W) &\mapsto \Ho(\C),
    \end{align*}
    where $\Loc$ is the $2$\nbd{}category of localizers.
  \end{example}
  \begin{paragr}\label{paragr:canmorphismcolimit}
    Let $\sD$ and $\sD'$ be op\nbd{}prederivators that admit left Kan extensions
    and let $F : \sD \to \sD'$ be a morphism of op\nbd{}prederivators. For every $u :
    A \to B$, there is a canonical natural transformation
    \[
      u_!\, F_A \Rightarrow F_B\, u_!
    \]
    defined as
    \[
      \begin{tikzcd}
        \sD(A) \ar[dr,"\mathrm{id}"',""{name=A,above}] \ar[r,"u_!"] &\sD(B) \ar[d,"u^*"] \ar[r,"F_B"]  & \sD'(B)\ar[d,"u^*"]  \ar[dr,"\mathrm{id}",""{name=B,below}]\\
        & \sD(A) \ar[r,"F_A"'] & \sD'(A) \ar[r,"u_!"'] & \sD'(B).
        \ar[from=2-2,to=1-3,Rightarrow,"F_u"',"\sim"]
        \ar[from=A,to=1-2,Rightarrow,"\eta"]
        \ar[from=2-3,to=B,Rightarrow,"\epsilon"]
      \end{tikzcd}
    \]
    For example, when $\sD$ is the homotopy op\nbd{}prederivator of a localizer
    and $B$ is the terminal category $e$, for every $X$ object of $\sD(A)$ the
    previous canonical morphism reads
    \[
      \hocolim_{A}(F_A(X))\to F_e(\hocolim_A(X)).
    \]
  \end{paragr}
  \begin{definition}\label{def:cocontinuous}
    Let $F : \sD \to \sD'$ be a morphism of op\nbd{}prederivators and suppose that
    $\sD$ and $\sD'$ both admit left Kan extensions. We say that $F$ is
    \emph{cocontinuous}\footnote{Some authors also say \emph{left exact}.} if for every $u: A \to B$, the
    canonical morphism
    \[
      u_! \, F_A \Rightarrow F_B \, u_!
    \]
    is an isomorphism.
  \end{definition}
  \begin{remark}
    When $\sD$ and $\sD'$ are homotopy op\nbd{}prederivators we will often say
    that a morphism $F : \sD \to \sD'$ is \emph{homotopy cocontinuous} instead
    of \emph{cocontinuous} to emphasize the fact that it preserves homotopy Kan
    extensions.
  \end{remark}
  \begin{example}
    Let $F : \C \to \C'$ be a functor and suppose that $\C$ and $\C'$ are
    cocomplete. The morphism induced by $F$ at the level of represented
    op\nbd{}prederivators is cocontinuous if and only if $F$ is cocontinuous in
    the usual sense.
  \end{example}
  \begin{paragr}\label{paragr:prederequivadjun}
    As in any $2$\nbd{}category, the notions of equivalence and adjunction make
    sense in $\PPder$. Precisely, we have that:
    \begin{itemize}
    \item[-] A morphism of op\nbd{}prederivators $F : \sD \to \sD'$ is an
      equivalence when there exists a morphism $G : \sD' \to \sD$ such that $FG$
      is isomorphic to $\mathrm{id}_{\sD'}$ and $GF$ is isomorphic to
      $\mathrm{id}_{\sD}$; the morphism $G$ is a \emph{quasi-inverse} of $F$.
    \item[-] A morphism of op\nbd{}prederivators $F : \sD \to \sD'$ is left
      adjoint to $G : \sD' \to \sD$ (and $G$ is right adjoint to $F$) if there
      exist $2$\nbd{}morphisms $\eta : \mathrm{id}_{\sD'} \Rightarrow GF$ and
      $\epsilon : FG \Rightarrow \mathrm{id}_{\sD}$ that satisfy the usual
      triangle identities.
    \end{itemize}
  \end{paragr}
  The following three lemmas are easy $2$\nbd{}categorical routine and are left
  to the reader.
  \begin{lemma}\label{lemma:dereq}
    Let $F : \sD \to \sD'$ be a morphism of op\nbd{}prederivators. If $F$ is an
    equivalence then $\sD$ is a right op-derivator (resp.\ left op-derivator,
    resp.\ op-derivator) if and only if $\sD'$ is one.
  \end{lemma}
  \begin{lemma}\label{lemma:eqisadj}
    Let $F : \sD \to \sD'$ be an equivalence and $G : \sD' \to \sD$ be a
    quasi-inverse of $G$. Then, $F$ is left adjoint to $G$.
  \end{lemma}
  \begin{lemma}\label{lemma:ladjcocontinuous}
    Let $\sD$ and $\sD'$ be op\nbd{}prederivators that admit left Kan extensions
    and let $F : \sD \to \sD'$ be a morphism of op\nbd{}prederivators. If $F$ is a left
    adjoint, then it is cocontinuous.
  \end{lemma}
  We end this section with a generalization of the notion of localization in the
  context of op\nbd{}prederivators.
  \begin{paragr}
    Let $(\C,\W)$ be a localizer. For every small category $A$, let
    \[
      \gamma_A : \C(A) \to \ho(\C(A))
    \]
    be the localization functor. The correspondence $A \mapsto \gamma_A$ is
    natural in $A$ and defines a strict morphism of
    op\nbd{}prederivators
    \[
      \gamma : \C \to \Ho(\C).
    \]
  \end{paragr}
  \begin{definition}\label{def:strnglyder}
    Let $(\C,\W)$ and $(\C',\W')$ be two localizers and $F : \C \to \C'$ a
    functor. We say that $F$ is \emph{strongly left derivable} if there exists a
    morphism of op\nbd{}prederivators
    \[
      \LL F : \Ho(\C) \to \Ho(\C')
    \]
    and a $2$\nbd{}morphism of op\nbd{}prederivators
    \[
      \begin{tikzcd}
        \C \ar[r,"F"] \ar[d,"\gamma"] & \C' \ar[d,"\gamma'"]\\
        \Ho(\C) \ar[r,"\LL F"'] & \Ho(\C'). \arrow[from=2-1,
        to=1-2,"\alpha",Rightarrow]
      \end{tikzcd}
    \]
    such that for every small category $A$, $((\LL F)_A,\alpha_A)$ is the
    \emph{absolute} total left derived functor of $F_A : \C(A) \to \C'(A)$. The
    pair $(\LL F, \alpha)$ is unique up to a unique isomorphism and is referred
    to as the \emph{left derived morphism of op\nbd{}prederivators of $F$}.
    Often, we will discard $\alpha$ and simply refer to $\LL F$ as the left
    derived morphism of $F$. The notion of \emph{strongly right derivable
      functor} is defined dually and the notation $\mathbb{R}F$ is used.
  \end{definition}
  \begin{example}
    Let $(\C,\W)$ and $(\C,\W')$ be localizers and $F : \C \to \C'$ a functor.
    If $F$ preserves weak equivalences (i.e.\ it is a morphism of localizers),
    then it is strongly left and right derivable and \[\overline{F} \simeq \LL F
      \simeq \RR G.\]
  \end{example}
  Gonzalez' criterion (Proposition \ref{prop:gonz}) admits the following
  generalization.
  \begin{proposition}\label{prop:gonzalezcritder}
    Let $(\C,\W)$ and $(\C',\W')$ be two localizers and
    \[
      \begin{tikzcd}
        F : \C \ar[r,shift left] & \C' : G \ar[l,shift left]
      \end{tikzcd}
    \]
    be an adjunction. If $G$ is strongly right derivable and if $\mathbb{R}G$
    has a left adjoint $F'$
    \[\begin{tikzcd} F' : \Ho(\C) \ar[r,shift left] & \Ho(\C') : \RR G,
        \ar[l,shift left] \end{tikzcd} \] then $F$ is strongly left derivable
    and
    \[
      \LL F \simeq F'.
    \]
  \end{proposition}
  \begin{proof}
    Let $\alpha : F' \circ \gamma \Rightarrow F\circ \gamma$ be the
    $2$\nbd{}morphism of op\nbd{}prederivators defined \emph{mutatis mutandis}
    as in \ref{paragr:prelimgonzalez} but at the level of op\nbd{}prederivators.
    Proposition \ref{prop:gonz} gives us that for every small category $A$, the
    functor $F_A$ is absolutely totally left derivable with $(F'_A,\alpha_A)$
    its total left derived functor. This means exactly that $F'$ is strongly
    left derivable and $(F',\alpha)$ is the left derived morphism of
    op\nbd{}prederivators of $F$.
  \end{proof}
  \section{Homotopy cocartesian squares}
  \begin{paragr}
    Let $\Delta_1$ be the ordered set $\{0 <1\}$ seen as category. We use the
    notation $\square$ for the category $\Delta_1\times \Delta_1$, which can be
    pictured as the \emph{commutative} square
    \[
      \begin{tikzcd}
        (0,0) \ar[r] \ar[d] & (0,1)\ar[d] \\
        (1,0) \ar[r] & (1,1)
      \end{tikzcd}
    \]
    and we use the notation $\ulcorner$ for the full subcategory of $\square$
    spanned by $(0,0)$, $(0,1)$ and $(1,0)$, which can be pictured as
    \[
      \begin{tikzcd}
        (0,0) \ar[d] \ar[r] & (0,1) \\
        (1,0) &.
      \end{tikzcd}
    \]
    Finally, we write $i_{\ulcorner} : \ulcorner \to \square$ for the canonical
    inclusion functor.
  \end{paragr}
  \begin{definition}\label{def:cocartesiansquare}
    Let $\sD$ be an op\nbd{}prederivator. An object $X$ of $\sD(\square)$ is
    \emph{cocartesian} if for every object $Y$ of $\sD(\square)$, the canonical
    map
    \[
      \Hom_{\sD(\square)}(X,Y) \to
      \Hom_{\sD(\ulcorner)}(i_{\ulcorner}^*(X),i_{\ulcorner}^*(Y))
    \]
    induced by the functor $i_{\ulcorner}^* : \sD(\square) \to \sD(\ulcorner)$,
    is a bijection.
  \end{definition}
  \begin{example}
    Let $\C$ be a category. An object of $\C(\square)$ is nothing but a
    commutative square in $\C$ and it is cocartesian in the sense of the
    previous definition if and only if it is cocartesian in the usual sense.
  \end{example}
  For the following definition to make sense, recall that for a localizer
  $(\C,\W)$ and a small category $A$, the objects of $\Ho(\C)(A)=\ho(\C(A))$ are
  identified with the objects of $\C(A)$ via the localization functor. In
  particular, an object of $\Ho(\C)(\square)$ is a commutative square of $\C$
  (up to weak equivalence).
  \begin{definition}\label{def:hmtpycocartesiansquare}
    Let $(\C,\W)$ be a localizer. A commutative square of $\C$ is said to be
    \emph{homotopy cocartesian} if it is cocartesian in $\Ho(\C)$ in the sense
    of Definition \ref{def:cocartesiansquare}.
  \end{definition}
  \begin{paragr}
    Let $\sD$ be an op\nbd{}prederivator. The object $(1,1)$ of $\square$ can be
    considered as a morphism of $\Cat$
    \[
      (1,1) : e \to \square
    \]
    and thus induces a functor $(1,1)^* : \sD(\square)\to \sD(e)$. For an object
    $X$ of $\sD(\square)$, we use the notation
    \[
      X_{(1,1)} := (1,1)^*(X).
    \]
    Now, since $(1,1)$ is the terminal object of $\square$, we have a canonical
    $2$\nbd{}triangle
    \[
      \begin{tikzcd}
        \ulcorner \ar[rr,"i_{\ulcorner}",""{name=A,below}] \ar[rd,"p"']  && \square\\
        &e\ar[ru,"{(1,1)}"']&, \ar[from=A,to=2-2,Rightarrow,"\alpha"]
      \end{tikzcd}
    \]
    where we wrote $p$ instead of $p_{\ulcorner}$ for short and where $\alpha$
    is the unique such natural transformation. Hence, we have a
    $2$\nbd{}triangle
    \[
      \begin{tikzcd}
        \sD(\ulcorner) && \sD(\square) \ar[ll,"i_{\ulcorner}^*"',""{name=A,below}] \ar[dl,"{(1,1)}^*"] \\
        & \sD(e) \ar[ul,"p^*"]&. \ar[from=A,to=2-2,Rightarrow,"\alpha^*"]
      \end{tikzcd}
    \]
    Suppose now that $\sD$ has left Kan extensions. For $X$ an object of
    $\sD(\square)$, we have a canonical morphism $p_!(i_{\ulcorner}^*(X)) \to
    X_{(1,1)}$ defined as the composition
    \[
      p_{!}(i_{\ulcorner}^*(X)) \to p_{!} p^* (X_{(1,1)}) \to X_{(1,1)},
    \]
    where the arrow on the left is induced by $\alpha^*$ and the arrow on the
    right is induced by the co-unit of the adjunction $p_{!} \dashv p^*$.

    When $\sD$ is the homotopy op\nbd{}prederivator of a localizer, and $X$ is a
    commutative square of $\C$
    \[
      X=
      \begin{tikzcd}
        A \ar[r,"u"] \ar[d,"f"]& B \ar[d,"g"] \\
        C \ar[r,"v"]&D,
      \end{tikzcd}
    \]
    this previous morphism reads
    \[
      \hocolim\left(\begin{tikzcd}[column sep=tiny, row sep=tiny] A \ar[d,"f"]
          \ar[r,"u"] & B \\ C & \end{tikzcd}\right) \to D.
    \]
  \end{paragr}
  \begin{proposition}
    Let $\sD$ be a right op\nbd{}prederivator. An object $X$ of $\sD(\square)$
    is cocartesian if and only if the canonical map $ p_!(i^*_{\ulcorner}(X))
    \to X_{(1,1)}$ is an isomorphism.
  \end{proposition}
  \begin{proof}
    Let $Y$ be another object of $\sD(\square)$. Using the adjunction
    $i_{\ulcorner!} \dashv i_{\ulcorner}^*$, the canonical map
    $\Hom_{\sD(\square)}(X,Y) \to
    \Hom_{\sD(\ulcorner)}(i_{\ulcorner}^*(X),i_{\ulcorner}^*(Y))$ can be
    identified with
    \[
      \Hom_{\sD(\square)}(X,Y) \to \Hom_{\sD(\ulcorner)}(i_{\ulcorner
        !}(i_{\ulcorner}^*(X)),Y).
    \]
    Hence, $X$ is cocartesian if and only if the co-unit map $X \to i_{\ulcorner
      !}(i_{\ulcorner}^*(X))$ is an isomorphism. The rest follows then from
    \cite[Lemma 9.2.2(i)]{groth2013book}.
  \end{proof}
  Hence, for a homotopy cocomplete localizer $(\C,\W)$, a commutative square of
  $\C$ is homotopy
  cocartesian if and only if the bottom right apex of the square is the homotopy
  colimit of the upper left corner of the square. This hopefully justifies the
  terminology of ``cocartesian square''.

  The previous proposition admits the following immediate corollary.
  \begin{corollary}
    Let $\sD$ and $\sD'$ be right op-derivators and $F : \sD \to \sD'$ a
    morphism of op\nbd{}prederivators. If $F$ is cocontinuous, then it preserves
    cocartesian squares. This means that if an object $X$ of $\sD(\square)$ is
    cocartesian, then $F_{\square}(X)$ is a cocartesian square of $\sD'$.
  \end{corollary}
  We end this section with two useful lemmas which show that homotopy
  cocartesian squares behave much like ``classical'' cocartesian squares.
  \begin{lemma}\label{lemma:hmtpycocartsquarewe}
    Let $(\C,\W)$ be a homotopy cocomplete localizer and let
    \[
      \begin{tikzcd}
        A \ar[r,"u"] \ar[d,"f"]& B \ar[d,"g"] \\
        C \ar[r,"v"]&D
      \end{tikzcd}
    \]
    be a commutative square in $\C$. If $f$ and $g$ are weak equivalences then
    the previous square is homotopy cocartesian.
  \end{lemma}
  \begin{proof}
    Using \textbf{Der 2}, one can show that the previous square is isomorphic in
    $\Ho(\C)(\square)$ to the square
    \[
      \begin{tikzcd}
        A \ar[r,"u"] \ar[d,"1_A"]& B \ar[d,"1_B"] \\
        A \ar[r,"u"]&B.
      \end{tikzcd}
    \]
    The result follows then from \cite[Proposition
    3.12(2)]{groth2013derivators}.
  \end{proof}
  \begin{lemma}[Pasting lemma]\label{lemma:pastinghmtpycocartesian}
    Let $(\C,\W)$ be a homotopy cocomplete localizer and let
    \[
      \begin{tikzcd}
        A \ar[r,"u"] \ar[d,"f"]& B \ar[d,"g"] \ar[r,"v"] & C \ar[d,"h"]\\
        D \ar[r,"w"]&E \ar[r,"x"] & F
      \end{tikzcd}
    \]
    be a commutative diagram in $\C$. If the square on the left is cocartesian,
    then the outer square is cocartesian if and only if the right square is
    cocartesian.
  \end{lemma}
  \begin{proof}
    This is a particular case of \cite[Proposition
    3.13(1)]{groth2013derivators}.
  \end{proof}
  \section{Model categories}
  In this section, we quickly review some aspects of the relation between
  Quillen's theory of model categories and Grothendieck's theory of derivators.
  We suppose that the reader is familiar with the former one and refer to the
  standard textbooks on the subject (such as \cite{hovey2007model,hirschhorn2009model,dwyer1995homotopy}) for basic definitions
  and results.

  For a model category $\M = (\M,\W,\Cof,\Fib)$, the homotopy
  op\nbd{}prederivator of $\M$, denoted by $\Ho(\M)$, is the homotopy
  op\nbd{}prederivator of the localizer $(\M,\W)$.

  The following theorems are due to Cisinski \cite{cisinski2003images} and can
  be summed up by the slogan:
  \begin{center}
    Model categories are homotopy cocomplete and left Quillen functors are
    homotopy cocontinuous.
  \end{center}
  \begin{theorem}[Cisinski]\label{thm:cisinskiI}
    Let $(\M,\W,\Cof,\Fib)$ be a model category. The localizer $(\M,\W)$ is
    homotopy cocomplete.
  \end{theorem}
  
  \begin{theorem}[Cisinski]\label{thm:cisinskiII}
    Let $\M$ and $\M'$ be model categories. Let $F \colon \M \to \M'$ be a left
    Quillen functor (i.e.\ the left adjoint in a Quillen adjunction). The
    functor $F$ is strongly left derivable and the morphism of
    op\nbd{}prederivators \[\LL F \colon \Ho(\M) \to \Ho(\M')\] is homotopy
    cocontinuous.
  \end{theorem}
  \begin{remark}
    The obvious duals of the two above theorems are also true. The reason we put
    emphasis on cocompleteness rather that completeness is because we will make
    no use whatsoever of homotopy limits in this dissertation.
  \end{remark}
  \begin{remark}
    Note that for a model category
    $(\M,\W,\Cof,\Fib)$, its homotopy op\nbd{}prederivator only depends on its underlying localizer. Hence, the existence
    of the classes $\Cof$ and $\Fib$ with the usual properties defining a model
    structure ought to be thought of as a \emph{property} of the localizer
    $(\M,\W)$, which is sufficient to define a ``homotopy theory''. For
    example, Theorem \ref{thm:cisinskiI} should have been stated by saying that
    if a localizer $(\M,\W)$ can be extended to a model category
    $(\M,\W,\Cof,\Fib)$, then it is homotopy cocomplete.
  \end{remark}

  Even if Theorem \ref{thm:cisinskiI} tells us that (the homotopy
  op\nbd{}prederivator of) a model category $(\M,\W,\Cof,\Fib)$ has homotopy
  left Kan extensions, it is not generally true that for a small category $A$
  the category of diagrams $\M(A)$ admits a model structure with the pointwise
  weak equivalences as its weak equivalences. Hence, in general we cannot use
  the theory of Quillen functors to compute homotopy left Kan extensions (and in
  particular homotopy colimits). However, in practice all the model categories
  that we shall encounter are \emph{cofibrantly generated}, in which case the
  theory is much simpler because $\M(A)$ does admit a model structure with the
  pointwise weak equivalences as its weak equivalences.
  \begin{paragr}\label{paragr:cofprojms}
    Let $\C$ be a category with coproducts and $A$ a small category. For every
    object $X$ of $\C$ and every object $a$ of $A$, we define $X\otimes a$ as
    the functor
    \[
      \begin{aligned}
        X \otimes a : A &\to \C \\
        b &\mapsto \coprod_{\Hom_A(a,b)}X.
      \end{aligned}
    \]
    For every object $a$ of $A$, this gives rise to a functor
    \[
      \begin{aligned}
        \shortminus \otimes a : \C &\to \C(A)\\
        X &\mapsto X \otimes a.
      \end{aligned}
    \]
  \end{paragr}
    
  \begin{proposition}\label{prop:modprs}
    Let $(\M,\W,\Cof,\Fib)$ be a cofibrantly generated model category with $I$
    (resp.\ $J$) as a set of generating cofibrations (resp. trivial
    cofibrations). For every small category $A$, there exists a model structure
    on $\M(A)$ such that:
    \begin{itemize}[label=-]
    \item the weak equivalences are the pointwise weak equivalences,
    \item the fibrations are the pointwise fibrations,
    \item the cofibrations are those morphisms which have the left lifting property
      to trivial fibrations. \end{itemize} Moreover, this model structure is
    cofibrantly generated and a set of generating cofibrations (resp.\ trivial
    cofibrations) is given by
    \[
      \{ f \otimes a : X \otimes a \to Y \otimes a \quad \vert \quad a \in
      \Ob(A), f \in I\}
    \]
    (resp.
    \[
      \{ f \otimes a : X \otimes a \to Y \otimes a \quad \vert \quad a \in
      \Ob(A), f \in J\}).
    \]
    
  \end{proposition}
  \begin{proof}
    See for example \cite[Proposition 3.4]{cisinski2003images}.
  \end{proof}
  \begin{paragr}\label{paragr:projmod}
    The model structure of the previous proposition is referred to as the
    \emph{projective model structure on $\M(A)$}.
  \end{paragr}
  \begin{proposition}
    Let $(\M,\W,\Cof,\Fib)$ be a cofibrantly generated model category. For
    every $u : A \to B$, the adjunction
    \[
      \begin{tikzcd}
        u_! : \M(A) \ar[r,shift left]& \ar[l,shift left] \M(B) : u^*
      \end{tikzcd}
    \]
    is a Quillen adjunction with respect to the projective model structures on
    $\M(A)$ and on $\M(B)$.
  \end{proposition}
  \begin{proof}
    By definition of the projective model structure, $u^*$ preserve weak
    equivalences and fibrations.
  \end{proof}
  \begin{paragr}\label{paragr:hocolimms}
    In particular, in the case that $B$ is the terminal category $e$, we have
    that
    \[
      \begin{tikzcd}
        \colim_A : \M(A) \ar[r,shift left] & \ar[l,shift left] \M(e) \simeq \M :
        p_A^*
      \end{tikzcd}
    \]
    is a Quillen adjunction. Since $\hocolim_A$ is the left derived functor of
    $\colim_A$, we obtain the following immediate corollary of the previous
    proposition.
  \end{paragr}
  \begin{corollary}\label{cor:cofprojms}
    Let $(\M,\W,\Cof,\Fib)$ be a cofibrantly generated model category, $A$ a
    small category and $X : A \to \M$ a diagram. If $X$ is cofibrant for the
    projective model structure on $\C(A)$, then the canonical morphism of
    $\ho(\C)$
    \[
      \hocolim_A(X) \to \colim_A(X)
    \]
    is an isomorphism.
  \end{corollary}
  \begin{proof}
    This is simply a particular case of the general fact that for a left Quillen
    functor $F : \M \to \M'$ and a cofibrant object $X$ of $\M$, the canonical
    map of $\ho(\M')$
    \[
      \LL F(X) \to F(X)
    \]
    is an isomorphism.
  \end{proof}

  Below is a particular case for which the previous corollary applies.
  \begin{proposition}\label{prop:sequentialhmtpycolimit}
    Let $(\M,\W,\Cof,\Fib)$ be a cofibrantly generated model category and let
    $X$ be a sequential diagram in $\M$
    \[
      X_0 \to X_1 \to X_2 \to \cdots
    \]
    (i.e.\ a diagram $X : (\mathbb{N},\leq) \to \M$). If $X_0$ is cofibrant and
    each $X_i \to X_{i+1}$ is a cofibration, then $X$ is cofibrant for the
    projective model structure on $\M((\mathbb{N},\leq))$.
  \end{proposition}
  \begin{proof}
    It is an easy exercise that uses only the fact that cofibrations of the
    projective model structure are, by definition, the morphisms with left
    lifting property to pointwise fibrations. For details see \cite[Example
    2.3.16]{schreiber2013differential}.
  \end{proof}
  Another setting for which a category of diagrams $\M(A)$ can be equipped with
  a model structure whose weak equivalences are the pointwise equivalences and for
  which the $A$-colimit functor is left Quillen is when the category $A$ is a
  \emph{Reedy category}. Rather that recalling this theory, we simply put here
  the only practical result that we shall need in the sequel.
    
    \begin{lemma}\label{lemma:hmtpycocartesianreedy}
      Let $(\M,\W,\Cof,\Fib)$ be a model category and let
      \[
        \begin{tikzcd}
          A \ar[r,"u"] \ar[d,"f"]& B \ar[d,"g"] \\
          C \ar[r,"v"]&D \ar[from=1-1,to=2-2,phantom,"\ulcorner",very near end]
        \end{tikzcd}
      \]
      be a \emph{cocartesian} square in $\M$. If either $u$ or $f$ is a
      cofibration and if $A$, $B$ and $C$ are cofibrant objects, then this
      square is \emph{homotopy cocartesian}.
    \end{lemma}
    \begin{proof}
      See for example \cite[Proposition A.2.4.4(i)]{lurie2009higher}.
    \end{proof}


\chapter{Homotopy theory of \texorpdfstring{$\oo$}{ω}-categories}
\section{Nerve}
\begin{paragr}\label{paragr:simpset}
  We denote by $\Delta$ the category whose objects are the finite non-empty totally ordered sets $[n]=\{0<\cdots<n\}$ and whose morphisms are the non-decreasing maps. For $n > 0$ and $0\leq i\leq n$, we denote by
  \[
  \delta^i : [n-1] \to [n] 
  \]
  the only injective increasing map whose image does not contain $i$, and for $n\geq 0$ and $0 \leq i \leq n$, we denote by
  \[
  \sigma^i : [n+1] \to [n]
  \]
  the only surjective non-decreasing map such that the pre-image of $i \in [n]$ contains exactly two elements.
  The category $\Psh{\Delta}$ of \emph{simplicial sets} is the category of presheaves on $\Delta$. For a simplicial set $X$, we use the notations
  \[
  \begin{aligned}
    X_n &:= X([n]) \\
    \partial_i &:= X(\delta^i): X_n \to X_{n\shortminus 1}\\
    s_i &:= X(\sigma^i): X_{n} \to X_{n+1}.
  \end{aligned}
  \]
  Elements of $X_n$ are referred to as \emph{$n$\nbd{}simplices of $X$}, the maps $\partial_i$ are the \emph{face maps} and the maps $s_i$ are the \emph{degeneracy maps}. 
\end{paragr}
\begin{paragr}\label{paragr:orientals}
  We denote by $\Or : \Delta \to \omega\Cat $ the cosimplicial object introduced by Street in \cite{street1987algebra}. The $\omega$-category $\Or_n$ is called the \emph{$n$\nbd{}oriental}. There are various ways to give a precise definition of the orientals, but each of them requires some machinery that we don't want to introduce here. Instead, we only recall some important facts on orientals that we shall need in the sequel and refer to the literature on the subject (such as \cite{street1987algebra,street1991parity,street1994parity,steiner2004omega,buckley2016orientals} or \cite[chapitre 7]{ara2016joint}) for details.

  The two main points to retain are:
  \begin{description}
  \item[(OR1)] Each $\Or_n$ is a free $\oo$\nbd{}category whose set of generating
    $k$\nbd{}cells is canonically isomorphic to the sets of increasing sequences
    \[
    0 \leq i_1 < i_2 < \cdots < i_k \leq n,
    \]
    or, which is equivalent, to injective increasing maps $[k] \to [n]$.
    \end{description}
  We use the notation $\langle i_1\, i_2\cdots i_k\rangle$ for such a cell. In particular, we have that:
  \begin{itemize}[label=-]
  \item There are no generating $k$-cells for $k>n$. Hence, $\Or_n$ is an $n$\nbd{}category.
  \item There is exactly one generating $n$\nbd{}cell of $\Or_n$, which is $\langle 0 \,1 \cdots n\rangle$. We refer to this cell as the \emph{principal cell of $\Or_n$}.
  \item There are exactly $n+1$ generating $(n-1)$-cells of $\Or_n$. They correspond to the maps
    \[
    \delta^i : [n-1] \to [n]
    \]
    for $i \in \{0,\cdots,n\}$.
  \end{itemize}
  \begin{description}
  \item[(OR2)] For $n>0$, the source (resp.\ target) of the principal cell of $\Or_n$ can be expressed as a composition of all the generating $(n-1)$\nbd{}cells corresponding to $\delta^i$ with $i$ odd (resp.\ even); each of these generating $(n-1)$\nbd{}cells appearing exactly once in the composite.
  \end{description}
  Another way of formulating \textbf{(OR2)} is: for $n>0$ the weight (\ref{paragr:weight}) of the $(n-1)$\nbd{}cell corresponding to $\delta_i$ in the \emph{source} of the principal cell of $\Or_n$ is $1$ if $i$ is odd and $0$ if $i$ is even and the other way around for the \emph{target} of the principal cell of $\Or_n$.
  Here are some pictures in low dimension:
  \[
  \Or_0 = \langle 0 \rangle,
  \]
  \[
  \Or_1=\begin{tikzcd}
    \langle 0 \rangle \ar[r,"\langle 01 \rangle"] &\langle 1 \rangle,
    \end{tikzcd}
  \]
  \[
  \Or_2=
  \begin{tikzcd}
    &\langle 1 \rangle \ar[rd,"\langle 12 \rangle"]& \\
    \langle 0 \rangle \ar[ru,"\langle 01 \rangle"]\ar[rr,"\langle 02 \rangle"',""{name=A,above}]&&\langle 2 \rangle,
    \ar[Rightarrow,from=A,to=1-2,"\langle 012 \rangle"]
    \end{tikzcd}
  \]
  \[
  \Or_3=
  \begin{tikzcd}
    & \langle 1 \rangle \ar[rd,"\langle 12 \rangle"]& \\
    \langle 0 \rangle \ar[ru,"\langle 01 \rangle"] \ar[rd,"\langle 03 \rangle"',""{name=B,above}] \ar[rr,"\langle 02 \rangle" description,""{name=A,above}]& & \langle 2 \rangle  \ar[ld,"\langle 23 \rangle"]\\
    & \langle 3 \rangle &
    \ar[from=A,to=1-2,Rightarrow,"\langle 012 \rangle", shorten <= 0.25em, shorten >= 0.25em]
    \ar[from=B,to=2-3,Rightarrow,"\langle 023 \rangle"', near start, shorten <= 1.1em, shorten >= 1.5em]
  \end{tikzcd}
  \overset{\langle 0123 \rangle}{\Rrightarrow}
    \begin{tikzcd}
    & \langle 1 \rangle \ar[rd,"\langle 12 \rangle"] \ar[dd,"\langle 13 \rangle"' description,""{name=B,right}] & \\
    \langle 0 \rangle \ar[ru,"\langle 01 \rangle"] \ar[rd,"\langle 03 \rangle"',""{name=A,above}] & & \langle 2 \rangle.  \ar[ld,"\langle 23 \rangle"]\\
    & \langle 3 \rangle &
    \ar[from=A,to=1-2,Rightarrow,"\langle 013 \rangle", near start, shorten <= 1em, shorten >= 1.5em]
    \ar[from=B,to=2-3,Rightarrow,"\langle 123 \rangle", shorten <= 0.75em, shorten >=0.75em]
    \end{tikzcd}
    \]
\end{paragr}
\begin{paragr}\label{paragr:nerve}
    For every $\omega$-category $C$, the \emph{nerve of $C$} is the simplicial set $N_{\omega}(C)$ defined as
    \[
    \begin{aligned}
     N_{\omega}(C) : \Delta^{\op} &\to \Set\\
      [n] &\mapsto \Hom_{\omega\Cat}(\Or_n,C).
      \end{aligned}
    \]
   By post-composition, this yields a functor
  \[
  \begin{aligned}
  N_{\omega} : \omega\Cat &\to \Psh{\Delta} \\
  C &\mapsto N_{\omega}(C),
  \end{aligned}
  \]
  which we refer to as the \emph{nerve functor for $\oo$\nbd{}categories}. Furthermore, for every $n \in \mathbb{N}$, we also define a nerve functor for $n$\nbd{}categories as the restriction of $N_{\oo}$ to $n\Cat$ (seen as a full subcategory of $\oo\Cat$)
  \[
  N_n := N_{\oo}{\big |}_{n\Cat} : n\Cat \to \Psh{\Delta}.
  \]
  By the usual Kan extension technique, we obtain for every $n \in \nbar$ a functor \[c_n : \Psh{\Delta} \to n\Cat,\] left adjoint to $N_n$.
\end{paragr}
  \begin{paragr}
    For $n=1$, the functor $N_1$ is the usual nerve of categories. Recall that for a (small) category $C$, an $m$-simplex $X$ of $N_1(C)$ is a sequence of composable arrows of $C$
    \[
    X_0 \overset{X_{0,1}}{\longrightarrow} X_1 \overset{X_{1,2}}{\longrightarrow} \cdots \longrightarrow X_{m-1}\overset{X_{m-1,m}}{\longrightarrow} X_m.
        \]
        For $m > 0$ and $0 \leq i \leq m$, the $(m-1)$-simplex $\partial_i(X)$ is obtained by composing arrows at $X_i$ (or simply deleting it for $i=0$ or $m$). For $m \geq 0$ and $0 \leq i \leq m$, the $(m+1)$-simplex $s_i(X)$ is obtained by inserting a unit map at $X_i$.
        
        For $n=2$, the functor $N_2$ is what is sometimes known as the \emph{Duskin nerve} \cite{duskin2002simplicial} (restricted from bicategories to $2$-categories). For a $2$-category $C$, an $m$\nbd{}simplex $X$ of $N_2(C)$ consists of:
  \begin{itemize}[label=-]
  \item for every $0\leq i \leq m$, an object $X_i$ of $C$,
  \item for all $0\leq i \leq  j \leq m$, an arrow $X_{i,j} : X_i \to X_j$ of $C$,
  \item for all $0 \leq i \leq  j \leq k \leq m$, a $2$-triangle
    \[
      \begin{tikzcd}
    & X_j \ar[rd,"{X_{j,k}}"]& \\
    X_i \ar[ru,"X_{i,j}"]\ar[rr,"X_{i,k}"',""{name=A,above}]&&X_k,
    \ar[Rightarrow,from=A,to=1-2,"X_{i,j,k}"]
    \end{tikzcd}
  \]
  \end{itemize}
  subject to the following axioms:
  \begin{itemize}[label=-]
  \item for all $0 \leq i \leq m$, we have
    \[
    X_{i,i}=1_{X_i},
    \]
  \item for all $0 \leq i \leq j \leq m$, we have
    \[
    X_{i,i,j}=X_{i,j,j}=1_{X_{i,j}},
    \]
  \item for all $0 \leq i < j < k < l \leq m$, we have the equality (known as the \emph{cocycle condition})
    \[
    (X_{k,l} \comp_0 X_{i,j,k})\comp_1 X_{i,k,l} = (X_{j,k,l} \comp_0 X_{i,j})\comp_1 X_{i,j,l}.
    \]
  \end{itemize}
  For $ m> 0$ and $0\leq l \leq m$, the $(m-1)$-simplex $\partial_l(X)$ is defined as
  \[
  \partial_l(X)_{i}=X_{\delta_l(j)}, \quad \partial_l(X)_{i,j}=X_{\delta_l(i),\delta_l(j)} \text{ and } \partial_l(X)_{i,j,k}=X_{\delta_l(i),\delta_l(j),\delta_l(k)}.
  \]
  And similarly, for $m \geq 0$ and $0\leq l \leq m$, the $(m+1)$-simplex $s_l(X)$ is defined as
  \[
  s_l(X)_{i}=X_{\sigma_l(j)}, \quad s_l(X)_{i,j}=X_{\sigma_l(i),\sigma_l(j)} \text{ and } s_l(X)_{i,j,k}=X_{\sigma_l(i),\sigma_l(j),\sigma_l(k)}.
  \]
\end{paragr}
  \section{Thomason equivalences}
  \begin{paragr}
From now on, we will consider that the category $\Psh{\Delta}$ is equipped with the model structure defined by Quillen in \cite{quillen1967homotopical}. A \emph{weak equivalence of simplicial sets} is a weak equivalence for this model structure. The cofibrations for this model structure are the monomorphisms.
  \end{paragr}
  \begin{definition}
   Let $n \in \nbar$. A morphism $f : X \to Y$ of $n\Cat$ is a \emph{Thomason equivalence} when ${N_n(f) : N_n(X) \to N_n(Y)}$ is a weak equivalence of simplicial sets. We denote by $\W_n^{\mathrm{Th}}$ the class of Thomason equivalences.
  \end{definition}
  \begin{paragr}\label{paragr:notationthom}
    We usually make reference to the name ``Thomason'' in the notations of homotopic constructions induced by Thomason equivalences. For example, we write $\Ho(n\Cat^{\Th})$ for the homotopy op\nbd{}prederivator of $(n\Cat,\W_n^{\Th})$ and
    \[
    \gamma^{\Th} : n\Cat \to \Ho(n\Cat^{\Th}) 
    \]
    for the localization morphism. The reason is to avoid confusion with other weak equivalences on $n\Cat$ that we will introduce later.
  \end{paragr}
  \begin{paragr}
    By definition, the nerve functor induces a morphism of localizers

    \[
      {N_n : (n\Cat,\W_n^{\Th}) \to (\Psh{\Delta},\W_{\Delta})}
    \]
    and hence a morphism of op\nbd{}prederivators
    \[
    \overline{N_n} : \Ho(n\Cat^{\Th}) \to \Ho(\Psh{\Delta}).
    \]
  \end{paragr}
  \begin{theorem}[Gagna]\label{thm:gagna}
    For every $1 \leq n \leq \oo$, the morphism \[{\overline{N}_n : \Ho(n\Cat^\Th) \to \Ho(\Psh{\Delta})}\] is an equivalence of op\nbd{}prederivators.
  \end{theorem}
  \begin{proof}
    Recall from \ref{paragr:nerve} that $c_n : \Psh{\Delta} \to n\Cat$ denotes
    the left adjoint of the nerve functor $N_n$. In \cite{gagna2018strict}, Gagna proves that there exists a functor $Q : \Psh{\Delta} \to \Psh{\Delta}$, as well as a zigzag of morphisms of functors
    \[
    N_{n}c_{n}Q \overset{\alpha}{\longleftarrow} Q \overset{\gamma}{\longrightarrow} \mathrm{id}_{\Psh{\Delta}}
    \]
    and a morphism of functors
    \[
    c_{n}Q N_{n} \overset{\beta}{\longrightarrow} \mathrm{id}_{\Psh{\Delta}},
    \]
    such that $c_{n}Q$ preserves weak equivalences and $\alpha$, $\beta$ and $\gamma$ are weak equivalences argument by argument. This easily implies that
       \[
    \overline{c_{n}Q} : \Ho(\Psh{\Delta})\to \Ho(n\Cat^{\Th})
    \]
    is a quasi-inverse (\ref{paragr:prederequivadjun}) of
    \[
    \overline{N_n} : \Ho(n\Cat^{\Th}) \to \Ho(\Psh{\Delta}).\qedhere
    \]
  \end{proof}
  From Lemma \ref{lemma:dereq}, we obtain the following corollary.
  \begin{corollary}\label{cor:thomhmtpycocomplete}
    For every $1 \leq n \leq \oo$, the localizer $(n\Cat^{\Th},\W_n^{\Th})$ is homotopy cocomplete (Definition \ref{def:cocompletelocalizer}).
  \end{corollary}
  We will speak of ``Thomason homotopy colimits'' and ``Thomason homotopy
  cocartesian squares'' for homotopy colimits and homotopy cocartesian squares in
  the localizer $(n\Cat^{\Th},\W_n^{\Th})$. (See also \ref{paragr:thomhmtpycol} below.)

  Another consequence of Gagna's theorem is the following
  corollary.
  \begin{corollary}\label{cor:thomsaturated}
    For every $1 \leq n \leq \oo$, the class $\W_n^{\Th}$ is saturated (\ref{paragr:loc}).
  \end{corollary}
  \begin{proof}
    This follows immediately from the fact that $\overline{N_n} : \ho(n\Cat^{\Th}) \to \ho(\Psh{\Delta})$ is an equivalence of categories and the fact that weak equivalences of simplicial sets are saturated (because they are the weak equivalences of a model structure). 
    \end{proof}
  \begin{remark}
    Corollaries \ref{cor:thomhmtpycocomplete} and \ref{cor:thomsaturated} would also follow from the existence of a model structure on $n\Cat$ with $\W^{\Th}_n$ as the weak equivalences. For $n=1$, this was established by Thomason \cite{thomason1980cat}, and for $n=2$, by Ara and Maltsiniotis \cite{ara2014vers}. For $n>3$, the existence of such a model structure is conjectured but not yet established.
  \end{remark}
   By definition, for all $1 \leq n \leq m \leq \omega$, the canonical inclusion \[n\Cat \hookrightarrow m\Cat\] sends the Thomason equivalences of $n\Cat$ to Thomason equivalences of $m\Cat$. Hence, it induces a morphism of localizers and then a morphism of op\nbd{}prederivators \[\Ho(n\Cat^\Th) \to \Ho(m\Cat^{\Th}).\]
  \begin{proposition}\label{prop:nthomeqder}
    For all $1 \leq n \leq m \leq \omega$, the canonical morphism
    \[
    \Ho(n\Cat^\Th) \to \Ho(m\Cat^{\Th})
    \]
  is an equivalence of op\nbd{}prederivators.
  \end{proposition}
  \begin{proof}
    This follows from Theorem \ref{thm:gagna} and the commutativity of the triangle
    \[
    \begin{tikzcd}[column sep=tiny]
      \Ho(n\Cat^{\Th}) \ar[rr] \ar[rd,"\overline{N_n}"'] & & \Ho(m\Cat^{\Th}) \ar[dl,"\overline{N_m}"] \\
      &\Ho(\Psh{\Delta})&.
    \end{tikzcd}
    \]
  \end{proof}
  \begin{paragr}\label{paragr:thomhmtpycol}
    It follows from the previous proposition that for all $1 \leq n \leq m \leq
    \omega$, the morphism $\Ho(n\Cat^\Th) \to \Ho(m\Cat^{\Th})$ of
    op\nbd{}prederivators is homotopy cocontinuous and reflects homotopy
    colimits (in an obvious sense). Hence, given a diagram $d : I
    \to n\Cat$ with $n>0$, we can harmlessly use the notation
    \[
      \hocolim^{\Th}_{i \in I}(d)
    \]
    for both the Thomason homotopy colimits in $n\Cat$ and in $\oo\Cat$ (or any
    $m\Cat$ with $n\leq m$). Similarly, a commutative square of
    $n\Cat$ is Thomason homotopy cocartesian in $n\Cat$ if and
    only if it is so in $\oo\Cat$. Hence, there is really no ambiguity when simply
    calling such a square \emph{Thomason homotopy cocartesian}.
  \end{paragr}
   \section{Tensor product and oplax transformations}
     Recall that $\oo\Cat$ can be equipped with a monoidal product $\otimes$, introduced by Al-Agl and Steiner in \cite{al1993nerves} and by Crans in \cite{crans1995combinatorial}, commonly referred to as the \emph{Gray tensor product}. The implicit reference for this section is \cite[Appendices A and B]{ara2016joint}. 
     \begin{paragr}
       The Gray tensor product makes $\oo\Cat$ into a monoidal category for which the unit is the $\oo$\nbd{}category $\sD_0$ (which is the terminal $\oo$\nbd{}category). This monoidal category is \emph{not} symmetric but it is biclosed \cite[Theorem A.15]{ara2016joint}, meaning that there exist two functors
       \[
       \underline{\hom}_{\mathrm{oplax}}(-,-),\, \underline{\hom}_{\mathrm{lax}}(-,-) : \oo\Cat^{\op}\times\oo\Cat \to \oo\Cat
       \]
       such that for all $\oo$\nbd{}categories $X,Y$ and $Z$, we have isomorphisms
       \begin{align*}
         \Hom_{\oo\Cat}(X\otimes Y , Z) &\simeq \Hom_{\oo\Cat}(X,  \underline{\hom}_{\mathrm{oplax}}(Y,Z))\\
         &\simeq  \Hom_{\oo\Cat}(Y,  \underline{\hom}_{\mathrm{lax}}(X,Z))
       \end{align*}
       natural in $X,Y$ and $Z$. When $X=\sD_0$, we have $\sD_0 \otimes Y \simeq Y$, and thus
       \[
       \Hom_{\oo\Cat}(Y,Z)\simeq \Hom_{\oo\Cat}(\sD_0,\underline{\hom}_{\mathrm{oplax}}(Y,Z)).
       \]
       Hence, the $0$-cells of the $\oo$\nbd{}category $\underline{\hom}_{\mathrm{oplax}}(Y,Z)$ are the $\oo$\nbd{}functors $Y \to Z$.
     \end{paragr}
     \begin{paragr}
       Let $u,v : X \to Y$ be two $\oo$\nbd{}functors. An \emph{oplax transformation} from $u$ to $v$ is a $1$-cell $\alpha$ of $\homoplax(X,Y)$ with source $u$ and target $v$. We usually use the double arrow notation \[
       \alpha : u \Rightarrow v
       \]
       for oplax transformations. By adjunction, we have
       \begin{align*}
         \Hom_{\oo\Cat}(\sD_1,\homoplax(X,Y)) &\simeq \Hom_{\oo\Cat}(\sD_1\otimes X , Y)\\
         &\simeq \Hom_{\oo\Cat}(X,\homlax(\sD_1,Y)).
       \end{align*}
       Hence, $\alpha : u \Rightarrow v$ can be encoded in the following two ways:
       \begin{itemize}[label=-]
         \item As an $\oo$\nbd{}functor $\alpha : \sD_1\otimes X \to Y$ such that the following diagram
       \[
       \begin{tikzcd}
         X\ar[rd,"u"] \ar[d,"i_0^X"']& \\
         \sD_1\otimes X \ar[r,"\alpha"] & Y, \\
         X \ar[ru,"v"'] \ar[u,"i_1^X"]&
       \end{tikzcd}
       \]
       where $i_0^X$ and $i_1^X$ are induced by the two $\oo$\nbd{}functors $\sD_0 \to \sD_1$ and where we implicitly used the isomorphism $\sD_0 \otimes X \simeq X$, is commutative.
     \item As an $\oo$\nbd{}functor $\alpha : X \to \homlax(\sD_1,Y)$ such that the following diagram
       \[
       \begin{tikzcd}
         & Y \\
         X \ar[ru,"u"] \ar[r,"\alpha"] \ar[rd,"v"']& \homlax(\sD_1,Y) \ar[u,"\pi_0^Y"'] \ar[d,"\pi_1^Y"] \\
         & Y,
       \end{tikzcd}
       \]
        where $\pi^Y_0$ and $\pi^Y_1$ are induced by the two $\oo$\nbd{}functors $\sD_0 \to \sD_1$ and where we implicitly used the isomorphism $\homlax(\sD_0,Y)\simeq Y$, is commutative.
       \end{itemize}
       The $\oo$\nbd{}category $\homlax(\sD_1,Y)$ is sometimes referred to as the $\oo$\nbd{}category of cylinders in $Y$. An explicit description of this $\oo$\nbd{}category can be found, for example, in \cite[Appendix A]{metayer2003resolutions}, \cite[Section 4]{lafont2009polygraphic} or \cite[Appendice B.1]{ara2016joint}. 
     \end{paragr}

     \begin{paragr}\label{paragr:formulasoplax}[Formulas for oplax
       transformations] We now give a third way of describing oplax
       transformations based on explicit formulas. The proof that this description is equivalent to those given in the previous paragraph can be found in \cite[Appendice B.2]{ara2016joint}.

       Let $u, v : X \to Y$ two $\oo$\nbd{}functors. An oplax transformation $\alpha : u \Rightarrow v$ is given by the data of:
       \begin{itemize}[label=-]
          \item for every $0$-cell $x$ of $X$, a $1$-cell of $Y$
      \[
      \alpha_x : u(x) \to v(x),
      \]
      \item for every $n$\nbd{}cell of $x$ of $X$ with $n>0$, an $(n+1)$-cell of $Y$
    \[
    \alpha_x : \alpha_{\trgt_{n-1}(x)}\comp_{n-1}\cdots\comp_1\alpha_{\trgt_0(x)}\comp_0u(x) \to v(x)\comp_0\alpha_{\src_0(x)}\comp_1\cdots\comp_{n-1}\alpha_{\src_{n-1}(x)}
    \]
    subject to the following axioms:
    \begin{enumerate}
    \item for every $n$\nbd{}cell $x$ of $X$,
      \[\alpha_{1_x}=1_{\alpha_x},\]
    \item for all $0\leq k < n$, for all $n$\nbd{}cells $x$ and $y$ of $X$ that are $k$-composable,
      \[
      \begin{multlined}
      \alpha_{x \comp_k y}={\left(v(\trgt_{k+1}(x))\comp_0\alpha_{\src_0(y)}\comp_1\cdots\comp_{k-1}\alpha_{\src_{k-1}(y)}\comp_k\alpha_y\right)}\\
            {\comp_{k+1}\left(\alpha_x \comp_k\alpha_{\trgt_{k-1}(x)}\comp_{k-1}\cdots\comp_1\alpha_{\trgt_0(x)}\comp_0u(\src_{k+1}(y))\right)}.
            \end{multlined}
      \]
      \end{enumerate}
   
    \end{itemize}
     Note that to read the formulas correctly, one has to remember the
     convention that for $k<n$, the composition $\comp_k$ has priority over $\comp_n$ (see \ref{paragr:defoomagma}).
     \end{paragr}
     \begin{example}\label{example:natisoplax}
       When $C$ and $D$ are $n$\nbd{}categories with $n$ finite and $u,v :C \to D$ are two $n$\nbd{}functors, an oplax transformation $\alpha : u \Rightarrow v$ amounts to the data of a $(k+1)$\nbd{}cell $\alpha_x$ of $D$ for each $k$\nbd{}cell $x$ of $C$ with $0 \leq k \leq n$, with source and target as in the previous paragraph. These data being subject to the axioms of the previous paragraph. Note that when $x$ is an $n$\nbd{}cell of $C$, $\alpha_x$ is necessarily a unit, which can be expressed as the equality
       \[
       \alpha_{t_{n-1}(x)}\comp_{n-1}\cdots\comp_1\alpha_{t_0(x)}\comp_0u(x) = v(x)\comp_0\alpha_{s_0(x)}\comp_1\cdots\comp_{n-1}\alpha_{s_{n-1}(x)}
       \]
       In particular, when $n=1$ and $C$ and $D$ are thus (small) categories, an oplax transformation $u \Rightarrow v$ is nothing but a natural transformation from $u$ to $v$.
     \end{example}
     \begin{paragr}
       Let $u : C \to D$ be an $\oo$\nbd{}functor. There is an oplax transformation from $u$ to $u$, denoted by $1_u$, which is defined as
       \[
       (1_u)_{x}:=1_{u(x)}
       \]
       for every cell $x$ of $C$. More abstractly, this oplax transformation corresponds to the $\oo$\nbd{}functor
       \[
       \sD_1 \otimes C \overset{p\otimes u}{\longrightarrow} \sD_0 \otimes D \simeq D,
       \]
       where $p$ is the only $\oo$\nbd{}functor $\sD_1\to \sD_0$.
       \end{paragr}
      \begin{paragr}
    Let
    \[
    \begin{tikzcd}
      B \ar[r,"f"] & C \ar[r,shift left,"u"]  \ar[r,shift right,"v"']&D \ar[r,"g"] &E
    \end{tikzcd}
    \]
    be a diagram in $\omega\Cat$ and $\alpha : u \Rightarrow v$ an oplax transformation. 
    The data of
    \[
    (g\star \alpha)_x := g(\alpha_x)
    \]
    for each cell $x$ of $C$ (resp. 
    \[
    (\alpha \star f)_x :=\alpha_{f(x)}
    \]
    for each cell $x$ of $B$) defines an oplax transformation from $g \circ u$ to $g \circ v$ (resp. $u \circ f$ to $v\circ f$) which we denote by $g\star \alpha$ (resp. $\alpha \star f$).

    More abstractly, if $\alpha$ is seen as an $\oo$\nbd{}functor $\sD_1 \otimes C \to D$, then $g \star \alpha$ (resp.\ $\alpha \star f)$ corresponds to the $\oo$\nbd{}functor obtained as the following composition 
    \[
    \sD_1 \otimes C \overset{\alpha}{\longrightarrow} D \overset{f}{\longrightarrow} E
    \]
    (resp.\
    \[
    \sD_1 \otimes B \overset{\sD_1 \otimes f}{\longrightarrow} \sD_1 \otimes C \overset{\alpha}{\longrightarrow} D).
    \]
  \end{paragr}
      \begin{remark}
        All the above descriptions of oplax transformations can be easily
        dualized for \emph{lax transformations} (i.e.\ the $1$\nbd{}cells of the $\oo$\nbd{}category $\underline{\hom}_{\mathrm{lax}}(X,Y)$ for some $\oo$\nbd{}categories $X$ and $Y$). Habit is the only reason why we put emphasis on oplax transformations rather than lax transformations.
      \end{remark}
      \section{Homotopy equivalences and deformation retracts}
      \begin{paragr}\label{paragr:hmtpyequiv}
    Let $C$ and $D$ be two $\oo$\nbd{}categories and consider the smallest equivalence relation on the set $\Hom_{\oo\Cat}(C,D)$ such that two $\oo$\nbd{}functors from $C$ to $D$ are equivalent if there is an oplax direction between them (in any direction). Let us say that two $\oo$\nbd{}functors $u, v : C \to D$ are \emph{oplax homotopic} if they are equivalent for this equivalence relation. 
      \end{paragr}
      \begin{definition}\label{def:oplaxhmtpyequiv}
        An $\oo$\nbd{}functor $u : C \to D$ is an \emph{oplax homotopy equivalence} if there exists an $\oo$\nbd{}functor $v : D \to C$ such that $u\circ v$ is oplax homotopic to $\mathrm{id}_D$ and $v\circ u$ is oplax homotopic to $\mathrm{id}_C$.
      \end{definition}
      Recall that we write $\gamma^{\Th} : \oo\Cat \to \ho(\oo\Cat^{\Th})$ for the localization functor with respect to the Thomason equivalences.
      \begin{lemma}\label{lemma:oplaxloc}
        Let $u, v : C \to D$ be two $\oo$\nbd{}functors. If there exists an oplax transformation $\alpha : u \Rightarrow v$, then $\gamma^{\Th}(u)=\gamma^{\Th}(v)$.
      \end{lemma}
      \begin{proof}
        This follows immediately from \cite[Théorème B.11]{ara2020theoreme}.
      \end{proof}
      From this lemma and the fact that the Thomason equivalences are saturated (Corollary \ref{cor:thomsaturated}), we deduce the following proposition.
      \begin{proposition}\label{prop:oplaxhmtpyisthom}
        Every oplax homotopy equivalence is a Thomason equivalence.
      \end{proposition}
      \begin{paragr}\label{paragr:defrtract}
        An $\oo$\nbd{}functor $i : C \to D$ is an \emph{oplax deformation retract} if there exists an $\oo$\nbd{}functor $r : D \to C$ such that:
        \begin{enumerate}[label=(\alph*)]
        \item $r\circ i=\mathrm{id}_C$,
        \item there exists an oplax transformation $\alpha : \mathrm{id}_D \Rightarrow i\circ r$.
        \end{enumerate}
        Furthermore, $i$ is a \emph{strong oplax deformation retract} if $\alpha$ can be chosen such that:
        \begin{enumerate}[label=(\alph*),resume]
          \item $\alpha \star i = 1_i$.
        \end{enumerate}
        An oplax deformation retract is a particular case of homotopy equivalence and thus of Thomason equivalence. 
      \end{paragr}
      \begin{lemma}\label{lemma:pushoutstrngdefrtract}
       The pushout of a strong oplax deformation retract is a strong oplax deformation retract.
      \end{lemma}
      \begin{proof}
        Let $i : A \to B$ be a strong oplax deformation retract and
        \begin{equation}\label{cocartsquareretract}\tag{i}
        \begin{tikzcd}
          A \ar[d,"i"] \ar[r,"u"] & A' \ar[d,"i'"] \\
          B \ar[r,"v"] & B'\ar[from=1-1,to=2-2,phantom,very near end,"\ulcorner"]
        \end{tikzcd}
        \end{equation}
        be a cocartesian square. We have to show that $i'$ is also a strong oplax deformation retract. By hypothesis there exist $r : B \to A$  such that $r \circ i = \mathrm{id}_A$ and $\alpha : \sD_1 \otimes B \to B$ such that the diagrams
        \begin{equation}\label{diagramtransf}\tag{ii}
       \begin{tikzcd}
         B\ar[rd,"\mathrm{id}_B"] \ar[d,"i_0^B"']& \\
         \sD_1\otimes B \ar[r,"\alpha"] & B, \\
         B \ar[ru,"i\circ r"'] \ar[u,"i_1^B"]&
       \end{tikzcd}
        \end{equation}
        and
        \begin{equation}\label{diagramstrong}\tag{iii}
          \begin{tikzcd}
            \sD_1 \otimes A \ar[rr, bend right,"p\otimes i"']\ar[r,"\sD_1 \otimes i"] & \sD_1 \otimes B \ar[r,"\alpha"] & B,
          \end{tikzcd}
        \end{equation}
        where $p$ is the unique morphism $\sD_1 \to \sD_0$, are commutative.

        From the commutativity of the following solid arrow diagram
        \[
        \begin{tikzcd}
          A \ar[r,"u"] \ar[d,"i"] & A' \ar[d,"i'"] \ar[dd,bend left=75,"\mathrm{id}_{A'}"] \\
          B \ar[d,"r"] \ar[r,"v"] & B' \ar[d,"r'",dashed ] \\
          A \ar[r,"u"] & A',
          \ar[from=1-1,to=2-2,phantom,"\ulcorner",very near end]
        \end{tikzcd}
        \]
        we deduce the existence of $r' : B' \to A'$ that makes the whole diagram commute. In particular, we have $r' \circ i' = \mathrm{id}_{A'}$.

        From the commutativity of (\ref{diagramstrong}), we easily deduce the commutativity of the following solid arrow diagram
        \[
        \begin{tikzcd}
          \sD_1\otimes A \ar[r,"\sD_1\otimes u"] \ar[d,"\sD_1\otimes i"] & \sD_1 \otimes A' \ar[d,"\sD_1 \otimes i'"] \ar[dd,bend left=75,"p\otimes i'"] \\
         \sD_1\otimes B \ar[d,"\alpha"] \ar[r,"\sD_1 \otimes v"] & \sD_1 \otimes B' \ar[d,"\alpha'",dashed ] \\
         B \ar[r,"v"] &  B'.
        \end{tikzcd}
        \]
        The existence of $\alpha' : \sD_1 \otimes B' \to B'$ that makes the whole diagram commutes follows from the fact that the functor $\sD_1 \otimes \shortminus$ preserves colimits. In particular, we have \[\alpha' \circ (\sD_1 \otimes i') = p \otimes i'.\]
        
        Now, notice that for every $\oo$\nbd{}category $C$, the maps
        \[
        i^C_0 : C \to \sD_1 \otimes C \text{ and } i^C_1 : C \to \sD_1 \otimes C
        \]
        are natural in $C$. Using this naturality and simple diagram chasing (left to the reader), we obtain the equalities
        \[
        \alpha ' \circ i_0^{B'} \circ v= v,\]
          \[
          \alpha' \circ i^{B'}_0 \circ i'=i',\]
          and the equalities
        \[
        \alpha ' \circ i_1^{B'} \circ v=  i' \circ r' \circ v \]
        \[
        \alpha' \circ i^{B'}_1 \circ i'=i' \circ r' \circ i'.
        \]
        Using the fact that square (\ref{cocartsquareretract}) is cocartesian, we deduce that $\alpha ' \circ i_0^{B'} = \mathrm{id}_{B'}$ and $\alpha' \circ i^{B'}_1 = i' \circ r'$. This proves that $i'$ is an oplax deformation retract, which is furthermore strong because of the equality $\alpha' \circ (\sD_1 \otimes i') = p \otimes i'$.
      \end{proof}
      In the following proposition, a \emph{co-universal Thomason equivalence} means a co-universal weak equivalence for the localizer $(\oo\Cat, \W^{\Th}_{\oo})$ (Definition \ref{def:couniversalwe}).
      \begin{proposition}
        Every strong oplax deformation retract is a co-universal Thomason equivalence.
      \end{proposition}
      \begin{proof}
        Immediate consequence of Lemma \ref{lemma:pushoutstrngdefrtract} and the fact that oplax transformation retracts are Thomason equivalences.
        \end{proof}
      \begin{remark}
        All the results we have seen in this section are still true if we replace ``oplax'' by ``lax'' everywhere.
      \end{remark}
\section[Equivalences of \texorpdfstring{$\oo$}{ω}-categories and the folk model
  structure]{Equivalences of \texorpdfstring{$\oo$}{ω}-categories and the folk model
  structure%
\sectionmark{The folk model structure}}
\sectionmark{The folk model structure}
\begin{paragr}\label{paragr:ooequivalence}
  Let $C$ be an $\omega$-category. We define the equivalence relation $\sim_{\omega}$ on the set $C_n$ by co-induction on $n \in \mathbb{N}$. For $x, y \in C_n$, we have $x \sim_{\omega} y $ when:
  \begin{itemize}
    \item[-] $x$ and $y$ are parallel,
    \item[-] there exist $r, s \in C_{n+1}$ such that $r : x \to y$, $s : y \to x$,
      \[
      r\ast_{n}s \sim_{\omega} 1_y
      \]
      and
      \[
      s\ast_nr \sim_{\omega} 1_x.
      \]

    \end{itemize}
    For details on this definition and the proof that it is an equivalence relation, see \cite[section 4.2]{lafont2010folk}.
\end{paragr}
\begin{example}
  Let $x$ and $y$ be two $0$-cells of an $n$\nbd{}category $C$.
  \begin{itemize}[label=-]
  \item When $n=1$, $x \sim_{\omega} y$ means that $x$ and $y$ are isomorphic.
  \item When $n=2$, $x \sim_{\omega} y$ means that $x$ and $y$ are equivalent, i.e.\ there exist $f : x \to y$ and $g : y \to x$ such that $fg$ is isomorphic to $1_y$ and $gf$ is isomorphic to $1_x$.
    \end{itemize}
\end{example}
For later reference, we put here the following trivial but important lemma, whose proof is omitted.
\begin{lemma}\label{lemma:ooequivalenceisfunctorial}
  Let $F : C \to D$ be an $\oo$\nbd{}functor, $n \geq 0$ and $x,y$ be $n$\nbd{}cells of $C$. If $x \sim_{\oo} y$, then $F(x) \sim_{\oo} F(y)$.
  \end{lemma}
  \begin{definition}\label{def:eqomegacat}
    An $\omega$-functor $F : C \to D$ is an \emph{equivalence of $\oo$\nbd{}categories} when:
    \begin{itemize}
      \item[-] for every $y \in D_0$, there exists $x \in C_0$ such that
      \[F(x)\sim_{\omega}y,\]
    \item[-] for every $n \geq 0$, for all $x,y \in C_n$ that are parallel and for every $\beta~\in~D_{n+1}$ such that \[\beta : F(x) \to F(y),\] there exists $\alpha \in C_{n+1}$ such that
      \[\alpha : x \to y
      \]
      and
      \[F(\alpha)\sim_{\omega}\beta.\]
      \end{itemize}
  \end{definition}
  \begin{example}\label{example:equivalencecategories}
    If $C$ and $D$ are (small) categories seen as $\oo$\nbd{}categories, then a functor $F : C \to D$ is an equivalence of $\oo$\nbd{}categories if and only if it is fully faithful and essentially surjective, i.e.\ an equivalence of categories. 
  \end{example}

  \begin{theorem}\label{thm:folkms}
    There exists a cofibrantly generated model structure on $\omega\Cat$ whose
    weak equivalences are the equivalences of $\oo$\nbd{}categories, and whose
    cofibrations are generated by the set $\{i_n : \sS_{n-1} \to \sD_n \vert n \in \mathbb{N}\}$ (see \ref{paragr:defglobe}).
  \end{theorem}
  \begin{proof}
    This is the main result of \cite{lafont2010folk}.
  \end{proof}
 
 \begin{paragr}\label{paragr:folkms}
   The model structure of the previous theorem is commonly referred to as \emph{folk model structure} on $\omega\Cat$.
   Data of this model structure will often be referred to by using the adjective folk, e.g.\ \emph{folk cofibration}. Consequently \emph{folk weak equivalence} and \emph{equivalence of $\oo$\nbd{}categories} mean the same thing.

   Furthermore, as in the Thomason case (see \ref{paragr:notationthom}), we
   usually make reference to the word ``folk'' in the notations of homotopic
   constructions induced by the folk weak equivalences. For example, we write
   $\W^{\folk}$ for the class of folk weak equivalences, $\Ho(\oo\Cat^{\folk})$ for the homotopy op\nbd{}prederivator of $(\oo\Cat,\W_{\oo}^{\folk})$ and
   \[
   \gamma^{\folk} : \oo\Cat \to \Ho(\oo\Cat^{\folk})
   \]
   for the localization morphism. It follows from the previous theorem and
   Theorem \ref{thm:cisinskiI} that the localizer $(\oo\Cat,\W_{\oo}^{\folk})$
   is homotopy cocomplete. We will speak of ``folk homotopy
   colimits'' and ``folk homotopy cocartesian squares'' for homotopy colimits
   and homotopy cocartesian squares in this localizer. 
 \end{paragr}
 \begin{paragr}\label{paragr:folktrivialfib}
   Using the set $\{i_n : \sS_{n-1} \to \sD_n \vert n \in \mathbb{N} \}$ of generating folk cofibrations, we obtain that an $\oo$\nbd{}functor $F : C \to D$ is a \emph{folk trivial fibration} when:
   \begin{itemize}[label=-]
   \item for every $y \in D_0$, there exists $x \in C_0$ such that
     \[
     F(x)=y
     \]
   \item for every $n\geq 0$, for all $x,y \in C_n$ that are parallel and for every $\beta~\in~D_{n+1}$ such that
     \[
     \beta : F(x) \to F(y)
     \]
     there exists $\alpha \in C_{n+1}$ such that
     \[
     \alpha : x \to y
     \]
     and
     \[
     F(\alpha)=\beta.
     \]
     This characterization of folk trivial fibrations is to be compared with Definition \ref{def:eqomegacat} of equivalences of $\oo$\nbd{}categories. 
     \end{itemize}
 \end{paragr}
  \begin{proposition}\label{prop:freeiscofibrant}
    An $\omega$-category is cofibrant for the folk model structure if and only if it is free.
  \end{proposition}
  \begin{proof}
    The fact that every free $\omega$-category is cofibrant follows immediately from the fact that the $i_n : \sS_{n-1} \to \sD_n$ are cofibrations and that every $\oo$\nbd{}category $C$ is the colimit of the canonical diagram (Lemma \ref{lemma:filtration})
    \[
     \sk_{0}(C) \to \sk_{1}(C) \to \cdots \to \sk_n(C) \to \sk_{n+1}(C) \to \cdots
    \]
    For the converse, see \cite{metayer2008cofibrant}.
  \end{proof}
  \section[Equivalences of \texorpdfstring{$\oo$}{ω}-categories vs Thomason
    equivalences ]{Equivalences of \texorpdfstring{$\oo$}{ω}-categories vs Thomason
      equivalences%
      \sectionmark{Folk vs Thomason}}
  \sectionmark{Folk vs Thomason}
    \begin{lemma}\label{lemma:nervehomotopical}
The nerve functor $N_{\omega} : \omega\Cat \to \Psh{\Delta}$ sends the equivalences of $\oo$\nbd{}categories to weak equivalences of simplicial sets.    
  \end{lemma}
  \begin{proof}
    Since every $\omega$-category is fibrant for the folk model structure
    \cite[Proposition 9]{lafont2010folk}, it follows from Ken Brown's Lemma
    \cite[Lemma 1.1.12]{hovey2007model} that it suffices to show that the nerve
    sends the folk trivial fibrations to weak equivalences of simplicial sets. In particular, it suffices to show the stronger condition that the nerve sends the folk trivial fibrations to trivial fibrations of simplicial sets.

    By adjunction, this is equivalent to showing that the functor $c_{\omega} : \Psh{\Delta} \to \omega\Cat$ sends the cofibrations of simplicial sets to folk cofibrations. Since $c_{\omega}$ is cocontinuous and the cofibrations of simplicial sets are generated by the inclusions
    \[
    \partial \Delta_n \to \Delta_n
    \]
    for $n \in \mathbb{N}$, it suffices to show that $c_{\omega}$ sends these inclusions to folk cofibrations.

    Now, it follows from \cite[Lemma 5.1]{street1987algebra} that the image of the inclusion $\partial \Delta_n \to \Delta_n$ by $c_{\omega}$ can be identified with the canonical inclusion
    \[
    \sk_{n-1}(\Or_n) \to \Or_n.
    \]
    Since $\Or_n$ is free, this last morphism is by definition a push-out of a coproduct of folk cofibrations (see Definition \ref{def:nbasis}), hence a folk cofibration.
  \end{proof}
  As an immediate consequence of the previous lemma, we have the following proposition.
  \begin{proposition}\label{prop:folkisthom}
    Every equivalence of $\oo$\nbd{}categories is a Thomason equivalence. 
  \end{proposition}
  \begin{remark}
    The converse of the above proposition is false. For example, the unique $\oo$\nbd{}functor
    \[
    \sD_1 \to \sD_0
    \]
    is a Thomason equivalence because its image by the nerve is the unique morphism of simplicial sets $\Delta_1 \to \Delta_0$ (which obviously is a weak equivalence), but it is \emph{not} an equivalence of $\oo$\nbd{}categories because $\sD_1$ and $\sD_0$ are not equivalent as categories (see Example \ref{example:equivalencecategories}).
    \end{remark}
  \begin{paragr}\label{paragr:compweakeq}
    Proposition \ref{prop:folkisthom} implies that the identity functor on $\oo\Cat$ induces a morphism of localizers $(\oo\Cat,\W^{\folk}) \to (\oo\Cat,\W^{\Th})$, which in turn induces a functor between localized categories
    \[
    \mathcal{J} : \ho(\oo\Cat^{\folk}) \to \ho(\oo\Cat^{\Th}).
    \]
    This functor cannot be an equivalence since this would imply that every Thomason equivalence is an equivalence of $\oo$\nbd{}categories. 
  \end{paragr}
   \section{{Slice \texorpdfstring{$\oo$}{ω}-categories and folk Theorem~A}}
  \begin{paragr}\label{paragr:slices}
    Let $A$ be an $\oo$\nbd{}category and $a_0$ an object of $A$. We define the slice $\oo$\nbd{}category $A/a_0$ as the following fibred product:
    \[
    \begin{tikzcd}
      A/a_0 \ar[d] \ar[r] & \homlax(\sD_1,A) \ar[d,"\pi_1^A"] \\
      \sD_0 \ar[r,"\langle a_0 \rangle"'] & A.
      \ar[from=1-1,to=2-2,phantom,very near start,"\lrcorner"]
    \end{tikzcd}
    \]
    We also define an $\oo$\nbd{}functor $\pi : A/a_0 \to A$ as the following composition
    \[
    \pi : A/a_0 \to \homlax(\sD_1,A) \overset{\pi^A_0}{\longrightarrow} A.
    \]
    Let us now give an alternative definition of the $\oo$\nbd{}category $A/a_0$ using explicit formulas. The equivalence with the previous definition follows from the dual of \cite[Proposition B.5.2]{ara2016joint}
    \begin{itemize}[label=-]
    \item An $n$\nbd{}cell of $A/a_0$ is a table 
      \[
      (x,a)=\begin{pmatrix}
      \begin{matrix}
      (x_0,a_1) & (x_1,a_2) & \cdots & (x_{n-1},a_n) \\[0.5em]
        (x_0',a_1') & (x_1',a_2') & \cdots & (x_{n-1}',a_n')
      \end{matrix}
      & (x_n,a_{n+1})
      \end{pmatrix}
      \]
      where $x_0$ and $x_0'$ are $0$-cells of $A$, and:
      
      \begin{tabular}{ll}
        $x_i : x_{i-1} \longrightarrow x'_{i-1}$, &for every $1 \leq i \leq n$,\\[0.75em]
         $x_i': x_{i-1} \longrightarrow x'_{i-1}$, &for every $1 \leq i \leq n-1$,\\[0.75em]
         $a_i : a'_{i-1}\comp_{i-2} a'_{i-2} \comp_{i-3} \cdots \comp_1 a'_1\comp_0 x_{i-1} \longrightarrow a_{i-1}$, &for every $1 \leq i \leq n+1$,\\[0.75em]
        $a'_i : a'_{i-1}\comp_{i-2} a'_{i-2} \comp_{i-3} \cdots \comp_1 a'_1 \comp_0 x'_{i-1} \longrightarrow a_{i-1}$, &for every $1 \leq i \leq n$\\
      \end{tabular}
      
      are $i$-cells of $A$. In low dimension, this gives:
      
      \begin{tabular}{ll}
     $(x_0,a_1)$ :&  {\begin{tikzcd}
          x_0 \ar[d,"a_1"] \\ a_0
        \end{tikzcd}} \\[2.75em]
        {$\begin{pmatrix}
      \begin{matrix}
      (x_0,a_1) \\[0.5em]
        (x_0',a_1')
      \end{matrix}
      & (x_1,a_{2})
          \end{pmatrix}$} :& {\begin{tikzcd}[column sep=small] x_0 \ar[rr,"x_1"] \ar[rd,"a_1"',""{name=A,left}] && x_0' \ar[ld,"a_1'"] \\ &a_0 & \ar[from=1-3,to=A,Rightarrow,"a_2", shorten <=1em, shorten >=1em]\end{tikzcd}} \\[2.75em]
        {$\begin{pmatrix}
      \begin{matrix}
      (x_0,a_1) & (x_1,a_2)  \\[0.5em]
        (x_0',a_1') & (x_1',a_2')
      \end{matrix}
      & (x_2,a_{3})
      \end{pmatrix}$}:&{$\begin{tikzcd}[column sep=small] x_0 \ar[rr,"x_1"] \ar[rd,"a_1"',""{name=A,left}] && x_0' \ar[ld,"a_1'"] \\ &a_0 & \ar[from=1-3,to=A,Rightarrow,"a_2", shorten <=1em, shorten >=1em]\end{tikzcd}\; \overset{a_3}{\Lleftarrow} \; \begin{tikzcd}[column sep=small] x_0\ar[rr,bend left=50,"x_1",pos=11/20,""{name=toto,below}] \ar[rr,"x_1'"description,""{name=titi,above}] \ar[rd,"a_1"',""{name=A,left}] && x_0' \ar[ld,"a_1'"] \\ &a_0 &. \ar[from=1-3,to=A,Rightarrow,"a_2'", shorten <=1em, shorten >=1em] \ar[from=toto,to=titi,Rightarrow,"x_2",pos=1/5]\end{tikzcd}$} 
        \end{tabular}
    \item The source and target of the $n$\nbd{}cell $(a,x)$ are given by the matrices:
      \[
      s(x,a)=\begin{pmatrix}
      \begin{matrix}
      (x_0,a_1) & (x_1,a_2) & \cdots & (x_{n-2},a_{n-1}) \\[0.5em]
        (x_0',a_1') & (x_1',a_2') & \cdots & (x_{n-2}',a_{n-1}')
      \end{matrix}
      & (x_{n-1},a_{n})
      \end{pmatrix}
      \]
      \[
      t(x,a)=\begin{pmatrix}
      \begin{matrix}
      (x_0,a_1) & (x_1,a_2) & \cdots & (x_{n-2},a_{n-1}) \\[0.5em]
        (x_0',a_1') & (x_1',a_2') & \cdots & (x_{n-2}',a_{n-1}')
      \end{matrix}
      & (x'_{n-1},a'_{n})
      \end{pmatrix}.
      \]
    \item The unit of the $n$\nbd{}cell $(a,x)$ is given by the table:
      \[
      1_{(x,a)}=\begin{pmatrix}
      \begin{matrix}
      (x_0,a_1) & (x_1,a_2) & \cdots & (x_{n-1},a_n) & (x_n,a_{n+1}) \\[0.5em]
        (x_0',a_1') & (x_1',a_2') & \cdots & (x_{n-1}',a_n') & (x_n,a_{n+1})
      \end{matrix}
      & (1_{x_n},1_{a_{n+1}})
      \end{pmatrix}.
      \]
    \item The composition of $n$\nbd{}cells $(x,a)$ and $(y,b)$ such that $\src_k(y,b)=\trgt_k(a,x)$, is given by the table:
      \[
      (y,b)\comp_k (x,a)=\begin{pmatrix}
      \begin{matrix}
        (x_0,a_1) & \cdots & (x_k,a_k) & (z_{k+1},c_{k+2}) & \cdots &(z_{n-1},c_n) \\[0.5em]
        (y_0',b_1') & \cdots & (y'_k,b'_k) &(z'_{k+1},c'_{k+2}) & \cdots & (z'_{n-1},c'_n) \\
      \end{matrix}
      & (z_n,c_{n+1})
      \end{pmatrix},
      \]
      where:
      
      \begin{tabular}{ll}
        $z_{i}=y_i\comp_k x_i$, & for every $k+1 \leq i \leq n$, \\[0.75em]
        $z'_i=y'_i \comp_k x'_i$, & for every $k+1 \leq i \leq n-1$, \\[0.75em]
        $c_i=a_i\comp_{k+1} b_i \comp_{k} a'_{k} \comp_{k-1} a'_{k-1} \comp_{k-2} \cdots \comp_{1} a'_1\comp_0 x_{k+1}$,&for every $k+2 \leq i \leq n+1$, \\[0.75em]
        $c'_i=a'_i\comp_{k+1} b'_i \comp_{k} a'_{k} \comp_{k-1} a'_{k-1} \comp_{k-2} \cdots \comp_{1} a'_1\comp_0 x'_{k+1}$,&for every $k+2 \leq i \leq n$.\\
      \end{tabular}
    \end{itemize}
    We leave it to the reader to check that the formulas are well defined and that the axioms for $\oo$\nbd{}categories are satisfied. The canonical forgetful $\oo$\nbd{}functor $\pi : A/a_0 \to A$ is simply expressed as:
        \begin{align*}
      A/a_0 &\to A \\
      (x,a) &\mapsto x_n.
        \end{align*}
        Notice that if $A$ is an $n$\nbd{}category, then so is $A/a_0$. In this case, for an $n$\nbd{}cell $(x,a)$, $a_{n+1}$ is a unit, hence
    \[
    a'_n \comp_{n-1} a'_{n-1} \comp_{n-2} \cdots \comp_1 a'_1 \comp_0 x_n = a_n.
    \]
  \end{paragr}
  \begin{example}\label{example:slicecategories}
    For a small category $A$ (considered as an $\oo$\nbd{}category) and an object $a_0$ of $A$, the category $A/a_0$ in the sense of the previous paragraph is nothing but the usual slice category of $A$ over $a_0$.
  \end{example}
  \begin{paragr}\label{paragr:comma}
    Let $u : A \to B$ be a morphism of $\oo\Cat$ and $b_0$ an object of $B$. We define the $\oo$\nbd{}category $A/b_0$ (also denoted by $u\downarrow b_0$) as the following fibred product:
    \[
    \begin{tikzcd}
    A/b_0 \ar[d,"u/b_0"'] \ar[r] & A \ar[d,"u"] \\
    B/b_0 \ar[r,"\pi"'] & B.
    \ar[from=1-1,to=2-2,phantom,description,very near start,"\lrcorner"]
    \end{tikzcd}
    \]
    More explicitly, an $n$\nbd{}cell $(x,b)$ of $A/b_0$ is a table
    \[
    (x,b)=\begin{pmatrix}
    \begin{matrix}
      (x_0,b_1) & (x_1,b_2) & \cdots & (x_{n-1},b_n) \\[0.5em]
      (x_0',b_1') & (x_1',b_2') & \cdots & (x_{n-1}',b_n')
    \end{matrix}
    & (x_n,b_{n+1})
    \end{pmatrix}
    \]
    where the $x_i$ and $x'_i$ are $i$-cells of $A$ such that
    
    \begin{tabular}{ll}
      $x_i : x_{i-1} \longrightarrow x'_{i-1}$, &for every $1 \leq i \leq n$,\\[0.75em]
      $x_i': x_{i-1} \longrightarrow x'_{i-1}$, &for every $1 \leq i \leq n-1$,\\[0.75em]
    \end{tabular}
    
     and the $b_i$ and $b'_i$ are $i$-cells of $B$ such that
    \[
    \begin{pmatrix}
    \begin{matrix}
      (u(x_0),b_1) & (u(x_1),b_2) & \cdots & (u(x_{n-1}),b_n) \\[0.5em]
      (u(x_0'),b_1') & (u(x'_1),b_2') & \cdots & (u(x'_{n-1}),b_n')
    \end{matrix}
    & (u(x_n),b_{n+1})
    \end{pmatrix}
    \]
    is an $n$\nbd{}cell of $B/b_0$.

    The canonical $\oo$\nbd{}functor $A/b_0 \to A$ is simply expressed as
    \begin{align*}
      A/b_0 &\to A\\
      (x,b) &\mapsto x_n,
    \end{align*}
    and the $\oo$\nbd{}functor $u/b_0$ as
    \begin{align*}
      u/b_0 : A/b_0 &\to B/b_0 \\
      (x,b) &\mapsto (u(x),b).
    \end{align*}
    More generally, if we have a commutative triangle in $\oo\Cat$
      \[
    \begin{tikzcd}[column sep=small]
    A \ar[rr,"u"] \ar[dr,"v"'] & &B \ar[dl,"w"] \\
    &C&,
    \end{tikzcd}
    \]
    then for every object $c_0$ of $C$, we have a functor $u/c_0 : A/c_0 \to B/c_0$ defined as
    \begin{align*}
      u/c_0 : A/c_0 &\to B/c_0 \\
      (x,c) &\mapsto (u(x),c).
    \end{align*}
  \end{paragr}

  \begin{theorem}\label{thm:folkthmA}(Folk Theorem $A$) Let
    \[
    \begin{tikzcd}[column sep=small]
    A \ar[rr,"u"] \ar[dr,"v"'] & &B \ar[dl,"w"] \\
    &C&
    \end{tikzcd}
    \]
    be a commutative triangle in $\oo\Cat$. If for every object $c_0$ of $C$ the induced morphism
    \[
    u/c_0 : A/c_0 \to B/c_0
    \]
    is an equivalence of $\oo$\nbd{}categories, then so is $u$.
  \end{theorem}
 \begin{proof}
   Before anything else, recall from Lemma \ref{lemma:ooequivalenceisfunctorial} that given an $\oo$\nbd{}functor $F : X \to Y$ and $n$\nbd{}cells $x$ and $y$ of $X$, if $x \sim_{\oo} y$, then $F(x) \sim_{\oo} F(y)$.
   \begin{enumerate}[label=(\roman*)]
   \item Let $b_0$ be a $0$\nbd{}cell of $B$ and set $c_0:=w(b_0)$. By definition, the pair $(b_0,1_{c_0})$ is a $0$-cell of $B/c_0$. By hypothesis, we know that there exists a $0$\nbd{}cell $(a_0,c_1)$ of $A/c_0$ such that $(u(a_0),c_1)\sim_{\oo} (b_0,1_{c_0})$. Hence, by applying the canonical functor $B/c_0 \to B$, we obtain that $u(a_0) \sim_{\oo} b_0$.

   \item Let $f$ and $f'$ be parallel $n$\nbd{}cells of $A$ and $\beta : u(f) \to u(f')$ an $(n+1)$\nbd{}cell of $B$. We need to show that there exists an $(n+1)$\nbd{}cell $\alpha : f \to f'$ of $A$ such that $u(\alpha) \sim_{\oo} \beta$.

     Let us use the notations:
     \begin{itemize}[label=-]
     \item $a_i := \src_i(f)=\src_i(f')$ for $0 \leq i <n$,
     \item $a_i' := \trgt_i(f)=\trgt_i(f')$ for $0 \leq i <n$,
       \item $a_n:=f$ and  $a_n'=f'$.
       \end{itemize}
     It is straightforward to check that  
       \[
       \begin{pmatrix}
    \begin{matrix}
      (a_0,v(a_1')) & (a_1,v(a_2')) & \cdots & (a_{n-1},v(a_n')) \\[0.5em]
      (a_0',1_{v(a_0')}) & (a_1',1_{v(a_1')}) & \cdots & (a_{n-1}',1_{(v(a_{n-1}'))})
    \end{matrix}
    & (a_n,w(\beta))
    \end{pmatrix}
       \]
       and
       \[
       \begin{pmatrix}
    \begin{matrix}
      (a_0,v(a_1')) & (a_1,v(a_2')) & \cdots & (a_{n-1},v(a_n')) \\[0.5em]
      (a_0',1_{v(a_0')}) & (a_1',1_{v(a_1')}) & \cdots & (a_{n-1}',1_{(v(a_{(n-1)}'))})
    \end{matrix}
    & (a_n',1_{v(a_n')})
    \end{pmatrix}
       \]
       are parallel $n$\nbd{}cells of $A/c_0$ where we set $c_0:=v(a_0')$. Similarly, we have an $(n+1)$\nbd{}cell of $B/c_0$
       \[
       \begin{pmatrix}
    \begin{matrix}
      (u(a_0),v(a_1')) & \cdots & (u(a_{n-1}),v(a_n')) &  (u(a_n),w(\beta)) \\[0.5em]
      (u(a_0'),1_{v(a_0')})  & \cdots & (u(a_{n-1}'),1_{(v(a_{n-1}'))}) & (u(a_n'),1_{v(a_n')})
    \end{matrix}
    & (\beta,1_{w(\beta)})
    \end{pmatrix}
       \]
       whose source and target respectively are the image by $u/c_0$ of the above two cells of $A/c_0$. By hypothesis,  there exists an $(n+1)$\nbd{}cell of $A/c_0$ of the form
       \[
       \begin{pmatrix}
         \begin{matrix}
           (a_0,v(a_1'))  & \cdots & (a_{n-1},v(a_n')) &  (a_n,w(\beta)) \\[0.5em]
           (a_0',1_{v(a_0')})  & \cdots & (a_{n-1}',1_{(v(a_{n-1}'))}) & (a_n',1_{v(a_n')})
         \end{matrix}
         & (\alpha,\Lambda)
       \end{pmatrix}
       \]
       whose image by $u/c_0$ is equivalent for the relation $\sim_{\oo}$ to the above ${(n+1)}$\nbd{}cell of $B/c_0$. In particular, the source and target of $\alpha$ are respectively $f$ and $f'$. Finally, we obtain that $u(\alpha) \sim_{\oo} \beta$ by applying the canonical $\oo$\nbd{}functor $B/c_0 \to B$.
   \end{enumerate}
 \end{proof}
 
 \begin{paragr} The name ``folk Theorem A'' is an explicit reference to Quillen's Theorem A \cite[Theorem A]{quillen1973higher} and its $\oo$\nbd{}categorical generalization by Ara and Maltsiniotis \cite{ara2018theoreme,ara2020theoreme}. For the sake of comparison we recall below the latter one.
   \end{paragr}
 \begin{theorem}[Ara and Maltsiniotis' Theorem A] Let
   \[
   \begin{tikzcd}[column sep=small]
     A \ar[rr,"u"] \ar[dr,"v"'] & &B \ar[dl,"w"] \\
     &C&
   \end{tikzcd}
   \]
   be a commutative triangle in $\oo\Cat$. If for every object $c_0$ of $C$ the induced morphism
   \[
   u/c_0 : A/c_0 \to B/c_0
   \]
   is a Thomason equivalence, then so is $u$.
 \end{theorem}


\chapter{Homology and abelianization of \texorpdfstring{$\oo$}{ω}-categories}
\chaptermark{Homology of $\omega$-categories}
\section{Homology via the nerve}
\begin{paragr}
We denote by $\Ch$ the category of non-negatively graded chain complexes of abelian groups. Recall that $\Ch$ can be equipped with a cofibrantly generated model structure, known as the \emph{projective model structure on $\Ch$}, where:
    \begin{itemize}
     \item[-] the weak equivalences are the quasi-isomorphisms, i.e.\ morphisms of chain complexes that induce an isomorphism on homology groups,
              \item[-] the cofibrations are the morphisms of chain complexes $f: X\to Y$ such that for every $n\geq 0$, $f_n : X_n \to Y_n$ is a monomorphism with projective cokernel,
     \item[-] the fibrations are the morphisms of chain complexes $f : X \to Y$ such that for every $n>0$, $f_n : X_n \to Y_n$ is an epimorphism.
     \end{itemize}
     (See for example \cite[Section 7]{dwyer1995homotopy}.)
    From now on, we will implicitly consider that the category $\Ch$ is equipped with this model structure. 
 \end{paragr}
\begin{paragr}
  Let $X$ be a simplicial set. We denote by $K_n(X)$ the abelian group of $n$\nbd{}chains of $X$, i.e.\ the free abelian group on the set $X_n$. For $n>0$, let $\partial : K_n(X) \to K_{n-1}(X)$ be the linear map defined for $x \in X_n$ by
  \[
  \partial(x):=\sum_{i=0}^n(-1)^i\partial_i(x).
  \]
  It follows from the simplicial identities (see \cite[section 2.1]{gabriel1967calculus}) that $\partial \circ \partial = 0$. Hence, the previous data defines a chain complex $K(X)$ and this defines a functor
  \begin{align*}
    K : \Psh{\Delta} &\to \Ch\\
    X &\mapsto K(X)
  \end{align*}
  in the expected way.
\end{paragr}
\begin{paragr}
  Recall that an $n$-simplex $x$ of a simplicial set $X$ is \emph{degenerate} if there exists an epimorphism $\varphi : [n] \to [m]$ with $m<n$ and an $m$-simplex $y$ such that $X(\varphi)(y)=x$. We denote by $D_n(X)$ the subgroup of $K_n(X)$ generated by the degenerate $n$-simplices and by $\kappa_n(X)$ the abelian group of \emph{normalized $n$\nbd{}chains}:
  \[
  \kappa_n(X)=K_n(X)/D_n(X).
  \]
  Using the simplicial identities, it can be shown that $\partial(D_n(X)) \subseteq D_{n-1}(X)$ for every $n>0$. Hence, there is an induced differential which we still denote by $\partial$:
  \[
  \partial : \kappa_n(X) \to \kappa_{n-1}(X).
  \]
This defines a chain complex $\kappa(X)$, which we call the \emph{normalized chain complex of $X$}. This yields a functor
  \begin{align*}
    \kappa : \Psh{\Delta} &\to \Ch \\
    X &\mapsto \kappa(X).
  \end{align*}
\end{paragr}
\begin{lemma}\label{lemma:normcompquil}
The functor $\kappa : \Psh{\Delta} \to \Ch$ is left Quillen and sends the weak equivalences of simplicial sets to quasi-isomorphisms.  
\end{lemma}
\begin{proof}
  Recall that the Quillen model structure on simplicial sets admits the set of inclusions
  \[
  \{\partial\Delta_n \hookrightarrow \Delta_n \vert n \in \mathbb{N} \}
  \]
  as generating cofibrations and the set of inclusions
  \[
  \{\Lambda^i_n \hookrightarrow \Delta_n \vert n \in \mathbb{N}, 0 \leq i \leq n\}
  \]
  as generating trivial cofibrations (see for example \cite[Section I.1]{goerss2009simplicial} for the notations). A quick computation, which we leave to the reader, shows that the image by $\kappa$ of $\partial\Delta_n \hookrightarrow \Delta_n$ is a monomorphism with projective cokernel and the image by $\kappa$ of $\Lambda^i_n \hookrightarrow \Delta_n$ is a quasi-isomorphism. This proves that $\kappa$ is left Quillen. Since all simplicial sets are cofibrant, it follows from Ken Brown's Lemma \cite[Lemma 1.1.12]{hovey2007model} that $\kappa$ also preserves weak equivalences.
  \end{proof}
\begin{remark}
  The previous lemma admits also a more conceptual proof as follows. From the Dold--Kan equivalence, we know that $\Ch$ is equivalent to the category $\Ab(\Delta)$ of simplicial abelian groups and with this identification the functor $\kappa : \Psh{\Delta} \to \Ch$ is left adjoint of the canonical forgetful functor
    \[
    U : \Ch \simeq \Ab(\Delta) \to \Psh{\Delta}
    \]
    induced by the forgetful functor from abelian groups to sets. The fact that $U$ is right Quillen follows then from \cite[Lemma 2.9 and Corollary 2.10]{goerss2009simplicial}.
\end{remark}
\begin{paragr}
  In particular, $\kappa$ induces a morphism of localizers \[\kappa : (\Psh{\Delta},\W_{\Delta}) \to (\Ch,\W_{\Ch}),\]
  where we wrote $\W_{\Ch}$ for the class of quasi-isomorphisms. 
  \end{paragr}
\begin{definition}\label{def:hmlgycat}
  The \emph{singular homology functor for $\oo$\nbd{}categories} $\sH^{\sing}$ is defined as the following composition
  \[
  \sH^{\sing} : \ho(\oo\Cat^{\Th}) \overset{\overline{N_{\omega}}}{\longrightarrow} \ho(\Psh{\Delta}) \overset{\overline{\kappa}}{\longrightarrow} \ho(\Ch).
    \]
    For an $\oo$\nbd{}category $C$, $\sH^{\sing}(C)$ is the \emph{singular homology of $C$}.
\end{definition}
\begin{paragr}\label{paragr:singularhmlgygroup}
  In other words, the singular homology of $C$ is the chain complex $\kappa(N_{\oo}(C))$ seen as an object of $\ho(\Ch)$ (see Remark \ref{remark:localizedfunctorobjects}). For $k \geq 0$, the $k$\nbd{}th singular homology group of an $\oo$\nbd{}category $C$ is defined as
  \[
  H_k^{\sing}(C):=H_k(\sH^{\sing}(C))=H_k(\kappa(N_{\oo}(C))),
  \]
  where $H_k : \ho(\Ch) \to \Ab$ is the usual functor that associates to an object of $\ho(\Ch)$ its $k$\nbd{}th homology group.
\end{paragr}


\begin{remark}\label{remark:singularhmlgyishmlgy}
  The adjective ``singular'' is there to avoid future confusion with another
  homological invariant for $\oo$\nbd{}categories that will be introduced later.
  As a matter of fact, the underlying point of view adopted in this thesis is
  that \emph{singular homology of $\oo$\nbd{}categories} ought to be simply called
  \emph{homology of $\oo$\nbd{}categories} as it is the only ``correct''
  definition of homology. This assertion will be justified in Remark \ref{remark:polhmlgyisnotinvariant}.
\end{remark}
\begin{remark}
   We could also have defined the singular homology of $\oo$\nbd{}categories using $K : \Psh{\Delta}\to \Ch$ instead of $\kappa : \Psh{\Delta} \to \Ch$ since these two functors are quasi-isomorphic (see \cite[Theorem 2.4]{goerss2009simplicial} for example). An advantage of the latter one is that it is left Quillen. 
\end{remark}
\begin{paragr}
  We will also denote by $\sH^{\sing}$ the morphism of op-prederivators defined as the following composition
  \[
  \sH^{\sing} : \Ho(\oo\Cat^{\Th}) \overset{\overline{N_{\omega}}}{\longrightarrow} \Ho(\Psh{\Delta}) \overset{\overline{\kappa}}{\longrightarrow} \Ho(\Ch).
  \]
\end{paragr}
\begin{proposition}\label{prop:singhmlgycocontinuous}
  The singular homology \[\sH^{\sing} : \Ho(\oo\Cat^{\Th}) \to \Ho(\Ch)\] is homotopy cocontinuous. 
\end{proposition}
\begin{proof}
  This follows from the fact that $\overline{N_{\oo}}$ and $\overline{\kappa}$ are both homotopy cocontinuous. For $\overline{N_{\oo}}$, this is because it is an equivalence of op\nbd{}prederivators and thus we can apply Lemma \ref{lemma:eqisadj} and Lemma \ref{lemma:ladjcocontinuous}. For $\overline{\kappa}$, this is because $\kappa$ is left Quillen and thus we can apply Theorem \ref{thm:cisinskiII}.
\end{proof}
\section{Abelianization}
We write $\Ab$ for the category of abelian groups and for an abelian group $G$, we write $\vert G \vert$ for the underlying set of $G$.
\begin{paragr}
  Let $C$ be an $\oo$\nbd{}category. For every $n\geq 0$, we define $\lambda_n(C)$ as the abelian group obtained by quotienting $\mathbb{Z}C_n$ (the free abelian group on $C_n$) by the congruence generated by the relations
  \[
  x \comp_k y \sim x+y
  \]
  for all $x,y \in C_n$ that are $k$\nbd{}composable for some $k<n$. For $n=0$,
  this means that $\lambda_0(C)=\mathbb{Z}C_0$. Now let $f : C \to D$ be an
  $\oo$\nbd{}functor. For every $n \geq 0$, the definition of $\oo$\nbd{}functor
  implies that the linear map
  \begin{align*}
    \mathbb{Z}C_n &\to \mathbb{Z}D_{n}\\
    x \in C_n &\mapsto f(x)
  \end{align*}
  induces a linear map
  \[
  \lambda_n(f) : \lambda_n(C) \to \lambda_n(D).
  \]
  This defines a functor $\lambda_n : \oo\Cat \to \Ab$.

  For $n>0$, consider the linear map
  \begin{align*}
  \mathbb{Z}C_n &\to \mathbb{Z}C_{n-1}\\
  x \in C_n &\mapsto t(x)-s(x).
  \end{align*}
  The axioms of $\oo$\nbd{}categories imply that it induces a map
  \[
  \partial : \lambda_{n}(C) \to \lambda_{n-1}(C)
  \]
  which is natural in $C$. Furthermore, it satisfies the equation $\partial \circ \partial = 0$. Thus, for every $\oo$\nbd{}category $C$, we have defined a chain complex $\lambda(C)$:
  \[
  \lambda_0(C) \overset{\partial}{\longleftarrow} \lambda_1(C) \overset{\partial}{\longleftarrow} \lambda_2(C) \overset{\partial}{\longleftarrow} \cdots
  \]
 and for every $f : C \to D$ a morphism of chain complexes
  \[
  \lambda(f) : \lambda(C) \to \lambda(D).
  \]
  Altogether, this defines a functor
  \[
  \lambda : \omega\Cat \to \Ch,
  \]
  which we call the \emph{abelianization functor}.
\end{paragr}
\begin{lemma}\label{lemma:adjlambda}
  The functor $\lambda$ is a left adjoint.
\end{lemma}
\begin{proof}
  The category $\Ch$ is equivalent to the category $\omega\Cat(\Ab)$ of $\oo$\nbd{}categories internal to abelian groups (see \cite[Theorem 3.3]{bourn1990another}) and with this identification, the functor $\lambda : \omega\Cat \to \omega\Cat(\Ab)$ is nothing but the left adjoint of the canonical forgetful functor $\omega\Cat(\Ab) \to \omega\Cat$.
\end{proof}
As we shall now see, when the $\oo$\nbd{}category $C$ is \emph{free} the chain complex $\lambda(C)$ admits a nice expression.
\begin{paragr}
  Let $n\geq 0$. Recall that for every monoid $M$ (supposed commutative if $n \geq 1$) we have defined in Section \ref{sec:suspmonoids} an $n$\nbd{}category $B^nM$ whose set of $n$\nbd{}cells is isomorphic to the underlying set of $M$. And the correspondence $M \mapsto B^nM$ defines a functor in the expected way. By considering abelian groups as particular cases of (commutative) monoids, we obtain a functor for each $n\geq 0$ 
  \begin{align*}
    B^n : \Ab &\to n\Cat \\
    G &\mapsto B^nG,
  \end{align*}
  where $\Ab$ is the category of abelian groups.
  
  Besides, let us write $\lambda_n$ again for the functor
  \begin{align*}
    \lambda_n : n\Cat &\to \Ab\\
    C&\mapsto \lambda_n(C).
  \end{align*}
  (That is the restriction of $\lambda_n : \oo\Cat \to \Ab$ to $n\Cat$.)
\end{paragr}
\begin{lemma}\label{lemma:adjlambdasusp}
  Let $n \geq 0$. The functor $\lambda_n : n\Cat \to \Ab$ is left adjoint to the functor $B^n : \Ab \to n\Cat$.
\end{lemma}
\begin{proof}
  The case $n=0$ is immediate since the functor $\lambda_0 : 0\Cat = \Set \to \Ab$ is the ``free abelian group'' functor and the functor $B^0 : \Ab \to 0\Cat=\Set$ is the ``underlying set'' functor.

  Suppose now that $n >0$. From Lemma \ref{lemma:nfunctortomonoid} we know that for every abelian group $G$ and every $n$\nbd{}category $C$, the map
  \begin{align*}
    \Hom_{n\Cat}(C,B^nG) &\to \Hom_{\Set}(C_n,\vert G \vert)\\
    F &\mapsto F_n,
  \end{align*}
 is injective and its image consists of those functions $f : C_n \to \vert G \vert$ such that:
  \begin{enumerate}[label=(\roman*)]
  \item\label{cond:comp} for every $0 \leq k <n $ and every pair $(x,y)$ of $k$\nbd{}composable $n$\nbd{}cells of $C$, we have
    \[
    f(x\comp_ky) = f(x)+f(y),
    \]
  \item\label{cond:unit} for every $x \in C_{n-1}$, we have
    \[
    f(1_x)=0.
    \]
  \end{enumerate}
  Let us see that because $G$ is an abelian group (recall that Lemma \ref{lemma:nfunctortomonoid} was stated for the general case of commutative monoids), condition \ref{cond:comp} imply condition \ref{cond:unit}. Let $f : C_n \to \vert G \vert$ be a function that satisfies condition \ref{cond:comp} and let $x \in C_{n-1}$. We have $1_x\comp_{n-1} 1_x = 1_x$, hence
  \[
  f(1_x)=f(1_x \comp_{n-1} 1_x)=f(1_x)+f(1_x),
  \]
  and then
  \[
  f(1_x)=0
  \]
  because every element of an (abelian) group has an inverse. Now, because of the adjunction morphism
  \[
  \Hom_{\Set}(C_n,\vert G \vert) \simeq \Hom_{\Ab}(\mathbb{Z}C_n,G),
  \]
  we have that $\Hom_{n\Cat}(C,B^nG)$ is naturally isomorphic to the set of morphisms of abelian groups $g : \mathbb{Z}C_n \to G$ such that for every pair $(x,y)$ of $k$\nbd{}composable elements of $C_n$ for some $k<n$, we have
  \[
  g(x\comp_ky)=g(x)+g(y).
  \]
  By definition, this set is naturally isomorphic to the set of morphisms of abelian groups from $\lambda_n(C)$ to $G$. In other words, we have
  \[
  \Hom_{n\Cat}(C,B^nG)\simeq \Hom_{\Ab}(\lambda_n(C),G).\qedhere
  \]
\end{proof}
\begin{paragr}\label{paragr:abelpolmap}
 Let $C$ be an $\oo$\nbd{}category, $n \in \mathbb{N}$ and $E \subseteq C_n$ a subset of the $n$-cells. We obtain a map $\mathbb{Z}E \to \lambda_n(C)$ defined as the composition
  \[
  \mathbb{Z}E \to \mathbb{Z}C_n \to \lambda_n(C),
  \]
  where the map on the left is induced by the canonical inclusion of $E$ in $C_n$ and the map on the right is the quotient map from the definition of $\lambda_n(C)$. 
\end{paragr}
\begin{lemma}\label{lemma:abelpol}
  Let $C$ be a \emph{free} $\oo$\nbd{}category and let $\Sigma=(\Sigma_n)_{n \in \mathbb{N}}$ be its basis. For every $n \in \mathbb{N}$, the map
  \[
  \mathbb{Z}\Sigma_n \to \lambda_n(C)
  \]
  from the previous paragraph, is an isomorphism.
\end{lemma}
\begin{proof}

    Notice first that for every $\oo$\nbd{}category $C$, we have $\lambda_n(\tau_{\leq n}^s(C))=\lambda_n(C)$. Suppose now that $C$ is free with basis $\Sigma=(\Sigma_n)_{n \in \mathbb{N}}$. Using Lemma \ref{lemma:adjlambdasusp} and Lemma \ref{lemma:freencattomonoid}, we obtain that for every abelian group $G$, we have
        \begin{align*}
      \Hom_{\Ab}(\lambda_n(C),G) &\simeq \Hom_{\Ab}(\lambda_n(\tau_{\leq n}^s(C)),G)\\
      &\simeq \Hom_{n\Cat}(\tau_{\leq n}^s(C),B^nG)\\
      &\simeq \Hom_{\Set}(\Sigma_n,\vert G \vert)\\
      &\simeq \Hom_{\Ab}(\mathbb{Z}\Sigma_n,G),
        \end{align*}
        and it is easily checked that this isomorphism is induced by
        precomposition with the map
        $\mathbb{Z}\Sigma_n \to \lambda_n(C)$ from the previous paragraph. The
        result follows then from the Yoneda Lemma.
\end{proof}
\begin{paragr}
  Let $C$ be a \emph{free} $\oo$\nbd{}category and write $\Sigma=(\Sigma_n)_{n \in \mathbb{N}}$ for its basis. For every $n \geq 0$ and every $\alpha \in \Sigma_n$, recall that we have proved in Proposition \ref{prop:countingfunction} the existence of a unique function $w_{\alpha} : C_n \to \mathbb{N}$ such that:
  \begin{enumerate}[label=(\alph*)]
  \item\label{cond:countingfunctionfirst} $w_{\alpha}(\alpha)=1$,
  \item\label{cond:countingfunctionsecond} $w_{\alpha}(\beta)=0$ for every $\beta \in \Sigma_n$ such that $\beta\neq \alpha$,
  \item\label{cond:countingfunctionthird} for every pair of $k$\nbd{}composable $n$\nbd{}cells of $C$ for some $k<n$, we have
    \[
    w_{\alpha}(x\comp_k y)=w_{\alpha}(x) + w_{\alpha}(y).
    \]
  \end{enumerate}
  We can then define for each $n \geq 0$, a map $w_n : C_n \to \mathbb{Z}\Sigma_n$ with the formula
  \[w_n(x)=\sum_{\alpha \in \Sigma_n}w_{\alpha}(x)\cdot \alpha\]
  for every $x \in C_n$.

  Condition \ref{cond:countingfunctionthird} implies that
  \[
  w_n(x\comp_k y)=w_n(x)+w_n(y)
  \]
  for every pair $(x,y)$ of $k$\nbd{}composable $n$\nbd{}cells of $C$, and conditions \ref{cond:countingfunctionfirst} and \ref{cond:countingfunctionsecond} imply that
  \[
  w_n(\alpha)=\alpha
  \]
  for every $\alpha \in \Sigma_n$. 
\end{paragr}
\begin{proposition}\label{prop:abelianizationfreeoocat}
  Let $C$ be a free $\oo$\nbd{}category and let $(\Sigma_n)_{n \in \mathbb{N}}$ be its basis. The chain complex $\lambda(C)$ is canonically isomorphic to the chain complex
  \[
  \mathbb{Z}\Sigma_0 \overset{\partial}{\longleftarrow} \mathbb{Z}\Sigma_1 \overset{\partial}{\longleftarrow} \mathbb{Z}\Sigma_2 \overset{\partial}{\longleftarrow} \cdots
  \]
  where $\partial : \mathbb{Z}\Sigma_n \to \mathbb{Z}\Sigma_{n-1}$ is the linear map defined by the formula
  \[
  \partial(x)=w_{n-1}(\trgt(x))-w_{n-1}(\src(x))
  \]
 for every $x \in \Sigma_n$.

 With this identification, if $C'$ is another free $\oo$\nbd{}category and if $F : C \to C'$ is an $\oo$\nbd{}functor (not necessarily rigid), then the map $\lambda_n(F) : \lambda_n(C) \to \lambda_{n}(C')$ reads
 \[
 \lambda_n(F)(x)=w'_n(F(x))
 \]
 for every $x \in \Sigma_n$.
\end{proposition}
\begin{proof}
  For $n \geq 0$, write $\phi_n : \mathbb{Z}\Sigma_n \to \lambda_n(C)$ for the map defined in \ref{paragr:abelpolmap} (which we know is an isomorphism from Lemma \ref{lemma:abelpol}).

  The map $w_n: C_n \to \mathbb{Z}\Sigma_n$ induces a map $\mathbb{Z}C_n \to \mathbb{Z}\Sigma_n$ by linearity, which in turn induces a map $\lambda_n(C) \to \mathbb{Z}\Sigma_n$ (because $w_n(x \comp_k y) = w_n(x)+w_n(y)$ for every pair $(x,y)$ of $k$\nbd{}composable $n$\nbd{}cells). Write $\psi_n$ for this last map. It is immediate to check that the composition
  \[
  \mathbb{Z}\Sigma_n \overset{\phi_n}{\longrightarrow} \lambda_n(C) \overset{\psi_n}{\longrightarrow} \mathbb{Z}\Sigma_n
  \]
  gives the identity on $\mathbb{Z}\Sigma_n$. Hence, $\psi_n$ is the inverse of $\phi_n$.

  Now, for $n>0$, notice that the map $\partial : \mathbb{Z}\Sigma_n \to \mathbb{Z}\Sigma_{n-1}$ given in the statement of the proposition is nothing but the composition
  \[
  \mathbb{Z}\Sigma_n \overset{\phi_n}{\longrightarrow} \lambda_n(C) \overset{\partial}{\longrightarrow} \lambda_{n-1}(C) \overset{\psi_{n-1}}{\longrightarrow} \mathbb{Z}\Sigma_{n-1}.
  \]
  The first part of the proposition follows then from Lemma \ref{lemma:abelpol}.

  As for the second part, it suffices to notice that if we identify $\lambda_n(C)$ with $\mathbb{Z}\Sigma_n$ via $\phi_n$ for every free $\oo$\nbd{}category $C$, then map $\mathbb{Z}\Sigma_n \to \mathbb{Z}\Sigma'_n$ (where $\Sigma'_n$ is the $n$-basis of $C'$) induced by $F$ is given by the composition
  \[
  \mathbb{Z}\Sigma_n \overset{\phi_n}{\longrightarrow} \lambda_n(C) \overset{\lambda_n(F)}{\longrightarrow} \lambda_n(C') \overset{\psi_n}{\longrightarrow} \mathbb{Z}\Sigma'_n.\qedhere
  \]
 \end{proof}
 \section{Polygraphic homology}\label{section:polygraphichmlgy}
 \begin{paragr}\label{paragr:chainhmtpy}
   Let $f,g : K \to K'$ be two morphisms of non-negatively graded chain
   complexes. Recall that a \emph{chain homotopy} from $f$ to $g$ consists of a
   sequence of linear maps $(h_n \colon K_n \to K'_{n+1})_{n \in \mathbb{N}}$
   such that
   \[
     \partial \circ h_0 = g_0-f_0
   \]
   and such that for every $n > 0$, we have
   \[
     \partial \circ h_n + h_{n-1} \circ \partial = g_n - f_n.
   \]
   Recall also that if there is a chain homotopy from $f$ to $g$, then the
   localization functor $\gamma^{\Ch} : \Ch \to \ho(\Ch)$ identifies $f$ and
   $g$, which means that \[\gamma^{\Ch}(f)=\gamma^{\Ch}(g).\]
   \end{paragr}
   \begin{lemma}\label{lemma:abeloplax}
    Let $u, v : C \to D$ be two $\oo$\nbd{}functors. If there is an oplax
    transformation $\alpha : u \Rightarrow v$, then there is a chain homotopy from $\lambda(u)$ to $\lambda(v)$.
   \end{lemma}
   \begin{proof}
     For an $n$-cell $x$ of $C$ (resp.\ $D$), let us use the notation $[x]$ for the image of $x$ in $\lambda_n(C)$ (resp.\ $\lambda_n(D)$).

  Let $h_n$ be the map
     \[
     \begin{aligned}
       h_n : \lambda_n(C) &\to \lambda_{n+1}(D)\\
       [x] & \mapsto [\alpha_x].
       \end{aligned}
     \]
     The formulas for oplax transformations from Paragraph
     \ref{paragr:formulasoplax} imply that $h_n$ is linear and that for every
     $n$-cell $x$ of $C$, if $n=0$, we have
     \[
       \partial(h_0(x))=[v(x)]-[u(x)],
     \]
     and if $n>0$, we have
     \[
     \partial (h_n(x)) + h_{n-1}(\partial(x)) = [v(x)] - [u(x)].
     \]
     Details are left to the reader. 
   \end{proof}
   \begin{proposition}
     The abelianization functor $\lambda : \oo\Cat \to \Ch$ is left Quillen with respect to the folk model structure on $\oo\Cat$.
   \end{proposition}
   \begin{proof}
          The fact that $\lambda$ is a left adjoint is Lemma \ref{lemma:adjlambda}.

     A simple computation using Lemma \ref{prop:abelianizationfreeoocat} shows that for every $n\in \mathbb{N}$,
     \[
     \lambda(i_n) : \lambda(\sS_{n-1}) \to \lambda(\sD_{n})
     \]
   is a monomorphism with projective cokernel. Hence $\lambda$ sends folk cofibrations to cofibrations of chain complexes.

     Then, we know from \cite[Sections 4.6 and 4.7]{lafont2010folk} and \cite[Remarque B.1.16]{ara2016joint} (see also \cite[Paragraph 3.11]{ara2019folk}) that there exists a set of generating trivial cofibrations $J$ of the folk model structure on $\omega\Cat$ such that every $j : X \to Y$ in $J$ is a deformation retract (see Paragraph \ref{paragr:defrtract}).
     From Lemma \ref{lemma:abeloplax}, we conclude that $\lambda$ sends folk trivial cofibrations to trivial cofibrations of chain complexes.
   \end{proof}
   In particular, $\lambda$ is totally left derivable (when $\oo\Cat$ is equipped with folk weak equivalences). This motivates the following definition.
   \begin{definition}\label{de:polhom}
     The \emph{polygraphic homology functor}
     \[
     \sH^{\pol} : \ho(\oo\Cat^{\folk}) \to \ho(\Ch)
     \]
     is the total left derived functor of $\lambda : \oo\Cat \to \Ch$ (where $\oo\Cat$ is equipped with folk weak equivalences). For an $\oo$\nbd{}category $C$, $\sH^{\pol}(C)$ is the \emph{polygraphic homology of $C$}. 
   \end{definition}
   \begin{paragr}
     Similarly to singular homology groups, for $k\geq0$ the $k$\nbd{}th polygraphic homology group of an $\oo$\nbd{}category $C$ is defined as
     \[
     H^{\pol}_k(C):=H_k(\sH^{\pol}(C))
     \]
     where $H_k : \ho(\Ch) \to \Ab$ is the usual functor that associate to an object of $\ho(\Ch)$ its $k$-th homology group. In practice, this means that one has to find a cofibrant replacement of $C$, that is to say a free $\oo$\nbd{}category $P$ and a folk trivial fibration
     \[
     P \to C,
     \]
     and then the polygraphic homology groups of $C$ are those of $\lambda(P)$ which are computed using Proposition \ref{prop:abelianizationfreeoocat}.
   \end{paragr}
   \begin{paragr}
     For later reference, let us recall here that since $\sH^{\pol}$ is the left derived functor of $\lambda$, it comes equipped with a universal natural transformation (see \ref{paragr:defleftderived})
         \[
     \begin{tikzcd}
        \oo\Cat \ar[d,"\gamma^{\folk}"'] \ar[r,"\lambda"] & \Ch \ar[d,"\gamma^{\Ch}"] \\
      \ho(\oo\Cat^{\Th}) \ar[r,"\sH^{\pol}"'] & \ho(\Ch).
       \ar[from=2-1,to=1-2,"\alpha^{\pol}",shorten <= 1em, shorten >= 1em, Rightarrow]
     \end{tikzcd}
     \]
Moreover, the functor $\lambda$ being left Quillen, it is strongly derivable (Definition \ref{def:strnglyder}) and hence induces a morphism of op-prederivators, which we again denote by $\sH^{\pol}$:
     \[
     \sH^{\pol} : \Ho(\oo\Cat^{\folk}) \to \Ho(\Ch).
     \]
\end{paragr}
The following proposition is an immediate consequence of Theorem \ref{thm:cisinskiII}.
\begin{proposition}\label{prop:polhmlgycocontinuous}
  The polygraphic homology
  \[
  \sH^{\pol} : \Ho(\oo\Cat^{\folk}) \to \Ho(\Ch)
  \]
  is homotopy cocontinuous.
\end{proposition}
   As we shall now see, oplax homotopy equivalences (Definition \ref{def:oplaxhmtpyequiv}) induce isomorphisms in polygraphic homology. In order to prove that, we first need a couple of technical lemmas.
   \begin{lemma}\label{lemma:liftingoplax}
      Let
      \[
      \begin{tikzcd}
        C' \ar[r,"f_{\epsilon}'"] \ar[d,"u"] & D' \ar[d,"v"]\\
        C \ar[r,"f_{\epsilon}"] & D
      \end{tikzcd}
      \]
      be commutative squares in $\omega\Cat$ for $\epsilon\in\{0,1\}$.

      If $C'$ is a free $\omega$-category and $v$ a folk trivial fibration, then for every oplax transformation \[\alpha : f_0 \Rightarrow f_1,\] there exists an oplax transformation \[\alpha' : f_0' \Rightarrow f_1'\] such that
      \[
      v \star \alpha' = \alpha \star u.
      \]
\end{lemma}
   \begin{proof}
     Notice first that because of the natural isomorphism \[(\sD_0\amalg \sD_0) \otimes C \simeq  C \amalg C,\] we have that $\alpha : f_0 \Rightarrow f_1$ can be encoded in a functor $\alpha : \sD_1 \otimes C \to D$ such that the diagram
      \[
      \begin{tikzcd}
            (\sD_0\amalg \sD_0) \otimes C \simeq  C \amalg C \ar[d,"i_1 \otimes C"'] \ar[dr,"{\langle u, v \rangle}"] &\\
        \sD_1 \otimes C \ar[r,"\alpha"'] & D
      \end{tikzcd}
      \]
      (where $i_1 : \sD_0 \amalg \sD_0 \simeq \sS_0 \to \sD_1$ is the morphism introduced in \ref{paragr:inclusionsphereglobe}) is commutative.

      Now, the hypotheses of the lemma yield the following commutative square
      \[
      \begin{tikzcd}
        (\sD_0 \amalg \sD_0)\otimes C' \ar[d,"{i_1\otimes C'}"'] \ar[rr,"{\langle f'_0, f_1' \rangle}"] && D' \ar[d,"v"] \\
        \sD_1\otimes C'\ar[r,"\sD_1 \otimes u"'] & \sD_1\otimes C \ar[r,"\alpha"] & D
        \end{tikzcd}
      \]
      and since $i_1$ is a folk cofibration and $C'$ is cofibrant, it follows
      that the left vertical morphism of the previous square is a folk
      cofibration (see \cite[Proposition 5.1.2.7]{lucas2017cubical} or
      \cite{ara2019folk}). By hypothesis, $v$ is a folk trivial fibration, and
      so the above square admits a lift
      \[
      \alpha' : \sD_1\otimes C' \to D'.
      \]
      The commutativity of the two induced triangles shows what we needed to prove.
   \end{proof}
From now on, for an $\oo$\nbd{}functor $u$, we write $\sH^{\pol}(u)$ instead of $\sH^{\pol}(\gamma^{\folk}(u))$ (where $\gamma^{\folk}$ is the localization functor $\oo\Cat \to \ho(\oo\Cat^{\folk})$) for the morphism induced by $u$ at the level of polygraphic homology.
\begin{lemma}\label{lemma:oplaxpolhmlgy}
  Let $u,v : C \to D$ be two $\oo$\nbd{}functors. If there exists an oplax transformation $u\Rightarrow v$, then
  \[
  \sH^{\pol}(u)=\sH^{\pol}(v).
  \]
\end{lemma}
\begin{proof}
  In the case that $C$ and $D$ are both folk cofibrant, this follows immediately
  from Lemma \ref{lemma:abeloplax} and the fact that the localization functor
  $\Ch \to \ho(\Ch)$ identifies chain homotopic maps (\ref{paragr:chainhmtpy}).

  In the general case, let
  \[
  p : C' \to C
  \]
  and
  \[
  q : D' \to D
  \]
  be folk trivial fibrations with $C'$ and $D'$ cofibrant. Using that $q$ is a trivial fibration and $C'$ is cofibrant, we know that there exist $u' : C' \to D'$ and $v' : C' \to D'$ such that the squares
  \[
  \begin{tikzcd}
    C' \ar[d,"p"] \ar[r,"u'"] & D' \ar[d,"q"] \\
    C \ar[r,"u"] & D
  \end{tikzcd}
  \text{ and }
  \begin{tikzcd}
    C' \ar[d,"p"] \ar[r,"v'"] & D' \ar[d,"q"] \\
    C \ar[r,"v"] & D
  \end{tikzcd}  
  \]
  are commutative. From Lemma \ref{lemma:liftingoplax}, we deduce the existence of an oplax transformation $u' \Rightarrow v'$. Since $C'$ and $D'$ are cofibrant, we have already proved that
  \[\sH^{\pol}(u')=\sH^{\pol}(v').\]
  The commutativity of the two previous squares and the fact that $p$ and $q$ are folk weak equivalences imply the desired result.
\end{proof}
The following proposition is an immediate consequence of the previous lemma.
\begin{proposition}\label{prop:oplaxhmtpypolhmlgy}
  Let $u : C \to D$ be an $\oo$\nbd{}functor. If $u$ is an oplax homotopy equivalence, then the induced morphism
  \[
  \sH^{\pol}(u) : \sH^{\pol}(C) \to \sH^{\pol}(D)
  \]
  is an isomorphism.
\end{proposition}
\begin{paragr}\label{paragr:polhmlgythomeq}
  Oplax homotopy equivalences being particular cases of Thomason equivalences, one may wonder whether it is true that \emph{every} Thomason equivalence induce an isomorphism in polygraphic homology. As we shall see later (Proposition \ref{prop:polhmlgynotinvariant}),  it is not the case.
\end{paragr}
\begin{remark}
  Lemma \ref{lemma:liftingoplax}, Lemma \ref{lemma:oplaxpolhmlgy} and Proposition \ref{prop:oplaxhmtpypolhmlgy} are also true if we replace ``oplax'' by ``lax'' everywhere. 
\end{remark}

   \section{Singular homology as derived abelianization}\label{section:singhmlgyderived}
   We have seen in the previous section that the polygraphic homology functor is the total left derived functor of $\lambda : \oo\Cat \to \Ch$  when $\oo\Cat$ is equipped with the folk weak equivalences. As it turns out, the abelianization functor is also totally left derivable when $\oo\Cat$ is equipped with the Thomason equivalences and the total left derived functor is the singular homology functor. In order to prove this result, we first need a few technical lemmas.
   \begin{lemma}\label{lemma:nuhomotopical}
     Let $\nu : \Ch \to \oo\Cat$ be the right adjoint of the abelianization
     functor (see Lemma \ref{lemma:adjlambda}). This functor sends the quasi-isomorphisms to Thomason equivalences.
   \end{lemma}
   \begin{proof}
     We have already seen that $\lambda : \oo\Cat \to \Ch$ is left Quillen with respect to the folk model structure on $\oo\Cat$. By adjunction, this means that $\nu$ is right Quillen for this model structure. In particular, it sends trivial fibrations of chain complexes to folk trivial fibrations. From Ken Brown's Lemma \cite[Lemma 1.1.12]{hovey2007model} and the fact that all chain complexes are fibrant, it follows that $\nu$ sends weak equivalences of chain complexes to weak equivalences of the folk model structure, which are in particular Thomason equivalences (Lemma \ref{lemma:nervehomotopical}).
   \end{proof}
   \begin{remark}
     The proof of the previous lemma shows the stronger result that $\nu$ sends
     the quasi-ismorphisms to folk weak equivalences. This will be of no use in the sequel.
   \end{remark}
   Recall that we write $c_{\oo} : \Psh{\Delta} \to \oo\Cat$ for the left adjoint of the nerve functor ${N_{\oo} : \oo\Cat \to \Psh{\Delta}}$ (see Paragraph \ref{paragr:nerve}).
   \begin{lemma}\label{lemma:abelor}
    The triangle of functors
    \[
    \begin{tikzcd}
      \Psh{\Delta} \ar[r,"c_{\oo}"] \ar[dr,"\kappa"']& \oo\Cat\ar[d,"\lambda"]\\
      &\Ch
    \end{tikzcd}
    \]
    is commutative (up to a canonical isomorphism).
   \end{lemma}
   \begin{proof}
     All the functors involved are cocontinuous, hence it suffices to prove that the triangle is commutative  when pre-composed by the Yoneda embedding $\Delta \to \Psh{\Delta}$. This follows immediately from the description of the orientals in \cite{steiner2004omega}.
   \end{proof}
   Recall now that the notions of adjunction and equivalence are valid in every $2$\nbd{}category and in particular in the $2$\nbd{}category of pre-derivators (see \ref{paragr:prederequivadjun}). We omit the proof of the following lemma, which is the same as when the ambient $2$\nbd{}category is the $2$\nbd{}category of categories.
    \begin{lemma}\label{lemma:adjeq}
     Let $\begin{tikzcd} f : y \ar[r,shift left]&z :g\ar[l,shift left] \end{tikzcd}$ be an adjunction and $h : x \to y$ an equivalence with quasi-inverse $k : y \to x$. Then $fh$ is left adjoint to $kg$.
   \end{lemma}
   We can now state and prove the promised result.
   \begin{theorem}\label{thm:hmlgyderived}
     Consider that $\oo\Cat$ is equipped with the Thomason equivalences. The abelianization functor $\lambda : \oo\Cat \to \Ch$ is strongly left derivable and the left derived morphism of op\nbd{}prederivators
     \[
     \LL \lambda^{\Th} : \Ho(\oo\Cat^{\Th}) \to \Ho(\Ch)
     \]
     is isomorphic to the singular homology
     \[
     \sH^{\sing} : \Ho(\oo\Cat^{\Th}) \to \Ho(\Ch).
     \]
   \end{theorem}
   \begin{proof}
     Let $\nu$ be the right adjoint of the abelianization functor (see Lemma \ref{lemma:adjlambda}) and consider the following adjunctions
     \[
     \begin{tikzcd}
       \Psh{\Delta} \ar[r,shift left,"c_{\omega}"] &\oo\Cat \ar[l,shift left,"N_{\omega}"] \ar[r,shift left,"\lambda"]& \Ch \ar[l,shift left,"\nu"]
       \end{tikzcd}
     \]
     where the functors from left to right are the left adjoints.
      We know that:
     \begin{itemize}[label=-]
       \item The functor $\nu$ induces a morphism of localizers
     \[
     \nu : (\Ch,\W_{\Ch}) \to (\oo\Cat,\W^{\Th}),
     \]
      thanks to Lemma \ref{lemma:nuhomotopical}.
   \item The functor $N_{\omega}$ induces a morphism of localizers
     \[
     N_{\omega} : (\oo\Cat,\W^{\Th}) \to (\Psh{\Delta},\W_{\Delta}),
     \]
     by definition of Thomason equivalences.
   \item There is an isomorphism of functors $\lambda c_{\omega} \simeq \kappa$ (Lemma \ref{lemma:abelor}), hence an induced morphism of localizers
     \[
     (\lambda c_{\omega})\simeq \kappa : (\Psh{\Delta},\W_{\Delta}) \to (\Ch,\W_{\Ch}),
     \]
     thanks to Lemma \ref{lemma:normcompquil}.
     \end{itemize}
     It follows that there is an induced adjunction at the level of op-prederivators:
     \[
     \begin{tikzcd}
       \overline{\kappa} \simeq \overline{\lambda c_{\omega}} :  \Ho(\Psh{\Delta}) \ar[r,shift left] &  \ar[l,shift left] \Ho(\Ch) :\overline{N_{\omega}}\overline{\nu}.
       \end{tikzcd}
     \]
     Now, we know from Theorem \ref{thm:gagna} that $\overline{N_{\omega}}$ is an equivalence of op\nbd{}prederivators, and thus admits a quasi-inverse. Let $ M : \Ho(\Psh{\Delta}) \to \Ho(\oo\Cat)$ be such a quasi-inverse. From Lemma \ref{lemma:adjeq}, we deduce that we have an adjunction:
     \[
     \begin{tikzcd}
       \overline{\kappa} \overline{N_{\oo}} : \Ho(\oo\Cat^{\Th}) \ar[r,shift left]& \ar[l,shift left] \Ho(\Ch) : M \overline{N_{\oo}} \overline{\nu} \simeq \overline{\nu}.
       \end{tikzcd}
     \]
     From Proposition \ref{prop:gonzalezcritder}, we  conclude that $\lambda : \oo\Cat \to \Ch$ is strongly left derivable and that $ \LL\lambda^{\Th}  \simeq \overline{\kappa} \overline{N_{\oo}}$, which is, by definition, the singular homology. 
   \end{proof}
   \begin{remark}
     Beware that neither $c_{\oo} : \Psh{\Delta} \to \oo\Cat$ sends all weak
     equivalences of simplicial sets to Thomason equivalences nor $\lambda :
     \oo\Cat \to \Ch$ sends all Thomason equivalences to quasi-isomorphisms.
     However, this does not contradict the fact that $\lambda c_{\oo} : \Psh{\Delta} \to \Ch$ does send all weak equivalences of simplicial sets to quasi-isomorphisms.
   \end{remark}
   \begin{paragr}\label{paragr:univmor}
     Since $\sH^{\sing} : \ho(\oo\Cat^{\Th}) \to \ho(\Ch)$ is the left derived
     functor of the abelianization functor, it comes with a universal natural transformation
     \[
     \begin{tikzcd}
        \oo\Cat \ar[d,"\gamma^{\Th}"'] \ar[r,"\lambda"] & \Ch \ar[d,"\gamma^{\Ch}"] \\
      \ho(\oo\Cat^{\Th}) \ar[r,"\sH^{\sing}"'] & \ho(\Ch).
       \ar[from=2-1,to=1-2,"\alpha^{\sing}",shorten <= 1em, shorten >= 1em, Rightarrow]
     \end{tikzcd}
     \]
     A thorough reading of the proofs of Proposition \ref{prop:gonzalezcritder}
     and Theorem \ref{thm:hmlgyderived} enables us to give the following
     description of $\alpha^{\sing}$. By post-composing the co-unit of the
     adjunction $c_{\oo} \dashv N_{\oo}$ with the abelianization functor, we
     obtain a natural transformation
     \[
     \lambda c_{\oo} N_{\oo} \Rightarrow \lambda.
     \]
     Then $\alpha^{\sing}$ is nothing but the following composition of natural transformations
     \[
     \begin{tikzcd}[column sep=huge]
       \oo\Cat \ar[d,"\gamma^{\Th}"]\ar[r,bend left,"\lambda",""{name=A,below}] \ar[r,"\lambda c_{\oo} N_{\oo}"',""{name=B,above}] & \Ch \ar[d,"\gamma^{\Ch}"] \\
        \ho(\oo\Cat^{\Th}) \ar[r,"\sH^{\sing}"'] & \ho(\Ch),
        \ar[from=B,to=A,Rightarrow]\ar[from=1-1,to=2-2,phantom,"\simeq" description]
     \end{tikzcd}
     \]
     where the square is commutative (up to an isomorphism) because $\sH^{\sing}\simeq\overline{\lambda c_{\oo}} \overline{N_{\oo}}$.
   \end{paragr}
   \section{Comparing homologies}
   \begin{paragr}\label{paragr:cmparisonmap}
     Recall from Proposition \ref{prop:folkisthom} that the identity functor on $\oo\Cat$ induces a morphism of localizers
     \[(\oo\Cat,\W^{\folk}) \to (\oo\Cat,\W^{\Th}),\]
     which in turn induces a functor
     \[\J : \ho(\oo\Cat^{\folk}) \to \ho(\oo\Cat^{\Th}),\]
     such that
     \[
     \gamma^{\Th} = \J \circ \gamma^{\folk}.
     \]
     Now, consider the following triangle 
     \begin{equation}\label{cmprisontrngle}
     \begin{tikzcd}
       \ho(\oo\Cat^{\folk}) \ar[r,"\J"] \ar[rd,"\sH^{\pol}"'] & \ho(\oo\Cat^{\Th}) \ar[d,"\sH^{\sing}"] \\
      & \ho(\Ch).
     \end{tikzcd}
     \end{equation}
     A natural question to ask is whether this triangle is commutative (up to an isomorphism). Since $\J$ is the identity on objects, this amounts to ask whether for every $\oo$\nbd{}category $C$ we have an isomorphism (natural in $C$)
     \[
     \sH^{\pol}(C)\simeq \sH^{\sing}(C).
     \]
     As it happens, this is not possible as the following counter-example, due to Ara and Maltsiniotis, shows.
       \end{paragr}
\begin{paragr}[Ara and Maltsiniotis' counter-example]\label{paragr:bubble}
  Write $\mathbb{N}=(\mathbb{N},+,0)$ for the commutative monoid of non-negative integers  and let $C$ be the $2$\nbd{}category defined as
  \[
  C:=B^2\mathbb{N}
  \]
  (see \ref{paragr:suspmonoid}). As usual, we consider $C$ as an $\oo$\nbd{}category with only unit cells strictly above dimension $2$. This $\oo$\nbd{}category is free; namely its $k$\nbd{}basis is a singleton for $k=0$ and $k=2$, and the empty set otherwise. In particular $C$ is cofibrant for the folk model structure (Proposition \ref{prop:freeiscofibrant}) and it follows from Proposition \ref{prop:abelianizationfreeoocat} that $\sH^{\pol}(C)$ is given by the chain complex (seen as an object of $\ho(\Ch)$)
       \[
     \begin{tikzcd}[column sep=small]
       \mathbb{Z} & 0 \ar[l] & \ar[l] \mathbb{Z} & \ar[l] 0 & \ar[l] 0 & \ar[l] \cdots
       \end{tikzcd}
     \]
     Hence, the polygraphic homology groups of $B$ are given by
     \[
     H^{\pol}_k(C)=\begin{cases} \mathbb{Z} \text{ if } k=0,2\\ 0 \text{ in other cases.}\end{cases}  
     \]
     On the other hand, it is proven in \cite[Theorem 4.9 and Example 4.10]{ara2019quillen} that (the nerve of) $C$ is a $K(\mathbb{Z},2)$. In particular, it has non-trivial singular homology groups in every even dimension. This proves that $\sH^{\pol}(C)$ is \emph{not} isomorphic to $\sH^{\sing}(C)$; which means that triangle \eqref{cmprisontrngle} cannot be commutative (up to an isomorphism). 
\end{paragr}
Another consequence of the above counter-example is the following result, which we claimed in \ref{paragr:polhmlgythomeq}. Recall that given a morphism $u : C \to D$ of $\oo\Cat$, we write $\sH^{\pol}(u)$ instead of $\sH^{\pol}(\gamma^{\folk}(u))$.
\begin{proposition}\label{prop:polhmlgynotinvariant}
  There exists at least one Thomason equivalence \[u : C \to D\] such that the induced morphism
  \[
  \sH^{\pol}(u) : \sH^{\pol}(C) \to \sH^{\pol}(D)
  \]
  is not an isomorphism of $\ho(\Ch)$.
\end{proposition}
\begin{proof}
  Suppose the converse, which is that the functor
  \[
  \sH^{\pol} \circ \gamma^{\folk} : \oo\Cat \to \ho(\Ch)
  \]sends the Thomason equivalences to isomorphisms of $\ho(\Ch)$. Because of the inclusion $\W^{\folk} \subseteq \W^{\Th}_{\oo}$, the category $\ho(\oo\Cat^{\Th})$ may be identified with the localization of $\ho(\oo\Cat^{\folk})$ with respect to $\gamma^{\folk}(\W^{\Th}_{\oo})$ and then the localization functor is nothing but
  \[
  \J : \ho(\oo\Cat^{\folk}) \to \ho(\oo\Cat^{\Th}).
  \]
  From this observation and because of the hypothesis we made on Thomason equivalences inducing isomorphisms in polygraphic homology, we deduce the existence of a functor
  \[
  \overline{\sH^{\pol}} : \ho(\oo\Cat^{\Th}) \to \ho(\Ch)
  \]
  such that we have
  \[
   \overline{\sH^{\pol}}\circ \J = \sH^{\pol},
   \]
   and because of the equality $\gamma^{\Th} = \J \circ \gamma^{\folk}$, the universal natural transformation $\alpha^{\pol}$ now reads
   \[
   \begin{tikzcd}
     \oo\Cat \ar[r,"\lambda"] \ar[d,"\gamma^{\Th}"] & \Ch \ar[d,"\gamma^{\Ch}"] \\
     \ho(\oo\Cat^{\Th}) \ar[r,"\overline{\sH^{\pol}}"'] & \ho(\Ch).
     \ar[from=2-1,to=1-2,"\alpha^{\pol}",shorten <= 1em, shorten >=1em,Rightarrow]
   \end{tikzcd}
   \]
   Let us show that $(\overline{\sH^{\pol}},\alpha^{\pol})$ is the left derived functor of $\lambda$ when $\oo\Cat$ is equipped with the Thomason equivalences. Let $G$ and $\beta$ be as in the following $2$\nbd{}diagram
      \[
   \begin{tikzcd}
     \oo\Cat \ar[r,"\lambda"] \ar[d,"\gamma^{\Th}=\J\circ \gamma^{\folk}"'] & \Ch \ar[d,"\gamma^{\Ch}"] \\
     \ho(\oo\Cat^{\Th}) \ar[r,"G"'] & \ho(\Ch).
     \ar[from=2-1,to=1-2,"\beta",shorten <= 1em, shorten >=1em,Rightarrow]
   \end{tikzcd}
   \]
   Since $\sH^{\pol}$ is the left derived functor of $\lambda$ when $\oo\Cat$ is equipped with the folk weak equivalences, there exists a unique $\delta : G \circ \J \Rightarrow \sH^{\pol}$ that factorizes $\beta$ as
   \[
      \begin{tikzcd}
     \oo\Cat \ar[r,"\lambda"] \ar[d,"\gamma^{\folk}"] & \Ch \ar[d,"\gamma^{\Ch}"] \\
     \ho(\oo\Cat^{\folk}) \ar[d,"\J"] \ar[r,"\sH^{\pol}"'] & \ho(\Ch)\\
     \ho(\oo\Cat^{\Th}) \ar[ru,"G"',bend right] &.
     \ar[from=2-1,to=1-2,"\alpha^{\pol}",shorten <= 1em, shorten >=1em,Rightarrow]\\
     \ar[from=3-1,to=2-2,"\delta"',shorten <= 1em, shorten >= 1em,Rightarrow]
   \end{tikzcd}
      \]
      But since $\J$ acts as a localization functor, $\delta$ also factorizes uniquely as
      \[
      \begin{tikzcd}[column sep=small] \ho(\oo\Cat^{\folk}) \ar[r,"\J"] & \ho(\oo\Cat^{\Th}) \ar[r,bend left,"\overline{\sH^{\pol}}",""{name=A,below}] \ar[r,bend right, "G"',pos=16/30,""{name=B,above}] & \ho(\Ch). \ar[from=B,to=A,Rightarrow,"\delta'"]\end{tikzcd}
      \]
      Altogether we have that $\beta$ factorizes as
      \[
   \begin{tikzcd}
     \oo\Cat \ar[r,"\lambda"] \ar[d,"\gamma^{\Th}"] & \Ch \ar[d,"\gamma^{\Ch}"] \\
     \ho(\oo\Cat^{\Th}) \ar[r,"\overline{\sH^{\pol}}",""{name=B,below}] & \ho(\Ch).
     \ar[from=2-1,to=1-2,"\alpha^{\pol}",shorten <= 1em, shorten >=1em,Rightarrow]
     \ar[from=2-1,to=2-2,"G"',pos=16/30,bend right,""{name=A,above}]
     \ar[from=A,to=B,Rightarrow,"\delta'"]
   \end{tikzcd}
   \]
   The uniqueness of such a factorization follows from a similar argument which is left to the reader. This proves that $\overline{\sH^{\pol}}$ is the left derived functor of $\lambda$ when $\oo\Cat$ is equipped with the Thomason equivalences and in particular we have
   \[
   \sH^{\sing}\simeq \overline{\sH^{\pol}}.
   \]
   But since $\J$ is the identity on objects, this implies that for every $\oo$\nbd{}category $C$ we have
   \[
   \sH^{\sing}(C)\simeq \overline{\sH^{\pol}}(C)=\sH^{\pol}(C),
   \]
   which we know is impossible.
 \end{proof}
 \begin{remark}\label{remark:polhmlgyisnotinvariant}
   It follows from the previous result that if we think of $\oo$\nbd{}categories as a
   model for homotopy types (see Theorem \ref{thm:gagna}), then the polygraphic
   homology of an $\oo$\nbd{}category is \emph{not} a well defined invariant. This
   justifies what we said in remark \ref{remark:singularhmlgyishmlgy}, which is
   that \emph{singular homology} is the only ``correct'' homology of $\oo$\nbd{}categories.
 \end{remark}
\begin{paragr}\label{paragr:defcancompmap}
  Even though triangle \eqref{cmprisontrngle} is not commutative (even up to an isomorphism), it can be filled up with a $2$\nbd{}morphism. Indeed, consider the following $2$\nbd{}square
     \[
     \begin{tikzcd}
       \oo\Cat \ar[d,"\gamma^{\Th}"] \ar[r,"\lambda"] & \Ch \ar[d,"\gamma^{\Ch}"] \\
       \ho(\oo\Cat^{\Th}) \ar[r,"\sH^{\sing}"] & \ho(\Ch),
       \ar[from=2-1,to=1-2,"\alpha^{\sing}",shorten <= 1em, shorten >= 1em, Rightarrow]
     \end{tikzcd}
     \]
     Since $\gamma^{\Th}=\J\circ \gamma^{\folk}$ and the polygraphic homology is
     the total left derived functor of the abelianization functor when $\oo\Cat$
     is equipped with the folk weak equivalences, we obtain by universal property (see \ref{paragr:defleftderived}) a unique natural transformation
     \begin{equation}\label{cmparisonmapdiag}
     \begin{tikzcd}
       \ho(\oo\Cat^{\folk}) \ar[d,"\J"] \ar[rd,"\sH^{\pol}",""{name=A,below}] & \\
     \ho(\oo\Cat^{\Th}) \ar[r,"\sH^{\sing}"']  & \ho(\Ch)\ar[from=2-1,to=A,"\pi",Rightarrow]
     \end{tikzcd}
     \end{equation}
     such that $\alpha^{\sing}$ factorizes as 
          \[
     \begin{tikzcd}
       \oo\Cat \ar[d,"\gamma^{\folk}"] \ar[r,"\lambda"] & \Ch \ar[d,"\gamma^{\Ch}"] \\
       \ho(\oo\Cat^{\folk})\ar[d,"\J"] \ar[r,"\sH^{\pol}",""{name=B,below}] & \ho(\Ch)\\
       \ho(\oo\Cat^{\Th}) \ar[ru,"\sH^{\sing}"',bend right=15] &
       \ar[from=2-1,to=1-2,"\alpha^{\pol}",shorten <= 1em, shorten >= 1em, Rightarrow]
       \ar[from=3-1,to=B,Rightarrow,"\pi",shorten <= 1em, shorten >= 1em]
     \end{tikzcd}
     \]
   
     Since $\J$ is nothing but the identity on objects, for every $\oo$\nbd{}category $C$, the natural transformation $\pi$ yields a map 
     \[
     \pi_C : \sH^{\sing}(C) \to \sH^{\pol}(C),
     \]
     which we shall refer to as the \emph{canonical comparison map.}
\end{paragr}
   \begin{remark}
     When $C$ is free, it follows from the considerations in \ref{paragr:univmor} that the canonical comparison map $\pi_C$ can be identified with the image by $\gamma^{\Ch}$ of the morphism of $\Ch$
     \[
    \lambda c_{\oo}N_{\oo}(C) \to \lambda(C)
     \]
     induced by the co-unit of $c_{\oo} \dashv N_{\oo}$.
   \end{remark}
   \begin{definition}
     An $\oo$\nbd{}category $C$ is said to be \emph{\good{}} when the canonical comparison map
     \[
     \pi_C : \sH^{\sing}(C) \to \sH^{\pol}(C)
     \]
     is an isomorphism of $\ho(\Ch)$.
     \end{definition}

  \begin{paragr}
The rest of this dissertation is devoted to the study of \good{} $\oo$\nbd{}categories. Examples of such $\oo$\nbd{}categories will be presented later. Following the perspective of Remark \ref{remark:polhmlgyisnotinvariant}, polygraphic homology can be thought of as a way to compute singular homology of \good{} $\oo$\nbd{}categories. 
     \end{paragr}

\begin{paragr}\label{paragr:prelimcriteriongoodcat}
  Similarly to \ref{paragr:cmparisonmap}, the morphism of localizers
  \[
 (\oo\Cat,\W^{\folk}) \to (\oo\Cat,\W^{\Th})
 \]
  induces a morphism of op-prederivators
  \[
 \J : \Ho(\oo\Cat^{\folk}) \to \Ho(\oo\Cat^{\Th})
 \]
  such that the triangle in the category of op-prederivators
  \[
  \begin{tikzcd}
    \oo\Cat \ar[d,"\gamma^{\folk}"'] \ar[rd,"\gamma^{\Th}"] &\\
    \Ho(\oo\Cat^{\folk}) \ar[r,"\J"'] & \Ho(\oo\Cat^{\Th})
  \end{tikzcd}
  \]
  is commutative.
  It follows from what we said in \ref{paragr:compweakeq} that the morphism $\J$ cannot be an equivalence of op-prederivators. As we shall see later, $\J$ is not even homotopy cocontinuous. In particular, this implies that given a diagram $d : I \to \oo\Cat$, the canonical arrow of $\ho(\oo\Cat^{\Th})$
  \begin{equation}\label{equation:Jhocolim}
  \hocolim^{\Th}_{I}(\J_I(d)) \to   \J_e(\hocolim_{I}^{\folk}(d))
  \end{equation}
  induced by $\J$ (see \ref{paragr:canmorphismcolimit}) is generally \emph{not} an isomorphism. Note that since
  \[\J_I : \Ho(\oo\Cat^{\folk})(I) \to \Ho(\oo\Cat^{\Th})(I)\]
  is the identity on objects for every small category $I$, morphism \eqref{equation:Jhocolim} simply reads
  \[
  \hocolim_I^{\Th}(d) \to \hocolim_I^{\folk}(d).
  \]
  Even if this is not always true, there are some particular diagrams $d$ for which the above morphism is indeed an isomorphism. The criterion to find \good{} $\oo$\nbd{}categories given in the proposition below is based on this observation.
\end{paragr}
\begin{proposition}\label{prop:criteriongoodcat}
  Let $C$ be an $\oo$\nbd{}category. Suppose that there exists a diagram
  \[
  d : I \to \oo\Cat
  \]
  and a co-cone
  \[
  (\varphi_i : d(i) \to C)_{i \in \Ob(I)}
  \]
  such that:
  \begin{enumerate}[label=(\roman*)]
    \item For every $i \in \Ob(I)$, the $\oo$\nbd{}category $d(i)$ is \good{}.
  \item The canonical morphism
    \[
    \hocolim^{\folk}d \to C
    \]
    is an isomorphism of $\ho(\oo\Cat^{\folk})$.
    \item The canonical morphism
      \[
      \hocolim^{\Th}d \to C
      \]
      is an isomorphism of $\ho(\oo\Cat^{\Th})$.
  \end{enumerate}
  Then the $\oo$\nbd{}category $C$ is \good{}.
\end{proposition}

\begin{proof}
  Notice first that all the constructions from \ref{paragr:defcancompmap} may be reproduced \emph{mutatis mutandis} at the level of op-prederivators. In particular, we obtain a $2$\nbd{}morphism of op-prederivators
  \[
     \begin{tikzcd}
      \Ho(\oo\Cat^{\folk}) \ar[d,"\J"] \ar[rd,"\sH^{\pol}",""{name=A,below}] & \\
      \Ho(\oo\Cat^{\Th}) \ar[r,"\sH^{\sing}"']  & \Ho(\Ch)\ar[from=2-1,to=A,"\pi",Rightarrow].
    \end{tikzcd}
  \]
  Then, by naturality, we have a commutative diagram in $\ho(\Ch)$:
  \[
  \begin{tikzcd}
    \displaystyle\hocolim_{i\in I}\sH^{\sing}(d_i) \ar[d] \ar[r] & \displaystyle\hocolim_{i \in I}\sH^{\pol}(d_i) \ar[d] \\
    \displaystyle\sH^{\sing}(\hocolim_{i \in I}^{\Th}(d_i)) \ar[r] \ar[d] & \displaystyle\sH^{\pol}(\hocolim^{\folk}_{i \in I}(d_i))\ar[d] \\
    \sH^{\sing}(C) \ar[r] & \sH^{\pol}(C),
  \end{tikzcd}
  \]
  where:
  \begin{itemize}[label=-]
  \item the top and bottom horizontal arrows are induced by $\pi$,
    \item the middle horizontal arrow is induced by $\pi$ and the canonical morphism \[\hocolim_{i \in I}^{\Th}(d_i)\to \hocolim_{i \in I}^{\folk}(d_i)\] from \ref{paragr:prelimcriteriongoodcat},
  \item the top vertical arrows are the canonical morphisms induced by every morphism of op-prederivators (see \ref{paragr:canmorphismcolimit}),
    \item the bottom vertical arrows are induced by the co-cone  \[(\varphi_i : d(i) \to C)_{i \in \Ob(I)}.\]
  \end{itemize}
  Since $\sH^{\pol}$ and $\sH^{\sing}$ are both homotopy cocontinuous (Proposition \ref{prop:singhmlgycocontinuous} and Proposition \ref{prop:polhmlgycocontinuous} respectively), both top vertical arrows are isomorphisms. Hypotheses $(ii)$ and $(iii)$ imply that the bottom vertical arrows are isomorphisms and hypothesis $(i)$ imply that the top horizontal arrow is an isomorphism. By a 2-out-of-3 property, the bottom horizontal arrow is an isomorphism, which means exactly that $C$ is \good{}.
\end{proof}
The previous proposition admits the following corollary, which will be of great use in later chapters.
\begin{corollary}\label{cor:usefulcriterion}
  Let
  \[
  \begin{tikzcd}
  A \ar[r,"u"] \ar[d,"f"] & B \ar[d,"g"] \\
  C \ar[r,"v"] & D
  \ar[from=1-1,to=2-2,phantom,"\ulcorner",very near end]
  \end{tikzcd}
  \]
  be a cocartesian square of $\oo\Cat$ such that:
  \begin{enumerate}[label=(\alph*)]
  \item the $\oo$\nbd{}categories $A$,$B$ and $C$ are free and \good{},
  \item at least one of the morphisms $u : A \to B$ or $f : A \to C$ is a folk cofibration,
  \item the square is Thomason homotopy cocartesian.
  \end{enumerate}
  Then, the $\oo$\nbd{}category $D$ is \good{}.
\end{corollary}
\begin{proof}
  The fact that $A$,$B$ and $C$ are free and one of the morphism $u$ or $f$ is a folk cofibration ensure that the square is folk homotopy cocartesian (Lemma \ref{lemma:hmtpycocartesianreedy}). The conclusion follows then from Proposition \ref{prop:criteriongoodcat}.
\end{proof}
\section{Equivalence of homologies in low dimension}
\begin{paragr}
  Recall that for every $n \geq 0$ we have taken the habit of identifying $n\Cat$ as a full subcategory of $\oo\Cat$ via the canonical fully faithful functor $\iota_n : n\Cat \to \oo\Cat$ (defined in \ref{paragr:defncat}) that sends an $n$\nbd{}category $C$ to the $\oo$\nbd{}category with the same $k$\nbd{}cells as $C$ for $k\leq n$ and only unit cells for $k > n$. In particular, we abusively wrote
  \[
  C=\iota_n(C).
  \]
  Within this section, \emph{and only within this section}, we try not to make this abuse of notation and explicitly write $\iota_n$ whenever we should.
  
  We have already seen that $\iota_n$ has a left adjoint $\tau^{i}_{\leq n} : \oo\Cat \to n\Cat$, where for an $\oo$\nbd{}category $C$, $\tau_{\leq n }^{i}(C)$ is the $n$\nbd{}category whose set of $k$\nbd{}cells is $C_k$ for $k<n$ and whose set of $n$\nbd{}cells is the quotient of $C_n$ by the equivalence relation $\sim$ generated by
  \[
  x \sim y \text{ when there exists } z : x \to y \text{ in } C_{n+1}.
\]
  As it happens, we can use the adjunction $\tau_{\leq n}^{i} \dashv \iota_n$ to transport the folk model structure from $\oo\Cat$ to $n\Cat$.
\end{paragr}
\begin{proposition}\label{prop:fmsncat}
  There exists a model structure on $n\Cat$ such that:
  \begin{itemize}[label=-]
  \item the weak equivalences are exactly those morphisms $f : C \to D$ such that $\iota_n(f)$ is a weak equivalence for the folk model structure on $\oo\Cat$,
  \item the fibrations are exactly those morphisms $f : C \to D$ such that $\iota_n(f)$ is a fibrations for the folk model structure on $\oo\Cat$.
  \end{itemize}
  Moreover, there exists a set $I$ of generating cofibrations (resp.\ a set $J$ of generating trivial cofibrations) for the folk model structure on $\oo\Cat$ such that the image by $\tau^{i}_{\leq n}$ of $I$ (resp.\ $J$) is a set of generating cofibrations (resp.\ generating trivial cofibrations) of the above model structure on $n\Cat$.
\end{proposition}
\begin{proof}
This is \cite[Theorem 5]{lafont2010folk}. (Although the part concerning generating cofibrations is not made explicit in the statement of the theorem, it is contained in proof.)
\end{proof}
\begin{paragr}
  We refer to the model structure of the above proposition as the \emph{folk
    model structure on $n\Cat$}. By definition, the functor $\iota_n : n\Cat \to
  \oo\Cat$ preserves weak equivalences and fibrations when $\oo\Cat$ and $n\Cat$
  are equipped with the folk model structure. In particular, the adjunction
  $\tau^{i}_{\leq n } \dashv \iota_n$ is a Quillen adjunction. As it happens,
  the functor $\tau^{i}_{\leq n}$ also preserves weak equivalences. 
\end{paragr}
\begin{proposition}\label{prop:truncationhomotopical}
  The functor $\tau^{i}_{\leq n} : \oo\Cat \to n\Cat$ sends the weak equivalences of the folk model structure on $\oo\Cat$ to weak equivalences of the folk model structure on $n\Cat$.
\end{proposition}
\begin{proof}
  Since every $\oo$\nbd{}category is fibrant for the folk model structure on
  $\oo\Cat$ \cite[Proposition 9]{lafont2010folk}, it suffices to show that
  $\tau^{i}_{\leq n}$ sends the folk trivial fibrations of $\oo\Cat$ to weak
  equivalences of $n\Cat$ (in virtue of
  Ken Brown's Lemma \cite[Lemma 1.1.12]{hovey2007model}).
  
  For convenience, let us write $T$ for the functor $\iota_n\tau^{i}_{\leq n} : \oo\Cat \to \oo\Cat$.
  By definition of folk weak
  equivalences on $n\Cat$, we have to show that
  for every folk trivial fibration $f : C \to D$ of $\oo\Cat$,
  the $\oo$\nbd{}functor
  \[
  T(f) :   T(C) \to   T(D)
  \]
  is a folk weak equivalence on $\oo\Cat$. Consider the following commutative square
  \[
  \begin{tikzcd}
    C \ar[d,"\eta_C"] \ar[r,"f"] & D \ar[d,"\eta_D"] \\
    T(C) \ar[r,"T(f)"] & T(D),
  \end{tikzcd}
  \]
  where $\eta$ is the unit of the adjunction $\tau^{i}_{\leq n} \dashv \iota_n$.

  Let us first treat the case of $0$\nbd{}cells. Let $y$ be a $0$\nbd{}cell of $T(D)$. The map
    $\eta_D$ being surjective on $0$\nbd{}cells (even if $n=0$), there exists
    $y'$ such that $\eta_D(y')=y$. Since $f$ is a folk trivial fibration, there
    exists a $0$\nbd{}cell $x'$ of $C$ such that $f(x')=y'$ and then if we set $x:=\eta_C(x')$,
    we have $T(f)(x)=y$.

    Now let $x,y$ be parallel $k$\nbd{}cells of $T(C)$ and let  $\beta : f(x)
    \to f(y)$ be a $(k+1)$\nbd{}cell of $T(D)$. We have to distinguish
    several cases.    
  \begin{description}
  \item[Case $k<n-1$:] Since $\eta_C$ and $\eta_D$ are identities on $k$-cells for every $0\leq k<n$ and since $f$ is a folk trivial fibration, there exists a $(k+1)$\nbd{}cell $\alpha : x \to y$ of $T(C)$ such that
    \[
    T(f)(\alpha)=\beta.
  \]
  \item[Case $k=n-1$:] By definition of $T(D)$, there exists an $n$\nbd{}cell
    $\beta ' : f(x) \to f(y)$ of $D$ such that
    $\eta_{D}(\beta')=\beta$. Since $f$ is a folk trivial fibration, there
    exists an $n$\nbd{}cell $\alpha' : x \to y$ of $C$ such that
    $f(\alpha')=\beta'$. If we set $\alpha:=\eta_{C}(\alpha')$, we have
    $T(f)(\alpha)=\beta$.
    \item[Case $k=n$:] Since all $l$\nbd{}cells of $T(C)$ and $T(D)$ with $l>n$
      are units, we trivially have that $f(x)=f(y)$ and $\beta$ is the unit on
      $f(x)$. Now let $x'$ and $y'$ be parallel $n$\nbd{}cells of $C$ such that
      $\eta_C(x')=x$ and $\eta_C(y')=y$ (this is always possible by definition
      of $T(C)$). We have $\eta_{D}(f(x'))=f(x)=f(y)=\eta_{D}(f(y'))$. By
      definition of the functor $\tau^{i}_{\leq n}$, this means that there
      exists a zigzag of $(n+1)$\nbd{}cells of $D$ from $f(x')$ to $f(y')$.
      More precisely, this means that there exists a sequence
      \[
        (z_0,\beta_1,z_1,\cdots,z_{p-1},\beta_p,z_p)
      \]
      where the $z_i$ are all parallel $n$\nbd{}cells of $D$ with $z_0=f(x')$
      and $z_p=f(y')$, and each $\beta_i$ is $(n+1)$\nbd{}cell of $D$ either
      from $z_{i-1}$ to $z_i$ or from $z_{i}$ to $z_{i-1}$. Using the fact that
      $f$ is a folk trivial fibration, it is easy to prove the existence of a
      zigzag from $x'$to $y'$, which implies in particular that $x=\eta_C(x')=\eta_C(y')=y$.
    \item[Case $k>n$:] Since all $k$\nbd{}cells of $T(C)$ and $T(D)$ with $k>n$
    are  units, we trivially have $f(x)=f(y)$ (and $\beta$ is the unit on $f(x)$) and $x=y$.
  \end{description}
  Altogether, this proves that $T(f)$ is a folk trivial fibration, hence a folk
  weak equivalence.
\end{proof}
For later reference, we put here the following lemma.
\begin{lemma}\label{lemma:cofncatfolk}
  If an $n$\nbd{}category $C$ has a $k$\nbd{}basis for every $0 \leq k \leq n-1$, then it is cofibrant for the folk model structure on $n\Cat$.
\end{lemma}
\begin{proof}
Since $\tau^{i}_{\leq n}$ is a left Quillen functor, it suffices in virtue of
Proposition \ref{prop:freeiscofibrant} to show that
there exists a free $\oo$\nbd{}category $C'$ such that $\tau^{i}_{\leq
  n}(C')=C$. 
First, consider the $n$\nbd{}category $(U_{n-1}(C))^*$ (for the
notations, see \ref{paragr:freecext} and \ref{paragr:cextlowdimension}). This $n$\nbd{}category
has the same $k$\nbd{}cells as $C$ for $k<n$ and has exactly one
\emph{generating} $n$\nbd{}cell for each $n$\nbd{}cell of $C$. It is obviously free and we have a canonical $n$\nbd{}functor
\[
  \epsilon_C : (U_{n-1}(C))^* \to C, 
\]
given by the co-unit of the adjunction $(-)^* \dashv U_{n-1}(-)$. Now, let $C'$
be the $(n+1)$\nbd{}category (considered as an $\oo$\nbd{}category) that has the same $k$\nbd{}cells as
$(U_{n-1}(C))^*$ for $k\leq n$ and whose set of $(n+1)$\nbd{}cells is freely
generated be the set
  \[
    \{(x,y)\,\vert\, x \text{ and } y \text{ are parallel } n \text{-cells of }
    (U_{n-1}(C))^* \text{ such that }\epsilon_C(x)=\epsilon_C(y)\}.
  \]
  The $(n+1)$\nbd{}category $C'$ is obviously free and it is a harmless
  verification, which we leave to reader, to check that $\tau^i_{\leq n}(C')=C$.
\end{proof}
\begin{example}
  Every (small) category is cofibrant for the folk model structure on $\Cat$. 
\end{example}
We now turn to truncations of chain complexes.
\begin{paragr}
 Let $\Ch^{\leq n}$ be the category of chain complexes concentrated in degrees between $0$ and $n$. This means that an object $K$ of $\Ch^{\leq n}$ is a diagram of abelian groups of the form
 \[
 K_0 \overset{\partial}{\longleftarrow} K_1 \overset{\partial}{\longleftarrow} K_2 \overset{\partial}{\longleftarrow} \cdots \overset{\partial}{\longleftarrow} K_n,
 \]
 where $\partial \circ \partial =0$, and morphisms of $\Ch^{\leq n}$ are defined
 the expected way. We write $\iota_n : \Ch^{\leq n} \to \Ch$ for the canonical functor that sends an object $K$ of $\Ch^{\leq n}$ to the chain complex
 \[
 K_0 \overset{\partial}{\longleftarrow} K_1 \overset{\partial}{\longleftarrow} K_2 \overset{\partial}{\longleftarrow} \cdots \overset{\partial}{\longleftarrow} K_n \longleftarrow 0 \longleftarrow 0 \longleftarrow \cdots.
 \]
 This functor is fully faithful and $\Ch^{\leq n}$ may be identified with the full subcategory of $\Ch$ spanned by chain complexes $K$ such that $K_k = 0$ for every $k >n$.

  Similarly to the case of $n$\nbd{}categories, the functor $\iota_n : \Ch^{\leq n} \to \Ch$ has a left adjoint $\tau^{i}_{\leq n} : \Ch \to \Ch^{\leq}$, where for a chain complex $K$, $\tau^{i}_{\leq n}(K)$ is the object of $\Ch^{\leq n}$ defined as
  \[
  K_0 \overset{\partial}{\longleftarrow} K_1 \overset{\partial}{\longleftarrow}  \cdots \overset{\partial}{\longleftarrow} K_{n-1} \overset{\partial}{\longleftarrow} K_{n}/{\partial(K_{n+1})}.
  \]
  Again, as with $n$\nbd{}categories, we can use the adjunction
  \[
  \begin{tikzcd}
    \tau^{i}_{\leq n} : \Ch\ar[r,shift left] & \ar[l,shift left]\Ch^{\leq n} : \iota_n 
  \end{tikzcd}
  \]
  to create a model structure on $\Ch^{\leq n}$.
\end{paragr}
\begin{proposition}
  There exists a model structure on $\Ch^{\leq n}$ such that:
  \begin{itemize}[label=-]
  \item the weak equivalences are exactly those morphisms $f : K \to K'$ such that $\iota_n(f)$ is a weak equivalence for the projective model structure on $\Ch$,
    \item the fibrations are exactly those morphisms $f : K \to K'$ such that $\iota_n(f)$ is a fibration for the projective model structure on $\Ch$.
    \end{itemize}
\end{proposition}
\begin{proof}
  This is a typical example of a transfer of a cofibrantly generated model
  structure along a right adjoint as in \cite[Proposition
  2.3]{beke2001sheafifiableII}. Since the weak equivalences of the projective model
  structure on $\Ch$ are closed under filtered colimits \cite[Theorem 2.6.15]{weibel1995introduction}, the only \emph{a priori} non-trivial hypothesis to check is that there exists a set $J$ of generating trivial cofibrations of the projective model structure on $\Ch$ such that for every $j \colon A \to B$ in $J$ and every cocartesian square
  \[
  \begin{tikzcd}
    \tau^{i}_{\leq n}(A) \ar[r] \ar[d,"\tau^{i}_{\leq n}(j)"'] & X \ar[d,"g"] \\
    \tau^{i}_{\leq n}(B) \ar[r] & Y,
    \ar[from=1-1, to=2-2, phantom, "\ulcorner",very near end]
  \end{tikzcd}
  \]
  the morphism $\iota_n(g)$ is a weak equivalence of $\Ch$. As explained in \cite[Proposition 7.19]{dwyer1995homotopy}, there exists a set of generating trivial cofibrations of the projective model structure on $\Ch$ consisting of the maps
  \[
  0 \to D_k
  \]
  for each $k >0$, where $D_k$ is the following chain complex concentrated in degree $k$ and $k-1$
  \[
  0 \leftarrow \cdots \leftarrow 0 \leftarrow \mathbb{Z} \overset{\mathrm{id}}{\leftarrow} \mathbb{Z} \leftarrow 0 \leftarrow \cdots
  \]
  What is left to show then is that for every $k > 0$ and every object $X$ of $\Ch^{\leq n}$, the canonical inclusion map
  \[
  X \to X \oplus \tau^{i}_{\leq n}(D_k)
  \]is sent by $\iota_n$ to a weak equivalence of $\Ch$. This follows immediately from the fact that homology groups commute with direct sums.
\end{proof}
\begin{paragr}
  We refer to the model structure of the previous proposition as the \emph{projective model structure on $\Ch^{\leq n}$}.
\end{paragr}
\begin{lemma}\label{lemma:unitajdcomp}
  For every chain complex $K$, the unit map
  \[
  K \to \iota_n\tau^{i}_{\leq n}(K)
  \]
  induces isomorphisms
  \[
  H_k(K) \simeq H_k(\iota_n\tau^{i}_{\leq n}(K))
  \]
  for every $0 \leq k \leq n$.
\end{lemma}
\begin{proof}
  For $0 \leq k < n-1$, this is trivial. For $k = n-1$, this follows easily from the fact that the image of $\partial : K_k/{\partial(K_{k+1})}\to K_{k-1}$ is equal to the image of $\partial : K_k \to K_{k-1}$. Finally for $k = n$, it is straightforward to check that
  \[
  H_n(K)=\frac{\mathrm{Ker}(\partial : K_n \to K_{n-1})}{\mathrm{Im}(\partial : K_{n+1} \to K_n)}
  \]
  is isomorphic to
  \[
  H_n(\iota_n\tau^{i}_{\leq n}(K))=\mathrm{Ker}(\partial : K_n/{\partial(K_{n+1})} \to K_{n-1}).
  \]
  The isomorphism being obviously induced by the unit map $K \to \iota_n\tau^{i}_{\leq n}(K)$.
\end{proof}
As a consequence of this lemma, we have the analogous of Proposition \ref{prop:truncationhomotopical}.
\begin{proposition}
  The functor $\tau^{i}_{\leq n} : \Ch \to \Ch^{\leq n}$ sends the weak equivalences of the projective model structure on $\Ch$ to weak equivalences of the projective model structure on $\Ch^{\leq n}$.
\end{proposition}
\begin{proof}
  Let $f : K \to K'$ be a weak equivalence for the projective model structure on $\Ch$ and consider the naturality square
  \[
  \begin{tikzcd}[column sep=huge]
    K \ar[d,"\eta_K"] \ar[r,"f"] & K' \ar[d,"\eta_K'"] \\
    \iota_n\tau^{i}_{\leq n}(K) \ar[r,"\iota_n\tau^{i}_{\leq n}(f)"] & \iota_n\tau^{i}_{\leq n}(K'),
    \end{tikzcd}
  \]
  where $\eta$ is the unit map of the adjunction $\tau^{i}_{\leq n} \dashv \iota_n$. It follows from Lemma \ref{lemma:unitajdcomp} that
  \[
  H_k(\iota_n\tau^{i}_{\leq n}(f)) : H_k(\iota_n\tau^{i}_{\leq n}(K)) \to H_k(\iota_n\tau^{i}_{\leq n}(K'))
  \]
  is an isomorphism for every $k \leq n$. Since obviously $H_k(\iota_n\tau^{i}_{\leq n}(f))$ is also an isomorphism for $k > n$, this proves the result.
\end{proof}
We now investigate the relation between truncation and abelianization.
\begin{paragr}
  Let $C$ be $n$\nbd{}category. A straightforward computation shows that the chain complex $\lambda(\iota_n(C))$ is such that
 \[
 \lambda_k(\iota_n(C))=0
 \]
 for every $k > n$ and thus $\lambda(\iota_n(C))$ can be seen as an object of $\Ch^{\leq n}$. Hence, we can define a functor $\lambda_{\leq n } : n\Cat \to \Ch^{\leq n}$ as
 \begin{align*}
   \lambda_{\leq n} : n\Cat &\to \Ch^{\leq n}\\
   C&\mapsto \lambda(\iota_n(C)),
 \end{align*}
 and we tautologically have that the square
 \[
 \begin{tikzcd}
   n\Cat \ar[d,"\iota_n"] \ar[r,"\lambda_{\leq n}"] & \Ch^{\leq n} \ar[d,"\iota_n"] \\
   \oo\Cat \ar[r,"\lambda"] & \Ch
 \end{tikzcd}
 \]
is commutative.
  \end{paragr}
\begin{lemma}\label{lemma:abelianizationtruncation}
  The square
  \[
  \begin{tikzcd}
    \oo\Cat \ar[d,"\tau_{\leq n}^{i}"] \ar[r,"\lambda"] & \Ch \ar[d,"\tau^{i}_{\leq n}"] \\
    n\Cat \ar[r,"\lambda_{\leq n}"] & \Ch^{\leq n}
    \end{tikzcd}
  \]
  is commutative (up to a canonical isomorphism).
\end{lemma}
\begin{proof}
  Notice first that we have a natural transformation
  \[
  \beta : \tau^{i}_{\leq n}\circ \lambda \Rightarrow \lambda_{\leq n} \circ \tau^{i}_{\leq n}
  \]
  defined as
  \[
  \begin{tikzcd}
 \oo\Cat \ar[r,"\tau_{\leq n}^{i}"] \ar[rd,"\mathrm{id}"',""{name=A,right}] &  n\Cat \ar[d,"\iota_n"] \ar[r,"\lambda_{\leq n}"] & \Ch^{\leq n} \ar[d,"\iota_n"'] \ar[dr,"\mathrm{id}",""{name=B,left}] & \\
 &\oo\Cat \ar[r,"\lambda"] & \Ch \ar[r,"\tau^{i}_{\leq n}"'] & \Ch^{\leq n}.
 \ar[from=A, to=1-2,Rightarrow,"\eta"]
 \ar[from=2-3,to=B,Rightarrow,"\epsilon"]
   \end{tikzcd}
  \]
  Since for every $\oo$\nbd{}category $C$ and every $k <n$, we have
  \[
  C_k=\iota_n(\tau^{i}_{\leq n}(C))_k
  \]
  and for every chain complex $K$ and every $k<n$, we have
  \[
  \iota_n(\tau^{i}_{\leq n}(K))_k=K_k,
  \]
  it follows that for every $\oo$\nbd{}category $C$ and every $k<n$, the morphism $\beta_k$ is nothing but the equality
  \[
  \tau^{i}_{\leq n}(\lambda(C))_k=\lambda_{\leq n}(\tau^{i}_{\leq n}(C))_k.
  \]
  Hence, all we have to prove is that
  \[
  \beta_n :   \tau^{i}_{\leq n}(\lambda(C))_n \to \lambda_{\leq n}(\tau^{i}_{\leq n}(C))_n
  \]
  is an isomorphism for every $\oo$\nbd{}category $C$.

  Recall from Lemma \ref{lemma:adjlambdasusp} that $\lambda_n \circ \iota_n : n\Cat \to \Ab$ (which we abusively wrote as $\lambda_n$) is left adjoint to the functor $B^n : \Ab \to n\Cat$. In particular, for every $\oo$\nbd{}category $C$ and every abelian group $G$, we have 
  \[
  \Hom_{\Ab}(\lambda_n\iota_n\tau^{i}_{\leq n}(C),G) \simeq \Hom_{n\Cat}(\tau^{i}_{\leq n}(C),B^nG) \simeq \Hom_{\oo\Cat}(C,\iota_n(B^nG)).
  \]
  Then, it follows from an argument similar to the proof of Lemma \ref{lemma:nfunctortomonoid} (see also the proof of Lemma \ref{lemma:adjlambdasusp}) that this last set is naturally isomorphic to the set of functions $f_n : C_n \to G$ such that:
  \begin{itemize}[label=-]
  \item for every $0 \leq k <n $ and every pair $(x,y)$ of $k$-composable $n$\nbd{}cells of $C$, we have
    \[
    f(x \comp_k y) = f(x) + f(y),
    \]
  \item for every $(n+1)$\nbd{}cell $\alpha : x \to y$ of $C$, we have
    \[
    f(x)=f(y).
    \]
  \end{itemize}
  By definition of $\lambda_n(C)$ and of $\partial : \lambda_{n+1}(C) \to \lambda_n(C)$, we thus have
  \[
  \Hom_{\Ab}(\lambda_n\iota_n\tau^{i}_{\leq n}(C),G)\simeq \Hom_{\Ab}(\lambda_n(C)/{\partial(\lambda_{n+1}(C))},G).
  \]
  Hence, we have $\lambda_{\leq n}(\tau^{i}_{\leq n}(C))_n \simeq \tau^{i}_{\leq n}(\lambda(C))_n$ and a thorough analysis of naturality shows that this isomorphism is nothing but $\beta_n$.
\end{proof}
In the following lemma, $n\Cat$ is equipped with the folk model structure and $\Ch^{\geq n}$ with the projective model structure.
\begin{lemma}
  The functor $\lambda_{\leq n} : n\Cat \to \Ch^{\leq n}$ is left Quillen.
\end{lemma}
\begin{proof}
  Let $I$ and $J$ respectively be sets of generating cofibrations and generating
  trivial cofibrations of the folk model structure on $\oo\Cat$ such that
  $\tau^{i}_{\leq n}(I)$ and $\tau^{i}_{\leq n}(J)$ respectively are sets of
  generating cofibrations and generating trivial cofibrations of the projective
  model structure on $n\Cat$ (which we know exist by the second part of Proposition
  \ref{prop:fmsncat}). What we have to show is that for every $f$ in $I$ (resp.
  $J$), $\lambda_{\leq n}(\tau^{i}_{\leq n}(f))$ is a cofibration (resp.\
  generating cofibration) for the folk model structure on $\Ch^{\leq n}$.
  From Lemma \ref{lemma:abelianizationtruncation}, we have
  \[
  \lambda_{\leq n }(\tau^{i}_{\leq n}(f)) \simeq \tau^{i}_{\leq n}(\lambda(f)).
  \]
  Since $\lambda$ and $\tau^{i}_{\leq n}$ are both left Quillen functors, this proves the result.
\end{proof}
As an immediate consequence of the previous lemma, the functor $\lambda_{\leq n}$ is left derivable and we have the following key result.
\begin{proposition}\label{prop:polhmlgytruncation}
  The square
  \[
  \begin{tikzcd}
    \ho(\oo\Cat^{\folk}) \ar[d,"\overline{\tau_{\leq n}^{i}}"] \ar[r,"\LL \lambda"] & \ho(\Ch) \ar[d,"\overline{\tau^{i}_{\leq n}}"] \\
    \ho(n\Cat^{\folk}) \ar[r,"\LL \lambda_{\leq n}"] & \ho(\Ch^{\leq n})
    \end{tikzcd}
  \]
  is commutative (up to a canonical isomorphism).
\end{proposition}
\begin{proof}
  Straightforward consequence of Lemma \ref{lemma:abelianizationtruncation} and the fact that the left derived functor of a composition of left Quillen functors is the composition of the left derived functors (see for example \cite[Theorem 1.3.7]{hovey2007model}).
\end{proof}
\begin{remark}
  Beware that the square
   \[
 \begin{tikzcd}
   \ho(n\Cat^{\folk}) \ar[d,"\overline{\iota_n}"] \ar[r,"\LL \lambda_{\leq n}"] & \ho(\Ch^{\leq n}) \ar[d,"\overline{\iota_n}"] \\
   \ho(\oo\Cat^{\folk}) \ar[r,"\LL \lambda"] & \ho(\Ch)
 \end{tikzcd}
 \]
is \emph{not} commutative. If it were, then for every $n$\nbd{}category $C$ and every $k >n$, we would have $H_k^{\pol}(\iota_n(C))=0$ for every $k >n$, which is not even true for the case $n=1$ as we shall see in the following chapter.
\end{remark}
A useful consequence of Proposition \ref{prop:polhmlgytruncation} is the following corollary.
\begin{corollary}\label{cor:polhmlgycofibrant}
  Let $n \geq 0$ and $C$ be an $\oo$\nbd{}category. If $C$ has a $k$\nbd{}basis for every $ 0 \leq k  \leq n-1$, then the canonical map of $\ho(\Ch)$
  \[
  \alpha^{\pol}_C : \sH^{\pol}(C) \to \lambda(C)
  \]
  induces isomorphisms
  \[
  H_k^{\pol}(C) \simeq H_k(\lambda(C))
  \]
  for every $0 \leq k \leq n$.
\end{corollary}
\begin{proof}
  From Lemma \ref{lemma:abelianizationtruncation} and Proposition \ref{prop:polhmlgytruncation}, we deduce that the morphism $\overline{\tau^{i}_{\leq n}}(\alpha_C^{\pol})$ of $\ho(\Ch^{\leq n})$ can be identified with the canonical morphism
  \[
  \LL \lambda_{\leq n}(\tau^{i}_{\leq n}(C)) \to \lambda_{\leq n}(\tau^{i}_{\leq n}(C)).
  \]
  From Lemma \ref{lemma:cofncatfolk}, we have that $\tau^{i}_{\leq n}(C)$ is cofibrant for the folk model structure on $n\Cat$, and the result follows immediately from the fact that $\lambda_{\leq n}$ is left Quillen.
\end{proof}
\begin{paragr}\label{paragr:polhmlgylowdimension}
  Since every  $\oo$\nbd{}category trivially admits its set of $0$\nbd{}cells as a $0$\nbd{}base, it follows from the previous proposition that for every $\oo$\nbd{}category $C$ we have
  \[
  \sH^{\pol}_0(C)\simeq H_0(\lambda(C))
  \]
  and
  \[
  \sH^{\pol}_1(C) \simeq H_1(\lambda(C)).
  \]
  Intuitively speaking, this means that no cofibrant resolution of $C$ is needed to compute its first two polygraphic homology groups.
\end{paragr}
We now turn to the relation between truncation and singular homology of $\oo$\nbd{}categories. Recall that for every $n \geq 0$, the nerve functor $N_n : n\Cat \to \Psh{\Delta}$ is defined as the following composition
\[
N_n : n\Cat \overset{\iota_n}{\longrightarrow} \oo\Cat \overset{N_{\oo}}{\longrightarrow} \Psh{\Delta},
\]
and for $n \in \mathbb{N}\cup \{\oo\}$ we write $c_n : \Psh{\Delta} \to n\Cat$ for the left adjoint of $N_n$.
\begin{lemma}
  For every $n \in \mathbb{N}$, the following triangle of functors
  \[
  \begin{tikzcd}
    \Psh{\Delta} \ar[r,"c_{\oo}"] \ar[d,"c_n"] & \oo\Cat \ar[dl,"\tau^{i}_{\leq n}"] \\
    n\Cat&
  \end{tikzcd}
  \]
  is commutative (up to an isomorphism).
\end{lemma}
\begin{proof}
Straightforward consequence of the fact that $N_n = N_{\oo} \circ \iota_n$ and the fact that the composition of left adjoints is the left adjoint of the composition.  
\end{proof}
\begin{paragr}
  In particular, it follows from the previous lemma that the co-unit of the adjunction $c_{\oo} \dashv N_{\oo}$ induces for every $\oo$\nbd{}category $C$ and every $n \geq 0$, a canonical morphism of $n\Cat$
  \[
  c_nN_{\oo}(C) \simeq \tau^{i}_{\leq n}c_{\oo}N_{\oo}(C) \to \tau^{i}_{\leq n}(C),
  \]
  which is natural in $C$.
\end{paragr}
\begin{lemma}\label{lemma:truncationcounit}
  For every $\oo$\nbd{}category $C$, the canonical morphism of $\Cat$
  \[
  c_1N_{\oo}(C) \to \tau^{i}_{\leq 1}(C)
  \]
  is an isomorphism.
\end{lemma}
\begin{proof}
  Let $C$ be an $\oo$\nbd{}category and $D$ be a (small) category. By adjunction, we have
  \begin{equation}
  \Hom_{\Cat}(c_1N_{\oo}(C),D) \simeq \Hom_{\Psh{\Delta}}(N_{\oo}(C),N_1(D)).
  \end{equation}
  Now let $\Delta_{\leq 2}$ be the full subcategory of $\Delta$ spanned by
  $[0]$, $[1]$ and $[2]$ and let $i : \Delta_{\leq 2} \to \Delta$ be the
  canonical inclusion. This inclusion induces by pre-composition a functor $i^*
  : \Psh{\Delta} \to \Psh{\Delta_{\leq 2}}$ which has a right-adjoint $i_* :
    \Psh{\Delta_{\leq 2}} \to \Psh{\Delta}$. Recall that the nerve of a (small) category is $2$\nbd{}coskeletal (see for example \cite[Theorem 5.2]{street1987algebra}), which means that for every category $D$, the unit morphism $ N_1(D) \to i_* i^* (N_1(D))$ is an isomorphism of simplicial sets. In particular, we have
  \begin{align*}
    \Hom_{\Psh{\Delta}}(N_{\oo}(C),N_1(D)) &\simeq \Hom_{\Psh{\Delta}}(N_{\oo}(C),i_* i^* (N_1(D)))\\
    &\simeq \Hom_{\Psh{\Delta_{\leq 2}}}(i^*(N_{\oo}(C)),i^*(N_1(D))).
  \end{align*}
  Using the description of $\Or_0$, $\Or_1$ and $\Or_2$ from \ref{paragr:orientals}, we deduce that a morphism $F : i^*(N_{\oo}(C)) \to i^*(N_1(D))$ of $\Psh{\Delta_{\leq 2}}$ consists of a function $F_0 : C_0 \to D_0$ and a function $F_1 : C_1 \to D_1$ such that
  \begin{enumerate}[label=(\alph*)]
  \item for every $x \in C_0$, we have $F_1(1_x)=1_{F_0(x)}$,
  \item for every $x \in C_1$, we have
    \[\src(F_1(x))=F_0(\src(x))) \text{ and }\trgt(F_1(x))=F_0(\trgt(x))),\]
  \item for every $2$\nbd{}triangle
    \[
    \begin{tikzcd}
      & Y \ar[rd,"g"] & \\
      X \ar[ru,"f"] \ar[rr,"h"',""{name=A,above}] & & Z
      \ar[from=A,to=1-2,Rightarrow,"\alpha"]
    \end{tikzcd}
    \]
    in $C$, we have $F_1(g)\comp_0 F_1(f)=F_1(h)$. 
  \end{enumerate}
  In particular, it follows that $F_1$ is compatible with composition of
  $1$\nbd{}cells in an obvious sense and that for every $2$\nbd{}cell $\alpha :
  f \Rightarrow g$ of $C$, we have $F_1(f)=F_1(g)$. And conversely, this last
  condition implies condition (c) above. This means exactly that we have a natural isomorphism
  \[
 \Hom_{\Psh{\Delta_{\leq 2}}}(i^*(N_{\oo}(C)),i^*(N_1(D))) \simeq \Hom_{\Cat}(\tau_{\leq 1}^{i}(C),D).
 \]
 Altogether, we have
 \[
 \Hom_{\Cat}(c_1N_{\oo}(C),D) \simeq \Hom_{\Cat}(\tau_{\leq 1}^{i}(C),D),
 \]
 which proves that
 \[
 c_1N_{\oo}(C) \simeq \tau_{\leq 1}^{i}(C)
 \]
 and a thorough analysis of naturality shows that this isomorphism is nothing but the canonical morphism $c_1N_{\oo}(C) \to \tau_{\leq 1}^{i}(C)$.
\end{proof}
We can now prove the important following proposition.
\begin{proposition}\label{prop:singhmlgylowdimension}
  For every $\oo$\nbd{}category $C$, the canonical map of $\ho(\Ch)$
  \[
  \alpha_C^{\sing}: \sH^{\sing}(C) \to \lambda(C)
  \]
  induces isomorphisms
  \[
  H^{\sing}_k(C) \simeq H_k(\lambda(C))
  \]
  for $k \in \{0,1\}$.
\end{proposition}
\begin{proof}
  Let $C$ be an $\oo$\nbd{}category. Recall from \ref{paragr:univmor} that the canonical morphism $\alpha_C^{\sing} : \sH^{\sing}(C) \to \lambda(C)$ is nothing but the image by the localization functor $\Ch \to \ho(\Ch)$ of the morphism
  \[
  \lambda c_{\oo}N_{\oo}(C) \to \lambda(C)
  \]
  induced by the co-unit of the adjunction $c_{\oo} \dashv N_{\oo}$. From \ref{prop:polhmlgytruncation} we have that
  \[
  \tau^{i}_{\leq 1}\lambda c_{\oo}N_{\oo}(C) \simeq \lambda_{\leq 1} \tau_{\leq 1}^{i} c_{\oo} N_{\oo}(C)=\lambda_{\leq 1} c_1 N_{\oo}(C),
  \]
  and from Lemma \ref{lemma:truncationcounit} we obtain
  \[
  \tau_{\leq 1 }^{i} \lambda c_{\oo} N_{\oo}(C) \simeq \lambda_{\leq 1} \tau^{i}_{\leq 1}(C) \simeq \tau^{i}_{\leq 1}\lambda(C).
  \]
  This means exactly that the image by $\overline{\tau^{i}_{\leq 1}}$ of $\alpha_C^{\sing}$ is an isomorphism, which is what we wanted to prove.
\end{proof}
Finally, we obtain the result we were aiming for.
\begin{proposition}\label{prop:comphmlgylowdimension}
  For every $\oo$\nbd{}category $C$, the canonical comparison map
  \[
  \pi_C : \sH^{\sing}(C) \to \sH^{\pol}(C)
  \]
  induces isomorphisms
  \[
  H^{\sing}_k(C) \simeq H^{\pol}_k(C)
  \]
  for $k \in \{0,1\}$.
\end{proposition}
\begin{proof}
  Let $C$ be an $\oo$\nbd{}category and consider the following commutative triangle of $\ho(\Ch)$
  \[
  \begin{tikzcd}[column sep=tiny]
    \sH^{\sing}(C) \ar[rd,"\alpha_C^{\sing}"'] \ar[rr,"\pi_C"] & & \sH^{\pol}(C) \ar[dl,"\alpha_C^{\pol}"] \\
    &\lambda(C)&.
  \end{tikzcd}
  \]
  From Proposition \ref{prop:singhmlgylowdimension}, we know that $\alpha_C^{\sing}$ induces isomorphisms \[H_k^{\sing}(C) \simeq H_k(\lambda(C))\] for $k \in \{0,1\}$ and from Corollary \ref{cor:polhmlgycofibrant} and Paragraph \ref{paragr:polhmlgylowdimension} we know that $\alpha_C^{\pol}$ induces isomorphisms $H_k^{\pol}(C) \simeq H_k(\lambda(C))$ for $k \in \{0,1\}$. The result follows then from an immediate 2-out-of-3 property.
\end{proof}
\begin{paragr}\label{paragr:conjectureH2}
  A natural question following the above proposition is:
  \begin{center}
    For which $k \geq 0$ do we have $H_k^{\sing}(C) \simeq H_k^{\pol}(C)$ for every $\oo$\nbd{}category $C$ ?
  \end{center}
  We have already seen in \ref{paragr:bubble} that when $C = B^2\mathbb{N}$ we have
  \[
  H_{2p}^{\sing}(B^2\mathbb{N}) \not\simeq H^{\pol}_{2p}(B^2\mathbb{N})
  \]
  for every $p \geq 2$. Furthermore, with a similar argument to the one given in \ref{paragr:bubble}, we have that for every $k \geq 3$, the (nerve of the) $\oo$\nbd{}category $B^k\mathbb{N}$ is a $K(\mathbb{Z},k)$. In particular, we have
  \[
  H_{2p+3}^{\sing}(B^{2p +1}\mathbb{N})\simeq \mathbb{Z}/{2\mathbb{Z}}
  \]
  for every $p \geq 1$ (see \cite[Theorem 23.1]{eilenberg1954groups}). On the other hand, since $B^k\mathbb{N}$ is a free $k$\nbd{}category, we have $H_n^{\pol}(B^k\mathbb{N})=0$ for all $n \geq k$. All in all, we have proved that for every $k \geq 4$, there exists at least one $\oo$\nbd{}category $C$ such that
  \[
  H_k^{\sing}(C) \not\simeq H_k^{\pol}(C).
  \]
  However, it is still an open question to know whether for $k \in \{2,3\}$ we have
  \[
  H^{\sing}_k(C) \simeq H^{\pol}_k(C)
  \]
  for every $\oo$\nbd{}category $C$. The only missing part to adapt the proof of Proposition \ref{prop:comphmlgylowdimension} for these values of $k$ is the analogue of Lemma \ref{lemma:truncationcounit}. But contrary to the case $k=1$, it is not generally true that the canonical morphism $c_k N_{\oo}(C) \to \tau^{i}_{\leq k}(C)$ is an isomorphism when $k \geq 2$. However, what we really need is that the image by $\lambda$ of this morphism be a quasi-isomorphism. In the case $k=2$, it seems that this canonical morphism admits an oplax $2$\nbd{}functor as an inverse up to oplax transformation which could be an hint towards the conjecture that $H^{\sing}_2(C) \simeq H^{\pol}_2(C)$ for every $\oo$\nbd{}category $C$. 
\end{paragr}


\chapter{Homology of contractible \texorpdfstring{$\oo$}{ω}-categories and its
  consequences}
\chaptermark{Contractible $\omega$-categories and consequences}
\section{Contractible \texorpdfstring{$\oo$}{ω}-categories}
Recall that for every $\oo$\nbd{}category $C$, we write $p_C : C \to \sD_0$ for
the canonical morphism to the terminal object $\sD_0$ of $\oo\Cat$.
\begin{definition}\label{def:contractible}
  An $\oo$\nbd{}category $C$ is \emph{oplax contractible} when the canonical morphism $p_C : C \to \sD_0$ is an oplax homotopy equivalence (Definition \ref{def:oplaxhmtpyequiv}).

\end{definition}



\begin{proposition}\label{prop:contractibleisgood}
  Every oplax contractible $\oo$\nbd{}category $C$ is \good{} and we have
  \[
  \sH^{\pol}(C)\simeq \sH^{\sing}(C)\simeq \mathbb{Z}
  \]
  where $\mathbb{Z}$ is seen as an object of $\ho(\Ch)$ concentrated in degree $0$.
\end{proposition}
\begin{proof}
Consider the commutative square
  \[
  \begin{tikzcd}
    \sH^{\pol}(C) \ar[d,"\sH^{\pol}(p_C)"] \ar[r,"\pi_C"] & \sH^{\sing}(C) \ar[d,"\sH^{\sing}(p_C)"] \\
    \sH^{\pol}(\sD_0) \ar[r,"\pi_{\sD_0}"] & \sH^{\sing}(\sD_0).
  \end{tikzcd}
  \]
  It follows respectively from Proposition \ref{prop:oplaxhmtpyisthom} and
  Proposition \ref{prop:oplaxhmtpypolhmlgy} that the right and left vertical
  morphisms of the above square are isomorphisms. Then, an immediate computation
  left to the reader shows that $\sD_0$ is \good{} and that
  $\sH^{\pol}(\sD_0)\simeq \sH^{\sing}(\sD_0)\simeq \mathbb{Z}$. By a 2-out-of-3
  property, we deduce that $\pi_C : \sH^{\sing}(C)\to \sH^{\pol}(C)$ is an
  isomorphism and $\sH^{\pol}(C)\simeq \sH^{\sing}(C)\simeq \mathbb{Z}$.
\end{proof}
\begin{remark}
  Definition \ref{def:contractible} admits an obvious ``lax'' variation and Proposition \ref{prop:contractibleisgood} is also true for lax contractible $\oo$\nbd{}categories.
  \end{remark}
We end this section with an important result on slice $\oo$\nbd{}categories (Paragraph \ref{paragr:slices}).
\begin{proposition}\label{prop:slicecontractible}
  Let $A$ be an $\oo$\nbd{}category and $a_0$ an object of $A$. The $\oo$\nbd{}category $A/a_0$ is oplax contractible. 
   \end{proposition}
\begin{proof}
  This follows from the dual of \cite[Proposition 5.22]{ara2020theoreme}.
  \end{proof}
\section{Homology of globes and spheres}
\begin{lemma}\label{lemma:globescontractible}
  For every $n \in \mathbb{N}$, the $\oo$\nbd{}category $\sD_n$ is oplax contractible.
\end{lemma}
\begin{proof}
  Recall that we write $e_n$ for the unique non-trivial $n$\nbd{}cell of $\sD_n$ and that by definition $\sD_n$ has exactly two non-trivial $k$\nbd{}cells for every $k$ such that $0\leq k<n$. These two $k$\nbd{}cells are parallel and are given by $\src_k(e_n)$ and $\trgt_k(e_n)$.

  Let $r : \sD_0 \to \sD_n$ be the $\oo$\nbd{}functor that points to $\trgt_0(e_n)$ (which means that $r=\langle \trgt_0(e_n) \rangle$ with the notations of \ref{paragr:defglobe}). For every $k$\nbd{}cell $x$ of $\sD_n$, we have
  \[
  r(p(x))=\1^k_{\trgt_0(e_n)},
  \]
  where we write $p$ for the unique $\oo$\nbd{}functor $\sD_n \to \sD_0$.

  Now for $0 \leq k <n$, we define $\alpha_{\src_k(e_n)}$ and $\alpha_{\trgt_k(e_n)}$ as
  \[
  \alpha_{\src_k(e_n)}=\begin{cases}\trgt_{k+1}(e_n), \text{ if } k<n-1 \\ e_n, \text{ if } k=n-1\end{cases} \text{ and }  \alpha_{\trgt_k(e_n)}=\begin{cases}1_{\trgt_k(e_n)}, \text{ if } k<n-1 \\ e_n, \text{ if } k=n-1.\end{cases}
  \]
  It is straightforward to check that this data defines an oplax transformation \[\alpha : \mathrm{id}_{\sD_n} \Rightarrow r\circ p\] (see \ref{paragr:formulasoplax} and Example \ref{example:natisoplax}), which proves the result. 
\end{proof}
In particular, for every $n \in \mathbb{N}$, $\sD_n$ is \good{}. Recall from \ref{paragr:inclusionsphereglobe} that for every $n \geq 0$, we have a cocartesian square
\[
\begin{tikzcd}
     \sS_{n-1} \ar[r,"i_n"] \ar[d,"i_n"] & \sD_n \ar[d,"j_n^+"]\\
    \sD_n \ar[r,"j_n^-"'] & \sS_{n}.
    \ar[from=1-1,to=2-2,phantom,very near end,"\ulcorner"]
  \end{tikzcd}
\]
\begin{lemma}\label{lemma:squarenerve}
  For every $n \geq 0$, the commutative square of simplicial sets
  \[
  \begin{tikzcd}
    N_{\oo}(\sS_{n-1}) \ar[r,"N_{\oo}(i_n)"] \ar[d,"N_{\oo}(i_n)"] & N_{\oo}(\sD_{n}) \ar[d,"N_{\oo}(j_n^+)"] \\
    N_{\oo}(\sD_{n}) \ar[r,"N_{\oo}(j_n^-)"] & N_{\oo}(\sS_{n})
  \end{tikzcd}
  \]
  is cocartesian.
\end{lemma}
\begin{proof}
  Since colimits in presheaf categories are computed pointwise, what we need
  to show is that for every $k\geq 0$, the following commutative square is
  cocartesian
  \begin{equation}\label{squarenervesphere}
    \begin{tikzcd}[column sep=huge,row sep=huge]
    \Hom_{\oo\Cat}(\Or_k,\sS_{n-1}) \ar[r,"{\Hom_{\oo\Cat}(\Or_k,i_n)}"] \ar[d,"{\Hom_{\oo\Cat}(\Or_k,i_n)}"'] & \Hom_{\oo\Cat}(\Or_k,\sD_{n}) \ar[d,"{\Hom_{\oo\Cat}(\Or_k,j_n^+)}"] \\
    \Hom_{\oo\Cat}(\Or_k,\sD_{n}) \ar[r,"{\Hom_{\oo\Cat}(\Or_k,j_n^-)}"'] & \Hom_{\oo\Cat}(\Or_k,\sS_{n}).
    \end{tikzcd}
  \end{equation}
  Notice first that the square
  \[
    \begin{tikzcd}
      \sS_{n-1} \ar[r,"i_n"] \ar[d,"i_n"] & \sD_n \ar[d,"j_n^+"]\\
      \sD_n \ar[r,"j_n^-"'] & \sS_{n}.
    \end{tikzcd}
  \]
  is cartesian and all four morphisms are monomorphisms. Since the
  functor \[\Hom_{\oo\Cat}(\Or_k,-):\oo\Cat \to \Set \] preserves limits, the square
  \eqref{squarenervesphere} is a cartesian square of $\Set$ all of
  whose four morphisms are monomorphisms.  Hence, in order to prove
  that square \eqref{squarenervesphere} is cocartesian, we only need
  to show that for every $k \geq 0$ and every $\oo$\nbd{}functor
  $\varphi : \Or_k \to \sS_{n}$, there exists an $\oo$\nbd{}functor
  $\varphi' : \Or_k \to \sD_n$ such that either $j_n^+ \circ \varphi '
  = \varphi$ or $j_n^- \circ \varphi' = \varphi$.
  
  For convenience, let us write $h_n^+$ (resp.\ $h_n^-$) for the only generating $n$\nbd{}cell of $\sS_n$ contained in the image of $j^+_n$ (resp.\ $j_n^-$). The cells $h_n^+$ and $h_n^-$ are the only non-trivial $n$\nbd{}cells of $\sS_n$. We also write $\alpha_k$ for the principal cell of $\Or_k$ (see \ref{paragr:orientals}). This is the only non-trivial $k$\nbd{}cell of $\Or_k$.

  Now, let $\varphi : \Or_k \to \sS_n$ be an $\oo$\nbd{}functor. There are several cases to distinguish.
  \begin{description}
  \item[Case $k<n$:] Since every generating cell of $\gamma$ of $\Or_k$ is of
    dimension not greater than $k$, the cell $\varphi(\gamma)$ is of dimension strictly lower than $n$. Since all cells of dimension strictly lower than $n$ are both in the image of $j^+_n$ and in the image of $j^-_n$, $\varphi$ obviously factors through $j^+_n$ (and $j^-_n$).
  \item[Case $k=n$:] The image of $\alpha_n$ is either a non-trivial $n$\nbd{}cell of $\sS_n$ or a unit on a strictly lower dimensional cell. In the second situation, everything works like the case $k<n$. Now suppose for example that $\varphi(\alpha_n)$ is $h^+_n$, which is in the image of $j^+_n$. Since all of the other generating cells of $\Or_n$ are of dimension strictly lower than $n$, their images by $\varphi$ are also of dimension strictly lower than $n$ and hence, are all contained in the image of $j^+_n$. Altogether this proves that $\varphi$ factors through $j^+_n$. The case where $\varphi(\alpha_n)=h^-_n$ is symmetric.
  \item[Case $k>n$:] Since $\sS_n$ is an $n$\nbd{}category, the image of
    $\alpha_k$ is necessarily of the form $\varphi(\alpha_k)=\1^k_{x}$ with $x$
    a cell of $\sS_n$ of dimension non-greater than $n$. If $x$ is a unit on a
    cell whose dimension is strictly lower than $n$, then everything works like
    in the case $k<n$. If not, this means that $x$ is a non-trivial
    $n$\nbd{}cell of $\sS_n$. Suppose for example that $x=h^+_n$. Now let
    $\gamma$ be a generator of $\Or_k$ of dimension $k-1$. We have
    $\varphi(\gamma)=\1^{k-1}_y$ with $y$ which is either a unit on a cell of
    dimension strictly lower than $n$, or a non-degenerate $n$\nbd{}cell of
    $\sS_n$ (if $k-1=n$, we use the convention that $\1^{k-1}_y=y$). In the
    first situation, $y$ is in the image of $j^+_n$ as in the case $k<n$, and
    thus, so is $\1^{k-1}_y$. In the second situation, this means \emph{a
      priori} that either $y=h_n^+$ or $y=h_n^-$. But we know that $\gamma$ is
    part of a composition that is equal to either the source or the target of
    $\alpha_k$ (see \ref{paragr:orientals}) and thus, $y$ is part of a composition that is equal to either the source or the target of $x=h^+_n$. Since no composition involving $h^-_n$ can be equal to $h^+_n$ (one could invoke the function introduced in \ref{prop:countingfunction}), this implies that $y=h_n^+$ and hence, $f(\gamma)$ is in the image of $j^+_n$. This goes for all generating cells of dimension $k-1$ of $\Or_k$ and we can recursively apply the same reasoning for generating cells of dimension $k-2$, then $k-3$ and so forth. Altogether, this proves that $\varphi$ factorizes through $j^+_n$. The case where $x=h^-_n$ and $\varphi$ factorizes through $j^-_n$ is symmetric.
    \end{description}
\end{proof}
From these two lemmas, follows the important proposition below.
\begin{proposition}\label{prop:spheresaregood}
  For every $n \geq -1$, the $\oo$\nbd{}category $\sS_n$ is \good{}.
\end{proposition}
\begin{proof}
  Recall that the cofibrations of simplicial sets are exactly the monomorphisms, and in particular all simplicial sets are cofibrant. Since $i_n : \sS_{n-1} \to \sD_n$ is a monomorphism for every $n \geq 0$ and since $N_{\oo}$ preserves monomorphisms (as a right adjoint), it follows from Lemma \ref{lemma:squarenerve} and Lemma \ref{lemma:hmtpycocartesianreedy} that the square
\[
  \begin{tikzcd}
    N_{\oo}(\sS_{n-1}) \ar[r,"N_{\oo}(i_n)"] \ar[d,"N_{\oo}(i_n)"] & N_{\oo}(\sD_{n}) \ar[d,"N_{\oo}(j_n^+)"] \\
    N_{\oo}(\sD_{n}) \ar[r,"N_{\oo}(j_n^-)"] & N_{\oo}(\sS_{n})
  \end{tikzcd}
  \]
is a homotopy cocartesian square of simplicial sets. Since $N_{\oo}$ induces an
  equivalence of op-prederivators $\Ho(\oo\Cat^{\Th}) \to \Ho(\Psh{\Delta})$
  (Theorem \ref{thm:gagna}), it follows that the square of $\oo\Cat$
  \[
    \begin{tikzcd}
    \sS_{n-1} \ar[r,"i_n"] \ar[d,"i_n"] & \sD_{n} \ar[d,"j_n^+"] \\
    \sD_{n} \ar[r,"j_n^-"] & \sS_{n}
    \ar[from=1-1,to=2-2,phantom,"\ulcorner",very near end]
    \end{tikzcd}
    \]
  is Thomason homotopy cocartesian for every $n\geq 0$. Finally, since $i_n : \sS_{n-1} \to \sD_{n}$ is
  a folk cofibration and $\sS_{n-1}$ and $\sD_{n}$ are folk cofibrant for every $n\geq0$, we deduce the desired result from Corollary \ref{cor:usefulcriterion} and an immediate induction. The base case being simply that $\sS_{-1}=\emptyset$ is obviously \good{}.
\end{proof}
\begin{paragr}
  The previous proposition implies what we claimed in Paragraph \ref{paragr:prelimcriteriongoodcat}, which is that the morphism of op-prederivators
  \[
  \J : \Ho(\oo\Cat^{\folk}) \to \Ho(\oo\Cat^{\Th})
  \]
  induced by the identity functor of $\oo\Cat$ is \emph{not} homotopy cocontinuous.
  Indeed, recall from Paragraph \ref{paragr:bubble} that the category
  $B^2\mathbb{N}$ is \emph{not} \good{}; but on the other hand we have a
  cocartesian square
  \[
  \begin{tikzcd}
    \sS_1 \ar[r] \ar[d,"i_2"] & \sD_0 \ar[d] \\
    \sD_2 \ar[r] & B^2\mathbb{N},
    \ar[from=1-1,to=2-2,phantom,"\ulcorner",very near end]
  \end{tikzcd}
  \]
  where the map $\sD_2 \to B^2\mathbb{N}$ points the unique generating
  $2$-cell of $B^2\mathbb{N}$ and $\sD_0 \to B^2\mathbb{N}$ points to
  the only object of $B^2\mathbb{N}$. Since $\sS_1$, $\sD_0$ and $\sD_2$ are
  free and $i_2$ is a folk cofibration,
  the square is also folk homotopy cocartesian. If $\J$ was homotopy
  cocontinuous, then this square would also be Thomason homotopy
  cocartesian. Since we know that $\sS_1$, $\sD_0$ and $\sD_2$ are
  \good{}, this would imply that $B^2\mathbb{N}$ is \good{}.

  From Proposition \ref{prop:spheresaregood}, we also deduce the proposition
  below which gives a criterion to detect \good{} $\oo$\nbd{}categories when we
  already know that they are free. 
\end{paragr}
\begin{proposition}
  Let $C$ be a free $\oo$\nbd{}category and for every $k \in \mathbb{N}$ let
  $\Sigma_k$ be its $k$\nbd{}basis. If for every $k \in \mathbb{N}$, the following
  cocartesian square (see \ref{def:nbasis})
  \[
  \begin{tikzcd}[column sep=large]
   \displaystyle \coprod_{x \in \Sigma_k} \sS_{k-1} \ar[d,"\coprod i_k"'] \ar[r,"{\langle s(x),t(x) \rangle}"]& \sk_{k-1}(C) \ar[d,hook] \\
    \displaystyle \coprod_{x \in \Sigma_k} \sD_{k} \ar[r,"\langle x \rangle"] &
    \sk_k(C)
    \ar[from=1-1,to=2-2,phantom,"\ulcorner",very near end]
  \end{tikzcd}
  \]
  is Thomason homotopy cocartesian, then $C$ is \good{}.
\end{proposition}
\begin{proof}
  Since the morphisms $i_k$ are folk cofibrations and the
  $\oo$\nbd{}categories $\sS_{k-1}$ and $\sD_{k}$ are folk cofibrant
  and \good{}, it follows from Corollary \ref{cor:usefulcriterion} and
  an immediate induction that all $\sk_k(C)$ are \good{}. The result
  follows then from Lemma \ref{lemma:filtration}, Corollary \ref{cor:cofprojms} and Proposition
  \ref{prop:sequentialhmtpycolimit}.
\end{proof}
\section{The miraculous case of 1-categories}
Recall that the terms \emph{1-category} and \emph{(small) category} are
synonymous. While we have used the latter one more often so far, in this section
we will mostly use the former one. As usual, the canonical functor $\iota_1 :
\Cat \to \oo\Cat$ is treated as an inclusion functor and hence we always consider
$1$\nbd{}categories as particular cases of $\oo$\nbd{}categories.

The goal of what follows is to show that every $1$\nbd{}category is \good{}. In order to do that, we will prove that every 1-category
is a canonical colimit of contractible $1$\nbd{}categories and that this colimit is
homotopic both
with respect to the folk weak equivalences and with respect to the Thomason equivalences.
We call the reader's attention to an important subtlety here: even though the
desired result only refers to $1$\nbd{}categories, we have to work in the setting
of $\oo$\nbd{}categories. This can be explained from the fact that if we take a
cofibrant resolution of a $1$\nbd{}category $C$ in the folk model structure on $\oo\Cat$
\[
P \to C,
\]
then $P$ is not necessarily a $1$\nbd{}category. In particular, polygraphic
homology groups of a $1$\nbd{}category need \emph{not} be trivial in dimension
higher than $1$.
\begin{paragr}
 Let $A$ be a $1$\nbd{}category and $a$ an object of $A$. Recall that we write $A/a$ for the slice $1$\nbd{}category of $A$ over $a$, that is the $1$\nbd{}category whose description is as follows:
  \begin{itemize}[label=-]
  \item an object of $A/a$ is a pair $(a', p : a' \to a)$ where $a'$ is an object of $A$ and $p$ is an arrow of $A$,
  \item an arrow of $A/a$ is a pair $(q,p : a' \to a)$ where $p$ is an arrow of
    $A$ and $q$ is an arrow of $A$ of the form $q : a'' \to a'$. The target of
    $(q,p)$ is given by $(a',p)$ and the source by $(a'',p\circ q)$. 
  \end{itemize}
  We write $\pi_a$ for the canonical forgetful functor
   \[
  \begin{aligned}
    \pi_{a} : A/a &\to A \\
    (a',p) &\mapsto a'.
  \end{aligned}
  \]
  This is a special case of the more general notion of slice $\oo$\nbd{}category introduced in \ref{paragr:slices}. In particular, given an $\oo$\nbd{}category $X$ and an $\oo$\nbd{}functor $f : X \to A$, we have defined the $\oo$\nbd{}category $X/a$ and the $\oo$\nbd{}functor
  \[
  f/a : X/a \to A/a
  \]
  as the following pullback
  \[
  \begin{tikzcd}
    X/a \ar[r]\ar[dr, phantom, "\lrcorner", very near start] \ar[d,"f/a"'] & X \ar[d,"f"] \\
    A/a \ar[r,"\pi_{a}"] &A.
  \end{tikzcd}
  \]
  More explicitly, the $n$\nbd{}cells of $X/a$ can be described as pairs $(x,p)$
  where $x$ is an $n$\nbd{}cell of $X$ and $p$ is an arrow of $A$ of the form
  \[
    p : f(x)\to a \text{ if }n=0
  \]
  and
  \[
    p : f(\trgt_0(x)) \to a \text{ if }n>0.
  \]
  \emph{From now on, let us use the convention that $\trgt_0(x)=x$ when $x$ is a
    $0$\nbd{}cell of $X$}.

  When $n>0$, the source and target of an $n$\nbd{}cell $(x,p)$ of $X/a$ are given by
  \[
  \src((x,p))=(\src(x),p) \text{ and } \trgt((x,p))=(\trgt(x),p).
\]
  Moreover, the $\oo$\nbd{}functor $f/a$ is described as 
  \[
  (x,p) \mapsto (f(x),p),
  \]
  and the canonical $\oo$\nbd{}functor $X/a \to X$ as
  \[
  (x,p) \mapsto x.
  \]
\end{paragr}
\begin{paragr}\label{paragr:unfolding}
  Let $f : X \to A$ be an $\oo$\nbd{}functor with $A$ a $1$\nbd{}category. Every arrow $\beta : a \to a'$ of $A$ induces an $\oo$\nbd{}functor
  \begin{align*}
    X/\beta : X/a &\to X/{a'} \\
    (x,p) & \mapsto (x,\beta \circ p),
  \end{align*}
  which takes part in a commutative triangle
  \[
  \begin{tikzcd}[column sep=tiny]
    X/{a} \ar[rr,"X/{\beta}"] \ar[dr] && X/{a'} \ar[dl] \\
    &X&.
  \end{tikzcd}
  \]
  This defines a functor
  \begin{align*}
    X/{-} : A &\to \oo\Cat\\
    a &\mapsto X/a
  \end{align*}
  and a canonical $\oo$\nbd{}functor
  \[
  \colim_{a \in A} (X/{a}) \to X.
  \]
  Let $f' : X' \to A$ be another $\oo$\nbd{}functor and let
  \[
  \begin{tikzcd}
    X \ar[rr,"g"] \ar[dr,"f"'] && X' \ar[dl,"f'"] \\
    &A&
  \end{tikzcd}
  \]
  be a commutative triangle in $\oo\Cat$. Recall from \ref{paragr:slices} that for every object $a$ of $A$, there is an $\oo$\nbd{}functor $g/a$ defined as
  \begin{align*}
    g/a : X/a &\to X'/a \\
    (x,p) &\mapsto (g(x),p).
  \end{align*}
  This defines a natural transformation
  \[g/- : X/- \Rightarrow X'/-,\]
  and thus induces an $\oo$\nbd{}functor
  \[
  \colim_{a \in A}(X/a) \to \colim_{a \in A}(X'/a).
  \]
  Furthermore, it is immediate to check that the square
  \[
  \begin{tikzcd}
  \displaystyle \colim_{a \in A}(X/a) \ar[d] \ar[r] & X \ar[d,"g"] \\
  \displaystyle\colim_{a \in A}(X'/a) \ar[r] & X',
  \end{tikzcd}
  \]
  is commutative.
\end{paragr}
\begin{lemma}\label{lemma:colimslice}
  Let $f : X \to A$ be an $\oo$\nbd{}functor such that $A$ is a $1$\nbd{}category. The canonical $\oo$\nbd{}functor
  \[
  \colim_{a \in A}(X/a) \to X
  \]
  is an isomorphism.
\end{lemma}
\begin{proof}
 We have to show that the cocone
 \[
   (X/a \to X)_{a \in \Ob(A)}
 \]
 is colimiting.
 Let
 \[
   (\phi_{a} : X/a \to C)_{a \in \Ob(A)}
 \]
 be another cocone and let $x$ be a $n$\nbd{}arrow of $X$. Notice that the pair
 \[
   (x,1_{f(\trgt_0(x))})
 \]
 is a $n$\nbd{}arrow of $X/f(\trgt_0(x))$. We leave it to the reader to check that the formula
\begin{align*}
     \phi : X &\to C \\
     x &\mapsto \phi_{f(\trgt_0(x))}(x,1_{f(\trgt_0(x))}).
   \end{align*}
 defines an $\oo$\nbd{}functor and it is straightforward to check that for every object $a$ of $A$ the triangle
 \[
   \begin{tikzcd}
     X/a\ar[dr,"\phi_{a}"']\ar[r] & X \ar[d,"\phi"] \\
     & C
   \end{tikzcd}
 \]
 is commutative. This proves the existence part of the universal property.
 
  Now let $\phi' : X \to C$ be another $\oo$\nbd{}functor that makes the previous triangle commute for every object $a$ of $A$ and let $x$ be an $n$\nbd{}cell of $X$. Since the triangle
 \[
   \begin{tikzcd}
     X/f(\trgt_0(x)) \ar[dr,"\phi_{f(\trgt_0(x))}"']\ar[r] & X \ar[d,"\phi'"] \\
     & C
   \end{tikzcd}
 \]
 is commutative, we necessarily have
 \[
   \phi'(x)=\phi_{f(\trgt_0(x))}(x,1_{f(\trgt_0(x))})
 \]
 which proves that $\phi'=\phi$.
\end{proof}
\begin{paragr}
  In particular, when we apply the previous lemma to $\mathrm{id}_A : A \to A$,
  we obtain that every $1$\nbd{}category $A$ is (canonically isomorphic to) the colimit
  \[
  \colim_{a \in A} (A/a).
\]
  We now proceed to prove that this colimit is homotopic with respect to
  the folk weak equivalences.
\end{paragr}
Up to Lemma \ref{lemma:basisofslice}, we fix once and for all an $\oo$\nbd{}functor $f : X \to A$ with $A$ a $1$\nbd{}category.
\begin{lemma}\label{lemma:sliceisfree}
  If $X$ is free, then for every object $a$ of $A$, the $\oo$\nbd{}category
  $X/a$ is free. More precisely, if $\Sigma^X_n$ is the $n$\nbd{}basis of $X$,
  then the $n$\nbd{}basis of $X/a$ is the set
  \[
  \Sigma^{X/a}_n := \{(x,p) \in (X/a)_n \vert x \in \Sigma^X_n \}.
  \]
\end{lemma}
\begin{proof}
  It is immediate to check that for every object $a$ of $A$, the canonical
  forgetful functor $\pi_{a} : A/a \to A$ is a discrete Conduché functor (see Section
  \ref{section:conduche}). Hence, from Lemma \ref{lemma:pullbackconduche} we
  know that $X/a \to X$ is a discrete Conduché $\oo$\nbd{}functor. The result follows
  then from Theorem \ref{thm:conduche}.
\end{proof}
\begin{paragr}
  When $X$ is free, every arrow $\beta : a \to a'$ of $A$ induces a map
  \begin{align*}
    \Sigma^{X/a}_n &\to \Sigma^{X/a'}_n \\
    (x,p) &\mapsto (x,\beta\circ p).
  \end{align*}
  This defines a functor
  \begin{align*}
    \Sigma^{X/{\shortminus}}_n : A &\to \Set \\
    a &\mapsto \Sigma^{X/a}_n.
  \end{align*}
\end{paragr}
\begin{lemma}\label{lemma:basisofslice}
  If $X$ is free, then there is an isomorphism of functors 
  \[
  \Sigma^{X/\shortminus}_n \simeq \coprod_{x \in \Sigma^X_n}\Hom_A\left(f(\trgt_0(x)),\shortminus\right).
  \]
\end{lemma}
\begin{proof}
  For every object $a$ of $A$ and every $x \in \Sigma_n^X$, we have a canonical map
  \begin{align*}
    \Hom_A\left(f(\trgt_0(x)),a\right) &\to \Sigma^{X/a}_n \\
    p &\mapsto (x,p).
  \end{align*}
  By universal property, this induces a map
  \[
  \coprod_{x \in \Sigma^X_n}\Hom_A\left(f(\trgt_0(x)),a\right) \to  \Sigma^{X/a}_n,
  \]
  which is natural in $a$. A simple verification shows that it is a bijection.
\end{proof}
\begin{proposition}\label{prop:sliceiscofibrant}
  Let $A$ be a $1$\nbd{}category, $X$ be a free $\oo$\nbd{}category and $f : X \to A$ be an $\oo$\nbd{}functor. The functor
  \begin{align*}
    A &\to \oo\Cat \\
    a &\mapsto X/a
  \end{align*}
  is a cofibrant object for the projective model structure on $\oo\Cat(A)$
  induced by the folk model structure on $\oo\Cat$ (\ref{paragr:projmod}).
\end{proposition}
\begin{proof}
  Recall that the set
  \[
    \{i_n: \sS_{n-1} \to \sD_n \vert n \in \mathbb{N} \}
  \]
  is a set of generating folk cofibrations.
  From Lemmas \ref{lemma:sliceisfree} and \ref{lemma:basisofslice} we deduce
  that for every object $a$ of $A$ and every $n \in \mathbb{N}$, the canonical square
  \[
    \begin{tikzcd}
      \displaystyle\coprod_{x \in \Sigma^X_n}\coprod_{\Hom_A(f(\trgt_0(x)),a)}\sS_{n-1} \ar[r] \ar[d] & \sk_{n-1}(X/a) \ar[d]\\
      \displaystyle\coprod_{x \in \Sigma^X_n}\coprod_{\Hom_A(f(\trgt_0(x)),a)}\sD_n \ar[r]& \sk_n{(X/a)}
    \end{tikzcd}
  \]
  is cocartesian. It is straightforward to check that this square is natural in
  $a$ in an obvious sense, which means that we have a cocartesian square in $\oo\Cat(A)$:
  \[
      \begin{tikzcd}
      \displaystyle\coprod_{x \in \Sigma^X_n}\sS_{n-1}\otimes f(\trgt_0(x)) \ar[r] \ar[d] & \sk_{n-1}(X/-) \ar[d]\\
      \displaystyle\coprod_{x \in \Sigma^X_n}\sD_n\otimes f(\trgt_0(x)) \ar[r]& \sk_n{(X/-)}
      \ar[from=1-1,to=2-2,phantom,"\ulcorner",very near end]
    \end{tikzcd}
  \]
  (see \ref{paragr:cofprojms} for notations). From the second part of Proposition \ref{prop:modprs}, we deduce that for every $n\geq 0$, \[\sk_{n-1}(X/-) \to \sk_{n}(X/-)\] is a cofibration for the
  projective model structure on $\oo\Cat(A)$. Thus, the transfinite composition
  \[
  \emptyset \to \sk_{0}(X/-) \to \sk_{1}(X/) \to \cdots \to \sk_{n}(X/-) \to \cdots,
  \]
  which is canonically isomorphic to $\emptyset \to X/-$ (see Lemma \ref{lemma:filtration}), is also a cofibration for the projective model structure. 
\end{proof}
\begin{corollary}\label{cor:folkhmtpycol}
  Let $A$ be a $1$\nbd{}category and $f : X \to A$ be an $\oo$\nbd{}functor. The canonical arrow of $\ho(\oo\Cat^{\folk})$
  \[
  \hocolim^{\folk}_{a \in A}(X/a) \to X,
  \]
  induced by the co-cone $(X/a \to X)_{a \in \Ob(A)}$, is an isomorphism.
\end{corollary}
Beware that in the previous corollary, we did \emph{not} suppose that $X$ was free.
\begin{proof}
  Let $P$ be a free $\omega$-category and $g : P \to X$ a folk trivial fibration
  and consider the following commutative diagram of $\ho(\oo\Cat^{\folk})$ 
  \begin{equation}\label{comsquare}
    \begin{tikzcd}
      \displaystyle\hocolim^{\folk}_{a \in A}(P/a) \ar[d] \ar[r] & \displaystyle\colim_{a \in A}(P/a) \ar[d] \ar[r] & P \ar[d]\\
      \displaystyle\hocolim^{\folk}_{a \in A}(X/a)  \ar[r] & \displaystyle\colim_{a \in A}(X/a) \ar[r] & X
    \end{tikzcd}
  \end{equation}
  where the middle and most left vertical arrows are induced by the arrows
  \[
    g/a : P/a \to X/a,
  \]
  and the most right vertical arrow is induced by $g$.
  Since trivial fibrations are stable by pullback, $g/a$ is a trivial fibration.
  This proves that the most left vertical arrow of diagram \eqref{comsquare} is an isomorphism.
  
  Now, from Proposition \ref{prop:sliceiscofibrant} and Corollary
  \ref{cor:cofprojms}, we deduce that the arrow \[\hocolim_{a \in
    A}^{\folk}(P/a)\to \colim_{a \in A}(P/a)\] is an isomorphism. Moreover, from Lemma \ref{lemma:colimslice}, we know that the arrows
  \[\colim_{a \in A}(P/a)\to P\] and \[\colim_{a \in A}(X/a)\to X\] are
  isomorphisms. 

  Finally, since $g$ is a folk weak equivalence, the most right vertical arrow of diagram \eqref{comsquare} is an
  isomorphism and by an immediate 2-out-of-3 property this proves
  that all arrows of \eqref{comsquare} are isomorphisms. In particular, so is
  the composition of the two bottom horizontal arrows, which is what we desired
  to show.
\end{proof}
We now move on to the next step needed to prove that every $1$\nbd{}category is \good{}. For that purpose, let us recall a construction commonly referred to as the ``Grothendieck construction''.
\begin{paragr}
  Let $A$ be a $1$\nbd{}category and $F : A \to \Cat$ a functor. We denote by $\int F$ or $\int_{a \in A}F(a)$ the category such that:
  \begin{itemize}[label=-]
  \item An object of $\int F$ is a pair $(a,x)$ where $a$ is an object of $A$ and $x$ is an object of $F(a)$.
  \item An arrow $(a,x) \to (a',x')$ of $\int F$ is a pair $(f,k)$ where
    \[
    f : a \to a'
    \]
    is an arrow of $A$, and
    \[
    k : F(f)(x) \to x'.
    \]
  \end{itemize}
  The unit on $(a,x)$ is the pair $(1_a,1_x)$ and the composition of $(f,k) : (a,x) \to (a',x')$ and $(f',k') : (a',x') \to (a'',x'')$ is given by:
  \[
  (f',k')\circ(f,k)=(f'\circ f,k'\circ F(f')(k)).
  \]
  Every natural transformation
  \[
  \begin{tikzcd}
    A \ar[r,bend left,"F",""{name=A,below},pos=19/30]\ar[r,bend right,"G"',""{name=B, above},pos=11/20] & \Cat \ar[from=A,to=B,Rightarrow,"\alpha",pos=9/20]
  \end{tikzcd}
  \]
  induces a functor
  \begin{align*}
    \int \alpha : \int F &\to \int G\\
    (a,x) &\mapsto (a,\alpha_a(x)).
  \end{align*}
  Altogether, this defines a functor
  \begin{align*}
    \int : \Cat(A)&\to \Cat \\
   F&\mapsto \int F,
  \end{align*}
  where $\Cat(A)$ is the category of functors from $A$ to $\Cat$. 
\end{paragr}
We now recall an important theorem due to Thomason.
\begin{theorem}[Thomason]\label{thm:Thomason}
  The functor $\int : \Cat(A) \to \Cat$ sends the pointwise Thomason equivalences (\ref{paragr:homder}) to Thomason equivalences and the induced functor
  \[
  \overline{\int} : \ho(\Cat^{\Th}(A)) \to \ho(\Cat^{\Th})
  \]
  is canonically isomorphic to the homotopy colimit functor
  \[
  \hocolim^{\Th}_A : \ho(\Cat^{\Th}(A)) \to \ho(\Cat^{\Th}).
  \]
\end{theorem}
\begin{proof}
  The original source for this Theorem is \cite{thomason1979homotopy}. However, the definition of homotopy colimit used by Thomason, albeit equivalent, is not the same as the one we used in this dissertation and is slightly outdated. A more modern proof of the theorem can be found in \cite[Proposition 2.3.1 and Théorème 1.3.7]{maltsiniotis2005theorie}.
\end{proof}
\begin{corollary}\label{cor:thomhmtpycol}
  Let $A$ be a $1$\nbd{}category. The canonical map
  \[
  \hocolim^{\Th}_{a \in A}(A/a) \to A
  \]
  induced by the co-cone $(A/a \to A)_{a \in \Ob(A)}$, is an isomorphism of $\ho(\Cat^{\Th})$. 
\end{corollary}
\begin{proof}
  For every object $a$ of $A$, the canonical map to the terminal category
  \[
  A/a \to \sD_0
  \]
  is a Thomason equivalence. This comes from the fact that $A/a$ is oplax contractible (Proposition \ref{prop:slicecontractible}), or from \cite[Section 1, Corollary 2]{quillen1973higher} and the fact that $A/a$ has a terminal object.

  In particular, the morphism of functors
  \[
  A/(-) \Rightarrow k_{\sD_0},
  \]
  where $k_{\sD_0}$ is the constant functor $A \to \Cat$ with value the terminal category $\sD_0$, is a pointwise Thomason equivalence. It follows from the first part of Theorem \ref{thm:Thomason} that
  \[
  \int_{a \in A}A/a \to \int_{a \in A}k_{\sD_0}
  \]
  is a Thomason equivalence and an immediate computation shows that \[\int_{a \in A}k_{\sD_0} \simeq A.\] From the second part of Theorem \ref{thm:Thomason}, we have that
  \[
  \hocolim^{\Th}_{a \in A}(A/a) \simeq A.
  \]
  A thorough analysis of all the isomorphisms involved shows that this last isomorphism is indeed induced by the co-cone $(A/a \to A)_{a \in \Ob(A)}$.
\end{proof}
\begin{remark}
  It is possible to extend the previous corollary to prove that for every functor $f : X \to A$ ($X$ and $A$ being $1$\nbd{}categories), we have \[\hocolim^{\Th}_{a \in A} (X/a) \simeq X.\] However, to prove that it is also the case when $X$ is an $\oo$\nbd{}category and $f$ an $\oo$\nbd{}functor, as in Corollary \ref{cor:folkhmtpycol}, one would need to extend the Grothendieck construction to functors with value in $\oo\Cat$ and to prove an $\oo$\nbd{}categorical analogue of Theorem \ref{thm:Thomason}. Such results, while being highly plausible, go beyond the scope of this dissertation.
\end{remark}
Putting all the pieces together, we are now able to prove the awaited theorem.
\begin{theorem}\label{thm:categoriesaregood}
  Every $1$\nbd{}category is \good{}.
\end{theorem}
\begin{proof}
  All the arguments of the proof have already been given and we sum them up here essentially for the sake of clarity.
  Let $A$ be a $1$\nbd{}category. Consider the diagram
  \begin{align*}
    A &\to \oo\Cat\\
    a &\mapsto A/a
  \end{align*}
  and the co-cone
  \[
  (A/a \to A)_{a \in \Ob(A)}.
  \]
  \begin{itemize}[label=-]
  \item The canonical map of $\ho(\oo\Cat^{\folk})$
    \[
    \hocolim_{a \in A}^{\folk} (A/a) \to A
    \]
    is an isomorphism thanks to Corollary \ref{cor:folkhmtpycol} applied to $\mathrm{id}_A : A \to A$.
  \item The canonical map of $\ho(\oo\Cat^{\Th})$
    \[
    \hocolim_{a \in A}^{\Th} (A/a) \to A
    \]
    is an isomorphism thanks to Corollary \ref{cor:thomhmtpycol} and the fact
    that the canonical morphisms of op-prederivators $\Ho(\Cat^{\Th}) \to
    \Ho(\oo\Cat^{\Th})$ is homotopy cocontinuous (see \ref{paragr:thomhmtpycol}).
    \item Every $A/a$ is \good{} thanks to Proposition \ref{prop:contractibleisgood} and Proposition \ref{prop:slicecontractible}.
  \end{itemize}
  Thus, Proposition \ref{prop:criteriongoodcat} applies and this proves that $A$ is \good{}.
\end{proof}


\chapter{Homotopy and homology type of free 2-categories}
\chaptermark{Homology of free $2$-categories}
\section{Preliminaries: the case of free 1-categories}\label{section:prelimfreecat}
In this section, we review some homotopical results on free
($1$-)categories that will be of great help in the sequel.
\begin{paragr}
  A \emph{reflexive graph} $G$ consists of the data of two sets $G_0$ and $G_1$
  together with
  \begin{itemize}[label=-]
  \item a ``source'' map $\src : G_1 \to G_0$,
  \item a ``target'' map $\trgt : G_1 \to G_0$,
  \item a ``unit'' map $1_{(-)} : G_0 \to G_1$,
  \end{itemize}
  such that for every $x \in G_0$,
  \[
    \src(1_{x}) = \trgt (1_{x}) = x.
  \]
  The vocabulary of categories is used: elements of $G_0$ are \emph{objects} or
  \emph{$0$-cells}, elements of $G_1$ are \emph{arrows} or \emph{$1$-cells},
  arrows of the form $1_{x}$ with $x$ an object are \emph{units}, etc. A
  \emph{morphism of reflexive graphs} $ f : G \to G'$ consists of maps $f_0 \colon
  G_0 \to G'_0$ and  $f_1 : G_1 \to G'_1$ that commute with sources, targets and
  units in an obvious sense. This defines the category $\Rgrph$ of reflexive
  graphs. Later we will make use of monomorphisms in the category $\Rgrph$; they
  are the morphisms $f : G \to G'$ that are injective on objects and on arrows,
  i.e. such that $f_0 : G_0 \to G_0'$ and $f_1 : G_1 \to G'_1$ are injective.
  
  There is an ``underlying reflexive graph'' functor
  \[
    U : \Cat \to \Rgrph,
  \]
  which has a left adjoint
  \[
    L : \Rgrph \to \Cat.
  \]
  For a reflexive graph $G$, the objects of $L(G)$ are exactly the objects of
  $G$ and an arrow $f$ of $L(G)$ is a chain
  \[
    \begin{tikzcd}
      X_0 \ar[r,"f_1"]& X_1\ar[r,"f_2"]& X_2 \ar[r]& \cdots \ar[r] &X_{n-1}
      \ar[r,"f_n"]& X_{n}
    \end{tikzcd}
  \]
  of arrows of $G$, such that \emph{none} of the $f_k$ are units. The integer
  $n$ is referred to as the \emph{length} of $f$ and is denoted by $\ell(f)$.
  Composition is given by concatenation of chains.
\end{paragr}
\begin{lemma}
  A category $C$ is free in the sense of \ref{def:freeoocat} if and only if
  there exists a reflexive graph $G$ such that
  \[
    C \simeq L(G).
  \]
\end{lemma}
\begin{proof}
  If $C$ is free, consider the reflexive graph $G$ such that $G_0 = C_0$ and
  $G_1$ is the subset of $C_1$ whose elements are either generating $1$-cells of
  $C$ or units. It is straightforward to check that $C\simeq L(G)$.

  Conversely, if $C \simeq L(G)$ for some reflexive graph $G$, then the
  description of the arrows of $L(G)$ given in the previous paragraph shows that
  $C$ is free and that its set of generating $1$-cells is (isomorphic to) the
  set of non unital $1$-cells of $G$.
\end{proof}
\begin{remark}
  In other words, a category is free on a graph if and only if it is free on a
  reflexive graph. The difference between these two notions is at the level of
  morphisms: there are more morphisms of reflexive graphs because (generating)
  $1$\nbd{}cells may be sent to units. Hence, for a morphism of reflexive graphs
  $f : G \to G'$, the induced functor $L(f)$ is not necessarily rigid in the
  sense of Definition \ref{def:rigidmorphism}.
\end{remark}
\begin{paragr}
  There is another important description of the category $\Rgrph$. Write
  $\Delta_{\leq 1}$ for the full subcategory of $\Delta$ spanned by $[0]$ and $[1]$.
  The category $\Rgrph$ is nothing but $\Psh{\Delta_{\leq 1}}$, the
  category of pre-sheaves on $\Delta_{\leq 1}$. In particular, the canonical
  inclusion $i : \Delta_{\leq 1} \rightarrow \Delta$ induces by pre-composition
  a functor
  \[
    i^* : \Psh{\Delta} \to \Rgrph,
  \]
  which, by the usual technique of Kan extensions, has a left adjoint
  \[
    i_! : \Rgrph \to \Psh{\Delta}.
  \]
  For a graph $G$, the simplicial set $i_!(G)$ has $G_0$ as its set of
  $0$-simplices, $G_1$ as its set of $1$-simplices and all $k$-simplices are
  degenerate for $k>1$. For future reference, we put here the following lemma.
\end{paragr}
\begin{lemma}\label{lemma:monopreserved}
  The functor $i_! : \Rgrph \to \Psh{\Delta}$ preserves monomorphisms.
\end{lemma}
\begin{proof}
  What we need to show is that, given a morphism of simplicial sets
  \[
    f : X \to Y,
  \]
  if $f_0 : X_0 \to Y_0$ and $f_1 : X_1 \to Y_1$ are monomorphisms and if all
  $n$\nbd{}simplices of $X$ are degenerate for $n\geq 2$, then $f$ is a
  monomorphism. A proof of this assertion is contained in \cite[Paragraph
  3.4]{gabriel1967calculus}. The key argument is the Eilenberg--Zilber Lemma
  (Proposition 3.1 of op. cit.).
\end{proof}
\begin{paragr}
  Let us denote by $N : \Psh{\Delta} \to \Cat$ (instead of $N_1$ as in Paragraph
  \ref{paragr:nerve}) the usual nerve of categories and by $c : \Cat \to
  \Psh{\Delta}$ its left adjoint. Recall that for a (small) category $C$, an
  $n$-simplex of $N(C)$ is a chain
  \[
    \begin{tikzcd}
      X_0 \ar[r,"f_1"]& X_1\ar[r,"f_2"]& X_2 \ar[r]& \cdots \ar[r] &X_{n-1}
      \ar[r,"f_n"]& X_{n}
    \end{tikzcd}
  \]
  of arrows of $C$. Such an $n$-simplex is degenerate if and only if at least
  one of the $f_k$ is a unit. It is straightforward to check that the composite
  of
  \[
    \Cat \overset{N}{\rightarrow} \Psh{\Delta} \overset{i^*}{\rightarrow} \Rgrph
  \]
  is nothing but the forgetful functor $U : \Cat \to \Rgrph$. Thus, the functor
  $L : \Rgrph \to \Cat$ is (isomorphic to) the composite of
  \[
    \Rgrph \overset{i_!}{\rightarrow} \Psh{\Delta} \overset{c}{\rightarrow}
    \Cat.
  \]

  We now review a construction due to Dwyer and Kan
  (\cite{dwyer1980simplicial}). Let $G$ be a reflexive graph. For
  every $k\geq 1$, we define the simplicial set $N^k(G)$ as the
  sub-simplicial set of $N(L(G))$ whose $n$-simplices are chains
  \[
    \begin{tikzcd}
      X_0 \ar[r,"f_1"]& X_1\ar[r,"f_2"]& X_2 \ar[r]& \cdots \ar[r] &X_{n-1}
      \ar[r,"f_n"]& X_{n}
    \end{tikzcd}
  \]
  of arrows of $L(G)$ such that
  \[
    \sum_{1 \leq i \leq n}\ell(f_i) \leq k.
  \] In particular, we have
  \[
    N^1(G)=i_!(G)
  \]
  and the transfinite composition of
  \[
    i_!(G) = N^1(G) \hookrightarrow N^2(G) \hookrightarrow \cdots
    \hookrightarrow N^{k}(G) \hookrightarrow N^{k+1}(G) \hookrightarrow \cdots
  \]
  is easily seen to be the map
  \[
    \eta_{i_!(G)} : i_!(G) \to Nci_!(G),
  \]
  where $\eta$ is the unit of the adjunction $c \dashv N$.
\end{paragr}
\begin{lemma}[Dwyer--Kan]\label{lemma:dwyerkan}
  For every $k\geq 1$, the canonical inclusion map
  \[
    N^{k}(G) \to N^{k+1}(G)
  \]
  is a trivial cofibration of simplicial sets.
\end{lemma}
\begin{proof}
  Let $A_{k+1}=\mathrm{Im}(\partial_0)\cup\mathrm{Im}(\partial_{k+1})$ be the
  union of the first and last face of the standard $(k+1)$-simplex
  $\Delta_{k+1}$. Notice that the canonical inclusion
  \[
    A_{k+1} \hookrightarrow \Delta_{k+1}
  \]
  is a trivial cofibration. Let $I_{k+1}$ be the set of chains
  \[
    \begin{tikzcd}
      f = X_0 \ar[r,"f_1"]& X_1\ar[r,"f_2"]& X_2 \ar[r]& \cdots \ar[r] &X_{k-1}
      \ar[r,"f_k"]& X_{k}\ar[r,"f_{k+1}"]& X_{k+1}
    \end{tikzcd}
  \]
  of arrows of $L(G)$ such that for every $1 \leq i \leq k+1$
  \[
    \ell(f_i)=1,
  \]
  i.e.\ each $f_i$ is a non-unit arrow of $G$. For every $f \in I_{k+1}$, we
  define a morphism $\varphi_f : A_{k+1} \to N^{k}(G)$ in the following fashion:
  \begin{itemize}
  \item[-]$\varphi_{f}\vert_{\mathrm{Im}(\partial_0)}$ is the $k$-simplex of
    $N^{k}(G)$
    \[
      \begin{tikzcd}
        X_1\ar[r,"f_2"]& X_2 \ar[r]& \cdots \ar[r] &X_{k} \ar[r,"f_{k+1}"]&
        X_{k+1},
      \end{tikzcd}
    \]
  \item[-] $\varphi_{f}\vert_{\mathrm{Im}(\partial_{k+1})}$ is the $k$-simplex
    of $N^{k}(G)$
    \[
      \begin{tikzcd}
        X_0 \ar[r,"f_1"]& X_1\ar[r,"f_2"]& X_2 \ar[r]& \cdots \ar[r] &X_{k-1}
        \ar[r,"f_k"]& X_{k}.
      \end{tikzcd}
    \]
  \end{itemize}
  All in all, we have a cocartesian square
  \[
    \begin{tikzcd}
      \displaystyle \coprod_{f \in I_{k+1}}A_{k+1} \ar[d] \ar[r,"(\varphi_f)_f"] & N^{k}(G)\ar[d] \\
      \displaystyle \coprod_{f \in I_{k+1}}\Delta_{k+1} \ar[r] & N^{k+1}(G),
      \ar[from=1-1,to=2-2,phantom,"\ulcorner",very near end]
    \end{tikzcd}
  \]
  which proves that the right vertical arrow is a trivial cofibration.
\end{proof}
From this lemma, we deduce the following proposition.
\begin{proposition}
  Let $G$ be a reflexive graph. The map
  \[
    \eta_{i_!(G)} : i_!(G) \to Nci_!(G),
  \]
  where $\eta$ is the unit of the adjunction $c \dashv N$, is a trivial
  cofibration of simplicial sets.
\end{proposition}
\begin{proof}
  This follows from the fact that trivial cofibrations are stable by transfinite
  composition.
\end{proof}
From the previous proposition, we deduce the following very useful corollary.
\begin{corollary}\label{cor:hmtpysquaregraph}
  Let
  \[
    \begin{tikzcd}
      A \ar[d,"\alpha"] \ar[r,"\beta"] &B \ar[d,"\delta"] \\
      C \ar[r,"\gamma"]& D
    \end{tikzcd}
  \]
  be a cocartesian square in $\Rgrph$. If either $\alpha$ or $\beta$ is a
  monomorphism, then the induced square of $\Cat$
  \[
    \begin{tikzcd}
      L(A) \ar[d,"L(\alpha)"] \ar[r,"L(\beta)"]& L(B) \ar[d,"L(\delta)"] \\
      L(C) \ar[r,"L(\gamma)"]& L(D)
    \end{tikzcd}
  \]
  is Thomason homotopy cocartesian.
\end{corollary}
\begin{proof}
  Since the nerve $N$ induces an equivalence of op-prederivators
  \[
    \Ho(\Cat^{\Th}) \to \Ho(\Psh{\Delta}),
  \]
  it suffices to prove that the induced square of simplicial sets
  \[
    \begin{tikzcd}
      NL(A) \ar[d,"NL(\alpha)"] \ar[r,"NL(\beta)"]& NL(B) \ar[d,"NL(\delta)"] \\
      NL(C) \ar[r,"NL(\gamma)"]& NL(D)
    \end{tikzcd}
  \]
  is homotopy cocartesian. But, since $L \simeq c\circ i_!$, it follows from
  Lemma \ref{cor:hmtpysquaregraph} that this last square is weakly equivalent to
  the square of simplicial sets
  \[
    \begin{tikzcd}
      i_!(A) \ar[d,"i_!(\alpha)"] \ar[r,"i_!(\beta)"] &i_!(B) \ar[d,"i_!(\delta)"] \\
      i_!(C) \ar[r,"i_!(\gamma)"]& i_!(D).
    \end{tikzcd}
  \]
  This square is cocartesian because $i_!$ is a left adjoint. Since
  $i_!$ preserves monomorphisms (Lemma \ref{lemma:monopreserved}), the
  result follows from Lemma \ref{lemma:hmtpycocartesianreedy} and the fact that
  the monomorphisms are the cofibrations of the standard Quillen model structure on simplicial
  sets.
\end{proof}
\begin{paragr}
  By working a little more, we obtain the more general result stated
  in the proposition below. Let us say that a morphism of reflexive
  graphs $\alpha : A \to B$ is \emph{quasi-injective on arrows} when
  for all  arrows $f$ and $g$ of $A$, if
  \[
    \alpha(f)=\alpha(g),
  \]
  then either $f=g$, or $f$ and $g$ are both units. In other words, $\alpha$
  never sends a non-unit arrow to a unit arrow and $\alpha$ never identifies two
  non-unit arrows. It follows that if $\alpha$ is quasi-injective on arrows and
  injective on objects, then it is also injective on arrows and hence, a
  monomorphism of $\Rgrph$.
\end{paragr}
\begin{proposition}\label{prop:hmtpysquaregraphbetter}
  Let
  \[
    \begin{tikzcd}
      A \ar[d,"\alpha"] \ar[r,"\beta"] &B \ar[d,"\delta"] \\
      C \ar[r,"\gamma"]& D
      \ar[from=1-1,to=2-2,phantom,"\ulcorner",very near end]
    \end{tikzcd}
  \]
  be a cocartesian square in $\Rgrph$. Suppose that the following two conditions
  are satisfied
  \begin{enumerate}[label=\alph*)]
  \item Either $\alpha$ or $\beta$ is injective on objects.
  \item Either $\alpha$ or $\beta$ is quasi-injective on arrows.
  \end{enumerate}
  Then, the induced square of $\Cat$
  \[
    \begin{tikzcd}
      L(A) \ar[d,"L(\alpha)"] \ar[r,"L(\beta)"] &L(B) \ar[d,"L(\delta)"] \\
      L(C) \ar[r,"L(\gamma)"] &L(D)
    \end{tikzcd}
  \]
  is Thomason homotopy cocartesian.
\end{proposition}
\begin{proof}
  The case where $\alpha$ or $\beta$ is both injective on objects and
  quasi-injective on arrows is Corollary \ref{cor:hmtpysquaregraph}. Hence, we
  only have to treat the case when $\alpha$ is injective on objects and $\beta$
  is quasi-injective on arrows; the remaining case being symmetric.

  Let use denote by $E$ the set of objects of $B$ that are in the image of
  $\beta$. We consider this set as well as the set $A_0$ of objects of $A$ as discrete reflexive graphs, i.e.\ reflexive graphs
  with no non-unit arrows. Now, let $G$ be the reflexive graph defined by the
  following cocartesian square
  \[
    \begin{tikzcd}
      A_0\ar[r] \ar[d] & E \ar[d]\\
      A \ar[r] & G, \ar[from=1-1,to=2-2,phantom,very near end,"\ulcorner"]
    \end{tikzcd}
  \]
  where the morphism \[ A_0 \to A\] is the canonical inclusion, and the
  morphism \[A_0 \to E\] is induced by the restriction of $\beta$ on objects. In other words, $G$ is
  obtained from $A$ by collapsing the objects that are identified through
  $\beta$. It admits the following explicit description: $G_0$ is (isomorphic
  to) $E$ and the set of non-unit arrows of $G$ is (isomorphic to) the set of
  non-unit arrows of $A$; the source (resp.\ target) of a non-unit arrow $f$ of
  $G$ is the source (resp.\ target) of $\beta(f)$. This completely describes $G$.

  Now, we have the following solid arrow commutative diagram
  \[
    \begin{tikzcd}
      A_0 \ar[r] \ar[d] & E  \ar[ddr,bend left]\ar[d]&\\
      A  \ar[drr,bend right,"\beta"'] \ar[r] & G \ar[dr, dotted]&\\
      &&B, \ar[from=1-1,to=2-2,phantom,very near end,"\ulcorner"]
    \end{tikzcd}
  \]
  where the arrow $E \to B$ is the canonical inclusion. Hence, by universal
  property, the dotted arrow exists and makes the whole diagram commute. A
  thorough verification easily shows that, because $\beta$ is quasi-injective on
  arrows, the morphism $G \to B$ is a
  monomorphism of $\Rgrph$.

  By forming successive cocartesian squares and combining with the square
  obtained earlier, we obtain a diagram of three cocartesian squares:
  \[
    \begin{tikzcd}[row sep = large]
      A_0\ar[r] \ar[d] & E \ar[d]&\\
      A \ar[d,"\alpha"] \ar[r] & G \ar[d] \ar[r] & B \ar[d,"\delta"]\\
      C \ar[r] & H \ar[r] & D. \ar[from=1-1,to=2-2,phantom,"\ulcorner" very near
      end,"\text{\textcircled{\tiny \textbf{1}}}" near start, description]
      \ar[from=2-1,to=3-2,phantom,"\ulcorner" very near
      end,"\text{\textcircled{\tiny \textbf{2}}}", description]
      \ar[from=2-2,to=3-3,phantom,"\ulcorner" very near
      end,"\text{\textcircled{\tiny \textbf{3}}}", description]
    \end{tikzcd}
  \]
  What we want to prove is that the image by the functor $L$ of the pasting of
  squares \textcircled{\tiny \textbf{2}} and \textcircled{\tiny \textbf{3}} is
  homotopy cocartesian. Since the morphism $G \to B$ is a monomorphism, we
  deduce from Corollary \ref{cor:hmtpysquaregraph} that the image by the functor
  $L$ of square \textcircled{\tiny \textbf{3}} is homotopy cocartesian. Hence,
  in virtue of Lemma \ref{lemma:pastinghmtpycocartesian}, all we have to show is
  that the image by $L$ of square \textcircled{\tiny \textbf{2}} is homotopy
  cocartesian. On the other hand, the morphisms
  \[
    A_0 \to A 
  \]
  and
  \[
    A_0 \to C
  \]
  are monomorphisms and thus, using Corollary
  \ref{cor:hmtpysquaregraph}, we deduce that the image by $L$ of square
  \textcircled{\tiny \textbf{1}} and of the pasting of squares
  \textcircled{\tiny \textbf{1}} and \textcircled{\tiny \textbf{2}}
  are homotopy cocartesian. By Lemma \ref{lemma:pastinghmtpycocartesian} again, this proves that the image by $L$ of
  square \textcircled{\tiny \textbf{2}} is homotopy cocartesian.
\end{proof}
We now apply Corollary \ref{cor:hmtpysquaregraph} and Proposition
\ref{prop:hmtpysquaregraphbetter} to a few examples.
\begin{example}[Identifying two objects]\label{example:identifyingobjects}
  Let $C$ be a free category, $A$ and $B$ be two objects of $C$ with $A\neq B$ and let $C'$ be
  the category obtained from $C$ by identifying $A$ and $B$, i.e.\ defined by
  the following cocartesian square
  \[
    \begin{tikzcd}
      \sS_0 \ar[d] \ar[r,"{\langle A,B \rangle}"] & C \ar[d] \\
      \sD_0 \ar[r] & C'.
      \ar[from=1-1,to=2-2,phantom,"\ulcorner",very near end]
    \end{tikzcd}
  \]
  Then, this square is Thomason homotopy cocartesian. Indeed, it is obviously
  the image by the functor $L$ of a cocartesian square of $\Rgrph$ and the top
  morphism is a monomorphism. Hence, we can apply Corollary \ref{cor:hmtpysquaregraph}.
  \end{example}
\begin{example}[Adding a generator]
  Let $C$ be a free category, $A$ and $B$ two objects of $C$ (possibly equal)
  and let $C'$ be the category obtained from $C$ by adding a generator $A \to
  B$, i.e.\ defined by the following cocartesian square:
  \[
    \begin{tikzcd}
      \sS_0 \ar[d,"i_1"] \ar[r,"{\langle A, B \rangle}"] & C \ar[d] \\
      \sD_1 \ar[r] & C'.
      \ar[from=1-1,to=2-2,phantom,"\ulcorner",very near end]
    \end{tikzcd}
  \]
  Then, this square is Thomason homotopy cocartesian. Indeed, it obviously is the image of a square of $\Rgrph$ by
  the functor $L$ and the morphism $i_1 : \sS_0 \to \sD_1$ comes from
  a monomorphism of $\Rgrph$. Hence, we can apply Corollary
  \ref{cor:hmtpysquaregraph}.
\end{example}
\begin{remark}
  Since $i_1 : \sS_0 \to \sD_1$ is a folk cofibration 
  and since every free category is obtained by recursively adding generators
  starting from a set of objects (seen as a $0$-category), the previous example
  yields another proof that \emph{free} (1\nbd{})categories are \good{} (which we
  already knew since we have seen that \emph{all} (1-)categories are \good{}).
\end{remark}
\begin{example}[Identifying two generators]
  Let $C$ be a free category and let $f,g : A \to B$ be parallel generating arrows of
  $C$ such that $f\neq g$. Now consider the category $C'$ obtained from $C$ by
  ``identifying'' $f$ and $g$, i.e. defined by the following cocartesian
  square
  \[
    \begin{tikzcd}
      \sS_1\ar[d] \ar[r,"{\langle f, g \rangle}"] &C \ar[d] \\
      \sD_1 \ar[r] & C',
      \ar[from=1-1,to=2-2,phantom,"\ulcorner",very near end]
    \end{tikzcd}
  \]
  where the morphism $\sS_1 \to \sD_1$ is the one that sends the two generating
  arrows of $\sS_1$ to the unique generating arrow of $\sD_1$. Then this square
  is Thomason homotopy cocartesian.
  Indeed, it is the image by the functor $L$ of a cocartesian square in
  $\Rgrph$, the morphism $\sS_1 \to \sD_1$ is injective on objects and the
  morphism $\sS_1 \to C$ is quasi-injective on arrows. Hence, we can apply
  Proposition \ref{prop:hmtpysquaregraphbetter}. Note that since we did
  \emph{not} suppose that $A\neq B$, the top morphism of the previous square is
  not necessarily a monomorphism and we cannot always apply Corollary
  \ref{cor:hmtpysquaregraph}.
\end{example}
\begin{example}[Killing a generator]\label{example:killinggenerator}
  Let $C$ be a free category and let $f : A \to B$ be one of its generating arrows
  such that $A \neq B$. Now consider the category $C'$ obtained from $C$ by
  ``killing'' $f$, i.e. defined by the following cocartesian square:
  \[
    \begin{tikzcd}
      \sD_1 \ar[d] \ar[r,"\langle f \rangle"] & C \ar[d] \\
      \sD_0 \ar[r] & C'.
      \ar[from=1-1,to=2-2,phantom,"\ulcorner",very near end]
    \end{tikzcd}
  \]
  Then, this square is Thomason homotopy cocartesian. Indeed, it
  obviously is the image of a cocartesian square in $\Rgrph$ by the
  functor $L$ and since the source and target of $f$ are different,
  the top map comes from a monomorphism of $\Rgrph$. Hence, we can
  apply Corollary \ref{cor:hmtpysquaregraph}.
\end{example}
\begin{remark}
  Note that in the previous example, we see that it was useful to consider the
  category of reflexive graphs and not only the category of graphs because the
  map $\sD_1 \to \sD_0$ does not come from a morphism in the category of graphs.

  Note also that the hypothesis that $A\neq B$ was fundamental in the previous
  example as for $A=B$ the square is \emph{not} Thomason homotopy cocartesian.
\end{remark}

\section{Preliminaries: bisimplicial sets}
\begin{paragr}
  A \emph{bisimplicial set} is a presheaf over the category $\Delta \times
  \Delta$,
  \[
    X : \Delta^{\op} \times \Delta^{\op} \to \Set.
  \]
  In a similar fashion as for simplicial sets (\ref{paragr:simpset}), for $n,m
  \geq 0$, we use the notations
  \begin{align*}
    X_{n,m} &:= X([n],[m]) \\
    \partial_i^v &:=X(\delta^i,\mathrm{id}) : X_{n+1,m} \to X_{n,m}\\
    \partial_j^h &:=X(\mathrm{id},\delta^j) : X_{n,m+1} \to X_{n,m}\\
    s_i^v &:=X(\sigma^i,\mathrm{id}): X_{n,m} \to X_{n+1,m}\\
    s_j^h&:=X(\mathrm{id},\sigma^j) : X_{n,m} \to X_{n,m+1}.
  \end{align*}
  The maps $\partial_i^h$ and $s_i^h$ will be referred to as the
  \emph{horizontal} face and degeneracy operators; and $\partial_i^v$ and
  $s_i^v$ as the \emph{vertical} face and degeneracy operators.
 
  Note that for every $n\geq 0$, we have simplicial sets
  \begin{align*}
    X_{\bullet,n} : \Delta^{\op} &\to \Set \\
    [k] &\mapsto X_{k,n}
  \end{align*}
  and
  \begin{align*}
    X_{n,\bullet} : \Delta^{\op} &\to \Set \\
    [k] &\mapsto X_{n,k}.
  \end{align*}
  The category of bisimplicial sets is denoted by $\Psh{\Delta\times\Delta}$.

\end{paragr}
\begin{paragr}
  The functor
  \begin{align*}
    \delta : \Delta &\to \Delta\times\Delta \\
    [n] &\mapsto ([n],[n])
  \end{align*}
  induces by pre-composition a functor
  \[
    \delta^* : \Psh{\Delta\times\Delta} \to \Psh{\Delta}.
  \]
  By the usual calculus of Kan extensions, $\delta^*$ admits a left adjoint
  $\delta_!$ and a right adjoint $\delta_*$
  \[
    \delta_! \dashv \delta^* \dashv \delta_*.
  \]
  We say that a morphism $f : X \to Y$ of bisimplicial sets is a \emph{diagonal
    weak equivalence} (resp.\ \emph{diagonal fibration}) when $\delta^*(f)$ is a
  weak equivalence (resp.\ fibration) of simplicial sets. By
  definition, $\delta^*$ induces a morphism of op-prederivators
  \[
    \overline{\delta^*} : \Ho(\Psh{\Delta\times \Delta}^{\mathrm{diag}}) \to
    \Ho(\Psh{\Delta}),
  \]
  where $\Ho(\Psh{\Delta\times \Delta}^{\mathrm{diag}})$ is the homotopy
  op-prederivator of $\Psh{\Delta\times \Delta}$ equipped with diagonal weak equivalences.
  Recall from \cite[Proposition 1.2]{moerdijk1989bisimplicial} that the category
  of bisimplicial sets can be equipped with a model structure whose weak
  equivalences are the diagonal weak equivalences and whose fibrations are the
  diagonal fibrations. We shall refer to this model
  structure as the \emph{diagonal model structure}.
\end{paragr}
\begin{proposition}\label{prop:diageqderivator}
  Consider that $\Psh{\Delta\times\Delta}$ is equipped with the diagonal model
  structure. Then, the adjunction
  \[
    \begin{tikzcd}
      \delta_! : \Psh{\Delta} \ar[r,shift left] & \Psh{\Delta\times\Delta}
      \ar[l,shift left]: \delta^*,
    \end{tikzcd}
  \]
  is a Quillen equivalence.
\end{proposition}
\begin{proof}
  By definition $\delta^*$ preserves weak equivalences and fibrations and thus,
  the adjunction is a Quillen adjunction. The fact that $\delta^*$ induces an
  equivalence at the level of homotopy categories is \cite[Proposition
  1.2]{moerdijk1989bisimplicial}.
\end{proof}
\begin{paragr}
  In particular, the morphism of op-prederivators
  \[
    \overline{\delta^*} : \Ho(\Psh{\Delta\times \Delta}^{\mathrm{diag}}) \to
    \Ho(\Psh{\Delta})
  \]
  is actually an equivalence of op-prederivators.
\end{paragr}
Diagonal weak equivalences are not the only interesting weak equivalences for
bisimplicial sets.
\begin{paragr}
  A morphism $f : X \to Y$ of bisimplicial sets is a \emph{vertical (resp.\
    horizontal) weak equivalence} when for every $n \geq 0$, the induced
  morphism of simplicial sets
  \[
    f_{\bullet,n} : X_{\bullet,n} \to Y_{\bullet,n}
  \]
  (resp.
  \[
    f_{n,\bullet} : X_{n,\bullet} \to Y_{n,\bullet})
  \]
  is a weak equivalence of simplicial sets. Recall now a very useful lemma.
\end{paragr}
\begin{lemma}\label{bisimpliciallemma}
  Let $f : X \to Y$ be a morphism of bisimplicial sets. If $f$ is a vertical or
  horizontal weak equivalence then it is a diagonal weak equivalence.
\end{lemma}
\begin{proof}
  See for example \cite[Chapter XII,4.3]{bousfield1972homotopy} or
  \cite[Proposition 2.1.7]{cisinski2004localisateur}.
\end{proof}
\begin{paragr}
  In particular, the identity functor of the category of bisimplicial sets
  induces the morphisms of op-prederivators:
  \[
    \Ho(\Psh{\Delta\times\Delta}^{\mathrm{vert}}) \to
    \Ho(\Psh{\Delta\times\Delta}^{\mathrm{diag}})
  \]
  and
  \[
    \Ho(\Psh{\Delta\times\Delta}^{\mathrm{hor}}) \to
    \Ho(\Psh{\Delta\times\Delta}^{\mathrm{diag}}).
  \]
\end{paragr}
\begin{proposition}\label{prop:bisimplicialcocontinuous}
  The morphisms of op-prederivators
  \[
    \Ho(\Psh{\Delta\times\Delta}^{\mathrm{vert}}) \to
    \Ho(\Psh{\Delta\times\Delta}^{\mathrm{diag}})
  \]
  and
  \[
    \Ho(\Psh{\Delta\times\Delta}^{\mathrm{hor}}) \to
    \Ho(\Psh{\Delta\times\Delta}^{\mathrm{diag}})
  \]
  are homotopy cocontinuous.
\end{proposition}
\begin{proof}
  Recall that the category of bisimplicial sets can be equipped with a model
  structure where the weak equivalences are the vertical (resp.\ horizontal)
  weak equivalences and the cofibrations are the monomorphisms (see for example
  \cite[Chapter IV]{goerss2009simplicial} or \cite{cisinski2004localisateur}).
  We respectively refer to these model structures as the \emph{vertical model
    structure} and \emph{horizontal model structure}. Since the functor
  $\delta^* : \Psh{\Delta\times\Delta} \to \Psh{\Delta}$ preserves
  monomorphisms, it follows from Lemma \ref{bisimpliciallemma} that the
  adjunction
  \[
    \begin{tikzcd}
      \delta^* : \Psh{\Delta\times\Delta} \ar[r,shift left] & \ar[l,shift left]
      \Psh{\Delta} : \delta_*
    \end{tikzcd}
  \]
  is a Quillen adjunction when $\Psh{\Delta\times\Delta}$ is equipped with
  either the vertical model structure or the horizontal model structure. In
  particular, the induced morphisms of op-prederivators
  \[
    \overline{\delta^*} : \Ho(\Psh{\Delta\times\Delta}^{\mathrm{vert}}) \to
    \Ho(\Psh{\Delta})
  \]
  and
  \[
    \overline{\delta^*} : \Ho(\Psh{\Delta\times\Delta}^{\mathrm{hor}}) \to
    \Ho(\Psh{\Delta})
  \]
  are homotopy cocontinuous. Now, the obvious identity
  $\delta^*=\delta^* \circ \mathrm{id}_{\Psh{\Delta\times\Delta}}$ implies that
  we have commutative triangles
  \[
    \begin{tikzcd}
      \Ho(\Psh{\Delta\times\Delta}^{\mathrm{vert}}) \ar[r]
      \ar[rd,"\overline{\delta^*}"']&
      \Ho(\Psh{\Delta\times\Delta}^{\mathrm{diag}})
      \ar[d,"\overline{\delta^*}"] \\
      &\Ho(\Psh{\Delta})
    \end{tikzcd}
  \]
  and
  \[
    \begin{tikzcd}
      \Ho(\Psh{\Delta\times\Delta}^{\mathrm{hor}}) \ar[r]
      \ar[rd,"\overline{\delta^*}"']&
      \Ho(\Psh{\Delta\times\Delta}^{\mathrm{diag}})
      \ar[d,"\overline{\delta^*}"] \\
      &\Ho(\Psh{\Delta}).
    \end{tikzcd}
  \]
  The result follows then from the fact that $\overline{\delta^*} :
  \Ho(\Psh{\Delta\times\Delta}^{\mathrm{diag}}) \to \Ho(\Psh{\Delta})$ is an
  equivalence of op-prederivators.
\end{proof}
In practice, we will use the following corollary.
\begin{corollary}\label{cor:bisimplicialsquare}
  Let
  \[
    \begin{tikzcd}
      A \ar[r,"u"]\ar[d,"f"] & B \ar[d,"g"] \\
      C \ar[r,"v"] & D
    \end{tikzcd}
  \]
  be a commutative square in the category of bisimplicial sets satisfying at least one of the two
  following conditions:
  \begin{enumerate}[label=(\alph*)]
  \item For every $n\geq 0$, the square of simplicial sets
    \[
      \begin{tikzcd}
        A_{\bullet,n} \ar[r,"{u_{\bullet,n}}"]\ar[d,"{f_{\bullet,n}}"] & B_{\bullet,n} \ar[d,"{g_{\bullet,n}}"] \\
        C_{\bullet,n} \ar[r,"{v_{\bullet,n}}"] & D_{\bullet,n}
      \end{tikzcd}
    \]
    is homotopy cocartesian.
  \item For every $n\geq 0$, the square of simplicial sets
    \[
      \begin{tikzcd}
        A_{n,\bullet} \ar[r,"{u_{n,\bullet}}"]\ar[d,"{f_{n,\bullet}}"] & B_{n,\bullet} \ar[d,"{g_{n,\bullet}}"] \\
        C_{n,\bullet} \ar[r,"{v_{n,\bullet}}"] & D_{n,\bullet}
      \end{tikzcd}
    \]
    is homotopy cocartesian.
  \end{enumerate}
  Then, the square
  \[
    \begin{tikzcd}
      \delta^*(A) \ar[r,"\delta^*(u)"]\ar[d,"\delta^*(f)"] & \delta^*(B) \ar[d,"\delta^*(g)"] \\
      \delta^*(C) \ar[r,"\delta^*(v)"] & \delta^*(D)
    \end{tikzcd}
  \]
  is a homotopy cocartesian square of simplicial sets.
\end{corollary}
\begin{proof}
  From \cite[Corollary 10.3.10(i)]{groth2013book} we know that the square of
  bisimplicial sets
  \[
    \begin{tikzcd}
      A \ar[r,"u"]\ar[d,"f"] & B \ar[d,"g"] \\
      C \ar[r,"v"] & D
    \end{tikzcd}
  \]
  is homotopy cocartesian with respect to the vertical weak equivalences if and
  only if for every $n\geq 0$, the square
  \[
    \begin{tikzcd}
      A_{\bullet,n} \ar[r,"{u_{\bullet,n}}"]\ar[d,"{f_{\bullet,n}}"] & B_{\bullet,n} \ar[d,"{g_{\bullet,n}}"] \\
      C_{\bullet,n} \ar[r,"{v_{\bullet,n}}"] & D_{\bullet,n}
    \end{tikzcd}
  \]
  is a homotopy cocartesian square of simplicial sets and similarly for
  horizontal weak equivalences. The result follows then from Proposition
  \ref{prop:bisimplicialcocontinuous}.
\end{proof}
\section{Bisimplicial nerve for 2-categories}\label{section:bisimplicialnerve}
We shall now describe a ``nerve'' for $2$\nbd{}categories with values in bisimplicial
sets and recall a few results that shows that this nerve is, in some sense,
equivalent to the nerve defined in \ref{paragr:nerve}.
\begin{notation}
  \begin{itemize}
  \item[-] Once again, we write $N : \Cat \to \Psh{\Delta}$ instead of $N_1$ for
    the usual nerve of categories. Moreover, using the usual notation for the
    set of $k$-simplices of a simplicial set, if $C$ is a (small) category, then
    \[
      N(C)_k
    \]
    is the set of $k$-simplices of the nerve of $C$.
  \item[-] Similarly, we write $N : 2\Cat \to \Psh{\Delta}$ instead of $N_2$ for
    the nerve of $2$\nbd{}categories. This makes sense since the nerve for categories
    is the restriction of the nerve for $2$\nbd{}categories.
  \item[-] For $2$\nbd{}categories, we refer to the $\comp_0$-composition of
    $2$\nbd{}cells as the \emph{horizontal composition} and the $\comp_1$-composition
    of $2$\nbd{}cells as the \emph{vertical composition}.
  \item[-] For a $2$\nbd{}category $C$ and $x$ and $y$ objects of $C$, we denote by
    \[
      C(x,y)
    \]
    the category whose objects are the $1$-cells of $C$ with $x$ as source and
    $y$ as target, and whose arrows are the $2$\nbd{}cells of $C$ with $x$ as
    $0$-source and $y$ as $0$-target. Composition is induced by vertical
    composition in $C$.
  \end{itemize}
\end{notation}
\begin{paragr}
  Every $2$\nbd{}category $C$ defines a simplicial object in $\Cat$,
  \[S(C): \Delta^{\op} \to \Cat,\] where, for each $n \geq 0$,
  \[
    S_n(C):= \coprod_{(x_0,\cdots,x_n)\in \Ob(C)^{\times (n+1)}}C(x_0,x_1)
    \times \cdots \times C(x_{n-1},x_n).
  \]
  Note that for $n=0$, the above formula reads $S_0(C)=C_0$. For $n>0$, the face operators  $\partial_i : S_{n}(C) \to S_{n-1}(C)$ are induced by horizontal
  composition for $0 < i <n$ and by projection for $i=0$ or $n$. The degeneracy operators $s_i \colon S_{n}(C) \to S_{n+1}(C)$ are
  induced by the units for the horizontal composition.

  Post-composing $S(C)$ with the nerve functor $N : \Cat \to \Psh{\Delta}$, we
  obtain a functor
  \[
    NS(C) : \Delta^{\op} \to \Psh{\Delta}.
  \]
\end{paragr}
\begin{remark}
  When $C$ is a $1$-category, the simplicial object $S(C)$ is nothing but the
  usual nerve of $C$ where, for each $n\geq 0$, $S_n(C)$ is seen as a discrete
  category.
\end{remark}
\begin{definition}
  The \emph{bisimplicial nerve} of a $2$\nbd{}category $C$ is the bisimplicial set
  $\binerve(C)$ defined as
  \[
    \binerve(C)_{n,m}:=N(S_n(C))_m,
  \]
  for all $n,m \geq 0$.
\end{definition}
\begin{paragr}\label{paragr:formulabisimplicialnerve}
  In other words, the bisimplicial nerve of $C$ is obtained by ``un-currying''
  the functor $NS(C) : \Delta^{\op} \to \Psh{\Delta}$.
  
  Since the nerve $N$ commutes with products and sums, we obtain the formula
  \begin{equation}\label{fomulabinerve}
    \binerve(C)_{n,m} = \coprod_{(x_0,\cdots,x_n)\in \Ob(C)^{\times (n+1)}}N(C(x_0,x_1))_m \times \cdots \times N(C(x_{n-1},x_n))_m.
  \end{equation}
  More intuitively, an element of $\binerve(C)_{n,m}$ consists of a ``pasting
  scheme'' in $C$ that looks like
  \[
    m \underbrace{\left\{\begin{tikzcd}[column sep=huge,ampersand
          replacement=\&] \bullet \ar[r,"\vdots"]\ar[r,bend
          left=90,looseness=1.4,""{name=A,below}] \ar[r,bend
          left=35,""{name=B,above}] \ar[r,bend
          right=35,"\vdots",""{name=G,below}]\ar[r,bend
          right=90,looseness=1.4,""{name=H,above}] \&
          \bullet\ar[r,"\vdots"]\ar[r,bend
          left=90,looseness=1.4,""{name=C,below}] \ar[r,bend
          left=35,""{name=D,above}] \ar[r,bend
          right=35,"\vdots",""{name=I,below}]\ar[r,bend
          right=90,looseness=1.4,""{name=J,above}]
          \&\bullet\ar[r,phantom,description,"\cdots"]\&\bullet\ar[r,"\vdots"]\ar[r,bend
          left=90,looseness=1.4,""{name=E,below}] \ar[r,bend
          left=35,""{name=F,above}] \ar[r,bend
          right=35,"\vdots",""{name=K,below}]\ar[r,bend
          right=90,looseness=1.4,""{name=L,above}] \&\bullet
          \ar[from=A,to=B,Rightarrow] \ar[from=C,to=D,Rightarrow]
          \ar[from=E,to=F,Rightarrow] \ar[from=G,to=H,Rightarrow]
          \ar[from=I,to=J,Rightarrow] \ar[from=K,to=L,Rightarrow]
        \end{tikzcd}\right.}_{ n }.
  \]
\end{paragr}
In the definition of the bisimplicial nerve of a $2$\nbd{}category we gave, one
direction of the bisimplicial set is privileged over the other. We now give
an equivalent definition of the bisimplicial nerve which uses the other direction.
\begin{paragr}
  Let $C$ be a $2$\nbd{}category. For every $k \geq 1$, we define a
  $1$\nbd{}category $V_k(C)$ in the following fashion:
  \begin{itemize}[label=-]
  \item The objects of $V_k(C)$ are the objects of $C$.
  \item A morphism $\alpha$ is a sequence
    \[
      \alpha=(\alpha_1,\alpha_2,\cdots,\alpha_k)
    \]
    of vertically composable $2$\nbd{}cells of $C$, i.e.\ such that for
    every $1 \leq i \leq k-1$, we have
    \[
      \src(\alpha_i)=\trgt(\alpha_{i+1}).
    \]
    The source and target of $\alpha$ are given by
    \[
      \src(\alpha):=\src_0(\alpha_1)\text{ and
      }\trgt(\alpha):=\trgt_0(\alpha_1).
    \]
    (Note that we could have used any of the $\alpha_i$ instead of $\alpha_1$
    since they all have the same $0$\nbd{}source and $0$\nbd{}target.)
  \item Composition is given by
    \[
      (\alpha_1,\alpha_2,\cdots,\alpha_k)\circ(\beta_1,\beta_2,\cdots,\beta_k):=(\alpha_1\comp_0\beta_1,\alpha_2\comp_0\beta_2,\cdots,\alpha_k\comp_0\beta_k)
    \]
    and the unit on an object $x$ is the sequence
    \[
      (1^2_x,\cdots, 1^2_x).
    \]
  \end{itemize}
  For $k=0$, we define $V_0(C)$ to be the category obtained from $C$ by simply
  forgetting the $2$\nbd{}cells (which is nothing but $\tau^{s}_{\leq 1}(C)$
  with the notations of \ref{paragr:defncat}). The correspondence $n \mapsto
  V_n(C)$ defines a simplicial object in $\Cat$
  \[
    V(C) : \Delta^{\op} \to \Cat,
  \]
  where the face operators are induced by the vertical composition and the
  degeneracy operators are induced by the units for the vertical composition.
\end{paragr}
\begin{lemma}\label{lemma:binervehorizontal}
  Let $C$ be a $2$\nbd{}category. For every $n \geq 0$, we have
  \[
    N(V_m(C))_n=(\binerve(C))_{n,m}.
  \]
\end{lemma}
\begin{proof}
  This is simply a reformulation of the formula given in Paragraph
  \ref{paragr:formulabisimplicialnerve}.
\end{proof}
\begin{paragr}
  The bisimplicial nerve canonically defines a functor
  \[
    \binerve : 2\Cat \to \Psh{\Delta\times\Delta}
  \]
  which enables us to compare the homotopy theory of $2\Cat$ with the homotopy
  theory of bisimplicial sets.
\end{paragr}
\begin{lemma}\label{lemma:binervthom}
  A $2$\nbd{}functor $F : C \to D$ is a Thomason equivalence if and only if
  $\binerve(F)$ is a diagonal weak equivalence of bisimplicial sets.
\end{lemma}
\begin{proof}
  It follows from what is shown in \cite[Section 2]{bullejos2003geometry} that
  there is a zigzag of weak equivalence of simplicial sets
  \[
    \delta^*(\binerve(C)) \leftarrow \cdots \rightarrow N(C)
  \]
  which is natural in $C$. This implies what we wanted to show. See also
  \cite[Théorème 3.13]{ara2020comparaison}.
\end{proof}
From this lemma, we deduce two useful criteria to detect Thomason equivalences
of $2$\nbd{}categories.
\begin{corollary}\label{cor:criterionThomeqI}
  Let $F : C \to D$ be a $2$\nbd{}functor. If
  \begin{enumerate}[label=\alph*)]
  \item $F_0 : C_0 \to D_0$ is a bijection,
  \end{enumerate}
  and
  \begin{enumerate}[resume*]
  \item for all objects $x,y$ of $C$, the functor
    \[
      C(x,y) \to D(F(x),F(y))
    \]
    induced by $F$ is a Thomason equivalence of $1$-categories,
  \end{enumerate}
  then $F$ is a Thomason equivalence of $2$\nbd{}categories.
\end{corollary}
\begin{proof}
  By definition, for every $2$\nbd{}category $C$ and every $m \geq 0$, we have
  \[
    (\binerve(C))_{\bullet,m} = NS(C).
  \]
  The result then follows from Lemma \ref{bisimpliciallemma} and the fact that the
  weak equivalences of simplicial sets are stable by coproducts and finite
  products.
\end{proof}
\begin{corollary}\label{cor:criterionThomeqII}
  Let $F : C \to D$ be a $2$\nbd{}functor. If for every $k \geq 0$,
  \[V_k(F) : V_k(C) \to V_k(D)\] is a Thomason equivalence of $1$-categories,
  then $F$ is a Thomason equivalence of $2$\nbd{}categories.
\end{corollary}
\begin{proof}
  From Lemma \ref{lemma:binervehorizontal}, we now that for every $m \geq 0$,
  \[
    \binerve(C)_{\bullet,m}=N(V_m(C)).
  \]
  The result then follows from Lemma \ref{bisimpliciallemma}.
\end{proof}
\begin{paragr}
  It also follows from Lemma \ref{lemma:binervthom} that the bisimplicial nerve
  induces a morphism of op\nbd{}prederivators
  \[
    \overline{\binerve} : \Ho(2\Cat^{\Th}) \to \Ho(\Psh{\Delta\times\Delta}^{\mathrm{diag}}).
  \]
\end{paragr}
 \begin{proposition}\label{prop:streetvsbisimplicial}
   The morphism of op\nbd{}prederivators
     \[
    \overline{\binerve} : \Ho(2\Cat^{\Th}) \to \Ho(\Psh{\Delta\times\Delta}^{\mathrm{diag}})
  \]
  is an equivalence of op\nbd{}prederivators.
\end{proposition}
\begin{proof}
   Consider the following triangle of functors
  \[
    \begin{tikzcd}
      2\Cat \ar[rr,"\binerve"] \ar[dr,"N"'] & & \Psh{\Delta\times\Delta} \ar[ld,"\delta^*"] \\
      &\Psh{\Delta}.
    \end{tikzcd}
  \]
  Even though it is \emph{not} commutative (even up to an isomorphism), it
  follows from the results contained in \cite[Section2]{bullejos2003geometry}
  that the induced triangle
    \[
    \begin{tikzcd}
      \Ho(2\Cat^{\Th}) \ar[rr,"\overline{\binerve}"] \ar[dr,"\overline{N}"'] & & \Ho(\Psh{\Delta\times\Delta}^{\mathrm{diag}}) \ar[ld,"\overline{\delta^*}"] \\
      &\Ho(\Psh{\Delta})
    \end{tikzcd}
  \]
  is commutative up to a canonical isomorphism. The result follows then from the
  fact that $\overline{\delta^*}$ and $\overline{N}$ are equivalences of
  op-prederivators (Proposition \ref{prop:diageqderivator} and Theorem
  \ref{thm:gagna} respectively).
  \end{proof}
 
%
  From Proposition \ref{prop:streetvsbisimplicial}, we deduce the proposition
  below which contains two useful criteria to detect when a commutative square
  of $2\Cat$ is Thomason homotopy cocartesian.                   
\begin{proposition}\label{prop:critverthorThomhmtpysquare}
  Let
  \begin{equation}\tag{$\ast$}\label{coucou}\begin{tikzcd}
      A \ar[r,"u"]\ar[d,"f"] & B \ar[d,"g"] \\
      C \ar[r,"v"] & D
    \end{tikzcd}\end{equation}
  be a commutative square in $2\Cat$ satisfying at least one of the two following conditions:
  \begin{enumerate}[label=(\alph*)]
  \item For every $n\geq 0$, the commutative square of $\Cat$ 
    \[
      \begin{tikzcd}
        V_{n}(A) \ar[r,"V_{n}(u)"]\ar[d,"V_{n}(f)"'] & V_n(B) \ar[d,"V_n(g)"] \\
        V_n(C) \ar[r,"V_n(v)"] & V_n(D)
      \end{tikzcd}
    \]
    is Thomason homotopy cocartesian,

  \item For every $n\geq 0$, the commutative square of $\Cat$
    \[
      \begin{tikzcd}
        S_{n}(A) \ar[r,"S_{n}(u)"]\ar[d,"S_{n}(f)"'] & S_n(B) \ar[d,"S_n(g)"] \\
        S_n(C) \ar[r,"S_n(v)"] & S_n(D)
      \end{tikzcd}
    \] is Thomason homotopy cocartesian.
  \end{enumerate}
  Then, square \eqref{coucou} is Thomason homotopy cocartesian.
\end{proposition}
\begin{proof}
  This is an immediate consequence of Proposition
  \ref{prop:streetvsbisimplicial} and Corollary \ref{cor:bisimplicialsquare}.
\end{proof}

\section{Zoology of 2-categories: basic examples}
\begin{paragr}\label{paragr:criterion2cat}
  Before embarking on computations of homology and homotopy types of
  $2$\nbd{}categories, let us recall the following particular case of Corollary
  \ref{cor:usefulcriterion}. Suppose given a cocartesian square
  \[
    \begin{tikzcd}
      A \ar[r,"u"] \ar[d,"f"] & B \ar[d,"g"] \\
      C \ar[r,"v"] & D
      \ar[from=1-1,to=2-2,phantom,"\ulcorner",very near end]
    \end{tikzcd}
  \]
  of 2\nbd{}categories. If $A$, $B$ and $C$ are free and \good{}, if at least
  one of $u$
  or $f$ is a folk cofibration and if the square is Thomason homotopy
  cocartesian, then $D$ is \good{}.
\end{paragr}
\begin{paragr}
  Let $n,m \geq 0$. We denote by $A_{(m,n)}$ the free $2$\nbd{}category with only one
  generating $2$\nbd{}cell whose source is a chain of length $m$ and whose target is  a
  chain of length $n$:
  \[
    \underbrace{\overbrace{\begin{tikzcd}[column sep=small, ampersand
          replacement=\&] \&\bullet
          \ar[r,description,"\cdots",phantom,""{name=A,below}] \& \bullet
          \ar[rd] \& \\ \bullet \ar[ru] \ar[rd] \& \& \&\bullet \\ \&\bullet
          \ar[r,description,"\cdots",phantom,""{name=B,above}] \& \bullet.
          \ar[ru]\ar[from=A,to=B,shorten <= 2em, shorten >=
          2em,Rightarrow]\end{tikzcd}}^{m}}_{n}\] More formally, $A_{(m,n)}$ is
  described in the following way:
  \begin{itemize}[label=-]
  \item generating $0$-cells: $A_0,\cdots, A_m$, $B_1,\cdots,B_{n-1}$
  \end{itemize}
  (and for convenience, we also set $B_0:=A_0$ and $B_n:=A_m$)
  \begin{itemize}[label=-]
  \item generating $1$-cells: $\begin{cases} f_{i+1} : A_i \to A_{i+1} & \text{
        for } 0\leq i \leq m-1 \\ g_{j+1} : B_j \to B_{j+1} &\text{ for } 0 \leq
      j \leq n-1 
    \end{cases}$
  \item generating $2$\nbd{}cell: $ \alpha : f_{m}\circ \cdots \circ f_1 \Rightarrow
    g_n \circ \cdots \circ g_1$.
  \end{itemize}
  Notice that $A_{(1,1)}$ is nothing but $\sD_2$. We are going to prove that if
  $n\neq 0$ or $m\neq 0$, then $A_{(m,n)}$ is \good{} and has the homotopy type
  of a point. When $m\neq0$ \emph{and} $n\neq0$, this result is not surprising,
  but when $n=0$ or $m=0$ (but not both), it is \emph{a priori} less clear what
  the homotopy type of $A_{(m,n)}$ is and whether it is \good{} or not. For
  example, $A_{(1,0)}$ can be pictured as
  \[
    \begin{tikzcd}
      A \ar[loop,in=50,out=130,distance=1.5cm,"f",""{name=A,below}]
      \ar[from=A,to=1-1,Rightarrow,"\alpha"]
    \end{tikzcd}
  \]
  and has many non trivial $2$\nbd{}cells, such as $f\comp_0 \alpha \comp_0 f$.

  Note that when $m=0$ \emph{and} $n=0$, then the $2$\nbd{}category $A_{(0,0)}$ is
  nothing but the $2$\nbd{}category $B^2\mathbb{N}$ and we have already seen that it
  is \emph{not} \good{} (see \ref{paragr:bubble}).
\end{paragr}
\begin{paragr}
  For $n\geq 0$, we write $\Delta_n$ for the linear order $\{0 \leq \cdots \leq
    n\}$ seen as a small category. Let $i : \Delta_1 \to \Delta_n$ be the unique
  functor such that
  \[
    i(0)=0 \text{ and } i(1)=n.
  \]
\end{paragr}
\begin{lemma}\label{lemma:istrngdefrtract}
  For $n\neq0$, the functor $i : \Delta_1 \to \Delta_n$ is a strong oplax deformation
  retract (\ref{paragr:defrtract}).
\end{lemma}
\begin{proof}
  Let $r : \Delta_n \to \Delta_1$ the unique functor such that
  \[
    r(0)=0\text{ and } r(k)=1 \text{ for } k>0.
  \]
  By definition we have $r \circ i = 1_{\Delta_1}$. Now, the natural order on
  $\Delta_n$ induces a natural transformation
  \[
    \alpha :  \mathrm{id}_{\Delta_n} \Rightarrow  i\circ r,
  \]
  and it is straightforward to check that $\alpha \ast i = \mathrm{id}_i$.
\end{proof}

\begin{paragr}\label{paragr:Amn}
  For every $n \geq 0$, consider the following cocartesian square
  \[
    \begin{tikzcd}
      \Delta_1 \ar[r,"i"] \ar[d,"\tau"] & \Delta_n \ar[d] \\
      A_{(1,1)} \ar[r] & A_{(1,n)}, \ar[from=1-1,to=2-2,phantom,very near
      end,"\ulcorner"]
    \end{tikzcd}
  \]
  where $\tau : \Delta_1 \to A_{(1,1)}$ is the $2$\nbd{}functor that sends the unique
  non-trivial $1$\nbd{}cell of $\Delta_1$ to the target of the generating
  $2$\nbd{}cell of $A_{(1,1)}$. It is not hard to check that $\tau$ is strong oplax
  deformation retract and thus, a co-universal Thomason equivalence (Lemma
  \ref{lemma:pushoutstrngdefrtract}). Hence, the morphism $\Delta_n \to
  A_{(1,n)}$ is also a (co-universal) Thomason equivalence and the square is
  Thomason homotopy cocartesian (Lemma \ref{lemma:hmtpycocartsquarewe}). Now,
  the morphism $\tau : \Delta_1 \to A_{(1,1)}$ is also a folk cofibration and
  since $\Delta_1$, $\Delta_n$ and $A_{(1,1)}$ are \good{}, it follows from what
  we said in \ref{paragr:criterion2cat} that $A_{(1,n)}$ is \good{}. Finally,
  since $\Delta_1$, $\Delta_n$ and $A_{(1,1)}$ have the homotopy type of a
  point, the fact that the previous square is Thomason homotopy cocartesian
  implies that $A_{(1,n)}$ has the homotopy type of a point.

  Similarly, for every $m \geq 0$, by considering the cocartesian square
  \[
    \begin{tikzcd}
      \Delta_1 \ar[r,"i"] \ar[d,"\sigma"] & \Delta_m \ar[d] \\
      A_{(1,1)} \ar[r] & A_{(m,1)}, \ar[from=1-1,to=2-2,phantom,very near
      end,"\ulcorner"]
    \end{tikzcd}
  \]
  where $\sigma : \Delta_1 \to A_{(1,1)}$ is the $2$\nbd{}functor that sends the
  unique non trivial $1$\nbd{}cell of $\Delta_1$ to the source of the generating
  $2$\nbd{}cell of $A_{(1,1)}$, we can prove that $A_{(m,1)}$ is \good{} and has
  the homotopy type of a point.

  Now, let $m\geq 0$ and $n > 0$ and consider the cocartesian square
  \[
    \begin{tikzcd}
      \Delta_1 \ar[r,"i"] \ar[d,"\tau"] & \Delta_n \ar[d] \\
      A_{(m,1)} \ar[r] & A_{(m,n)},
      \ar[from=1-1,to=2-2,phantom,very near
      end,"\ulcorner"]
    \end{tikzcd}
  \]
  where $\tau$ is the $2$\nbd{}functor that sends the unique non-trivial $1$-cell of
  $\Delta_1$ to the target of the generating $2$\nbd{}cell of $A_{(m,1)}$. This
  $2$\nbd{}functor is once again a folk cofibration, but it is \emph{not} in general
  a co-universal Thomason equivalence (it would be if we had made the hypothesis that
  $m\neq 0$, but we did not). However, since we made the hypothesis that $n\neq
  0$, it follows from Lemma \ref{lemma:istrngdefrtract} that $i : \Delta_1 \to
  \Delta_n$ is a co-universal Thomason equivalence. Hence, the previous square
  is Thomason homotopy cocartesian and $A_{(m,n)}$ has the homotopy type of a
  point. Since $A_{(m,1)}$, $\Delta_1$ and $\Delta_n$ are \good{}, this shows
  that for $m \geq 0$ and $n >0$, $A_{(m,n)}$ is \good{}.

  Similarly, if $m >0$ and $ n\geq 0$, then $A_{(m,n)}$ has the homotopy type of
  a point and is \good{}.
\end{paragr}
Combined with the result of Paragraph \ref{paragr:bubble}, we have proved the
following proposition.
\begin{proposition}\label{prop:classificationAmn}
  Let $m,n \geq 0$ and consider the $2$\nbd{}category $A_{(m,n)}$. If $m\neq 0$ or
  $n\neq 0$, then $A_{(m,n)}$ is \good{} and has the homotopy type of a point.
  If $n=m=0$, then $A_{(0,0)}$ is not \good{} and has the homotopy type of a
  $K(\mathbb{Z},2)$.
\end{proposition}
\section{Zoology of 2-categories: more examples}
As a warm-up, let us begin with an example which is a direct consequence of the
results at the end of the previous section.
\begin{paragr}
  Let $P$ the free $2$\nbd{}category defined as follows:
  \begin{itemize}[label=-]
  \item generating $0$\nbd{}cell: $A$,
  \item generating $1$\nbd{}cells: $f,g : A \to A$,
  \item generating $2$\nbd{}cells: $\alpha : f \Rightarrow 1_A$, $\beta : g
    \Rightarrow 1_A$.
  \end{itemize}
  In pictures, this gives
  \[
    \begin{tikzcd}[column sep=huge]
      A \ar[r,bend left=75,"f",""{name=A,below}]\ar[r,bend
      right=75,"g"',""{name=B,above}] \ar[r,"1_A"
      pos=1/3,""{name=C,above},""{name=D,below}]& A
      \ar[from=A,to=C,Rightarrow,"\alpha"] \ar[from=B,to=D,"\beta"
      pos=9/20,Rightarrow]
    \end{tikzcd}
    \text{ or }
    \begin{tikzcd}
      A. \ar[loop,in=50,out=130,distance=1.5cm,"f",""{name=A,below}]
      \ar[loop,in=-50,out=-130,distance=1.5cm,"g"',""{name=B,above}]
      \ar[from=A,to=1-1,Rightarrow,"\alpha"]
      \ar[from=B,to=1-1,Rightarrow,"\beta"]
    \end{tikzcd}
  \]
  Notice that this category has many non-trivial $2$\nbd{}cells and it is not
  \emph{a priori} clear what its homotopy type is and whether or not it is
  \good{}. Observe that $P$ is obtained as the following amalgamated sum
  \begin{equation}\label{square:lemniscate}
    \begin{tikzcd}
      \sD_0 \ar[r] \ar[d] & A_{(1,0)} \ar[d]\\
      A_{(1,0)} \ar[r] & P. \ar[from=1-1,to=2-2,very near
      end,"\ulcorner",phantom]
    \end{tikzcd}
  \end{equation}
  Since $\sD_0$, $A_{(1,0)}$ are free and \good{} and since
  $\sD_0 \to A_{(1,0)}$ is a folk cofibration, all we
  have to show to prove that $P$ is \good{} is that the above square is Thomason
  homotopy cocartesian. Notice that the $2$\nbd{}category $A_{(1,0)}$ is
  obtained as the following amalgamated sum
  \[
    \begin{tikzcd}
      \Delta_1 \ar[r,"\tau"] \ar[d]& A_{(1,1)} \ar[d] \\
      \Delta_0 \ar[r] &A_{(1,0)}, \ar[from=1-1,to=2-2,very near
      end,"\ulcorner",phantom]
    \end{tikzcd}
  \]
  where $\tau : \Delta_1 \to A_{(1,1)}$ has already been defined in
  \ref{paragr:Amn}. We have seen that $\tau$ is a co\nbd{}universal Thomason
  equivalence and thus, so is $\sD_0 \to A_{(0,1)}$ (as $\sD_0$ and
  $\Delta_0$ are two different names for the same category). Hence, square
  \eqref{square:lemniscate} is Thomason homotopy cocartesian and this proves that
  $P$ is \good{} and has the homotopy type of a point.

  All the variations by reversing the direction of $\alpha$ or $\beta$ work
  exactly the same way.
\end{paragr}
Let us now get into more sophisticated examples.
\begin{paragr}[Variations of spheres]\label{paragr:variationsphere}
  Let $P$ the free $2$\nbd{}category defined as follows:
  \begin{itemize}[label=-]
  \item generating $0$\nbd{}cells: $A,B$,
  \item generating $1$\nbd{}cells: $f,g : A \to B$,
  \item generating $2$\nbd{}cells: $\alpha : f \Rightarrow g$, $\beta: g
    \Rightarrow f$.
  \end{itemize}
  In pictures, this gives
  \[
    \begin{tikzcd}
      A \ar[r,bend left=75,"f",""{name=A,below}] \ar[r,bend
      right=75,"g"',""{name=B,above}] & B. \ar[from=A,to=B,bend
      right,Rightarrow,"\alpha"'] \ar[from=B,to=A,bend
      right,Rightarrow,"\beta"']
    \end{tikzcd}
  \]
    
  Let $P'$ be the free $2$\nbd{}category defined as follows:
  \begin{itemize}[label=-]
  \item generating $0$\nbd{}cells: $A',B'$,
  \item generating $1$\nbd{}cell: $h:A'\to B'$,
  \item generating $2$\nbd{}cell: $\gamma : h \Rightarrow h$,
  \end{itemize}
  which can be pictured as
  \[
    \begin{tikzcd}
      A' \ar[r,bend left=75,"h",""{name=A,below}] \ar[r,bend right=75,"h"']& B'
      \ar[from=A,to=B,Rightarrow,"\gamma"]
    \end{tikzcd}
    \text{ or }
    \begin{tikzcd}
      A' \ar[r,"h"',""{name=A,above}] & B' \ar[from=A,to=A,loop, in=130,
      out=50,distance=1cm, Rightarrow,"\gamma"']
    \end{tikzcd}
  \]
  and let $F : P \to P'$ be the unique $2$\nbd{}functor such that
  \begin{itemize}[label=-]
  \item $F(A)=A'$ and $F(B)=B'$,
  \item $F(f)=F(g)=h$,
  \item $F(\alpha)=\gamma$ and $F(\beta)=1_h$.
  \end{itemize}
  We wish to prove that this $2$\nbd{}functor is a Thomason equivalence. Since
  it is an isomorphism on objects, it suffices to prove that the functors induced
  by $F$
  \[
    F_{A,A} : P(A,A) \to P'(F(A),F(A)),
  \]
  \[
    F_{B,B} : P(B,B) \to P'(F(B),F(B)),
  \]
  \[
    F_{B,A} : P(B,A) \to P'(F(B),F(A))
  \]
  and
  \[
    F_{A,B} : P(A,B) \to P'(F(A),F(B))
  \]
  are Thomason equivalences of categories (Corollary
  \ref{cor:criterionThomeqI}). For the first two ones, this follows trivially
  from the fact that the categories $P(A,A)$, $P'(A',A')$, $P(B,B)$ and
  $P'(B',B')$ are all isomorphic to $\sD_0$. For the third one, this follows
  trivially from the fact that the categories $P(B,A)$ and $P'(B',A')$ are the
  empty category. For the fourth one, this can be seen as follows. The category $P(A,B)$ is the free category on the
  graph
  \[
    \begin{tikzcd}
      f \ar[r,shift left,"\alpha"] & g \ar[l,shift left,"\beta"]
    \end{tikzcd}
  \]
  ($2$\nbd{}cells of $P$ become $1$\nbd{}cells of $P(A,B)$ and $1$\nbd{}cells of
  $P$ become $0$\nbd{}cells of $P(A,B)$) and the category
  $P'(F(A),F(B))=P'(A',B')$ is the free category on the graph
  \[
    \begin{tikzcd}
      h. \ar[loop above,"\gamma"]
    \end{tikzcd}
  \]
  The functor $F_{A,B}$ comes from a morphism of reflexive graphs and is
  obtained by ``killing the generator $\beta$'' (see Example
  \ref{example:killinggenerator}). In particular, the square
  \[
    \begin{tikzcd}
      \sD_1 \ar[r,"\langle \beta \rangle"] \ar[d]& P(A,B) \ar[d,"F_{A,B}"] \\
      \sD_0 \ar[r,"\langle h \rangle" ] & P'(A',B')
    \end{tikzcd}
  \]
  is Thomason homotopy cocartesian and thus, $F_{(A,B)}$ is a Thomason
  equivalence.

  Now consider (a copy of) $\sS_2$ labelled as follows:
  \[
    \begin{tikzcd}
      \overline{A} \ar[r,bend left=75,"i",""{name=A,below}] \ar[r,bend
      right=75,"j"',""{name=B,above}] & \overline{B} \ar[from=A,to=B,bend
      right,Rightarrow,"\delta"'] \ar[from=A,to=B,bend left
      ,Rightarrow,"\epsilon"]
    \end{tikzcd}
  \]
  and let $G : \sS_2 \to P'$ be the unique $2$\nbd{}functor such that
  \begin{itemize}[label=-]
  \item $G(\overline{A})=A'$ and $G(\overline{B})=B'$,
  \item $G(i)=G(j)=h$,
  \item $G(\delta)=\gamma$ and $G(\epsilon)=1_h$.
  \end{itemize}
  For similar reasons as for $F$, the $2$\nbd{}functor $G$ is a Thomason
  equivalence. This proves that both $P'$ and $P$ have the homotopy type of
  $\sS_2$.

  Now, let $P''$ be the free $2$\nbd{}category defined as follows:
  \begin{itemize}[label=-]
  \item generating $0$\nbd{}cell: $A''$,
  \item generating $1$\nbd{}cell: $l : A'' \to A''$,
  \item generating $2$\nbd{}cells: $\lambda : l \Rightarrow 1_{A''}$ and $\mu: l
    \Rightarrow 1_{A''}$.
  \end{itemize}
  In pictures, this gives
  \[
    \begin{tikzcd}[column sep=huge]
      A'' \ar[r,bend left=75,"l",""{name=A'',below}]\ar[r,bend
      right=75,"l"',""{name=B,above}] \ar[r,"1_{A''}"
      pos=1/3,""{name=C,above},""{name=D,below}]& A''
      \ar[from=A'',to=C,Rightarrow,"\lambda"] \ar[from=B,to=D,"\mu"
      pos=9/20,Rightarrow]
    \end{tikzcd}
    \qquad \quad \text{ or }
    \begin{tikzcd}
      A''. \ar[loop,in=30,out=150,distance=3cm,"l",""{name=A,below}]
      \ar[from=A,to=1-1,bend right,Rightarrow,"\lambda"'] \ar[from=A,to=1-1,bend
      left,Rightarrow,"\mu"]
    \end{tikzcd}
  \]
  Let $H : \sS_2 \to P''$ be the unique $2$\nbd{}functor such that:
  \begin{itemize}[label=-]
  \item $H(\overline{A})=H(\overline{B})=A''$,
  \item $H(i)=l$ and $H(j)=1_{A''}$,
  \item $H(\delta)=\lambda$ and $H(\epsilon)=\mu$.
  \end{itemize}
  Let us prove that $H$ is a Thomason equivalence using Corollary
  \ref{cor:criterionThomeqII}. In order to do so, we have to compute $V_k(H) :
  V_k(\sS_2) \to V_k(P'')$ for every $k\geq 0$. For $k=0$, the category
  $V_0(\sS_2)$ is the free category on the graph
  \[
    \begin{tikzcd}
      \overline{A} \ar[r,"i",shift left] \ar[r,"j"',shift right] & \overline{B},
    \end{tikzcd}
  \]
  the category $V_0(P'')$ is the free category on the graph
  \[
    \begin{tikzcd}
      A \ar[loop above,"l"]
    \end{tikzcd}
  \]
  and $V_0(H)$ comes from a morphism of reflexive graphs obtained by ``killing
  the generator $j$''. Hence, it is a Thomason equivalence of categories. For
  $k>0$, the category $V_k(\sS_2)$ has two objects $\overline{A}$ and
  $\overline{B}$ and an arrow $\overline{A} \to \overline{B}$ is a
  $k$\nbd{}tuple of one of the following forms:
  \begin{itemize}[label=-]
  \item $(1_i,\cdots,1_i,\delta,1_j,\cdots,1_j)$,
  \item $(1_i,\cdots,1_i,\epsilon,1_j,\cdots,1_j)$,
  \item $(1_i,\cdots,1_i)$,
  \item $(1_j,\cdots,1_j)$,
  \end{itemize}
  and these are the only non-trivial arrows. In other words, $V_k(\sS_2)$ is the
  free category on the graph with two objects and $2k+2$ parallel arrows between
  these two objects. In order to compute $V_k(P'')$, let us first notice that
  every $2$\nbd{}cell of $P''$ (except for $\1^2_{A''}$) is uniquely encoded as
  a finite word on the alphabet that has three symbols : $1_l$, $\lambda$ and
  $\mu$. Concatenation corresponding to the $0$\nbd{}composition of these cells.
  This means exactly that $V_1(P'')$ is free on the graph that has one object
  and three arrows. More generally, it is a tedious but harmless exercise to
  prove that for every $k>0$, the  category $V_k(P'')$ is the
  free category on the graph that has one object $A''$ and $2k+1$
  arrows which are of one of the following forms:
  \begin{itemize}[label=-]
  \item $(1_l,\cdots,1_l,\lambda,1^2_{A''},\cdots,1^2_{A''})$,
  \item $(1_l,\cdots,1_l,\mu,1^2_{A''},\cdots,1^2_{A''})$,
  \item $(1_l,\cdots,1_l)$.
  \end{itemize}
  Once again, the functor $V_k(H)$ comes from a morphism of reflexive graphs and
  is obtained by ``killing the generator $(1_j,\cdots,1_j)$''. Hence, it is a
  Thomason equivalence and thus, so is $H$. This proves that $P''$ has the
  homotopy type of $\sS_2$.
  
  Finally, consider the commutative diagram of $\ho(\Ch)$
  \[
    \begin{tikzcd}[column sep=huge]
      \sH^{\sing}(P) \ar[r,"\sH^{\sing}(F)"] \ar[d,"\pi_P"] & \sH^{\sing}(P')
      \ar[d,"\pi_{P'}"] & \sH^{\sing}(\sS_2) \ar[l,"\sH^{\sing}(G)"']
      \ar[d,"\pi_{\sS_2}"] \ar[r,"\sH^{\sing}(H)"]&\sH^{\sing}(P'') \ar[d,"\pi_{P''}"]\\
      \sH^{\pol}(P) \ar[r,"\sH^{\pol}(F)"'] & \sH^{\pol}(P') & \sH^{\pol}(\sS_2)
      \ar[l,"\sH^{\pol}(G)"] \ar[r,"\sH^{\pol}(H)"']&\sH^{\pol}(P'').
    \end{tikzcd}
  \]
  Since $F$, $G$ and $H$ are Thomason equivalences, the three top horizontal
  morphisms are isomorphisms. Besides, a simple computation using Proposition
  \ref{prop:abelianizationfreeoocat}, which we leave to the reader, shows that
  the three bottom horizontal morphisms are also isomorphisms. Since $\sS_2$ is
  \good{} (Proposition \ref{prop:spheresaregood}), the morphism $\pi_{\sS_2}$ is
  an isomorphism. This implies that $\pi_{P'}$, $\pi_{P'}$ and $\pi_{P''}$ are
  isomorphisms, which means by definition that $P$, $P'$ and $P''$ are \good{}.
\end{paragr}
\begin{paragr}\label{paragr:anothercounterexample}
  Let $P$ be the free $2$\nbd{}category defined as follows:
  \begin{itemize}[label=-]
  \item generating $0$\nbd{}cell: $A$,
  \item generating $1$\nbd{}cell: $f : A \to A$,
  \item generating $2$\nbd{}cells: $\alpha : f \Rightarrow 1_A$ and $\beta: 1_A
    \Rightarrow f$.
  \end{itemize}
  In pictures, this gives
  \[
    \begin{tikzcd}[column sep=huge]
      A \ar[r,bend left=75,"f",""{name=A,below}]\ar[r,bend
      right=75,"f"',""{name=B,above}] \ar[r,"1_A"
      pos=1/3,""{name=C,above},""{name=D,below}]& A
      \ar[from=A,to=C,Rightarrow,"\alpha"] \ar[from=D,to=B,"\beta"
      pos=9/20,Rightarrow]
    \end{tikzcd}
    \qquad \quad \text{ or }
    \begin{tikzcd}
      A. \ar[loop,in=30,out=150,distance=3cm,"f",""{name=A,below}]
      \ar[from=A,to=1-1,bend right,Rightarrow,"\alpha"'] \ar[from=1-1,to=A,bend
      right,Rightarrow,"\beta"']
    \end{tikzcd}
  \]
  Now, let $P'$ be a copy of $B^2\mathbb{N}$ labelled as follows
  \[
    \begin{tikzcd}
      A' \ar[r,bend left=75,"1_{A'}",""{name=A,below}] \ar[r,bend
      right=75,"1_{A'}"',""{name=B,above}] & A',
      \ar[from=A,to=B,"\gamma",Rightarrow]
    \end{tikzcd}
  \]
  which can be also pictured as
  \[
    \begin{tikzcd}
      A'. \ar[loop,in=120,out=60,distance=1.5cm,"\gamma"',Rightarrow]
    \end{tikzcd}
  \]
  Let $F : P \to P'$ be the unique $2$\nbd{}functor such that:
  \begin{itemize}[label=-]
  \item $F(A)=A'$,
  \item $F(f)=1_{A'}$
  \item $F(\alpha)=\1^2_{A'}$ and $F(\beta)=\gamma$,
  \end{itemize}
  and $G : P' \to P$ be the unique $2$\nbd{}functor such that:
  \begin{itemize}[label=-]
  \item $G(A')=A$,
  \item $G(\gamma)=\alpha\comp_1\beta$.
  \end{itemize}
  Notice that we have $F\circ G = \mathrm{id}_{P'}$ and that we have an oplax
  transformation \[h \colon \mathrm{id}_{P} \Rightarrow G \circ F\] defined as
  \begin{itemize}[label=-]
  \item $h_A:=1_A$,
    \item $h_f:=\alpha$.
    \end{itemize}
  Hence, $F$ is a Thomason equivalence (Proposition \ref{prop:oplaxhmtpyisthom}) and $P$ has the homotopy type of a
  $K(\mathbb{Z},2)$ (see \ref{paragr:bubble}). In particular, it has non-trivial
  singular homology groups in all even dimension; but since it is a
  free $2$\nbd{}category, all its polygraphic homology groups are trivial strictly above
  dimension $2$, which means that $P$ is \emph{not} \good{}.
\end{paragr}
\begin{paragr}\label{paragr:sumupsphere}
  All the results from \ref{paragr:variationsphere} and
  \ref{paragr:anothercounterexample} are summed up by the following table.
  
  \begin{center}
    \begin{tabular}{| l || c | c |}
      \hline
      $2$\nbd{}category & \good{}? & homotopy type \\ \hline \hline
      {
      $\begin{tikzcd}
        \bullet \ar[r,bend
        left=75,""{name=A,below,pos=9/20},""{name=C,below,pos=11/20}]
        \ar[r,bend
        right=75,""{name=B,above,pos=9/20},""{name=D,above,pos=11/20}] & \bullet
        \ar[from=C,to=D,bend left,Rightarrow] \ar[from=A,to=B,bend
        right,Rightarrow]
      \end{tikzcd}$
                                                                         } & yes & $\sS_2$\\
    
      \hline
      { $ \begin{tikzcd}
          \bullet \ar[r,bend left=75,""{name=A,below}] \ar[r,bend
          right=75,""{name=B,above}] & \bullet \ar[from=A,to=B,bend
          right,Rightarrow] \ar[from=B,to=A,bend
          right,Rightarrow]
        \end{tikzcd}$} & yes & $\sS_2$ \\ \hline {$ \begin{tikzcd} \bullet
          \ar[r,""{name=A,above}] & \bullet \ar[from=A,to=A,loop, in=130,
          out=50,distance=1cm, Rightarrow] \end{tikzcd}$} & yes &$\sS_2$ \\
      \hline
    
      {
      $\begin{tikzcd}
        \bullet \ar[loop,in=30,out=150,distance=2cm,""{name=A,below}]
        \ar[from=A,to=1-1,bend right,Rightarrow]
        \ar[from=A,to=1-1,bend left,Rightarrow]
      \end{tikzcd}$ } & yes & $\sS_2$ \\ \hline { $\begin{tikzcd} \bullet
        \ar[loop,in=30,out=150,distance=2cm,""{name=A,below}]
        \ar[from=A,to=1-1,bend right,Rightarrow] \ar[from=1-1,to=A,bend
        right,Rightarrow]
      \end{tikzcd}$ } & no & $K(\mathbb{Z},2)$ \\ \hline {$\begin{tikzcd}
        \bullet \ar[loop,in=120,out=60,distance=1.2cm,Rightarrow]
      \end{tikzcd}$} & no & $K(\mathbb{Z},2)$ \\ \hline
    \end{tabular}
  \end{center}
  Notice that the fourth and fifth entries of this table only differ by the
  direction of a generating $2$\nbd{}cell but the homotopy types are not the same. 
\end{paragr}
Let us now move on to bouquets of spheres.
\begin{paragr}
  Let $P$ be the free $2$\nbd{}category defined as follows:
  \begin{itemize}[label=-]
  \item generating $0$\nbd{}cell: $A$,
  \item generating $1$\nbd{}cells: $f,g: A \to A$,
  \item generating $2$\nbd{}cells: $\alpha,\beta : f \Rightarrow g$.
  \end{itemize}
  In pictures, this gives:
  \[
    \begin{tikzcd}[column sep=huge]
      A \ar[r,bend
      left=75,"f",""{name=A,below,pos=9/20},""{name=C,below,pos=11/20}]
      \ar[r,bend
      right=75,"g"',""{name=B,above,pos=9/20},""{name=D,above,pos=11/20}] & A.
      \ar[from=C,to=D,bend left,"\alpha",Rightarrow] \ar[from=A,to=B,bend
      right,"\beta"',Rightarrow]
    \end{tikzcd}
  \]
Now let $\sS_2$ be labelled as
 follows:
 \[
   \begin{tikzcd}
         C \ar[r,bend
      left=75,"h",""{name=A,below,pos=9/20},""{name=C,below,pos=11/20}]
      \ar[r,bend
      right=75,"i"',""{name=B,above,pos=9/20},""{name=D,above,pos=11/20}] & D.
      \ar[from=C,to=D,bend left,"\gamma",Rightarrow] \ar[from=A,to=B,bend
      right,"\delta"',Rightarrow]
   \end{tikzcd}
 \]
 Notice that $P$ is obtained as the following amalgamated sum:

   \begin{equation}
   \begin{tikzcd}
     \sS_0 \ar[d,"p"] \ar[r,"{\langle C,D \rangle}"] & \sS_2 \ar[d]\\
     \sD_0 \ar[r] & P.
     \ar[from=1-1,to=2-2,"\ulcorner",phantom,very near end]
   \end{tikzcd}\label{eq:squarebouquethybrid}
   \end{equation}

 Let us prove that this square is Thomason homotopy cocartesian using the first
 part of Proposition \ref{prop:critverthorThomhmtpysquare}. This means that we
 have to show that the induced square of $\Cat$
 \begin{equation}
   \begin{tikzcd}
     V_k(\sS_0) \ar[r,"{V_k(\langle C,D \rangle)}"] \ar[d,"V_k(p)"] & V_k(\sS_2) \ar[d]\\
     V_k(\sD_0) \ar[r] & V_k(P)
   \end{tikzcd}\label{eq:squarebouquethybridvertical}
 \end{equation}
 is Thomason homotopy cocartesian for every $k \geq 0$. Notice first that we
 trivially have $V_k(\sS_0)\simeq \sS_0$ and $V_k(\sD_0)\simeq \sD_0$ for every
 $k\geq 0$ and that $V_0(\sS_2)$ is the free category on the
 graph
 \[   
   \begin{tikzcd}
     C \ar[r,shift left,"f"] \ar[r,shift right,"g"'] & D.
   \end{tikzcd}
 \]
 Besides, $V_0(P)$ is the free category on the
 graph
 \[
   \begin{tikzcd}
     A. \ar[loop above,"f"] \ar[loop below,"g"]
   \end{tikzcd}
 \]
 In particular, square \eqref{eq:squarebouquethybridvertical} is
 cocartesian for $k=0$ and we are in the situation of identification
 of two objects of a free category (see Example
 \ref{example:identifyingobjects}). Hence, square
 \eqref{eq:squarebouquethybridvertical} is Thomason cocartesian for
 $k=0$.  Similarly, for $k>0$, we have already seen that $V_k(\sS_2)$
 is the free category on the graph that has $2$ objects and $2k+2$
 parallel arrows between these two objects and we leave as an easy
 exercise to the reader to check that the category $V_k(P)$ is the
 free category on the graph that has one object and $2k+2$ arrows, which are
 the $k$\nbd{}tuples of one of the following forms:
 \begin{itemize}[label=-]
 \item $(1_f,\cdots,\alpha,\cdots,1_g)$,
 \item $(1_f,\cdots,\beta,\cdots,1_g)$,
 \item $(1_f,\cdots,1_f)$,
 \item $(1_g,\cdots,1_g)$.
 \end{itemize}
 In particular, square \eqref{eq:squarebouquethybridvertical} is again a
 cocartesian square of identification of two objects of a free category, and
 thus, it is Thomason homotopy cocartesian. This implies that square
 \eqref{eq:squarebouquethybrid} is Thomason homotopy cocartesian. Since $\sS_0$,
 $\sD_0$ and $\sS_2$ are \good{} and since  $\langle C, D \rangle : \sS_0 \to
 \sS_2$ is a folk cofibration, this proves that $P$ is \good{}
 and has the homotopy type of the bouquet of a $1$\nbd{}sphere with a $2$\nbd{}sphere.
\end{paragr}

\begin{paragr}
  Let $P$ be the free $2$\nbd{}category defined as follows:
  \begin{itemize}[label=-]
  \item generating $0$\nbd{}cells: $A$ and $B$,
  \item generating $1$\nbd{}cells: $f,g : A \to B$,
  \item generating $2$\nbd{}cells: $\alpha,\beta,\gamma : f \Rightarrow g$.
  \end{itemize}
  In pictures, this gives
  \[
    \begin{tikzcd}[column sep=huge]
      A \ar[r,bend
      left=75,"f",""{name=A,below,pos=8/20},""{name=C,below,pos=1/2},""{name=E,below,pos=12/20}]
      \ar[r,bend
      right=75,"g"',""{name=B,above,pos=8/20},""{name=D,above,pos=1/2},""{name=F,above,pos=12/20}]&
      B. \ar[from=A,to=B,Rightarrow,"\alpha"',bend right]
      \ar[from=C,to=D,Rightarrow,"\beta"]
      \ar[from=E,to=F,Rightarrow,"\gamma",bend left]
    \end{tikzcd}
  \]
  Now let $P'$ be the sub-$2$\nbd{}category of $P$ spanned by $A$, $B$, $f$, $g$, $\alpha$
  and $\beta$, and let $P''$ be the sub-$2$\nbd{}category of $P$ spanned by $A$,
  $B$, $f$, $g$, $\beta$ and $\gamma$. These $2$\nbd{}categories are simply copies of
  $\sS_2$. Notice that we have a cocartesian square
  \begin{equation}\label{square:bouquet}
    \begin{tikzcd}
      \sD_2 \ar[r,"\langle \beta \rangle"] \ar[d,"\langle \beta \rangle"] &
      P' \ar[d] \\
      P'' \ar[r] & P, \ar[from=1-1,to=2-2,phantom,"\ulcorner",very near end]
    \end{tikzcd}
  \end{equation}
  and by reasoning as in the proof of Lemma \ref{lemma:squarenerve}, one can
  show that the square induced by the nerve
  \[
    \begin{tikzcd}
      N_{\oo}(\sD_2) \ar[r,"\langle \beta \rangle"] \ar[d,"\langle \beta
      \rangle"] &
      N_{\oo}(P') \ar[d] \\
      N_{\oo}(P'') \ar[r] & N_{\oo}(P)
    \end{tikzcd}
  \]
  is also cocartesian. Since $\langle \beta \rangle : \sD_2 \to P'$ and $\langle
  \beta \rangle : \sD_2 \to P''$ are monomorphisms and $N_{\oo}$ preserves
  monomorphisms, it follows from Lemma \ref{lemma:hmtpycocartesianreedy} that
  square \eqref{square:bouquet} is Thomason homotopy cocartesian and in
  particular that $P$ has the homotopy type of a bouquet of two
  $2$\nbd{}spheres. Since $\sD_2$, $P'$ and $P''$ are free and \good{} and since
  $\langle \beta \rangle : \sD_2 \to P'$ and $\langle \beta \rangle : \sD_2 \to
  P''$ are folk cofibrations, this also proves that $P$ is \good{} (see \ref{paragr:criterion2cat}).
\end{paragr}
\begin{paragr}
  Let $P$ be the free $2$\nbd{}category defined as follows:
  \begin{itemize}[label=-]
  \item generating $0$\nbd{}cells: $A$ and $B$,
  \item generating $1$\nbd{}cells: $f,g,h : A \to B$,
  \item generating $2$\nbd{}cells: $\alpha,\beta:f \Rightarrow g$ and
    $\delta,\gamma:g \Rightarrow h$.
  \end{itemize}
  In pictures, this gives:
  \[
    \begin{tikzcd}[column sep=huge]
      A \ar[r,bend
      left=75,"f",""{name=A,below,pos=8/20},""{name=E,below,pos=12/20}]
      \ar[r,"g",""{name=B,above,pos=8/20},""{name=C,below,pos=8/20},""{name=F,above,pos=12/20},""{name=G,below,pos=12/20}]
      \ar[r,bend
      right=75,"h"',""{name=D,above,pos=8/20},""{name=H,above,pos=12/20}] & B.
      \ar[from=A,to=B,Rightarrow,"\alpha"',bend right]
      \ar[from=C,to=D,Rightarrow,"\gamma"',bend right]
      \ar[from=E,to=F,Rightarrow,"\beta",bend left]
      \ar[from=G,to=H,Rightarrow,"\delta",bend left]
    \end{tikzcd}
  \]
  Let us prove that this $2$\nbd{}category is \good{}. Let $P_0$ be the
  sub-$1$\nbd{}category of $P$ spanned by $A$, $B$ and $g$, let $P_1$ be the
  sub-$2$\nbd{}category of $P$ spanned by $A$, $B$, $g$, $h$, $\gamma$ and
  $\delta$ and let $P_2$ be the sub-$2$\nbd{}category of $P$ spanned by $A$,
  $B$, $f$, $g$, $\alpha$ and $\beta$. The $2$\nbd{}categories $P_1$ and $P_2$
  are copies of $\sS_2$, and $P_0$ is a copy of $\sD_1$. Moreover, we have a
  cocartesian square of inclusions
  \begin{equation}\label{squarebouquet}
    \begin{tikzcd}
      P_0 \ar[d,hook] \ar[r,hook] & P_2 \ar[d,hook] \\\
      P_1 \ar[r,hook] & P.
    \end{tikzcd}
  \end{equation}
  Let us prove that this square is Thomason homotopy cocartesian using the
  second part of Corollary \ref{prop:critverthorThomhmtpysquare}. This means
  that we have to show that for every $k \geq 0$, the induced square of $\Cat$
  \begin{equation}\label{squarebouquethorizontal}
    \begin{tikzcd}
      S_k(P_0) \ar[d] \ar[r] & S_k(P_2) \ar[d] \\\
      S_k(P_1) \ar[r] & S_k(P)
    \end{tikzcd}
  \end{equation}
  is Thomason homotopy cocartesian. For $k=0$, this is obvious since all of the
  morphisms of square \eqref{squarebouquet} are isomorphisms at the level of
  objects and the functor $S_0$ is the functor that sends a $2$\nbd{}category to
  its set of objects (seen as a discrete category). Now, notice that the
  categories $P_i(A,A)$ and $P_i(B,B)$ for $0 \leq i \leq 2$ are all isomorphic
  to the terminal category $\sD_0$ and the categories $P_i(B,A)$ for $0 \leq i
  \leq 2$ are all the empty category. It follows that for $k>0$, we have
  \[
    S_k(P_i)\simeq \sD_0\amalg \left( \coprod_{E_k}P_i(A,B) \right)\amalg \sD_0
  \]
  where $E_k$ is the set of all $k$\nbd{}tuples of the form
  \[
    (A,\cdots,A,B,\cdots,B).
  \]
  The set $E_1$ is empty and thus all of the morphisms of square
  \eqref{squarebouquethorizontal} for the value $k=1$ are isomorphisms. This
  makes square \eqref{squarebouquethorizontal} Thomason homotopy cocartesian for
  $k=1$. For $k>1$, notice that the categories $P(A,B)$, $P_2(A,B)$ and
  $P_1(A,B)$ are respectively free on the graphs
  \[
    \begin{tikzcd}
      f \ar[r,shift left,"\alpha"] \ar[r, shift right, "\beta"'] & g \ar[r,shift
      left,"\gamma"] \ar[r,shift right,"\delta"'] & h,
    \end{tikzcd}
  \]
  \[
    \begin{tikzcd}
      f \ar[r,shift left,"\alpha"] \ar[r, shift right, "\beta"'] & g,
    \end{tikzcd}
  \]
  and
  \[
    \begin{tikzcd}
      g \ar[r,shift left,"\gamma"] \ar[r,shift right,"\delta"'] & h.
    \end{tikzcd}
  \]
  It is then straightforward to check that we are in a situation where Corollary
  \ref{cor:hmtpysquaregraph} applies and thus square
  \eqref{squarebouquethorizontal} is Thomason homotopy cocartesian for $k\geq
  1$. Altogether, this proves that square \eqref{squarebouquet} is Thomason homotopy
  cocartesian. Since the inclusions $P_0 \hookrightarrow P_1$ and $P_0 \hookrightarrow P_2$ are folk cofibrations and since $P_0$, $P_1$ and $P_2$ are \good{}, this proves that $P$ is \good{} and has the homotopy type of a bouquet of
  two $2$\nbd{}spheres.
\end{paragr}
\begin{paragr}
  Let $P$ be the free $2$\nbd{}category defined as follows:
  \begin{itemize}[label=-]
  \item generating $0$\nbd{}cells: $A$, $B$ and $C$,
  \item generating $1$\nbd{}cells: $f,g : A \to B$ and $h,i : B \to C$,
  \item generating $2$\nbd{}cells: $\alpha,\beta : f \Rightarrow g$ and
    $\gamma,\delta : h \Rightarrow i$.
  \end{itemize}
  In pictures, this gives:
  \[
    \begin{tikzcd}[column sep=huge]
      A \ar[r,bend
      left=75,"f",""{name=A,below,pos=8/20},""{name=C,below,pos=12/20}]
      \ar[r,bend
      right=75,"g"',""{name=B,above,pos=8/20},""{name=D,above,pos=12/20}] & B
      \ar[r,bend
      left=75,"h",""{name=E,below,pos=8/20},""{name=G,below,pos=12/20}]
      \ar[r,bend
      right=75,"i"',""{name=F,above,pos=8/20},""{name=H,above,pos=12/20}] & C.
      \ar[from=A,to=B,bend right,"\alpha",Rightarrow] \ar[from=C,to=D,bend
      left,"\beta",Rightarrow] \ar[from=E,to=F,bend right,"\gamma",Rightarrow]
      \ar[from=G,to=H,bend left,"\delta",Rightarrow]
    \end{tikzcd}
  \]
  Let us prove that $P$ is \good{}. Let $P'$ be the sub-$2$\nbd{}category of $P$
  spanned by $A$, $B$, $f$, $g$, $\alpha$ and $\beta$ and let $P''$ be the
  sub-$2$\nbd{}category of $P$ spanned by $B$, $C$, $h$, $i$, $\gamma$ and
  $\delta$. These two $2$\nbd{}categories are copies of $\sS_2$ and we have a
  cocartesian square
  \begin{equation}\label{squarebouquetbis}
    \begin{tikzcd}
      \sD_0 \ar[r,"\langle B \rangle"] \ar[d,"\langle B \rangle"] & P' \ar[d,hook]\\
      P'' \ar[r,hook] & P, \ar[from=1-1,to=2-2,phantom,very near end]
    \end{tikzcd}
  \end{equation}
  where the anonymous arrows are the canonical inclusions. Let us prove that
  this cocartesian square is Thomason homotopy cocartesian using the first part
  of Proposition \ref{prop:critverthorThomhmtpysquare}. This means that we need
  to prove that for every $k\geq 0$, the induced square of $\Cat$
  \begin{equation}\label{squarebouquetvertical}
    \begin{tikzcd}
      V_k(\sD_0) \ar[r] \ar[d] & V_k(P') \ar[d]\\
      V_k(P'') \ar[r] & V_k(P)
    \end{tikzcd}
  \end{equation}
  is Thomason homotopy cocartesian. For every $k\geq 0$, we have $V_k(\sD_0) \simeq \sD_0$
  and for $k=0$ the categories $V_0(P')$, $V_0(P'')$ and $V_0(P)$ are respectively free on
  the graphs
  \[
    \begin{tikzcd}
      A \ar[r,"f",shift left] \ar[r,"g"',shift right] & B,
    \end{tikzcd}
  \]
  \[
    \begin{tikzcd}
      B \ar[r,"h",shift left] \ar[r,"i"',shift right] & C,
      \end{tikzcd}
    \]
    and
  \[
    \begin{tikzcd}
      A \ar[r,"f",shift left] \ar[r,"g"',shift right] & B \ar[r,"h",shift left] \ar[r,"i"',shift right] & C.
    \end{tikzcd}
  \]
  This implies that square \eqref{squarebouquetvertical} is cocartesian for $k=0$ and in
  virtue of Corollary \ref{cor:hmtpysquaregraph} it is also Thomason homotopy
  cocartesian for this value of $k$. For $k>0$, since $P'$ and $P''$ are both
  isomorphic to $\sS_2$, we have already seen in \ref{paragr:variationsphere} that $V_k(P')$ and $V_k(P'')$ are
  (isomorphic) to the free category on the graph that has two objects and $2k+2$ parallel
  arrows between these two objects.
  Similarly, the category $V_k(P)$ is free on the graph
  that has three objects $A$, $B$, $C$, whose arrows from $A$ to $B$
  are $k$\nbd{}tuples of one of the following form
    \begin{itemize}[label=-]
  \item $(1_f,\cdots,1_f,\alpha,1_g,\cdots,1_g)$,
  \item $(1_f,\cdots,1_f,\beta,1_g,\cdots,1_g)$,
  \item $(1_f,\cdots,1_f)$,
  \item $(1_g,\cdots,1_g)$,
  \end{itemize}
  whose arrows from $B$ to $C$ are $k$\nbd{}tuple of one of the
  following form
      \begin{itemize}[label=-]
  \item $(1_h,\cdots,1_h,\gamma,1_i,\cdots,1_i)$,
  \item $(1_h,\cdots,1_h,\delta,1_i,\cdots,1_i)$,
  \item $(1_h,\cdots,1_h)$,
  \item $(1_i,\cdots,1_i)$,
  \end{itemize}
  and with no other arrows. This implies that square \eqref{squarebouquetvertical}
  is cocartesian for every $k>0$ and in virtue of Corollary
  \ref{cor:hmtpysquaregraph} it is also Thomason homotopy cocartesian for these values of $k$.
  Altogether, this proves that square \eqref{squarebouquetbis} is Thomason
  homotopy cocartesian. Hence, $P$ is \good{} and has the homotopy type of a
  bouquet of two $2$\nbd{}spheres.
\end{paragr}
Let us end this section with an example of a $2$\nbd{}category that has the
homotopy type of the torus.
\begin{paragr}
  Let $P$ be the free $2$\nbd{}category defined as follows:
  \begin{itemize}[label=-]
  \item generating $0$\nbd{}cell: $A$,
  \item generating $1$\nbd{}cells: $f , g : A \to A$,
  \item generating $2$\nbd{}cell: $\alpha : g\comp_0 f \Rightarrow f \comp_0
    g$.
  \end{itemize}
  In pictures, this gives:
  \[
    \begin{tikzcd}
      A \ar[r,"f"] \ar[d,"g"'] & A \ar[d,"g"] \\
      A \ar[r,"f"'] & A. \ar[from=2-1,to=1-2,Rightarrow,"\alpha"]
    \end{tikzcd}
  \]
  From now on, we will use concatenation instead of the symbol $\comp_0$ for the
  $0$\nbd{}composition. For example, $fg$ will stand for $f \comp_0 g$. With
  this notation, the set $1$\nbd{}cells of $P$ is canonically isomorphic to the
  set of finite words in the alphabet $\{f,g\}$ and the set of $2$\nbd{}cells of
  $P$ is canonically isomorphic to the set of finite words in the alphabet
  $\{f,g,\alpha\}$. For a $1$\nbd{}cell $w$ such that $f$ appears $n$ times in
  $w$ and $g$ appears $m$ times in $w$, it is a simple exercise left to the
  reader to show that there exists a unique $2$\nbd{}cell of $P$ from $w$ to the
  word
  \[
    f \cdots fg\cdots g
  \]
  where $f$ is repeated $n$ times and $g$ is repeated $m$ times.
  Recall that we write $B^1(\mathbb{N}\times\mathbb{N})$ for the 
  monoid $\mathbb{N}\times\mathbb{N}$ considered as a category with only one
  object, and let \[F : P \to B^1(\mathbb{N}\times\mathbb{N})\] be the unique
  $2$\nbd{}functor such that:
  \begin{itemize}[label=-]
  \item $F(f)=(1,0)$ and $F(g)=(0,1)$,
  \item $F(\alpha)=1_{(1,1)}$.
  \end{itemize}
  This last equation makes sense since $(1,1)=(0,1)+(1,0)=(1,0)+(0,1)$. For every
  $1$\nbd{}cell $w$ of $P$ (encoded as a finite words in the alphabet $\{f,g\}$)
  such that $f$ appears $n$ times and $g$ appears $m$ times, we have
  $F(w)=(n,m)$. Let us prove that $F$ is a Thomason equivalence using a dual of
  \cite[Corollaire 5.26]{ara2020theoreme} (see Remark 5.20 of op.\ cit.). If we
  write $\star$ for the only object of $B^1(\mathbb{N}\times\mathbb{N})$, what
  we need to show is that the canonical $2$\nbd{}functor from $P/{\star}$ (see
  \ref{paragr:comma}) to the terminal $2$\nbd{}category
  \[
    P/{\star} \to \sD_0
  \]
  is a Thomason equivalence. The $2$\nbd{}category $P/{\star}$ is described as
  follows:
  \begin{itemize}[label=-]
  \item A $0$\nbd{}cell is a $1$\nbd{}cell of $B^1(\mathbb{N}\times
    \mathbb{N})$.
  \item For $(n,m)$ and $(n',m')$ two $0$\nbd{}cells of $P/{\star}$, a
    $1$\nbd{}cell from $(n,m)$ to $(n',m')$ is a $1$\nbd{}cell $w$ of $P$ such
    that the triangle
    \[
      \begin{tikzcd}[column sep=small,row sep=small]
        \star \ar[rr,"F(w)"] \ar[rd,"{(n,m)}"']& & \star \ar[dl,"{(n',m')}"]\\
        &\star&
      \end{tikzcd}
    \]
    is commutative. More explicitly, if $F(w)=(n'',m'')$, the commutativity of
    the previous triangle means
    \[
      n'+n''=n \text{ and } m'+m''=m.
    \]
  \item Given two parallel $1$\nbd{}cells $w$ and $w'$ of $P/\star$, a
    $2$\nbd{}cell of $P/{\star}$ from $w$ to $w'$ is simply a $2$\nbd{}cell of
    $P$ from $w$ to $w'$ seen as $1$\nbd{}cells of $P$.
  \end{itemize}
  From what we said earlier on the $1$\nbd{}cells and $2$\nbd{}cells of $P$, it
  follows easily that for every $0$\nbd{}cell $(n,m)$ of $P/{\star}$, the
  category
  \[
    P/{\star}((m,n),(0,0))
  \]
  has a terminal object, which is given by
  \[
    f\cdots fg\cdots g
  \]
  where $f$ is repeated $n$ times and $g$ is repeated $m$ times. Then, it
  follows from \cite[Théorème 5.27 and Remarque 5.28]{ara2020theoreme} that
  $P/{\star} \to \sD_0$ is a Thomason equivalence of $2$\nbd{}categories and
  this proves that $F$ is a Thomason equivalence. Since
  $B^1(\mathbb{N}\times\mathbb{N})\simeq B^1(\mathbb{N})\times B^1(\mathbb{N})$
  and $B^1(\mathbb{N})$ has the homotopy type of $\sS_1$, we conclude that $P$
  has the homotopy type of $\sS_1 \times \sS_1$, i.e.\ the homotopy type of the
  torus.

  Consider now the commutative square
  \[
    \begin{tikzcd}
      \sH^{\sing}(P) \ar[r,"\sH^{\sing}(F)"]  \ar[d,"\pi_{P}"] & \sH^{\sing}(B^1(\mathbb{N}\times\mathbb{N})) \ar[d,"\pi_{B^1(\mathbb{N}\times\mathbb{N})}"] \\
      \sH^{\pol}(P) \ar[r,"\sH^{\pol}(F)"] &
      \sH^{\pol}(B^1(\mathbb{N}\times\mathbb{N}))
    \end{tikzcd}
  \]
  Since $F$ is a Thomason equivalence, the top horizontal arrow is an
  isomorphism and since $B^1(\mathbb{N}\times\mathbb{N})$ is a
  $1$\nbd{}category, it is \good{} (Theorem \ref{thm:categoriesaregood}), which
  means that the right vertical arrow is an isomorphism. The $1$\nbd{}category
  $B^1(\mathbb{N}\times \mathbb{N})$ is not free but since it has the homotopy
  type of the torus, we have \[H^{\sing}_k(B^1(\mathbb{N}\times \mathbb{N}))=0=H_k^{\pol}(B^1(\mathbb{N}\times \mathbb{N}))\]
  for $k\geq 2$ and it follows then from Corollary \ref{cor:polhmlgycofibrant}
  and Paragraph \ref{paragr:polhmlgylowdimension} that the map canonical map
  \[
    \alpha^{\pol}_{B^{1}(\mathbb{N}\times\mathbb{N})} :
    \sH^{\pol}(B^{1}(\mathbb{N}\times\mathbb{N})) \to \lambda(B^{1}(\mathbb{N}\times\mathbb{N}))
  \]
  is a quasi-isomorphism. Since $P$ is free, it follows that the map $\sH^{\pol}(F)$
  can be identified with the image in $\ho(\Ch)$ of the map
  \[
    \lambda(F) : \lambda(P) \to \lambda(B^1(\mathbb{N}\times\mathbb{N})),
  \]
  which is easily checked to be a quasi-isomorphism. 
  By a 2-out-of-3 property, we deduce that $\pi_P : \sH^{\pol}(P) \to
  \sH^{\sing}(P)$ is an isomorphism, which means by definition that $P$ is
  \good{}.
\end{paragr}
\section{The ``bubble-free'' conjecture}
\begin{definition}
  Let $C$ be a $2$\nbd{}category. A \emph{bubble} (in $C$) is a $2$\nbd{}cell
  $x$ of $C$ such that:
  \begin{itemize}[label=-]
  \item $x$ is not a unit,
  \item $\src_0(x)=\trgt_0(x)$,
  \item $\trgt_1(x)=\src_1(x)=1_{\src_0(x)}$.
  \end{itemize}
\end{definition}
\begin{paragr}
  In pictures, a bubble $x$ is represented as 
  \[
    \begin{tikzcd}
      A \ar[r,bend left=75,"1_A",""{name=A,below}] \ar[r,bend
      right=75,"1_A"',pos=21/40,""{name=B,above}] &A,
      \ar[from=A,to=B,"x",Rightarrow]
    \end{tikzcd}
    \text{ or }
        \begin{tikzcd}
      A \ar[loop,in=120,out=60,distance=1cm,"x"',Rightarrow]
    \end{tikzcd}
  \]
  where $A=\src_0(x)=\trgt_0(x)$.
\end{paragr}
\begin{definition}\label{def:bubblefree2cat}
  A $2$\nbd{}category is said to be \emph{bubble-free} if it has no bubbles.
\end{definition}
\begin{paragr}
  The archetypal example of a $2$\nbd{}category that is \emph{not} bubble-free
  is $B^2\mathbb{N}$. Another non-bubble $2$\nbd{}category is the one from Paragraph \ref{paragr:anothercounterexample}. It is
  remarkable that of all the free $2$\nbd{}categories we have seen so far, these
  are the only examples that are non-\good{}. This motivates the following conjecture.
\end{paragr}
\begin{conjecture}\label{conjecture:bubblefree}
  A free $2$\nbd{}category is \good{} if and only if it is bubble-free.
\end{conjecture}
\begin{paragr}
  At the time of writing, I do not have a real hint towards a proof of the above conjecture.
  Yet, in light of all the examples seen in the previous section, it seems very
  likely to be true. Note that we have also conjectured in Paragraph
  \ref{paragr:conjectureH2} that for every $\oo$\nbd{}category $C$, we have
  \[
    H_2^{\pol}(C)\simeq H_2^{\sing}(C).
  \]
  If this conjecture on the second homology group is true, then conjecture \ref{conjecture:bubblefree} may be reformulated as: A free $2$\nbd{}category $P$ has trivial
  singular homology groups strictly above dimension $2$ if and only if it is bubble-free.
\end{paragr}

\addcontentsline{toc}{chapter}{Bibliography}
\bibliographystyle{alpha}
\bibliography{memoire}
\end{document}